\documentclass{amsbook}

\usepackage{a4wide}
\usepackage[french]{babel}
\usepackage[all]{xy}
\usepackage{latexsym}
\usepackage{amssymb}
\usepackage{fancyhdr}

\makeindex

\theoremstyle{plain}
\newtheorem{pro}{Proposition}
\newtheorem{thm}[pro]{Th\'eor\`eme}
\newtheorem{lem}[pro]{Lemme}
\newtheorem{cor}[pro]{Corollaire}
\newtheorem*{conj}{Conjecture}
\newtheorem*{thmintro}{Th\'eor\`eme}
\theoremstyle{definition}
\newtheorem*{dei}{D\'efinition}
\newtheorem*{deo}{D\'emonstration}

\newcommand{\cone}{\textrm{c\^one}}
\newcommand{\B}{\mathcal{B}}
\newcommand{\oC}{\bar{\mathcal{C}}}

\newcommand{\oPo}{\overline{\Po}}
\newcommand{\Bi}{\mathcal{B}i}
\newcommand{\BLi}{\mathcal{B}i\mathcal{L}ie}
\newcommand{\IBi}{\varepsilon\mathcal{B}i}
\newcommand{\G}{\mathcal{G}}
\newcommand{\N}{\mathcal{N}}
\newcommand{\Sh}{\mathcal{S}}

\newcommand{\kj}{{\bar{k},\, \bar{\jmath}}}
\newcommand{\oi}{{\bar{\imath}}}
\newcommand{\oj}{{\bar{\jmath}}}
\newcommand{\ok}{{\bar{k}}}
\newcommand{\ol}{{\bar{l}}}
\newcommand{\coker}{\mathop{\mathrm{coker}}}
\newcommand{\Colim}{\mathop{\mathrm{Colim}}}
\newcommand{\FS}{\mathcal{FS}}
\newcommand{\Sy}{\mathbb{S}}
\newcommand{\CA}{\mathcal{A}}
\newcommand{\C}{\mathcal{C}om}
\newcommand{\A}{\mathcal{A}s}

\newcommand{\Po}{\mathcal{P}}
\newcommand{\F}{\mathcal{F}}
\newcommand{\ac}{\scriptstyle \textrm{!`}}
\newcommand{\K}{\mathcal{K}}
\newcommand{\II}{I}
\newcommand{\Qo}{\mathcal{Q}}
\newcommand{\Li}{\mathcal{L}ie}
\newcommand{\Sc}{\mathbb{S}^c}

\newcommand{\Co}{\mathcal{C}}
\newcommand{\cqfd}{\ \hfill \square}

\title{\bf Dualit\'e de Koszul des PROPs}
\author{\bf Bruno V.}
\date{}

\begin{document}


\thispagestyle{empty}
\begin{center}
\vbox to 20mm{  \Huge\textbf{\uppercase{Dualit\'e de Koszul des
PROPs\\}}}
\medskip
\Large
\textrm{\huge Bruno Vallette}\\
\medskip
\textrm{2003}
\end{center}
\vfill
\textbf{Adresse:} Bruno Vallette\\
\indent Institut de Recherche Math\'ematique Avanc\'ee \\
\indent Universit\'e Louis Pasteur et C.N.R.S. (UMR 7501)\\
\indent 7, Rue Ren\'e Descartes\\
\indent F-67084 STRASBOURG Cedex\\
\indent France\\
\textbf{Adresse \'electronique :} vallette@math.u-strasbg.fr\\
\textbf{Page personnelle :} {http://www-irma.u-strasbg.fr/\~{
}vallette/}

\cleardoublepage

\newpage

\thispagestyle{empty}

Un grand merci \`a tous ceux qui ont rendu possible l'accomplissement de cette 
th\`ese. \\

Je pense particuli\`erement \`a mon directeur, Jean-Louis 
Loday et aux nombreuses heures qu'il a pass\'ees \`a m'expliquer, avec une passion 
toujours vive, cette ``math\'ematique bleue'' que nous trouvons si belle.  
Je dois aussi beaucoup aux
 diff\'erents membres de mon jury, 
notamment pour avoir lu, relu et critiqu\'e cette th\`ese. \\

Deux personnes m'ont soutenu durant ce travail, comme elles l'ont fait depuis 
si longtemps. La pudeur m'interdit de les nommer ici, mais je sais qu'elles se 
reconna\^\i tront dans ces mots. Cette th\`ese est la v\^otre. 

\cleardoublepage

\newpage

\thispagestyle{empty}
$ \ $
\vspace{3cm}

\hspace{6cm}\begin{tabular}{l}
``C'est gr\^ace au progr\`es fantastique de la science \\
que l'on sait maintenant que,\\
quand on plonge un corps dans une baignoire, \\
le t\'el\'ephone sonne.''\\
\hspace{5cm}Pierre Desproges
\end{tabular}

\cleardoublepage

\pagestyle{fancy}
\setlength{\parindent}{0pt}

\fancyhead{}

\chapter*{Introduction}

\thispagestyle{empty}

\fancyhead{}
\fancyhead[RO]{INTRODUCTION}
\fancyhead[LE]{INTRODUCTION}

Le but de cette th\`ese est d'\'etablir une th\'eorie de dualit\'e
de Koszul pour les PROPs, c'est-\`a-dire les objets qui
mod\'elisent les op\'erations \`a plusieurs entr\'ees et sorties sur
diff\'erents types de structures alg\'ebriques, comme les alg\`ebres et les
big\`ebres par exemple.\\

La dualit\'e de Koszul des alg\`ebres associatives est une
th\'eorie qui a \'et\'e d\'e\-ve\-lo\-pp\'ee par S. Priddy
\cite{Priddy} dans les ann\'ees 1970 . Elle associe \`a toute
alg\`ebre quadratique $A$ une cog\`ebre duale ${A^!}^*$ et un
complexe de cha\^\i nes appel\'e complexe de Koszul. Lorsque ce
dernier est acyclique, on dit que l'alg\`ebre $A$ est de Koszul.
Une telle alg\`ebre, ainsi que ses repr\'esentations, ont de
nombreuses propri\'et\'es (\emph{cf.} les travaux de Beilinson,
Ginzburg et Soergel, entre autres \cite{BGS}).

Dans les ann\'ees 1990, une th\'eorie similaire a \'et\'e
d\'evelopp\'ee par V. Ginzburg et M. M. Kapranov \cite{GK} pour
les op\'erades alg\`ebriques. Une op\'erade est un objet qui
mod\`elise les op\'erations d'un type d'alg\`ebre donn\'e et les
compositions de celles-ci. La dualit\'e de Koszul des op\'erades a
de nombreuses applications : construction d'un ``petit'' complexe
pour le calcul des groupes d'homologie d'une alg\`ebre, notion
d'alg\`ebre \`a homotopie pr\`es, mod\`ele minimal d'une
op\'erade.

Les op\'erades ne tiennent compte que des op\'erations \`a $n$
entr\'ees et une seule sortie. Or, dans le cas
des big\`ebres, on a des op\'erations et aussi des coop\'erations
(\`a plusieurs sorties)
. On doit
alors enrichir la notion d'op\'erade, c'est-\`a-dire travailler avec des
PROPs.\\

Il est naturel d'essayer d'\'etendre la dualit\'e de Koszul des
op\'erades aux PROPs. Plusieurs travaux existent d\'ej\`a dans
cette direction  par W. L. Gan \cite{Gan}, M. Markl et A. A.
Voronov \cite{MV} mais, dans le premier cas par exemple, l'auteur ne traite
qu'une sous-cat\'egorie stricte de PROPs. \\

Pour tout PROP quadratique $\Po$, nous d\'efinissons ici un coPROP
dual, not\'e $\Po^{\ac}$, qui est une g\'en\'eralisation des notions
de cog\`ebre duale et de coop\'erade duale. En outre, nous
g\'en\'eralisons aux PROPs les notions de bar et de cobar
constructions, not\'ees respectivement $\B$ et $\B^c$. Rappelons
que dans le cadre des alg\`ebres et des op\'erades de Koszul,
la cobar construction
fournit une r\'esolution quasi-libre de l'alg\`ebre (ou de l'op\'erade)
de d\'epart.
Nous \'etendons ce th\'eor\`eme aux PROPs.\\

Le principal r\'esultat de cette th\`ese est le th\'eor\`eme
suivant qui donne un crit\`ere pour que la cobar construction sur le
coPROP dual fournisse une r\'esolution quasi-libre du PROP de
d\'epart.

\begin{thmintro}[Crit\`ere de Koszul des PROPs]
Soit $\Po$ un PROP diff\'erentiel quadratique. Les propositions
suivantes sont \'equivalentes :
\begin{enumerate}

\item Le complexe de Koszul $\Po^{\ac}\boxtimes \Po$ est acyclique.

\item Le morphisme naturel de PROPs diff\'erentiels gradu\'es par un poids
$\bar{\B}^c(\Po^{\ac})\to \Po $ est un quasi-isomorphisme.

\end{enumerate}
\end{thmintro}

 Lorsque c'est le cas, on dit que $\Po$ est un \emph{PROP de
Koszul} et la r\'esolution $(2)$ fournit le
\emph{mod\`ele minimal} de $\Po$.\\

Dans cette th\`ese, nous commen\c cons par g\'en\'eraliser les
notions relatives aux anneaux et aux alg\`ebres \`a toute
cat\'egorie mono\"\i dale. Par exemple, on d\'efinit les notions
de module, modules lin\'eaire et multilin\'eaire sur un mono\"\i
de, de produit relatif, et d'id\'eal d'un mono\"\i de. Notons que
ces g\'en\'eralisations ne sont pas imm\'ediates,
notamment lorsque le produit mono\"\i dal n'est pas bilin\'eaire.\\

Nous donnons aussi une construction du mono\"\i de libre. Dans le
cas o\`u le produit mono\"\i dal est biadditif (ou bilin\'eraire,
c'est-\`a-dire
lorsqu'il pr\'eserve les coproduits \`a gauche et \`a droite), on
sait que le mono\"\i de libre sur un objet $V$ est donn\'e par les
mots en $V$. Le cas g\'en\'eral est plus compliqu\'e et a \'et\'e
tr\`es peu \'etudi\'e. Nous
d\'ecrivons ici la construction du mono\"\i de libre dans le cas
o\`u la cat\'egorie mono\"\i dale pr\'eserve les co\'egalisateurs
r\'eflexifs. Notons que cette hypoth\`ese est assez peu
restrictive. Nous montrons que les foncteurs analytiques scind\'es
pr\'eservent les co\'egalisateurs r\'eflexifs. Comme tous les
produits mono\"\i daux, que nous \'etudions dans cette th\`ese,
induisent des foncteurs analytiques scind\'es, cette propri\'et\'e
 est v\'erifi\'ee par tous nos exemples. On peut donc
leur appliquer la construction propos\'ee ici. \\

Pour d\'emontrer le th\'eor\`eme \'enonc\'e pr\'ec\'edemment, on se
place dans la cat\'egorie des $\Sy$-bimodules. Un $\Sy$-bimodule
est une collection de $(\Sy_m,\, \Sy_n)$-bimodules, o\`u $\Sy_n$
est le groupe sym\'etrique. Les $\Sy$-bimodules servent \`a
repr\'esenter les op\'erations \`a $n$ entr\'ees et $m$ sorties
sur un certain type de \emph{g\`ebre} (alg\`ebre, cog\`ebre,
big\`ebre, etc ...) : $\Po(m,\, n)\otimes A^{\otimes n} \to
A^{\otimes m}$. Dans ce cadre, nous introduisons un analogue au
produit tensoriel $\otimes_k$ des espaces vectoriels sur un corps $k$
et au produit
$\circ$ de composition des op\'erades, que l'on note $\boxtimes$.
Pour deux $\Sy$-bimodules $\Qo$ et $\Po$, le produit $\Qo\boxtimes
\Po$ repr\'esente les compositions d'op\'erations de $\Po$ avec
celles de $\Qo$. Comme ce produit n'a pas d'unit\'e, on ne
consid\`ere que la partie de celui-ci qui s'\'ecrit \`a l'aide de
graphes connexes. Le produit engendr\'e est unitaire et on le
note $\boxtimes_c$.

Un PROP est d\'efini comme une ``alg\`ebre'' pour le produit
$\boxtimes$. Comme toute l'information des PROPs que nous
consid\'erons ici s'\'ecrit \`a l'aide du produit connexe
$\boxtimes_c$, nous d\'efinissons un analogue connexe aux PROPs
que nous appelons les \emph{pro\-p\'e\-ra\-des}. Une prop\'erade
est un mono\"\i de dans la cat\'egorie mono\"\i dale des
$\Sy$-bimodules munie du produit $\boxtimes_c$. D\`es lors, on
travaille au niveau des prop\'erades. Ce choix de pr\'esentation
permet d'obtenir des r\'esultats un peu plus fins et il n'est pas
r\'educteur. Les cat\'egories des PROPs et des prop\'erades sont
reli\'ees par une paire de foncteurs adjoints. A partir d'un PROP,
on d\'efinit une prop\'erade en oubliant les compositions non
connexes. Ce foncteur admet un adjoint \`a gauche $S_\otimes$
bas\'e sur le produit de concat\'enation des $\Sy$-bimodules
$\otimes$ qui est bien connu, notamment du
point de vue homologique. Pour \'etudier le produit $\boxtimes$,
il suffit de faire l'\'etude sur le produit mono\"\i dal
$\boxtimes_c$ et
d'utiliser cet adjoint \`a gauche. Ainsi tous les th\'eor\`emes donn\'es
dans cette th\`ese au niveau des prop\'erades ont un \'equivalent au niveau
des PROPs.\\

La difficult\'e pour g\'en\'eraliser la dualit\'e de Koszul des
alg\`ebres aux op\'erades vient du fait que le produit tensoriel
$\otimes_k$ est bilin\'eaire alors que le produit $\circ$ n'est
lin\'eaire qu'\`a gauche et qu'il s'exprime avec des actions du
groupe sym\'etrique. Le produit $\boxtimes_c$ (et $\boxtimes$)
introduit ici  n'est lin\'eaire ni \`a gauche ni \`a droite et il
s'\'ecrit lui aussi avec des actions du groupe sym\'etrique.
Afin de surmonter cette difficult\'e, nous \'etudions
les propri\'et\'es homologiques du produit $\boxtimes_c$,
en g\'en\'eralisant la d\'emarche conceptuelle de B. Fresse
\cite{Fresse} des op\'erades aux prop\'erades. Pour mener \`a bien cette
\'etude, on ne peut pas se contenter de reproduire simplement les
d\'emonstrations du cas op\'eradique. L'id\'ee principale que nous
ajoutons ici vient du fait que le produit mono\"\i dal
$\boxtimes_c$ induit des foncteurs analytiques. Et c'est
pr\'ecisement les graduations induites par ces foncteurs
analytiques qui nous permettent de d\'ecomposer les diff\'erents complexes
de cha\^\i nes en jeu, et ainsi de conclure les
d\'emonstrations.\\

Nous \'etendons le produit $\boxtimes_c$ (et $\boxtimes$) au cadre
des $\Sy$-bimodules diff\'erentiels gradu\'es par un poids.
Dans le m\^eme chapitre, nous
d\'efinissons les notions de $\Po$-module quasi-libre et de
prop\'erade quasi-libre et nous en \'etudions les propri\'et\'es.
Nous poursuivons avec la g\'en\'eralisation des notions de bar et
de cobar constructions aux prop\'erades (et aux PROPs). Les d\'efinitions
de ces constructions reposent sur des g\'en\'eralisations naturelles
des notions d'\emph{edge contraction} et de \emph{vertex expansion}
donn\'ees par M. Kontsevich dans
la cadre de la (co)homologie des graphes (\emph{cf.} \cite{Ko}).
Nous montrons un premier r\'esultat
significatif :

\begin{thmintro}[Acyclicit\'e de la bar construction augment\'ee]
Pour toute prop\'erade diff\'erentielle $\Po$, le morphisme
d'augmentation
$$\Po \boxtimes_c \bar{\B}(\Po) \to I$$
est un quasi-isomorphisme.
\end{thmintro}

Nous utiliserons ce r\'esultat homologique pour construire des
r\'esolutions au niveau des prop\'erades. Pour cela, on d\'emontre les deux
lemmes de comparaison suivants.

\begin{thmintro}[Lemme de comparaison des modules quasi-libres analytiques]
Dans la cat\'egorie mono\"\i dale des $\Sy$-bimodules
diff\'erentiels gra\-du\'es par un poids, on consid\`ere  $\Psi\,
:\, \Po \to \Po'$ un morphisme de prop\'erades augment\'ees et
$(L,\, \lambda)$ et $(L',\, \lambda')$ deux modules quasi-libres
analytiques sur $\Po$ et $\Po'$ de la forme $L=M\boxtimes \Po$ et
$L'=M'\boxtimes \Po'$. Soit $\Phi \, :\, L \to L'$ un morphisme de
$\Po$-modules analytiques, o\`u la structure de $\Po$-module sur
$L'$ est celle donn\'ee par le foncteur de restriction $\Psi^!$.
On pose $\bar{\Phi} \, :\, M\to
M'$ le morphisme de dg-$\Sy$-bimodules induit par $\Phi$.\\

Si deux des trois morphismes suivants
$\left\{ \begin{array}{l}
\Psi \ : \ \Po \to \Po' \\
\bar{\Phi} \ : \ M \to M' \\
\Phi \ : \ L \to L'
\end{array} \right. $
sont des quasi-isomorphismes, alors le troisi\`eme est aussi un
quasi-i\-so\-mor\-phi\-sme.
\end{thmintro}

\begin{thmintro}[Lemme de comparaison des prop\'erades quasi-libres]
Soient $M$ et $M'$ deux dg-$\Sy$-bimodules gradu\'es par un poids
et de degr\'e sup\'erieur \`a $1$. Soient $\Po$ et $\Po'$ deux
prop\'erades quasi-libres de la forme $\Po=\F(M)$ et
$\Po'=\F(M')$, munies de d\'erivations $d_\theta$ et $d_{\theta'}$
provenant de morphismes $\theta \, :\, M \to \bigoplus_{s\ge
2}\F_{(s)}(M)$ et
 $\theta' \, :\, M' \to \bigoplus_{s\ge 2}\F_{(s)}(M')$ qui pr\'eservent la
graduation totale venant de celle de $M$ et $M'$. Et soit, un
morphisme de dg-$\Sy$-bimodules $\Phi \, :\, \Po \to \Po'$  qui
respecte la graduation analytique de $\F$ et la graduation totale.
Alors, $\Phi$ induit un morphisme $\bar{\Phi} \, :\, M=\F_{(1)}(M)
\to M'=\
F_{(1)}(M')$.\\

Le morphisme $\Phi$ est un quasi-isomorphisme si et seulement si
$\bar{\Phi}$ est un quasi-isomorphisme.
\end{thmintro}

Ces deux lemmes sont des g\'en\'eralisations aux prop\'erades de
lemmes donn\'es par B. Fresse \cite{Fresse} dans le cadre des
op\'erades. Comme les objets li\'es au produit $\circ$ se
d\'ecrivent \`a l'aide des arbres, l'auteur utilise les
propri\'et\'es des arbres pour construire des suites spectrales
dont la convergence permet de d\'emontrer les deux lemmes. N'ayant
pas de telles propri\'et\'es dans le contexte des PROPs,
 nous avons raffin\'e le raisonnement en introduisant une
graduation suppl\'ementaire qui vient du nombre de sommets des
graphes consid\'er\'es dans le cas des prop\'erades quadratiques. Un
autre avantage de cette d\'emarche est qu'elle inclut le cas des alg\`ebres.\\

Le lemme de comparaison des modules quasi-libres analytiques joint \`a
l'acyclicit\'e de la bar construction augment\'ee permet de montrer le
th\'eor\`eme suivant :

\begin{thmintro}[R\'esolution bar-cobar]
Pour toute prop\'erade diff\'erentielle augment\'ee gradu\'ee par
un poids $\Po$, le morphisme naturel d'augmentation
$$\bar{\B}^c(\bar{\B}(\Po))\to \Po$$
est un quasi-isomorphisme.
\end{thmintro}

Enfin, on peut conclure la d\'emonstration du th\'eor\`eme de dualit\'e
de Koszul des prop\'erades avec le lemme de comparaison des prop\'erades
quasi-libres.

\begin{thmintro}[Dualit\'e de Koszul des prop\'erades]
Soit $\Po$ une prop\'erade diff\'e\-ren\-ti\-elle quadratique. Les
propositions suivantes sont \'equivalentes
\begin{enumerate}

\item Le complexe de Koszul $\Po^{\ac}\boxtimes_c \Po$ est
acyclique.

\item Le morphisme de prop\'erades diff\'erentielles gradu\'ees
par un poids $\bar{\B}^c(\Po^{\ac})\to \Po $ est un
quasi-isomorphisme.

\end{enumerate}
\end{thmintro}

Comme le foncteur $S_\otimes$ qui relie les prop\'erades aux PROPs
pr\'eserve l'homologie et comme, dans tous les complexes pr\'esents
ici, la diff\'erentielle agit composante connexe par composante
connexe, on g\'en\'eralise tous ces r\'esultats au cadre des
PROPs. Ceci permet d'achever la d\'emonstration du th\'eor\`eme
de dualit\'e de Koszul pour les PROPs. Enfin, on montre qu'un PROP
est de Koszul si et seulement si la prop\'erade qui lui est
associ\'ee par le foncteur oubli est de Koszul, ce qui justifie, \`a
nouveau, le
fait de travailler au niveau des prop\'erades.\\

Remarquons que les alg\`ebres associatives et les op\'erades sont des exemples de
prop\'erades. Ainsi, les th\'eor\`emes d\'emontr\'es ici s'appliquent \'a ces
deux cas particuliers. On retrouve exactement les r\'esultats de B. Fresse
\cite{Fresse} pour
les op\'erades. Par contre, les d\'emonstrations de B. Fresse n'incluent pas le cas
des alg\`ebres, alors que celles donn\'ees ici les incluent.\\

Pour pouvoir monter qu'une prop\'erade est de Koszul, il reste \`a montrer
que le complexe de Koszul est acyclique. En s'inspirant des travaux de T. Fox
et M. Markl (\emph{cf.} \cite{FM}), on introduit une large classe de prop\'erades
$\Po$ d\'efinies comme un ``m\'elange'' de deux prop\'erades $A$ et $B$
plus simples. Nous d\'emontrons que lorsque ces deux prop\'erades $A$ et $B$ sont
de Koszul, la prop\'erade $\Po$ est de Koszul. Ce r\'esultat permet
d'affirmer que la prop\'erade $\BLi$ des big\`ebres de Lie
(\emph{cf.} V. Drinfeld \cite{Drinfeld}) et la
prop\'erade $\IBi$ des big\`ebres de Hopf infinit\'esimales
(\emph{cf.} M. Aguiar \cite{Ag1} et \cite{Ag2}) sont de Koszul.
(Elles sont construites \`a partir
de deux op\'erades de Koszul \`a chaque fois). En interpr\'etant la
cobar construction sur les duales de telles prop\'erades, on peut calculer
la cohomologie de certains graphes
au sens de M. Kontsevitch \cite{Ko}. Dans les cas $\BLi$ et $\IBi$, on retrouve
les r\'esultats de \cite{MV} sur la cohomologie des graphes classiques
ainsi que des graphes \emph{ribbons}.\\

Dans une derni\`ere partie, nous g\'en\'eralisons les
d\'efinitions des \emph{s\'eries de Poincar\'e} des alg\`ebres et
des op\'erades aux prop\'erades.
Nous \'etablissons une \'equation fonctionnelle qui relie la
s\'erie de Poincar\'e d'une prop\'erade de Koszul \`a celle de sa
duale. Ceci nous permet de g\'en\'eraliser aux op\'erades quadratiques
quelconques les
r\'esultats obtenus par \cite{GK} au niveau des op\'erades binaires. En
appliquant cette \'equation fonctionnelle \`a  une op\'erade libre particuli\`ere,
nous red\'emontrons une formule v\'erif\'ee par la s\'erie g\'en\'eratrice
des polytopes de Stasheff.

\chapter*{Conventions}

\thispagestyle{empty}

\fancyhead{}
\fancyhead[RO]{CONVENTIONS}
\fancyhead[LE]{CONVENTIONS}

On travaille  sur un corps de base $k$ de caract\'eristique nulle
(sauf au chapitre $3$).

\subsection{Les diff\'erentes cat\'egories en jeu} La premi\`ere
cat\'egorie pr\'esente dans cette th\`ese est celle des modules
sur le corps $k$ (espaces vectoriels sur $k$) que l'on note
$k$-Mod. Munie du produit tensoriel $\otimes_k$, elle forme une
cat\'egorie mono\"\i dale.
(Lorsqu'il n'y a pas d'ambigu\"\i t\'e, on note ce produit $\otimes$).\\

Soient $V$, $W$ et $A$ des modules sur $k$ et $f\, :\, V\to W$ un morphisme
de $k$-modules. L'application $f\otimes_k id_A \, : \,
V \otimes_k A \to W\otimes_k A$ sera souvent not\'e $f\otimes_k A$. Et on fera
de m\^eme dans toutes les cat\'egories mono\"\i dales pr\'esentes ici.\\

On travaillera aussi dans la cat\'egorie des $k$-modules
gradu\'es, not\'ee g-Mod. La graduation est ici positive et un
$k$-module gradu\'e $V$ est un module sur $k$ qui admet une
d\'ecomposition de la forme $V=\bigoplus_{n\in \mathbb{N}} V_n$.
On dit qu'un morphisme de modules gradu\'es $f \, :\, V \to W$ est
homog\`ene de degr\'e $d$ si $f(V_n)\subset W_{n+d}$
pour tout $n$ dans $\mathbb{N}$.\\

On consid\`ere la cat\'egorie des $k$-modules diff\'erentiels
gradu\'es, not\'ee dg-Mod. Un $k$-module di\-ff\'e\-ren\-tiel
gradu\'e $V$ est un module sur $k$ qui admet une d\'ecomposition
de la forme $V=\bigoplus_{n\in \mathbb{N}} V_n$ et qui est muni
d'une diff\'erentielle $\delta \, : \, V_n \to V_{n-1}$,
c'est-\`a-dire un morphisme de degr\'e $-1$ qui v\'erifie
$\delta^2=0$. La cat\'egorie g-Mod des modules gradu\'es est une
sous-cat\'egorie pleine de la cat\'egorie dg-Mod des modules
diff\'erentiels gradu\'es. Elle correspond aux modules
diff\'erentiels gradu\'es dont la diff\'erentielle $\delta$ est
nulle. On note $H_*(V)$ l'homologie du complexe de cha\^\i nes
d\'efinie par $V$
 et $|v|=n$
repr\'esente le degr\'e homologique $n$ d'un \'el\'ement $v$ de
$V$. On dit qu'un morphisme de modules diff\'erentiels gradu\'es
$f \, :\, V \to W$ est homog\`ene de degr\'e $d$ si $f(V_n)\subset
W_{n+d}$ pour tout $n$ dans $\mathbb{N}$ et s'il commute avec les
diff\'erentielles respectives. On appelle quasi-isomorphisme
\index{quasi-isomorphisme} tout
morphisme homog\`ene de degr\'e $0$, $f\, : \, V \to W$, qui induit
un isomorphisme en homologie
$H_*(f) \, : \, H_*(V) \to H_*(W)$.\\

Un point crucial dans le pr\'esent travail est l'utilisation d'une graduation
suppl\'ementaire donn\'ee par un poids \index{graduation par un poids}.
Les $k$-modules $V$ qui admettent une
d\'ecomposition en fonction d'un poids
$V=\bigoplus_{\rho \in \mathbb{N}} V^{(\rho)}$ forment une cat\'egorie not\'ee
gr-Mod. Dans le m\^eme esprit, on appelle  module diff\'erentiel gradu\'e
par un poids tout module diff\'erentiel $V$ qui se d\'ecompose en une
somme directe de sous-modules diff\'erentiels not\'es $V^{(\rho)}$. On note
la cat\'egorie associ\'ee gr-dg-Mod. Ces modules
sont bigradu\'es par le degr\'e homologique d'une part
et par un poids d'autre part. On note la graduation homologique par $V_n$ et
celle donn\'ee par le poids par $V^{(\rho)}$ (voire $V_{(\rho)}$). On dit qu'un
foncteur est exact \index{foncteur exact} lorsqu'il pr\'eserve l'homologie.\\

Toutes les inclusions de cat\'egories sont r\'esum\'ees dans le
diagramme
$$\xymatrix{k\textrm{-Mod} \ar[r] \ar[d] &  \textrm{g-Mod} \ar[d] \ar[r]&
\textrm{dg-Mod}
\ar[d] \\
\textrm{gr-Mod} \ar[r]& \textrm{gr-g-Mod} \ar[r]& \textrm{gr-dg-Mod}.}$$

Toutes ces cat\'egories sont munies d'un produit tensoriel. A
partir de $V$ et $W$ deux modules diff\'erentiels, on associe la
produit $V\otimes_k W$ , d\'efini par $(V\otimes_k
W)_n=\bigoplus_{i+j=n} V_i\otimes_k W_j$ et la diff\'erentielle
$\delta$ d'un tenseur \'el\'ementaire homog\`ene $v\otimes_k w$ est
donn\'ee par $\delta(v\otimes_k w)=\delta(v)\otimes_k w + (-1)^{|v|}
v\otimes_k \delta(w)$. On utilise dans ce cadre les r\`egles de
signe de Koszul-Quillen : lorsque l'on doit commuter deux objets
(morphismes, \'el\'ements, etc ...) de degr\'e $d$ et $e$, on
introduit un signe $(-1)^{de}$.

A deux modules (\'eventuellement diff\'erentiels) gradu\'es par un poids $V$ et
$W$, on associe le produit $V\otimes_k W$ donn\'e par la formule analogue
$(V \otimes_k W)^{(\rho)} = \bigoplus_{s+t=\rho}V^{(s)}\otimes_k W^{(t)}$.

Ces produits mono\"\i daux transforment les inclusions pr\'ec\'edentes en
inclusions de cat\'egories mo\-no\-\"\i \-da\-les.

\subsection{$n$-uplets}

\index{$n$-uplet}
Pour simplifier les \'ecritures, un $n$-uplet $(i_1,\, \ldots,\,
i_n)$ sera not\'e $\bar{\imath}$. On aura affaire ici \`a des
$n$-uplets d'entiers strictement positifs. Et on repr\'esente la
quantit\'e $i_1+\cdots+i_n$ par $|\bar{\imath}|$. Lorsqu'il n'y a
pas d'ambigu\"\i t\'e sur les nombres de termes, on note $\bar{1}$
le $n$-uplet $(1,\, \ldots,\, 1)$.

On se sert de la notation $\bar{\imath}$ pour repr\'esenter des produits
d'\'el\'ements indic\'es par le $n$-uplet $(i_1,\, \ldots ,\, i_n)$.
Par exemple, dans le cadre des $k$-modules, $V_{\bar{\imath}}$ correspond
au produit $V_{i_1}\otimes_k \cdots \otimes_k V_{i_n}$.

\subsection{Groupe sym\'etrique}
\index{groupe sym\'etrique}
On note $\Sy_n$ le groupe des permutations de $\{1,\, \ldots,\, n\}$.
On re\-pr\'e\-sen\-te une permutation $\sigma$ de $\Sy_n$ par la
$n$-uplet $(\sigma(1),\, \ldots,\, \sigma(n))$. On prolonge la
remarque pr\'ec\'edente, pour tout $n$-uplet $(i_1,\, \ldots,\,
i_n)$, on note $\Sy_\oi$ le sous-groupe $\Sy_{i_1}\times \cdots
\times \Sy_{i_n}$ de $\Sy_{|\oi|}$. A partir de toute permutation
$\tau$ de $\Sy_n$ et de tout $n$-uplet $\oi=(i_1,\, \ldots,\,
i_n)$, on associe une permutation de $\Sy_{|\oi|}$ dite
permutation par blocs  d\'efinie par
\begin{eqnarray*}
\tau_\oi=\tau_{i_1,\ldots,\, i_n} &=& (i_1+\cdots+i_{\tau^{-1}(1)-1}+1,\, \ldots ,\, i_1+\cdots+i_{\tau^{-1}(1)},\,
\ldots ,\, \\
&& i_1+\cdots+i_{\tau^{-1}(n)-1}+1,\, \ldots ,\,
i_1+\cdots+i_{\tau^{-1}(n)}).
\end{eqnarray*}
\index{permutation par blocs}

\subsection{Suites spectrales associ\'ees \`a un bicomplexe}

Soit $(V,\, \delta_h,\, \delta_v)$ un bicomplexe. A ce bicomplexe, on associe
deux suites spectrales qui convergent vers l'homologie du complexe total
$\delta=\delta_h+\delta_v$.
$$I^r(V)\Rightarrow H_*(V,\,\delta) \quad \textrm{et} \quad II^r(V)\Rightarrow
H_*(V,\,\delta).$$ Plus pr\'ecisement, on a
$$(I^0_{s,\, t},\,
d^0)=(V_{s,\,t},\, \delta_v), \ (I^1_{s,\, t},\,
d^1)=(H_t(V_{s,\,* },\, \delta_v),\, \delta_h) \ \textrm{et} \
I^2_{s,\,t}=H_s(H_t(V_{*,\, *},\, \delta_v),\, \delta_h).$$ Et
pour la seconde suite spectrale, on a
$$(II^0_{s,\, t},\,
d^0)=(V_{s,\,t},\, \delta_h)\ , \ (II^1_{s,\, t},\,
d^1)=(H_s(V_{*,\,t },\, \delta_h),\, \delta_v) \ \textrm{et}
 \ II^2_{s,\,t}=H_t(H_s(V_{*,\, *},\, \delta_h),\, \delta_v).$$

\subsection{G\`ebre}

\index{g\`ebre}
On regroupe sous le terme g\'en\'erique de \emph{g\`ebre} touts
les diff\'erents types d'al\-g\`e\-bres, de cog\`ebres, de big\`ebres,
etc ... Cette terminologie a \'et\'e propos\'ee par Jean-Pierre
Serre (\emph{cf.} \cite{Serre}).

\chapter{Notions Mono\"\i dales}

\thispagestyle{empty}

\fancyhead{}
\fancyhead[RO]{\rightmark}
\fancyhead[LE]{\leftmark}

Le but de ce chapitre est d'abord de fixer les notions que l'on
rencontre dans diff\'erentes cat\'egories mono\"\i dales. La plus
utilis\'ee est probablement celle de mono\"\i de. La notion
de mono\"\i de inclut celles d'anneau, d'alg\`ebre et
d'op\'erade. L'utilisation de la notion de cat\'egorie mono\"\i dale
permet de g\'en\'eraliser les raisonnements effectu\'es dans le cadre des
anneaux et des alg\`ebres. On donne par exemple les d\'efinitions de module
sur un mono\"\i de, de produits relatifs, de mono\"\i de libre et
d'id\'eal d'un mono\"\i de.

Il faut cependant faire attention lorsque l'on g\'en\'eralise ces notions.
Elles viennent toutes d'un cadre tr\`es particulier o\`u le
produit mono\"\i dal est biadditif. Plusieurs
g\'en\'eralisations d'une m\^eme notion sont possibles, et certaines
reposent sur ce que nous appelons
les parties lin\'eaires et multilin\'eaires du produit mono\"\i dal.
La d\'efinition d'id\'eal que nous proposons ici entre dans ce cas de
figure. La g\'en\'eralisation stricto sensu de la notion d'id\'eal
ne permet pas de conserver la propri\'et\'e que le quotient d'un
mono\"\i de par un id\'eal est muni d'une structure naturelle de
mono\"\i de. Pour pallier cette difficult\'e, nous d\'efinissons les
id\'eaux \`a
partir de la notion plus fine de partie multilin\'eaire.

Il en va de m\^eme pour la construction du mono\"\i de libre. Le cas
biadditif est connu depuis longtemps (\emph{cf.} \cite{MacLane1}). Pour
avoir le mono\"\i de libre sur un objet $V$, il
suffit alors de prendre les mots en $V$. Le cas g\'en\'eral a
\'et\'e tr\`es
peu trait\'e. Nous donnons \`a la section $6$, une
construction dans le cas o\`u le produit mono\"\i dal pr\'eserve les
co\'egalisateurs r\'eflexifs. Cette hypoth\`ese est v\'erifi\'ee
par tous les produits
rencontr\'es dans cette th\`ese.

\section{Cat\'egorie mono\"\i dale}

On rappelle ici rapidement les d\'efinitions usuelles des
cat\'egories mono\"\i dales. Pour plus de d\'etails, on renvoie le
lecteur au livre de S. Mac Lane \cite{MacLane1} (chapitre VII). La
d\'efinition de cat\'egorie mono\"\i dale est inspir\'ee par la
cat\'egorie des $k$ modules munie du produit tensoriel sur $k$.

\subsection{Cat\'egorie mono\"\i dale stricte}

\index{cat\'egorie mono\"\i dale stricte}
\begin{dei}[Cat\'egorie mono\"\i dale stricte] On appelle \emph{cat\'egorie mono\"\i dale stricte} toute cat\'egorie $\CA$
munie d'un bifoncteur $\Box : \CA\times \CA \rightarrow \CA$
associatif, c'est-\`a-dire v\'erifiant l'identit\'e
$$ \Box(\Box\times id_{\CA})=\Box(id_{\CA} \times \Box)\ : \ \CA\times \CA\times \CA  \rightarrow \CA, $$
et d'un objet $\II$, unit\'e \`a gauche et \`a droite pour le
produit mono\"\i dal $\Box$, c'est-\`a-dire v\'erifiant
l'identit\'e
$$ \Box(\II \times id_\CA) = \Box(id_\CA\times \II) = id_\CA.$$
On la note $(\CA,\, \Box,\, \II)$.
\end{dei}
\textsc{Exemple :} La cat\'egorie des endofoncteurs d'une
cat\'egorie $\mathcal{C}$ munie de la composition des foncteurs et
du foncteur identit\'e $(\mathcal{C}^\mathcal{C},\, \circ, \,
id_\mathcal{C})$ est une cat\'egorie mono\"\i dale stricte (\`a
cause de la stricte associativit\'e de la composition).

\subsection{Cat\'egorie mono\"\i dale}
\index{cat\'egorie mono\"\i dale}
Comme nous venons de le voir, la d\'efinition pr\'ec\'edente est
trop rigide pour inclure tous les cas que nous aimerions traiter.
Pour pouvoir englober plus de cas, il faut relacher les
hypoth\`eses et consid\'erer l'associativit\'e et les unit\'es \`a
isomorphisme pr\`es.

\begin{dei}[Cat\'egorie mono\"\i dale] Une \emph{cat\'egorie mono\"\i
dale} est une cat\'egorie $\CA$ munie
\begin{itemize}
\item d'un bifoncteur $\Box\, : \, \CA\times\CA \rightarrow \CA$ et
d'une famille d'isomorphismes
$$\alpha_{a,\, b,\, c} \ : \ (a\Box b)\Box c
\cong a\Box (b\Box c)$$
 naturels en $a,\, b$ et $c$, tels que le
pentagone suivant commute pour tout $a,\, b,\, c,\, d$ dans $\CA$ :
$$\xymatrix@R=20pt@C=10pt{
 &(a\Box (b\Box c))\Box d \ar[rr]^-{\alpha_{a,b\Box c,d}}& & a\Box ((b\Box c)\Box d)
 \ar[rd]^-{a \Box \alpha_{b,c,d}}&  \\
((a\Box b)\Box c)\Box d \ar[ru]^-{\alpha_{a,b,c}\Box d}
\ar[rrd]^-{\alpha_{a\Box b,c,d}}
& & & & a\Box ((b\Box (c\Box d)) \\
 & &(a\Box b)\Box (c\Box d) \ar[rru]^-{\alpha_{a,b,c\Box d}}& & }$$

\item et d'un objet $\II$ ainsi que deux isomorphismes $\lambda_a
\, : \, \II \Box a \cong a$ et $\rho_a \, : \, a \Box \II \cong a$
naturels en $a$ tels que le diagramme triangulaire suivant commute
pour tout $a,\, c$ dans $\CA$ :
$$\xymatrix{
(a\Box \II) \Box c \ar[rr]^-{\alpha_{a,\II,c}} \ar[rd]_-{\rho_a
\Box c}& & a\Box (\II \Box c)
\ar[ld]^-{a \Box\lambda_c}\\
 & a \Box c & }$$
et tels que $$\lambda_\II = \rho_\II \ : \ \II \Box \II \rightarrow
\II.$$
\end{itemize}
On la note $(\CA,\, \Box,\, \II,\, \alpha ,\,\lambda,\, \rho)$
voire $(\CA,\, \Box,\, \II)$.
\end{dei}

\begin{dei}[Foncteurs mono\"\i daux]
Tout foncteur entre deux cat\'egories mono\"\i dales qui
pr\'eserve la structure mono\"\i dale est appel\'e \emph{foncteur
mono\"\i dal}.
\end{dei}

\subsection{Exemples}
\begin{enumerate}

\item La cat\'egorie des ensembles munie du produit cart\'esien et
d'un ensemble \`a un \'el\'ement pour unit\'e forme une
cat\'egorie mono\"\i dale not\'ee $(\textrm{Ens},\, \times,\,
\{*\})$.

\item  Sur le m\^eme mod\`ele, la cat\'egorie des espaces
topologiques forme une cat\'egorie mono\"\i dale pour le produit
et un ensemble r\'eduit \`a un point pour unit\'e
$(\textrm{Top},\, \times,\, \{*\})$.

\item Les groupes ab\'eliens avec le produit tensoriel sur
$\mathbb{Z}$ et le groupe $\mathbb{Z}$ lui-m\^eme forment une
cat\'egorie mono\"\i dale not\'ee $(\textrm{Ab},\,
\otimes_\mathbb{Z} ,\, \mathbb{Z})$.

\item Vient ensuite la famille d'exemples form\'ee par les
cat\'egories de $k$-modules. Le plus simple est celui de la
cat\'egorie des modules sur $k$ munie du produit tensoriel
classique $\otimes_k$ avec $k$ pour unit\'e. On la note  $(k\textrm{-Mod},\,
\otimes_k,\, k)$. Puis, en affinant la d\'efinition du produit
tensoriel (\emph{cf.} Conventions), on fournit \`a la cat\'egorie
des $k$-modules gradu\'es et \`a celle des $k$-modules
diff\'erentiels gradu\'es une structure de cat\'egorie mono\"\i
dale not\'ees respectivement $(\textrm{g-Mod},\, \otimes_k,\, k)$ et
$(\textrm{dg-Mod},\, \otimes_k,\, k)$. On peut aussi citer l'exemple
des repr\'esentations vectoriels sur diff\'erentes structures
alg\'ebriques (\emph{cf.} D. Calaque et P. Etingof
\cite{Etingof}).

\item Un exemple plus compliqu\'e, inspir\'e par la th\'eorie des
op\'erades (\emph{cf.} V. Ginzburg et M.M. Kapranov \cite{GK} et
J.-L. Loday \cite{Loday3}) , est donn\'e par la cat\'egorie des $\Sy$-modules.
Un $\Sy$-module est une collection
$(\Po(n))_{n \in \mathbb{N}^*}$ de modules sur $\Sy_n$. \index{$\Sy$-module}
On munit cette cat\'egorie
du produit $\circ$ d\'efini par
$$ \Po\circ\Qo\, (n) = \bigoplus_{1\leqslant k\leqslant n \atop i_1+\cdots +i_k=n}
\Po(k)\otimes_{\Sy_k}\Qo(i_1)\otimes_k \cdots \otimes_k \Qo(i_k) $$ et
de l'unit\'e $\II=(k,\, 0,\, \ldots)$. Elle est
not\'ee $(\Sy\textrm{-Mod},\, \circ,\, \II)$. Remarquons que dans
tous les exemples pr\'ec\'edents le produit mono\"\i dal
pr\'eserve les coproduits \`a gauche comme \`a droite, alors que cet
exemple-ci ne v\'erifie cette propri\'et\'e qu'\`a gauche.
\end{enumerate}

\textsc{Remarque :} Dans la th\'eorie des groupes quantiques, on
se sert de la notion de cat\'egorie tensorielle qui est une
cat\'egorie mono\"\i dale dont les morphismes (Homs) sont munis
d'une structure de $k$-modules de dimension finie. Ici, il nous
suffira de consid\'erer seulement une cat\'egorie mono\"\i dale
ab\'elienne.

\begin{dei}[Cat\'egorie mono\"\i dale sym\'etrique]
\index{cat\'egorie mono\"\i dale sym\'etrique}
Une cat\'egorie mono\"\i dale est dite \emph{sym\'etrique} si elle
poss\`ede des isomorphismes $\tau_{a,b} \, : \, a\Box b\ \to \
b\Box a$ naturels en $a,\, b$ tels que
$$\tau_{a,b}\circ \tau_{b,a}=id_{b\Box a}, \quad \rho_b=\lambda_b\circ\tau_{b,\II}, $$
et tels que le diagramme suivant commute :
$$\xymatrix@C=35pt{(a\Box b) \Box c \ar[r]^-{\alpha_{a,b,c}} \ar[d]^-{\tau_{a,b}\Box c}&
 a\Box (b \Box c)\ar[r]^-{\tau_{a,b\Box c}} & (b\Box c) \Box a\ar[d]^-{\alpha_{b,c,a}} \\
  (b\Box a) \Box c \ar[r]^-{\alpha_{b,a,c}}& b\Box (a \Box c) \ar[r]^-{b\Box\tau_{a,c}} & b\Box (c \Box a).}$$
\end{dei}

Les exemples de $(1)$ \`a $(4)$ sont des exemples de cat\'egories
mono\"\i dales sym\'etriques.

\subsection{Th\'eor\`eme de coh\'erence de Mac Lane}
Gr\^ace au th\'eor\`eme suivant, on peut souvent se passer des
isomorphismes $\alpha$, $\lambda$ et $\rho$ pour ne consid\'erer
que des cat\'egories mono\"\i dales strictes. C'est pourquoi, on
omet souvent en pratique d'\'ecrire ces trois isomorphismes.

\begin{thm}[Th\'eor\`eme de coh\'erence de Mac Lane, \emph{cf.}
\cite{MacLane1}] Toute
cat\'egorie mono\"\i dale est \'equivalente \`a une cat\'egorie
mono\"\i dale stricte. (La relation d'\'equivalence est celle des
cat\'egories mono\"idales).
\end{thm}

\begin{cor}
Tout diagramme construit avec les isomorphismes $\alpha$,
$\lambda$ et $\rho$, les identit\'es et les produits
mono\"\i daux est commutatif.
\end{cor}

\subsection{Cat\'egories mono\"\i dales ab\'eliennes}

\index{cat\'egorie mono\"\i dale ab\'elienne}
Soit $(\CA,\, \Box,\, I)$ une cat\'egorie mono\"\i dale a\-b\'e\-li\-enne.
Dans l'\'etude d'une telle cat\'egorie, un des enjeux fondamentaux
est de comprendre le comportement du produit mono\"\i dal
vis-\`a-vis du coproduit. Pour cela, on introduit les deux
foncteurs suivants :

\begin{dei}[Foncteurs de multiplication]
\index{foncteurs de multiplication}
Dans une cat\'egorie mono\"\i dale $(\CA$, $\Box$, $\II)$, pour
tout objet $A$, on appelle \emph{foncteur de multiplication (ou
de composition) \`a gauche} par $A$ (respectivement \emph{\`a
droite}), le foncteur d\'efini par $L_A \, : \, N \mapsto A\Box N$
(respectivement, $R_A \, : \, N \mapsto N\Box A$).
\end{dei}

\begin{dei}[Cat\'egorie biadditive]
\index{cat\'egorie biadditive}
On appelle \emph{cat\'egorie biadditive} toute ca\-t\'e\-go\-rie
mono\"\i dale ab\'elienne telle que les foncteurs de
multiplication $L_A$ et$ R_A$ soient additifs pour tout objet $A$,
c'est-\`a-dire qu'ils pr\'eservent le coproduit.
\end{dei}

Dans une cat\'egorie biadditive $\CA$, on sait construire certains
objets importants comme le mono\"\i de libre par exemple (\emph{cf.}
section~\ref{monolibre}). Par contre, le cas g\'en\'eral est plus
compliqu\'e et il faut souvent faire la distinction entre deux
types de notions. Pour cela, on d\'efinit l'objet suivant :

\begin{dei}[Partie multilin\'eaire]
\index{partie multilin\'eaire}
\label{partielin} Soient $A,\, B ,\, X$ et $Y$ des objets de
$\CA$. On appelle \emph{partie multilin\'eaire en $X$} le conoyau
de l'application
$$A\Box Y \Box B \xrightarrow{A\Box
i_Y \Box B} A\Box(X\oplus Y)\Box B,$$
que l'on note $A\Box
(\underline{X} \oplus Y)\Box B$.
\end{dei}

\textsc{Exemples :}
\begin{itemize}
\item Dans le cas d'une cat\'egorie mono\"\i dale biadditive, on a
toujours $A\Box (\underline{X} \oplus Y)\Box B=A\Box X \Box B$.
\item Dans le cas des $\Sy$-modules, $A\circ (\underline{X} \oplus
Y)$ correspond aux \'el\'ements de la forme
$A(n)\otimes_{\Sy_n}Z(i_1)\otimes_k \cdots \otimes_k Z(i_n)$ avec
$Z=X$ ou $Y$ mais avec globalement au moins un $X$, d'o\`u la notation.
\end{itemize}

\textsc{Remarque :} La partie multilin\'eaire en $X$ de l'expression
$A\Box(X\oplus Y)$ montre le d\'efaut pour le foncteur
$R_A\, :\, Z\to A\Box Z$ \`a pr\'eserver les conoyaux.
En effet, si ce foncteur pr\'eserve les conoyaux, alors la partie
multilin\'eaire $A\Box(\underline{X}\oplus Y)$  se r\'eduit \`a
$A\Box X$. On voit par exemple, que le foncteur $Z\to A\circ Z$ li\'e
aux $\Sy$-modules ne pr\'eseve en g\'en\'eral pas les conoyaux.

\begin{lem}
\label{firstlemme} Soient $A$ et $B$ deux objets d'une cat\'egorie
ab\'elienne $\CA$. Soient $\pi$ et $i$ deux morphismes
$$\xymatrix{B \ar@<0.5ex>[r]^{i}& \ar@<0.5ex>[l]^{\pi} A} $$ tels que
$i$ soit une section de $\pi$, c'est-\`a-dire $\pi\circ i=id_A$.
Alors le noyau de $\pi$ est naturellement isomorphe au conoyau de
$i$.
\end{lem}

\begin{deo}
Comme $i$ est une section de $\pi$, l'objet $B$ se d\'ecompose
sous la forme $B\simeq \textrm{im}\, A \oplus \ker \pi$. Ce qui
montre que $\coker i\simeq \ker \pi$. $\cqfd$
\end{deo}

Gr\^ace \`a ce lemme, on peut \'ecrire la partie multilin\'eaire
en $X$ comme un sous-objet de $A\Box (X\oplus Y)\Box B$.

\begin{cor}
La partie multilin\'eaire en $X$ correspond aussi au noyau de
l'application
$$A\Box (X\oplus Y)\Box B \xrightarrow{A\Box \pi_Y \Box B} A\Box Y\Box B,$$
o\`u $\pi_Y$ est la projection $X\oplus Y \to Y$.
\end{cor}

Nous aurons besoin plus loin du lemme technique suivant.

\begin{lem}
\label{lemmeconoyau}
 On se place dans une cat\'egorie mono\"\i
dale ab\'elienne $(\CA,\, \Box,\, \II)$. Soit $C=A\oplus B$ un
objet de $\CA$. Alors, le conoyau de $i_A\Box i_A \, : \, A\Box A
\hookrightarrow C\Box C$ est donn\'e par $C\Box (A\oplus
\underline{B})+(A\oplus \underline{B})\Box C$.
\end{lem}

\begin{deo}
On a le diagramme commutatif suivant :
$$\xymatrix{ A\Box(A\oplus \underline{B})\ \ar@{^{(}->}[r] \ar@/^/[d]&
C \Box (A \oplus \underline{B}) \ar[r] \ar@/^/[d]& (C\Box C) /(A\Box A) \\
A\Box C\ \ar@{^{(}->}[r]^-{i_A\Box C} \ar@{>>}[u] & C\Box C\ar@{>>}[u]
\ar@{>>}[ur] \ar@{>>}[r] & \ar[u] (A\oplus \underline{B})\Box C \ar@/^/[l]\\
A\Box A \  \ar[ur]^-{i_A\Box i_A} \ar@{^{(}->}[u]^-{A\Box i_A}
\ar@{^{(}->}[r]^-{i_A\Box A}&  C\Box A  \ar@{^{(}->}[u]_-{C\Box
i_A} \ar@{>>}[r] &(A\oplus \underline{B})\Box A
\ar@{^{(}->}[u]\ar@/^/[l].}
$$
D'o\`u,
\begin{eqnarray*}
C\Box C &=& C\ \Box A \oplus C\Box (A\oplus \underline{B}) \\
&=& A\Box A \oplus \underbrace{(A\oplus \underline{B})\Box A
\oplus C\Box (A\oplus \underline{B} )}_{\coker (i_A\Box i_A).}
\end{eqnarray*}
De la m\^eme mani\`ere, $$ C\Box C = A\Box A \oplus \underbrace{A
\Box (A\oplus \underline{B}) \oplus (A\oplus \underline{B} )\Box
C}_{\coker (i_A\Box i_A)}.$$

Ainsi, $C\Box (A\oplus \underline{B}) + (A\oplus
\underline{B})\Box C \hookrightarrow (C\Box C)/(A\Box A)$ et
l'\'egalit\'e vient de la pr\'ec\'edente combin\'ee \`a $A\Box (A
\oplus \underline{B})\hookrightarrow C\Box (A\oplus
\underline{B})$. $\cqfd$
\end{deo}

\textsc{Remarque :} Le corollaire pr\'ec\'edent montre que le conoyau de
$i_A\Box i_A$ peut \^etre vu comme un sous-objet de $C\Box C$.
Plus pr\'ecisement, le lemme~\ref{firstlemme} montre que le conoyau
de $i_A\Box i_A$ correspond au
noyau de $\pi_A\Box
\pi_A$ (c'est-\`a-dire \`a $C\Box (A\oplus
\underline{B})+(A\oplus \underline{B})\Box C$).

\section{Mono\"\i de}

La notion de mono\"\i de est la g\'en\'eralisation naturelle de celle
d'alg\`ebre.

\subsection{D\'efinition}

\begin{dei}[Mono\"\i de]
\index{mono\"\i de}
Dans une cat\'egorie mono\"\i dale $(\CA,\, \Box,\, \II)$, un \emph{mono\"\i
de} est un objet $M$ muni de deux morphismes :

\begin{itemize}

\item une \emph{composition (\textrm{ou} multiplication)} $\mu \ :
\ M\Box M \rightarrow M$,

\item une \emph{unit\'e} $\eta \ : \ \II
\rightarrow M$
\end{itemize}

tels que les deux diagrammes suivants soient commutatifs
$$
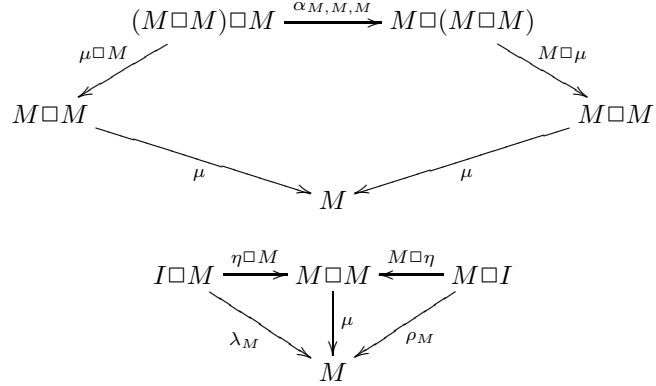
\xymatrix@R=20pt@C=10pt{
& (M\Box M) \Box M \ar[ld]_-{\mu \Box M}  \ar[rr]^-{\alpha_{M,M,M}} & &M \Box (M \Box M) \ar[rd]^-{M \Box \mu} & \\
M \Box M \ar[rrd]_-{\mu}& & & &M \Box M \ar[lld]^-{\mu}\\
& &M & &} $$

$$\xymatrix{
\II \Box M  \ar[r]^-{\eta \Box M} \ar[rd]_-{\lambda_M} &  M\Box M
\ar[d]^-{\mu} & M\Box \II \ar[l]_-{M\Box \eta}
\ar[ld]^-{\rho_M}\\
& M& }$$ On le note souvent $(M,\, \mu ,\, \eta)$.
\end{dei}

\subsection{Exemples}

Dans le cas strict, un mono\"\i de dans la cat\'egorie
$(\mathcal{C}^\mathcal{C},\, \circ, \, id_\mathcal{C}))$ des
endofoncteurs d'une cat\'egorie $\mathcal{C}$ n'est autre qu'une
 monade. \index{monade}

Dans les exemples de cat\'egories mono\"\i dales donn\'es
pr\'ec\'edemment, la notion de mono\"\i de correspond aux
d\'efinitions suivantes :

\begin{enumerate}

\item Dans la cat\'egorie des ensembles, on retrouve la notion
classique de mono\"\i de.

\item Dans le cas des espaces topologiques, on parle de mono\"\i
de topologique.

\item La d\'efinition d'un anneau est exactement celle d'un
mono\"ide dans la cat\'egorie $(\textrm{Ab}$, $\otimes_\mathbb{Z}$
, $\mathbb{Z})$.

\item Pour les cat\'egories construites \`a partir de $k$-modules,
la d\'efinition de mono\"\i de est celle de $k$-alg\`ebre
($k$-alg\`ebre gradu\'ee et $k$-alg\`ebre diff\'erentielle
gradu\'ee).

\item La donn\'ee d'un mono\"\i de dans la cat\'egorie des
$\Sy$-modules est celle d'une op\'erade.\index{op\'erade}
\end{enumerate}

\begin{dei}[Morphismes de mono\"\i des]
Un \emph{morphisme de mono\"\i des} est un morphisme qui commute
avec les compositions et les unit\'es des mono\"\i des source et
but.
\end{dei}
Ainsi, les mono\"\i des de $\CA$ munis de leurs morphismes forment
une
cat\'egorie not\'ee $\textrm{Mon}_{\CA}$.\\

 Dualement, on a la
notion de comono\"\i de.

\begin{dei}[Comono\"\i de]
\index{comono\"\i de}
On appelle \emph{comono\"\i de} $(\Co,\, \Delta,\,\varepsilon)$ de
la cat\'egorie $(\CA,\, \Box,\, \II)$, tout mono\"\i de dans la
cat\'egorie oppos\'ee $(\CA^{op},\, \Box^{op},\, \II)$.
\end{dei}

De mani\`ere e\'quivalente un comono\"\i de $C$ est la donn\'ee
\begin{itemize}
\item d'un morphisme $\Delta \, : \, \ \Co \to \Co \Box
\Co$ appel\'e comultiplication,

\item et d'un morphisme $\varepsilon \, : \,  \Co \to I$ appel\'e
counit\'e.
\end{itemize}

La comultiplication est coassociative, ce qui se repr\'esente (en
omettant les isomorphismes naturels) par le diagramme commutatif

$$\xymatrix{ \Co \ar[r]^{\Delta} \ar[d]^{\Delta}& \Co \Box \Co
\ar[d]^{\Co \Box \Delta}\\
\Co \Box \Co \ar[r]^{\Delta \Box \Co}& \Co \Box \Co \Box \Co.} $$

Et la relation de counit\'e s'\'ecrit
$$\xymatrix{ \Co\Box I \ar[r]^{\rho_\Co} &\Co \ar[d]^{\Delta}  &
\ar[l]_{\lambda_\Co}I \Box \Co \\
 & \ar[ul]^{\Co\Box \varepsilon}
\Co \Box \Co. \ar[ur]_{\varepsilon \Box \Co} &}$$

\subsection{Mono\"\i de augment\'e}

Enfin, dans de nombreux cas, nous aurons affaire \`a des mono\"\i
des munis d'une counit\'e.

\begin{dei}[Mono\"\i de augment\'e]
\index{mono\"\i de augment\'e}
\label{monoaugm}
 Un mono\"ide $(M,\, \mu,\, \eta)$ est dit
\emph{augment\'e} s'il poss\`ede un morphisme de mono\"\i des
$\varepsilon \ : \ M \to \II$. Cela signifie que le diagramme
suivant commute :

$$\xymatrix{ M\Box M \ar[r]^-{\varepsilon \Box \varepsilon} \ar[d]^-{\mu}& \II \Box \II \ar[d]^-{\lambda_\II=\rho_\II} \\
M \ar[r]^-{\varepsilon}& \II \\} $$ et que
$$\xymatrix{\II \ar[r]^-{\eta} & M \ar[r]^-{\varepsilon} & \II}=id_\II.$$
\end{dei}

Si, de plus,  $\CA$ est une cat\'egorie ab\'elienne, on pose
$\overline{M}=ker \ \varepsilon$, que l'on appelle \emph{id\'eal
d'augmentation} de $M$.\index{id\'eal
d'augmentation}

\begin{pro}
Dans une cat\'egorie mono\"\i dale ab\'elienne, tout mono\"\i de
augment\'e est isomorphe \`a $\overline{M}\oplus \II$.
\end{pro}

\begin{deo}
Le morphisme $\eta$ est un rel\`evement de la suite exacte
$$ \xymatrix@C=40pt{ \ar@{^{(}->}[r]^-{ker \ \varepsilon} \overline{M}\ & M
\ar@<1ex>[r]^-{\varepsilon} &
\ar@<1ex>[l]^-{\eta}\II.  }$$ $\cqfd$
\end{deo}

\section{Modules sur un mono\"\i de}

\label{modules} Plusieurs g\'en\'eralisations de la notion de module sur
une alg\`ebre sont
propos\'ees ici. Elles sont toutes \'equivalentes dans le cas
d'une cat\'egorie biadditive.

\subsection{D\'efinition de module}

\index{module}
\begin{dei}[Module sur un mono\"\i de] Un \emph{module \`a gauche} $R$ sur un
mono\"\i de $(M,\, \mu,\,\eta)$ est la donn\'ee d'un objet $R$ de
$\CA$ avec un morphisme $r\, : \, M\Box R \rightarrow R$ tels que
les diagrammes suivants commutent

$$\xymatrix@R=20pt@C=10pt{
& (M\Box M) \Box R \ar[ld]_-{\mu \Box R}  \ar[rr]^-{\alpha_{M,M,R}} & &M \Box (M \Box R) \ar[rd]^-{M \Box r} & \\
M \Box R \ar[rrd]_-{r}& & & &M \Box R \ar[lld]^-{r}\\
& &R & &} $$

$$\xymatrix{
\II \Box R  \ar[r]^-{\eta \Box R} \ar[rd]_-{\lambda_R} &  M\Box R
\ar[d]^-{r} \\
& R. }$$ On note la cat\'egorie des modules \`a gauche sur $M$
 par $M$-Mod.

La d\'efinition de \emph{module \`a droite} $(L, \,
l)$ est sym\'etrique et on note la cat\'egorie associ\'ee Mod-$M$.
Enfin, on appellera \emph{bimodule} \index{bimodule} tout objet $B$ qui est \`a la
fois un module \`a gauche et \`a droite et tels que les deux
actions commutent. On note cette cat\'egorie $M$-biMod.
\end{dei}

\textsc{Exemple :} Un mono\"\i de $(M,\, \mu,\, \eta)$ agit sur
lui-m\^eme par multiplication $\mu$. On parle alors de
repr\'esentation r\'eguli\`ere.

\textsc{Remarque :} Comme les d\'efinitions entre les modules \`a
gauche et les modules \`a droite sont similaires, on se contentera
dans la suite de ne d\'etailler qu'un des deux cas.\\


En dualisant la d\'efinition de module sur un mono\"\i de, on
aboutit \`a celle de comodule sur un comono\"\i de.

\begin{dei}[Comodule sur un comono\"\i de]
\index{comodule}
Soit $(C,\, \Delta,\, \varepsilon)$ un comono\"\i de. Un objet $R$
de $\CA$ muni d'un morphisme $r\, :\, R\to C\Box R$ est appel\'e
\emph{comodule sur $C$ \`a gauche} si $(R, r)$ est un module \`a
gauche dans la cat\'egorie oppos\'ee $(\CA^{op},\, \Box^{op},\,
I)$.
\end{dei}

\subsection{Module libre}

L'oubli de la structure de module d\'efinit un foncteur de la cat\'egorie des
modules (\`a gauche) sur un mono\"\i de $M$ vers la cat\'egorie
$\CA$.
$$U \ : \ (R,\, r) \mapsto R.$$
Ce foncteur admet un adjoint \`a gauche qui est donn\'e dans la
proposition suivante.

\begin{pro}
Pour tout objet $A$ de $\CA$, l'objet $M\Box A$ avec le morphisme
$$\xymatrix{M\Box ( M \Box A) \ar[r]^-{\alpha^{-1}_{M,M,A}}& (M\Box M) \Box A
\ar[r]^-{\mu \Box A}&  M\Box A}$$
forment le module \`a gauche libre sur $A$.
\end{pro}

De la m\^eme mani\`ere, on a :

\begin{pro}
Pour tout objet $A$ de $\CA$, l'objet $C\Box A$ avec le morphisme
$$\xymatrix{C \Box A \ar[r]^-{\Delta \Box A}& (C\Box C) \Box A \ar[r]^-{\alpha_{C,C,A}}
&  C\Box (C \Box A)}$$
forment le comodule \`a gauche colibre sur $A$.
\end{pro}

\subsection{Modules multilin\'eaire et lin\'eaire}
\label{modlin}

Dans le cadre d'une cat\'egorie mono\"\i dale a\-b\'e\-li\-enne, on
peut proposer deux autres g\'en\'eralisations de la notion
classique de module sur un anneau, celles de module
multilin\'eaire et de module lin\'eaire.

\begin{dei}[Module multilin\'eaire]
\index{module  multilin\'eaire}

Un objet $R$ est un \emph{module multilin\'eaire \`a gauche} s'il
admet un morphisme $r\, :\, M\Box (M\oplus \underline{R}) \to R$
v\'erifiant les diagrammes commutatifs suivants :
$$\xymatrix@R=20pt@C=20pt{
& (M\Box M) \Box(M\oplus \underline{R}) \ar[ld]_-{\mu \Box (M\oplus R)}  \ar[rr]^-{\alpha_{M,M,(M\oplus R)}} & &M \Box (M \Box (M \oplus \underline{R}) \oplus M\Box M) \ar[rd]^-{M \Box (r+\mu)} & \\
M \Box (M \oplus \underline{R}) \ar[rrd]_-{r}& & & &M \Box R \ar[lld]^-{r}\\
& &R & &} $$
$$\xymatrix{\II \Box R  \ar[r]^-{\eta \Box R} \ar[rrd]_-{\lambda_R} &  M\Box R
\ar[r] & M \Box (M \oplus \underline{R})
\ar[d]^-{r} \\
& & R. }$$
\end{dei}

On a alors le m\^eme type de propoposition pour le module
multilin\'eaire libre.

\begin{pro}
Pour tout objet $A$ de $\CA$, l'objet $M\Box (M \oplus
\underline{A})$ muni du morphisme d\'efini par la composition

$$M\Box (M \oplus \underline{M \Box (M\oplus
\underline{A})}) \xrightarrow{\alpha^{-1}_{M,M,M\oplus A} \circ
M\Box (M\Box \eta + M\Box(M\oplus A))}  (M\Box M) \Box (M\oplus
\underline{A}) \xrightarrow{\mu \Box (M\oplus A)}  M\Box (M \oplus
\underline{A})$$
 forme le module lin\'eaire \`a gauche libre sur
$A$.
\end{pro}

La terminologie de ``modules multilin\'eaires'' vient du fait que
l'on fait agir $M$ sur des \'el\'ements de $M$ et de $R$ mais avec
au moins un \'el\'ement de $R$. \\

Une autre notion, celle de
\emph{module lin\'eaire}, a \'et\'e introduite dans le cadre de la
cohomologie de Quillen des mono\"\i des par H. J. Baues, M.
Jibladze et A. Tonks dans \cite{BJT}. La d\'efinition de module
lin\'eaire se r\'esume en disant que $M$ agit sur des \'el\'ements
de $M$ et
de $R$ mais avec un seul \'el\'ement de $R$.\\

Pour pouvoir d\'efinir la partie de $M\Box(M\oplus R)$ qui
s'\'ecrit avec un seul \'el\'ement de $R$ \`a droite, on
lin\'earise le foncteur $\mathcal{R}\, :\, X \to M\Box(M\oplus
X)$.

\begin{dei}[Effet crois\'e]
\index{effet crois\'e}

Soit $\Gamma \,  : \, \CA \to \CA$ un endofoncteur de la
cat\'egorie ab\'elienne $\CA$. Pour deux objets $X$ et $Y$ de
$\CA$, on d\'efinit l'\emph{effet crois\'e} $\Gamma(X|Y)$ par le
noyau de l'application
$$\Gamma(X|Y)=\ker \left( \Gamma(X \oplus Y) \xrightarrow{\Gamma(\pi_X)\oplus
\Gamma(\pi_Y)} \Gamma(X)\oplus \Gamma(Y)\right).$$
\end{dei}

L'image de l'effet crois\'e $\Gamma(X|X)$ via l'application
$$\Gamma(X|X) \hookrightarrow \Gamma(X\oplus X) \xrightarrow{\Gamma(+)}
 \Gamma(X),$$
correspond \`a la partie non additive du foncteur $\Gamma$. Il
suffit donc de quotienter $\Gamma(X)$ par cet objet pour obtenir
un foncteur additif.

\begin{pro}
Le foncteur $\Gamma^{\textrm{add}}$ d\'efini par
$$\Gamma^{\textrm{add}}(X)=\coker \left(\Gamma(X|X) \hookrightarrow
\Gamma(X\oplus X) \xrightarrow{\Gamma(+)}
 \Gamma(X) \right)$$
est un foncteur additif.

De plus, la transformation de foncteurs $\textrm{add} \, : \,
\Gamma \to \Gamma^{\textrm{add}}$ factorise de mani\`ere unique
toute transformation naturelle $\Gamma \to \Pi$, o\`u $\Pi$ est un
foncteur additif.
\end{pro}

Lorsque l'on lin\'earise le foncteur $\mathcal{R}\, :\, X \to
M\Box(M\oplus X)$, on obtient la partie lin\'eaire en $X$ de
l'expression $M\Box(M\oplus X)$.

\begin{dei}[Partie lin\'eaire]
\index{partie lin\'eaire}

La \emph{partie multilin\'eaire} en $X$ de l'expression
$M\Box(M\oplus X)$ correspond \`a l'image $\mathcal{R}^{add}(X)$.
\end{dei}

On d\'efinit la notion de module lin\'eaire gr\^ace \`a cet objet
de $\CA$.

\begin{dei}[Module lin\'eaire]
\index{module lin\'eaire}
On appelle \emph{module lin\'eaire \`a
gauche}, tout objet $R$ de $\CA$ muni d'un morphisme $r\, : \,
\mathcal{R}^{\textrm{add}}(R) \to R$ qui v\'erifie le m\^eme type
de diagrammes commutatifs que ceux des deux d\'efinitions de
modules pr\'ec\'edentes.
\end{dei}

\begin{cor}
Lorsque la cat\'egorie $\CA$ est biadditive, les notions de
module, module lin\'eaire et module multilin\'eaire se confondent.
\end{cor}

\begin{deo}
Lorsque la cat\'egorie $\CA$ est biadditive, la partie lin\'eaire
et multilin\'eaire en $R$ de $M\Box(M\oplus R)$ correspond \`a
$M\Box R$. $\cqfd$
\end{deo}

\section{Produits mono\"\i daux relatifs}

La d\'efinition du produit mono\"\i dal relatif est la
g\'en\'eralisation naturelle de la notion de produit tensoriel
relatif sur une alg\`ebre.

A partir de maintenant, nous nous placerons toujours dans une
cat\'egorie mono\"\i dale ab\'elienne $(\CA,\, \Box,\, \II)$. Soit
$(M,\, \mu,\, \eta)$ un mono\"\i de de $\CA$.

\subsection{D\'efinition et premi\`eres propri\'et\'es}

\begin{dei}[Produit mono\"\i dal relatif]
\index{produit mono\"\i dal relatif}
Soient $(L,\, l)$ un
$M$-module \`a droite et $(R,\, r)$ un $M$-module \`a gauche. On
d\'efinit le \emph{produit mono\"\i dal relatif} $L\Box_M R$ par
le conoyau de
$$ \xymatrix@C=35pt{ L\Box M \Box R \ar@<0.5ex>[r]^-{l\Box R}
\ar@<-0.5ex>[r]_-{L \Box r} & L\Box R \ar@{>>}[r]^-{\coker} & L\Box_M R.}$$
\end{dei}

On a les premi\`eres propri\'et\'es suivantes.

\begin{pro}
Pour $R$ un $M$-module \`a gauche, on a $M\Box_M R = R$.

Et pour un $M$-module libre \`a droite $A\Box M$, on a $(A\Box
M)\Box_M R = A\Box R$.
\end{pro}

\begin{deo}
On montre que $r$ correspond au conoyau voulu.
D\'ej\`a, la composition
$$ \xymatrix@C=35pt{ M\Box M \Box R \ar@<0.5ex>[r]^-{\mu \Box R} \ar@<-0.5ex>[r]_-{M \Box r} & M\Box R
\ar@{>>}[r]^-{r}& R}$$ est nulle par d\'efinition de l'action $r$.
Ensuite, consid\'erons $f \ : \ M\Box R \to A$ telle que
$$ \xymatrix@C=35pt{M\Box M \Box R \ar[r]^-{\mu \Box R}_{- M\Box r} & M\Box R \ar[r]^-{f}& A }$$
soit nulle. L'application $(\mu \Box R) \circ (\eta \Box M \Box
R)$ est un isomorphisme de $I\Box M \Box R$ vers $M\Box R$. On
regarde alors la composition
\begin{eqnarray*}
f\circ (\mu\Box R)\circ (\eta \Box M \Box R) &=& f\circ (M\Box
r)\circ (\eta \Box M \Box R) \\
&=& f\circ (\eta \Box R)\circ (\II \Box r).
\end{eqnarray*}
Posons $\bar{f}=f\circ(\eta \Box R)$. Alors, aux isomorphismes
li\'es \`a $\II$ pr\`es, $f$ se factorise en $f=\bar{f}\circ r$.
Comme $r$ est un \'epimorphisme ($r\circ \eta = \lambda$), on en
conclut qu'une telle application $\bar{f}$ est unique. $\cqfd$
\end{deo}

\subsection{Foncteurs de restriction et d'extension}

Soit un morphisme de mo\-no\-\"\i des $\Phi \, : \, M \to M'$. On
construit \`a partir de $\Phi$ deux foncteurs entre les modules
sur $M$ et ceux sur $M'$.

\begin{dei}[Foncteur de restriction]
\index{foncteur de restriction}
 On appelle \emph{foncteur de
restriction} induit par $\Phi$, le foncteur $\Phi^!\ :\
M'\textrm{-Mod} \to M\textrm{-Mod}$ d\'efini par
$$ \xymatrix@C=35pt{M\Box R' \ar[r]^-{\Phi \Box R'}& M'\Box R' \ar[r]^-{r'}& R'.}$$
\end{dei}

Le morphisme $\Phi$ induit sur $M'$ une action \`a droite par $M$.
On peut d\'efinir un foncteur reciproque \`a $\Phi^!$.

\begin{dei}[Foncteur d'extension]
\index{foncteur d'extension}
Le \emph{foncteur d'extension} issu de
$\Phi$ est le foncteur $\Phi_! \ : \ M\textrm{-Mod} \to
M'\textrm{-Mod}$ donn\'e par  le produit de composition relatif
$$ \Phi_!(R)=M'\Box_M R.$$
\end{dei}

\begin{pro}
Les deux foncteurs $\xymatrix{\Phi_! \ : \ M\textrm{-Mod}
\ar@<0.5ex>@^{->}[r] & \ar@<0.5ex>@^{->}[l] M'\textrm{-Mod}\ : \
\Phi^!}$ sont adjoints.
\end{pro}

\begin{deo}
Commen\c cons par d\'efinir les deux transformations suivantes :
\begin{itemize}
\item L'unit\'e d'adjonction $u \ : \ id_{M\textrm{-Mod}}
\Rightarrow \Phi^! \circ \Phi_!$ est donn\'ee par
$$ u_R \ : \ R=M\Box_M R \mapsto M'\Box_M R \ = \Phi\Box_M R .$$
\item Quant \`a la counit\'e $c \ : \ \Phi_! \circ \Phi^!
\Rightarrow id_{M'\textrm{-Mod}} $, elle correspond au passage au
quotient, pour le produit relatif sur $M$, du conoyau
d\'efinissant le produit relatif sur $M'$ de $M'\Box R'$. Ce qui
se r\'esume ainsi :
$$\xymatrix{M'\Box_M \Phi^!(R') \ar[rr]^-{c_{R'}} & & M'\Box_{M'} R'=R'\\
 &M'\Box R' \ar@{>>}[ul]^-{\coker} \ar@{>>}[ur]_-{\coker}& \\
M'\Box M' \Box R' \ar@<0.5ex>[ur]^-{\mu' \Box R'}
\ar@<-0.5ex>[ur]_-{M' \Box r'} & & M'\Box M \Box R'
\ar@<0.5ex>[ul]^-{(\mu' \Box R')  \circ \atop (M'\Box \Phi \Box
R')} \ar@<-0.5ex>[ul]_{(M' \Box r')\circ \atop (M'\Box \Phi \Box
R')} \ar@/_-1pc/[ll]^-{M' \Box \Phi \Box R'}.}$$
\end{itemize}
On v\'erifie ensuite les relations voulues :
\begin{enumerate}
\item La composition
$$ \xymatrix{ \Phi_!(R) \ar[r]^-{\Phi_!(u)} &\Phi_!\circ \Phi^!
\circ \Phi_!(R) \ar[r]^-{c_{\Phi_!(R)}} &\Phi_!(R)}$$
correspond au morphisme identit\'e $id_{\Phi_!(R)}$.

En effet, on a le diagramme commutatif suivant :

$$\xymatrix{ M'\Box_M R \ar[r]^-{\Phi_!(u)} &  M'\Box_M (M' \Box_M R)
\ar[r]^-{c_{\Phi_!(R)}} & M' \Box_M R  \\
M' \Box R = M' \Box (M\Box_M R) \ar@{>>}[u]-^{\coker} \ar[r]
\ar@/_1pc/[rr]_(0.7){\coker}& \ar@{>>}[u]^-{\coker}M' \Box (M'
\Box_M R) \ar@{>>}[r]^-{\coker}&
M' \Box_{M'}(M'\Box_M R)=M' \Box_M R  \ar@{=}[u]\\
M' \Box M \Box R \ar[u]^-{M' \Box \coker} \ar[r]^-{M' \Box \Phi
\Box R} & M' \Box M' \Box R \ar[r]^-{\coker \Box R} \ar[u]^-{M'
\Box \coker} & M' \Box R. \ar@{>>}[u]^-{\coker}}$$

\item D'autre part, on a
\begin{eqnarray*}
\xymatrix{\Phi^!(R') \ar[r]^-{u(\Phi^!)} & \Phi^! \circ \Phi_! \circ
\Phi^!(R') \ar[r]^-{\Phi^!(c_{R'})}& \Phi^!(R')} &=&\\
\xymatrix{R' \ar[r]& M'\Box_M R'\ar[r]& R'}&=& id_{R'}.
\end{eqnarray*}$\cqfd$
\end{enumerate}
\end{deo}

\subsection{Quotient ind\'ecomposable}

Lorsque $(M,\, \mu,\, \II, \, \varepsilon)$ est un mono\"\i de
augment\'e (\emph{cf.} sous-section~\ref{monoaugm}), la ocunit\'e
$\varepsilon$ induit un foncteur $\varepsilon^! \, : \,
\II\textrm{-Mod}=\CA \to M\textrm{-Mod}$ qui fournit une structure
de $M$-module \`a tout objet de $\CA$. On parle alors
d'\emph{action triviale} ou d'\emph{action scalaire} \index{action
scalaire}. En outre, on a un foncteur $\varepsilon_!\, : \,
M\textrm{-Mod} \to \II\textrm{-Mod}=\CA$ tel que
$\varepsilon_!(R)=\II\Box_M R$.

\begin{dei}[Quotient ind\'ecomposable]
\index{quotient ind\'ecomposable} Le produit relatif
$\overline{R}=\II\Box_M R$ est appel\'e \emph{quotient
in\-d\'e\-com\-po\-sa\-ble} de $R$.
\end{dei}

La bijection naturelle d'adjonction s'\'ecrit alors
$$\textrm{Hom}_\CA(\overline{R},\, A)\cong  \textrm{Hom}_{M\textrm{-Mod}}
(R ,\,\varepsilon^!(A)).$$
Et, l'unit\'e d'adjonction devient
$$\xymatrix@C=45pt{R=M\Box_M R \ar@{>>}[r]^-{u_R=\varepsilon\Box_M R}
& I\Box_M R = \overline{R}.}$$ \textsc{Remarque :} Le nom
(quotient ind\'ecomposable) vient du fait que, dans un cadre
ensembliste, il repr\'esente les \'el\'ements du module qui ne
peuvent \^etre obtenus comme image
d'autres \'el\'ements par l'action de $M$.\\

Rappelons que lorsqu'un mono\"\i de est augment\'e, alors on a
$$ \xymatrix{A=A\Box \II \ar[r]^-{A\Box \eta}& A\Box M \ar[r]^-{A \Box \varepsilon}& A}\quad =id_A.$$
Ceci donne la proposition suivante.

\begin{pro}
\label{quotientindecomposable}
 Lorsque le mono\"\i de $M$ est
augment\'e, on peut identifier le quotient
in\-d\'e\-com\-po\-sa\-ble de la repr\'esentation r\'eguli\`ere
avec $\II$ et celui du module libre $M\Box A$ avec $A$ lui-m\^eme.
\end{pro}

\section{Co\'egalisateurs r\'eflexifs}

Nous avons vu que, dans l'\'etude d'une cat\'egorie mono\"\i dale
ab\'elienne, on cherchait \`a comprendre le comportement du
produit mono\"\i dal vis-\`a-vis du coproduit. Plus fort encore,
on peut aussi chercher \`a savoir si le produit mono\"\i dal
pr\'eserve les conoyaux. Dans le cas du produit tensoriel
classique $\otimes_k$, comme les foncteurs de multiplication $R_A$
et $L_A$ sont des adjoints \`a gauche d'autres foncteurs, ils sont
exacts \`a droite. Et comme ils sont additifs, ils pr\'eservent
les conoyaux. Plus g\'en\'eralement, un foncteur pr\'eserve les
conoyaux si et seulement si il pr\'eserve les coproduits et qu'il
est exact \`a droite. On voit donc que, dans une cat\'egorie
mono\"\i dale qui n'est pas biadditive, le produit mono\"\i dal
n'a aucune chance de pr\'eserver les conoyaux. C'est pour cela que
nous avons introduit la notion de partie multilin\'eaire qui sert
\`a mesurer le d\'efaut pour les foncteurs de multiplication \`a
pr\'eserver les conoyaux (\emph{cf.} section $1$).

On sait que dans un cat\'egorie ab\'elienne, la notion de
co\'egalisateur correspond \`a celle de conoyau. Par contre, il
existe une notion plus fine, celle de \emph{co\'egalisateur
r\'eflexif}. Cette notion est plus facilement pr\'eserv\'ee par
les foncteurs, notamment les foncteurs de multiplications. Pour
preuve, le produit mono\"\i dal $\circ$ pr\'eserve les
co\'egalisateurs r\'eflexifs (\emph{cf.} \cite{GH}). Nous verrons
\`a la section $8$ que tout foncteur analytique scind\'e
pr\'eserve les co\'egalisateurs r\'eflexifs et que produits
mono\"\i daux, que nous consid\'erons dans cette th\`ese,
induisent tous des foncteurs analytiques scind\'es.

\subsection{D\'efinition et premi\`eres propri\'et\'es}

\begin{dei}[Paire r\'eflexive]
\index{paire r\'eflexive} Une paire de morphismes $\xymatrix{X_1
\ar@<-0.5ex>[r]_{d_0} \ar@<0.5ex>[r]^{d_1} & X_0}$ est dite
\emph{r\'eflexive} s'il existe un morphisme $s_0\, :\, X_0 \to
X_1$ tel que $d_0 \circ s_0=d_1 \circ s_0=id_{X_0}$.
\end{dei}

\begin{dei}[Co\'egalisateur r\'eflexif]
\index{co\'egalisateur r\'eflexif}
On appelle
\emph{co\'egalisateur r\'eflexif} tout co\'egalisateur provenant
d'une paire r\'eflexive.
\end{dei}

\begin{pro}
\label{reflexifepi} Soit $\Gamma \, : \, \CA \to \CA$ un
endofoncteur d'une cat\'egorie ab\'elienne $\CA$. Si $\Gamma$
pr\'eserve les co\'egalisateurs r\'eflexifs alors $\Gamma$
pr\'eserve les \'epimorphismes.
\end{pro}

\begin{deo}
 Soit $\xymatrix{B \ar[r]^{\pi} & C}$ un \'epimorphisme.
Comme on s'est plac\'e dans une cat\'egorie ab\'elienne, on sait
que $\pi$ correspond au conoyau de son noyau (que l'on note $i$)
$$\xymatrix{{A\ } \ar@{^{(}->}[r]^{i}& {B\ } \ar@{>>}[r]^{\pi}& C.} $$
Le conoyau $\pi$ peut s'\'ecrire comme le co\'egalisateur
r\'eflexif de la paire suivante
$$\xymatrix{A\oplus B\ar@<-1ex>[r]_{d_0}
\ar[r]^{d_1} & \ar@/_1pc/[l]_{s_0}B \ar@{>>}[r]^{\pi}& C,}$$ o\`u
$d_0=i+id_B$, $d_1=id_B$ et $s_0=i_B$.

Comme $\Gamma$ pr\'eserve les co\'egalisateurs r\'eflexifs, on
obtient que $\Gamma(\pi)$ est le co\'egalisateur de
$(\Gamma(d_0)$, $\, \Gamma(d_1))$. En tant que co\'egalisateur,
$\Gamma(\pi)$ est un \'epimorphisme. $\cqfd$
\end{deo}

\subsection{Lien avec le produit mono\"\i dal}

Nous allons \'etudier certaines propri\'et\'es v\'erifi\'ees par
un produit mono\"\i dal lorsque ce dernier pr\'eserve les
co\'egalisateurs r\'eflexifs.

\begin{dei}
On dit d'un produit mono\"\i dal $(\CA, \, \Box,\, I)$ qu'il
\emph{pr\'eserve les co\'egalisateurs r\'eflexifs} si pour tout
objet $A$ de $\CA$, les foncteurs de multiplications $R_A$ et
$L_A$ pr\'eserve les co\'egalisateurs r\'eflexifs.
\end{dei}

\begin{pro}
\label{coegalisateurreflexifproduit}
 Soit $\CA$ est une
cat\'egorie mono\"\i dale ab\'elienne dont le produit pr\'eserve
les co\'egalisateurs r\'eflexifs. Soient
$$\xymatrix{M_1 \ar@<-1ex>[r]_{d_0} \ar[r]^{d_1}
 &\ar@/_1pc/[l]_{s_0} M_0 \ar@{>>}[r]^{\pi_M}&
M} \quad \textrm{et} \quad \xymatrix{N_1 \ar@<-1ex>[r]_{d_0}
\ar[r]^{d_1}  & \ar@/_1pc/[l]_{s_0}N_0 \ar@{>>}[r]^{\pi_N}& N} $$
deux co\'egalisateurs r\'eflexifs. Alors $M\Box N$ est le
co\'egalisateur de
$$\xymatrix@C=40pt{M_1\Box N_1 \ar@<-1ex>[r]_{d_0\Box
d_0} \ar[r]^{d_1\Box d_1}  &\ar@/_1.5pc/[l]_{s_0\Box s_0} M_0\Box
N_0 \ar@{>>}[r]^{\pi_M\Box \pi_N}& M\Box N}.$$
\end{pro}

\begin{deo}
Soit $\phi \, : \, M_0 \Box N_0 \to A$ un morphisme tel que
$\phi(d_0\Box d_0)=\phi(d_1\Box d_1)$.

L'hypoth\`ese que le produit mono\"\i dal $\Box$ pr\'eserve les
co\'egalisateurs r\'eflexifs donne que $M_0 \Box \pi_N$ est le
co\'egalisateur de $(M_0\Box d_0,\ M_0\Box d_1)$. Comme
\begin{eqnarray*}
\phi(M_0 \Box d_0) &=& \phi(d_0\Box d_0) \circ (s_0\Box
N_1) \\
&=& \phi(d_1\Box d_1) \circ (s_0\Box
N_1) \\
&=& \phi(M_0 \Box d_1),
\end{eqnarray*}
on a que $\phi$ se factorise de mani\`ere unique sous la forme
$\phi=\phi_1\circ(M_0 \Box \pi_N)$, o\`u $\phi_1 \, : \, M_0\Box N
\to A$.

On veut montrer que $\phi_1(d_0\Box N)=\phi_1(d_1\Box N)$. Pour
cela, il suffit de montrer que $\phi_1(d_0\Box
\pi_N)=\phi_1(d_1\Box \pi_N)$ parce que $M_1\Box \pi_N$ est un
co\'egalisateur (r\'eflexif) et donc un \'epimorphisme. On a
\begin{eqnarray*}
\phi_1(d_0\Box \pi_N ) &=& \phi(d_0\Box N_0) \\
&=& \phi(d_0\Box d_0) \circ (M_1\Box s_0) \\
&=& \phi(d_1\Box d_1) \circ (M_1\Box s_0) \\
&=& \phi(d_1 \Box N_0) \\
&=& \phi_1(d_1\Box \pi_N).
\end{eqnarray*}

Comme $\pi_M \Box N$ est le co\'egalisateur de $(d_0 \Box N ,\,
d_1 \Box N)$, on peut factoriser de mani\`ere unique $\phi_1$ en
$\phi_1 = \widetilde{\phi} \circ (\pi_M \Box N)$ avec
$\widetilde{\phi} \, : \, M \Box  N \to A$. Ainsi, on a pu
factoriser $\phi$ par $\pi_M\Box \pi_N$, puisque
$\phi=\widetilde{\phi}\circ (\pi_M \Box \pi_N)$.

Soit $\phi=\widetilde{\psi}\circ (\pi_M \Box \pi_N)$ une autre
factorisation , on a alors que $(\widetilde{\psi}
-\widetilde{\phi})\circ (\pi_M \Box \pi_N)=0$. D'apr\`es la
proposition~\ref{reflexifepi}, $\pi_M \Box \pi_N$ apparait comme
une composition de deux \'epimorphismes, il s'agit donc d'un
\'epimorphisme, d'o\`u $\widetilde{\psi}=\widetilde{\phi}$.$\cqfd$
\end{deo}

\section{Mono\"\i de libre}
\label{monolibre} \index{mono\"\i de libre} Dans une cat\'egorie
mono\"\i dale ab\'elienne $(\CA,\, \Box,\, \II)$, pour tout objet
$A$, on peut consid\'erer les deux foncteurs de multiplication
(composition) \`a gauche et \`a droite par $A$ ($L_A \, : \, N
\mapsto A\Box N$ et $R_A \, : \, N \mapsto N\Box A$). Lorsque ces
foncteurs pr\'eservent les coproduits, c'est-\`a-dire qu'ils sont
additifs, alors la description du mono\"\i de libre sur $V$ est
assez simple et est donn\'ee par les mots en $V$ (\emph{cf.}
\cite{MacLane1} VII.3).

Dans le cas o\`u seul un de ces deux foncteurs est additif, la
construction du mono\"\i de libre peut s'\'ecrire \`a l'aide d'une
colimite astucieusement choisie. (\emph{cf.} \cite{BJT}  Appendix
B). L'exemple des op\'erades entre dans ce cadre. Ainsi,
l'op\'erade libre, d\'ej\`a explicit\'ee en terme d'arbres par V.
Ginzburg et M.M. Kapranov dans \cite{GK} correspond \`a cette
colimite.

Par contre, le cas g\'en\'eral a \'et\'e tr\`es peu trait\'e. Seul
E. Dubuc propose une solution dans (\cite{Dubuc}) \`a l'aide d'un
raffinement transfini de la construction usuelle. Nous proposons
ici, une construction diff\'erente et plus concr\`ete qui
s'applique aux exemples que nous \'etudierons par la suite.

\subsection{Construction du mono\"\i de libre}

On se place dans une cat\'egorie mo\-no\-\"\i \-dale ab\'elienne
$(\CA ,\, \Box,\, \II)$ telle que, pour tout objet $A$ de $\CA$,
les foncteurs de multiplications $L_A$ et $R_A$ pr\'eservent les
colimites s\'equentielles ainsi que les co\'egalisateurs
r\'eflexifs.\\

Soit $V$ un objet de $\CA$. On consid\`ere l'objet augment\'e
$V_+=\II\oplus V$, et on pose $\eta \, : \, \II \hookrightarrow
V_+$ l'injection de $\II$ dans $V_+$ et $\varepsilon \, : \, V_+
\twoheadrightarrow \II\,$ la projection de $V_+$ sur $\II$. On
note $V_n=(V_+)^{\Box n}$ (par convention $V_0=(V_+)^{\Box
0}=\II$) et on appelle $\FS(V)$ l'objet d\'efini par
$\bigoplus_{n\geqslant 0}
V_n$.\\

\textsc{Remarque :} Cette objet est muni de d\'eg\'en\'erescences
$$\xymatrix@C=55pt{\eta_i \, :\, V_{n}=(V_+)^{\Box i}\Box \II
\Box (V_+)^{\Box (n-i)} \ar[r]^-{V_{i}\Box \eta \Box V_{n-i}} &
(V_+)^{\Box i}\Box V_+ \Box (V_+)^{\Box (n-i)}=V_{n+1}}.$$ Et,
lorsque $V$ est un mono\"\i de augment\'e, $\FS(\bar{V})$ est muni
de faces pour donner la bar construction simpliciale sur $V$
(\emph{cf.} chapitre $6$).

\begin{dei}[Les cat\'egories simpliciales $\Delta$ et $\Delta_{\rm face}$]
\index{cat\'egories simpliciales $\Delta$ et $\Delta_{\rm face}$}
La \emph{cat\'egorie $\Delta_{\rm face}$} est une sous-cat\'egorie
de la \emph{cat\'egorie simpliciale} $\Delta$. Dans les deux cas,
les objets correspondent aux ensembles finis ordonn\'es
$[n]=\{0<1<\cdots <n\}$, pour $n\in \mathbb{N}$. Les ensembles de
morphismes $Hom_{\Delta}([n],\, [m])$ sont form\'es des
applications croissantes de $[n]$ vers $[m]$. Pour $i=0, \ldots,\,
n$, on d\'efinit les applications \emph{faces} $\varepsilon_i$ par
$$ \varepsilon_i(j)=\left\{ \begin{array}{ll}
j & \textrm{si}\quad j<i, \\
j+1 & \textrm{si}\quad j\ge i.
\end{array}
\right.$$ Les ensembles de morphismes de $\Delta_{\rm face}$ sont
form\'es des seules compositions d'applications faces (et des
identit\'es $id_{[n]}$).
\end{dei}

\textsc{Remarque :} La cat\'egorie $\Delta_{\textrm{face}}$ est
parfois not\'ee
$\Delta^+$ dans la litt\'erature.\\

Dans le cas o\`u le produit mono\"\i dal pr\'eserve les coproduits
\`a gauche comme \`a droite, la colimite de $\FS(V)$ sur la petite
cat\'egorie $\Delta_{\rm face}$ correspond aux mots en $V$ et
fournit ainsi le mono\"\i de libre sur $V$ not\'e $\F(V)$. C'est
le cas dans les exemples num\'erot\'es de (1) \`a (4)
pr\'ec\'edemment (\emph{cf.} section $1$). Ceci explique pourquoi
la construction du mono\"\i de classique libre (cat\'egorie
\textbf{Ens}, exemple (1)), de l'anneau libre (cat\'egorie
\textbf{Ab}, exemple (3)) et de l'alg\`ebre libre (cat\'egorie
\textbf{$k$-Mod}, exemple (4)) se fait selon le m\^eme sch\'ema.

Dans le cas contraire, la colimite $\Colim_{\Delta_{\rm
face}}\FS(V)$ est un objet trop gros pour \^etre un bon candidat
au mono\"\i de libre. Dit autrement, le morphisme de
concat\'enation $V_n\Box V_m \to V_{n+m}$ ne passe pas \`a la
colimite. On quotiente donc les $V_n$ avant de passer \`a
la colimite sur $\Delta_{\rm face}$. \\

Posons $\tau\, :\, V \to V_2$ le morphisme d\'efini par la
composition
$$\xymatrix@C=50pt{ V \ar[r]^(.4){\lambda_V^{-1}\oplus \rho_V^{-1}}& \II\Box
V \oplus V \Box \II \ar[r]^-{i_\II\Box i_V - i_V \Box i_\II }&
(I\oplus V)\Box (I \oplus V)=V_2}.$$

Pour $A, B$ deux objets de $\CA$, $A\Box V_2 \Box B$ s'injecte
dans $A\Box (V_2 \oplus V)\Box B$ via le monomorphisme $A\Box
i_{V_2} \Box B$. En consid\`erant la partie multilin\'eaire en $V$
de l'objet $A\Box (V_2 \oplus V)\Box B$, on a
$$A\Box (V \oplus V_2)\Box B =A\Box V_2 \Box B \, \oplus
 \, A\Box (V_2\oplus \underline{V})\Box B.$$
 On d\'efinit
$$ R_{A,B}=\textrm{im} \left( A\Box (V_2\oplus \underline{V})\Box B
\hookrightarrow
  A\Box (V \oplus V_2)\Box B
\xrightarrow{A\Box (\tau + id_{V_2}) \Box B} A\Box V_2\Box B
\right).$$

\begin{dei}[$\widetilde{V}_n$]
On d\'efinit $\widetilde{V}_n$ par
$$\widetilde{V}_n=\coker\left(\bigoplus_{i=0}^{n-2} R_{V_i,\, V_{n-i-2}}
 \to V_n     \right) .$$
On le note aussi $V_n/{\textstyle (\sum_{i=0}^{n-2}R_{i,n-i-2})}$.
\end{dei}

\textsc{Remarque :} Dans le cas des op\'erades, le $\Sy$-module
$V_n$ correspond aux arbres \`a $n$ niveaux dont les sommets sont
indic\'es par des \'el\'ements de $V$. Alors que $\bigoplus_n
\widetilde{V}_n$ correspond aux arbres sans niveaux. En effet,
quotienter par le $\Sy$-module $\sum_{i=0}^{n-2}R_{i,n-i-2}$
revient \`a identifier un arbre avec pour sous-arbre $\vcenter{
\xymatrix@M=2pt@R=5pt@C=5pt{ \ar@{-}[dr]& &\ar@{-}[dl] \\
\ar@{.}[r]& {\scriptstyle V}\ar@{-}[d]\ar@{.}[r]& \\
 \ar@{.}[r]& { \scriptstyle \II}\ar@{.}[r]&}} $ au m\^eme arbre avec
 le sous-arbre $\vcenter{
\xymatrix@M=2pt@R=5pt@C=5pt{ \ar@{-}[d]& &\ar@{-}[d] \\
{ \scriptstyle \II}\ar@{-}[dr] \ar@{.}[rr]& &{ \scriptstyle \II} \ar@{-}[dl]\\
\ar@{.}[r] & { \scriptstyle V}\ar@{.}[r]&}} $ \`a la place.

\begin{lem}
$\ $
\begin{enumerate}
\item Les morphismes $\eta_i$ entre $V_n$ et $V_{n+1}$ passent au
quotient pour donner des applications $\widetilde{\eta}_i$ de
$\widetilde{V}_n$ vers $\widetilde{V}_{n+1}$. \item Pour tout
couple $i,\, j$, les applications $\widetilde{\eta}_i$ et
$\widetilde{\eta}_j$ sont \'egales.
\end{enumerate}
\end{lem}

\begin{deo}
$\ $
\begin{enumerate}
\item  Il suffit de voir que l'on a
$$\left\{ \begin{array}{ll}
\eta_i(R_{j,\, n-j-2})\subset R_{j,\,n-j-1} & \textrm{si} \quad  j\le i-2, \\
\eta_i(R_{j,\, n-j-2})\subset R_{j+1,\,n-j-2} &
\textrm{si}\quad j\ge i,  \\
\eta_i(R_{i-1,\, n-i-1})\subset R_{i-1,\,n-i}+R_{i,\,n-i-1}&
\textrm{lorsque} \quad i=1,\ldots,\, n-1.
\end{array} \right.$$
\item Comme
$(\widetilde{\eta}_i-\widetilde{\eta}_{i+1})(V_n)\subset R_{i,\,
n-i-1}$, on a $\widetilde{\eta}_i=\widetilde{\eta}_{i+1}$.$\cqfd$
\end{enumerate}
\end{deo}

On pose alors,

\begin{dei}[$\F(V)$]
On d\'efinit l'objet $\F(V)$ par la colimite s\'equentielle
suivante :
$$ \xymatrix@M=10pt{\II \ar[r]^-{\widetilde{\eta}} \ar[rd]_(0.35){j_0} &
\widetilde{V}_1=V_1=V_+ \ar[r]^-{\widetilde{\eta}}
 \ar[d]^(0.45){j_1} &
\widetilde{V}_2  \ar[r]^-{\widetilde{\eta}}\ar[dl]_(0.55){j_2}&
\widetilde{V}_3\ar[r]^-{\widetilde{\eta}}
 \ar[dll]_(0.53){j_3}&
\widetilde{V}_4 \,
\cdots \ar[dlll]^(0.35){j_4} \\
& \F(V)=\Colim_{\mathbb{N}} \widetilde{V}_n .& & &}$$
\end{dei}

Le fait d'avoir quotient\'e les $V_n$ en $\widetilde{V}_n$ a
permis de transformer une colimite sur la cat\'egorie $\Delta_{\rm
face}$ en une colimite s\'equentielle. L'hypoth\`ese que le
produit mono\"\i dal $\Box$ pr\'eserve ce type de colimite donne
ici la propri\'et\'e suivante :

\begin{lem}
\label{lemColim} Pour tout objet $A$ de $\CA$, les foncteurs de
multiplication \`a gauche et \`a droite par $A$, $L_A$ et $R_A$,
pr\'eserve la colimite pr\'ec\'edente $\F(V)$. De mani\`ere
explicite, on a
$$A\Box\Colim_{\mathbb{N}} \widetilde{V}_n \cong \Colim_{\mathbb{N}}
 (A\Box \widetilde{V}_n)\quad \textrm{et} \quad
\Colim_{\mathbb{N}} \widetilde{V}_n\Box A \cong
\Colim_{\mathbb{N}} (\widetilde{V}_n \Box A). $$
\end{lem}

On va maintenant chercher \`a munir l'objet $\F(V)$ d'une
structure de mono\"\i de.

L'unit\'e, not\'ee $\bar{\eta}$ correspond au morphisme $j_0\, :
\, \II \to \F(V)$. Quant au produit, on le d\'efinit \`a partir de
la concat\'enation $V_n\Box V_m \to V_{n+m}$. Posons
$$\mu_{n,\, m} \ : \ \xymatrix{ V_n\Box V_m \ar[r]^-{\sim}& V_{n+m} \ar@{>>}[r]&
 \widetilde{V}_{n+m} \ar[r]^-{j_{n+m}}& \F(V)  }.$$

Posons $R_n=\sum_{i=0}^{n-2} R_{i,\, n-i-2}$, on a alors la suite
exacte courte
$$\xymatrix{0 \ar[r] & {R_n\ } \ar@{^{(}->}[r]^{i_n} & {V_n\ }
\ar@{>>}[r]^{\pi_n} & \widetilde{V}_n \ar[r] & 0.}$$

\begin{pro}
Il existe une unique application $\widetilde{\mu}_{n,\, m}\, : \,
\widetilde{V}_n\Box \widetilde{V}_m \to \F(V)$ qui factorise
$\mu_{n,\, m}$ de la mani\`ere suivante
$$\xymatrix{V_n\Box V_m \ar@{>>}[r]^{\pi_n \Box \pi_m}
\ar[rd]_{\mu_{n,\, m}}&\widetilde{V}_n\Box \widetilde{V}_m
\ar[d]^{\widetilde{\mu}_{n,\, m}} \\
& \F(V).} $$
\end{pro}

\begin{deo}
On \'ecrit les conoyaux $\widetilde{V}_n$ comme des
co\'egalisateurs r\'eflexifs sous la forme
$$\xymatrix{R_n\oplus V_n \ar@<-1ex>[r]_{d_0} \ar[r]^{d_1}
  & \ar@/_1pc/[l]_{s_0} V_n \ar@{>>}[r]^{\pi_n}&
\widetilde{V}_n,}$$ avec $d_0=i_n+id_{V_n}$, $d_1=id_{V_n}$ et
$s_0=i_{V_n}$. L'hypoth\`ese que le produit mono\"\i dal $\Box$
pr\'eserve les co\'egalisateurs r\'eflexifs permet d'affirmer,
gr\^ace \`a la proposition~\ref{coegalisateurreflexifproduit}, que
$\pi_n\Box \pi_m$ est le co\'egalisateur (r\'eflexif) de $(d_0\Box
d_0,\, d_1 \Box d_1)$. La proposition d\'ecoule de la
propri\'et\'e universelle v\'erifi\'ee par les co\'egalisateurs.
Il suffit pour cela de montrer que $\mu_{n,\, m}(d_0\Box d_0)
=\mu_{n,\, m}(d_1\Box d_1)$. Cette \'egalit\'e vient du diagramme
suivant
$$\xymatrix@C=60pt{(R_n\oplus V_n)\Box(R_m \oplus V_m) \ar[d]^{id\Box id}
\ar[r]^(0.6){(i_n+id)\Box(i_m+id)} & V_n\Box V_m \ar[r] & V_{n+m}
\ar[d]^{\pi_{n+m}}\\
V_n\Box V_m  \ar[r] &V_{n+m} \ar[r]^{\pi_{n+m}} &
\widetilde{V}_{n+m},}$$ qui est commutatif en vertu des inclusions
$(\underline{R_n}\oplus V_n)\Box V_m \hookrightarrow R_{n,\, m}$
et $V_n \Box (\underline{R_m}\oplus V_m) \hookrightarrow R_{n,\,
m}$.$\cqfd$
\end{deo}

\begin{lem}
Il existe un unique morphisme $\widetilde{\mu}_{n,\, *}$ rendant
le diagramme suivant commutatif
$$\xymatrix@R=30pt@C=40pt{ \widetilde{V}_n \Box \II \ar[r]^-{\widetilde{V}_n\Box
\widetilde{\eta}}
 \ar[d]_(0.57){\widetilde{\mu}_{n,\, 0}}
 \ar[drr] |\hole  &
  \widetilde{V}_n \Box \widetilde{V}_1
  \ar[r]^-{\widetilde{V}_n\Box \widetilde{\eta}}
\ar[dl]_(0.7){\widetilde{\mu}_{n,\, 1}} |(0.25)\hole \ar[dr]
|(.4)\hole & \widetilde{V}_n \Box \widetilde{V}_2 \ \cdots
\ar[dll]^(0.65){\widetilde{\mu}_{n,\, 2}}
\ar[d]\\
\F(V) & &\ar[ll]^-{\exists !\,  \widetilde{\mu}_{n,\, *}}
\widetilde{V}_n \Box \F(V)=\Colim_{\mathbb{N}}(\widetilde{V}_n\Box
\widetilde{V}_m) .}$$
\end{lem}

\begin{deo}
D\'ej\`a, les applications $\widetilde{\mu}_{n,\, m}$ commutent
avec les  $\widetilde{V}_n\Box \widetilde{\eta}$
$$\xymatrix@C=30pt{ \widetilde{V}_n\Box \widetilde{V}_m
\ar[r]^-{\widetilde{V}_n \Box \widetilde{\eta}}
\ar[d]_-{\widetilde{\mu}_{n,\, m}}
& \widetilde{V}_n\Box \widetilde{V}_{m+1} \ar[dl]^-{\widetilde{\mu}_{n,\, m+1}   }\\
 \F(V). & }$$
Par d\'efinition de la colimite, elles engendrent donc une unique
application
$$\widetilde{\mu}_{n,\, *}\, : \,
\Colim_{\mathbb{N}}(\widetilde{V}_n\Box \widetilde{V}_m) \to
\F(V)$$
 rendant le diagramme de l'\'enonc\'e commutatif. On
conclut en utilisant le lemme \ref{lemColim} pour justifier que $
\Colim_{\mathbb{N}}(\widetilde{V}_n\Box
\widetilde{V}_m)=\widetilde{V}_n \Box \F(V).\cqfd$
\end{deo}

De la m\^eme mani\`ere, on a le lemme suivant.

\begin{lem}
Il existe un unique morphisme $\bar{\mu}$ rendant le diagramme
suivant commutatif
$$\xymatrix@R=30pt@C=40pt{ \II \Box \F(V) \ar[r]^-{\widetilde{\eta} \Box \F(V)}
 \ar[d]_(0.57){\widetilde{\mu}_{0,\, *}}
 \ar[drr] |\hole  &
 \widetilde{V}_1 \Box \F(V)
   \ar[r]^-{\widetilde{\eta} \Box \F(V)}
\ar[dl]_(0.7){\widetilde{\mu}_{1,\, *}} |(0.25)\hole \ar[dr]
|(.4)\hole & \widetilde{V}_2\Box \F(V) \ \cdots
\ar[dll]^(0.65){\widetilde{\mu}_{2,\, *}}
\ar[d]\\
\F(V) & &\ar[ll]^-{\exists !\,  \bar{\mu}} \F(V) \Box
\F(V)=\Colim_{\mathbb{N}}(\widetilde{V}_n\Box \F(V)) .}$$
\end{lem}

\begin{deo}
Les arguments sont les m\^emes. $\cqfd$
\end{deo}

\textsc{Remarque :} La construction de $\bar{\mu}$ en commen\c
cant \`a passer \`a la colimite par la gauche puis par la droite
donnerait le m\^eme morphisme.

\begin{pro}
L'objet $\F(V)$ muni de la multiplication $\bar{\mu}$ et de
l'unit\'e $\bar{\eta}$ forme un mono\"ide dans la cat\'egorie
$(\CA,\, \Box ,\, \II)$.

De plus, ce mono\"\i de est augment\'e et on note
$\overline{\F}(V)$ son id\'eal d'augmentation.
\end{pro}

\begin{deo}
La relation v\'erifi\'ee par l'identit\'e est \'evidente.
L'associativit\'e de $\bar{\mu}$ vient de celle des $\mu_{n,\,
m}$.

On d\'efinit la counit\'e par passage \`a la colimite des
applications
$$ \xymatrix@C=35pt{R_n=\sum_{i=0}^{n-2}R_{i,\, n-2-i} \
\ar@{^{(}->}[r] & V_n=(V \oplus \II)^{\Box n}\ar@{>>}[r]^-{\coker}
\ar[d]_-{\varepsilon^{\Box n}}
 & \widetilde{V}_n \ar@{-->}[dl]^-{\exists !\,
 \widetilde{\varepsilon^{\Box n}}}   \\
& \II^{\Box n}=\II.  &  }$$ qui, apr\`es passage \`a la colimite,
donne la counit\'e $\varepsilon$ voulue. $\cqfd$
\end{deo}

\begin{thm}[Mono\"\i de libre]
\label{monoidelibre} Dans une cat\'egorie mono\"\i dale
ab\'elienne qui admet des colimites s\'equentielles
et telle que le produit mono\"\i dal pr\'eserve ce
type de colimite ainsi que les co\'egalisateurs
r\'eflexifs, le mono\"\i de $(\F(V),\, \bar{\mu},\,
\bar{\eta})$ est libre sur $V$.
\end{thm}

\begin{deo}
L'unit\'e d'adjonction est d\'efinie par
$$u_V\ : \ \xymatrix{V \ \ar@{^{(}->}[r] & V\oplus \II \ar[r]^-{j_1}& \F(V).}$$

Quant \`a la counit\'e $c_M \, :\, \F(M) \to M$, pour un mono\"\i
de $(M,\,  \nu,\, \zeta)$ de $\CA$, on la d\'efinit par passage
\`a la colimite des applications $\widetilde{\nu^n}$ suivantes :
$$ \xymatrix@C=35pt{R_n=\sum_{i=0}^{n-2}R_{i,\, n-2-i}\    \ar@{^{(}->}[r] &
M_n=(M\oplus \II)^{\Box n}\ar@{>>}[r]^-{\coker}
\ar[d]_-{\nu^n\circ (M + \zeta)^{\Box n}}
 & \widetilde{M}_n \ar@{-->}[dl]^-{\exists !\, \widetilde{\nu^n}} \\
& M,  &  }$$ o\`u les morphismes $\nu^n$ repr\'esentent $n-1$
compositions de $\nu$ : $M^{\Box n} \xrightarrow{{\nu^n}} M$. Les
$ \widetilde{\nu^n}$ sont bien d\'efinis parce que $\nu^n\circ (M
+ \zeta)^{\Box n} (R_{i,\, n-2-i})=0$, pour tout $i$.

On a alors imm\'ediatement les deux relations d'adjonction
\begin{eqnarray*}
\xymatrix{ \F(V) \ar[r]^-{\F(u_V)}& \F(\F(V))\ar[r]^-{c_{\F(V)}}
&\F(V)}&= &id_{\F(V)}
\quad \textrm{et} \\
\xymatrix{ M \ar[r]^-{u_M}& \F(M) \ar[r]^-{c_{\F(M)}}& M }& =&
id_M.
\end{eqnarray*} $\cqfd$
\end{deo}

\textsc{Remarque :} Dans le cas des $\Sy$-modules, l'objet $\F(V)$
correspond \`a la somme directe sur les arbres des $\Sy$-modules
obtenus en indi\c cant les sommets des arbres par des \'el\'ements
de $V$. On retrouve l'op\'erade libre donn\'ee par \cite{GK} en
termes d'arbres (sans niveau) ainsi que la construction de
\cite{BJT}.

\subsection{Comono\"\i de colibre}
On a une autre d\'efinition \'equivalente du mono\"\i de libre.
Lorsqu'il existe, le mono\"\i de libre sur $V$ est l'unique objet
(\`a isomorphisme pr\`es) qui v\'erifie la propri\'et\'e suivante
: pour tout morphisme $f\, :\, V \to M$ o\`u $M$ est un mono\"\i
de, il existe un unique morphisme de mono\"\i des $\widetilde{f}
\, : \, \F(V) \to M$ tel que le diagramme suivant soit commutatif
$$\xymatrix{ V \ar[r]^(0.4){\widetilde{\eta}} \ar[dr]_{f}& \F(V)
\ar[d]^{\widetilde{f}} \\ & M . }$$

On peut dualiser cette d\'efinition pour obtenir celle de
comono\"\i de colibre.

\begin{dei}[Comono\"\i de colibre]
\index{comono\"\i de colibre}
 Soit $V$ un objet de $\CA$.
Lorsqu'il existe, le \emph{comono\"\i de colibre} est l'unique
objet $\F^c(V)$ tel que pour tout morphisme $f\, : \, C \to V$,
o\`u $C$ est un comono\"\i de, il existe un unique morphisme de
comono\"\i des $\widetilde{f}\, :\, C \to \F^c(V)$ rendant le
diagramme suivant commutatif
$$\xymatrix{ V  & \ar[l]_{\widetilde{\varepsilon}}\F^c(V) \\ & C
\ar[u]_{\widetilde{f}}. \ar[ul]^{f} }$$
\end{dei}

\section{Id\'eal}

Rappelons qu'un id\'eal d'une alg\`ebre $A$ est un sous-module  de
$J$ de $A$ tel que $A\otimes J \xrightarrow{\mu} J$ et $J \otimes
A \xrightarrow{\mu} J$. Lorsque $J$ est un id\'eal d'une alg\`ebre
$A$, le module quotient $A/J$ est naturellement muni d'une
structure d'alg\`ebre.

Nous avons vu pr\'ec\'edemment (\emph{cf.} section~\ref{modules})
que le cas biadditif ne permettait pas de faire la diff\'erence
entre les diff\'erentes notions de modules. De la m\^eme mani\`ere
lorsque l'on veut g\'en\'eraliser la notion d'id\'eal dans une
cat\'egorie mono\"\i dale quelconque, si on veut conserver la
propri\'et\'e que l'objet quotient est naturellement muni d'une
structure de mono\"\i de, il ne faut pas prendre pour d\'efinition
d'id\'eal la g\'en\'eralisation stricto sensu $A\Box J
\xrightarrow{\mu} J$ et $J \Box A \xrightarrow{\mu} J$. La
d\'efinition que nous proposons ici ne repose pas sur le produit
$A\Box J$ mais sur la partie multilin\'eaire
$A\Box(A\oplus\underline{J})$.

\subsection{D\'efinition et mono\"\i de quotient}

Pour $J\hookrightarrow M$ un sous-objet de $M$ dans $\CA$, on note
$M/J$ le conoyau (quotient) de $M$ par $J£$, soit
$$\xymatrix@C=40pt{J\ \ar@{^{(}->}[r]^-{i}& M \ar@{>>}[r]^-{\pi =\coker i}& M/J.} $$

\begin{dei}[Id\'eal]
\index{id\'eal}
 Un sous-objet $J$ d'un mono\"\i de $(M,\,
\mu,\, \eta)$ est appel\'e \emph{id\'eal} de $M$ si la composition
$\pi \circ \mu \circ \ker (\pi\Box \pi)$ est nulle
$$\xymatrix@C=40pt{K_J\  \ar@{^{(}->}[r]^-{\ker \pi\Box\pi}&
M\Box M \ar[r]^-{\mu}& M \ar@{>>}[r]^-{\pi} & M/J}.$$
\end{dei}

La d\'efinition d'id\'eal est faite pour avoir la proposition
suivante.
\begin{pro}
Dans une cat\'egorie mono\"\i dale ab\'elienne $(\CA,\, \Box,\,
\II)$ telle que le produit mono\"\i dal pr\'eserve les
\'epimorphismes, le quotient $M/J$ est muni d'une structure
naturelle de mono\"\i de.
\end{pro}

\begin{deo}
D'apr\`es la condition de l'\'enonc\'e sur les
\'epi\-mor\-phi\-smes, $\pi \Box \pi$ est un \'epi\-mor\-phi\-sme.
Comme $\CA$ est une cat\'egorie ab\'elienne, $\pi\Box \pi = \coker
(\ker (\pi \Box \pi))$. Par d\'efinition du conoyau, il existe un
unique morphisme $\overline{\mu}$ rendant le diagramme suivant
commutatif
$$\xymatrix@M=5pt{ M/J\Box M/J \ar@{-->}[r]^-{\bar{\mu}} & M/J \\
 M\Box M \ar[r]^-{\mu} \ar@{>>}[u]^-{\pi \Box \pi}& M \ar@{>>}[u]^-{\pi} \\
 K_J \ar@{^{(}->}[u]^-{\ker(\pi \Box \pi)} & J.  \ar@{^{(}->}[u]^-{i}}$$
On d\'efinit l'unit\'e par $\bar{\eta}=\pi\circ \eta \, : \,
\xymatrix{\II \ar[r]^-{\eta} & M \ar[r]^-{\pi}& M/J}.$ On peut
montrer l'associativit\'e de $\bar{\mu}$ \`a partir de celle de
$\mu$ :
$$\xymatrix@C=45pt{ M\Box M\Box M \ar[rrr]^-{M\Box \mu} \ar@{>>}[dr]_-{\pi \Box \pi \Box
\pi}
\ar[ddd]_-{\mu \Box M}& & &M \Box M \ar[dl]^-{\pi \Box \pi} \ar[ddd]^-{\mu}\\
& M/J \Box M/J \Box  M/J \ar[r]^{M/J\Box \bar{\mu}}
\ar[d]^-{\bar{\mu} \Box M/J}
& M/J \Box M/J \ar[d]^-{\bar{\mu}}& \\
 & M/J \Box M/J  \ar[r]^-{\bar{\mu}}& M/J & \\
M\Box M \ar[ru]^-{\pi \Box \pi} \ar[rrr]^-{\mu}& & & M,
\ar[lu]_-{\pi} }$$ parce que $\pi \Box \pi \Box \pi$ est un
\'epimorphisme. On proc\`ede de la m\^eme mani\`ere pour montrer
la relation v\'erifi\'ee par l'unit\'e. $\cqfd$
\end{deo}

\textsc{Remarque :} Gra\^ce \`a la proposition~\ref{reflexifepi},
on a que cette proposition est vraie dans toute cat\'egorie
mono\"\i dale ab\'elienne qui pr\'eserve les conoyaux r\'eflexifs.\\

On justifie la terminologie utilis\'ee pour les mono\"\i des
augment\'es par la proposition suivante.

\begin{pro}
Soit $(M,\, \mu,\, \eta,\, \varepsilon)$ un mono\"\i de
augment\'e. Alors l'id\'eal d'augmentation $\overline{M}$ est bien
un id\'eal au sens pr\'ec\'edent.
\end{pro}

\begin{deo}
Par d\'efinition de $\varepsilon$ morphisme de mono\"\i des, on a
$ \varepsilon \circ \mu \circ \ker(\varepsilon \Box \varepsilon) =
\mu_\II \circ \varepsilon \Box \varepsilon \circ \ker(\varepsilon
\Box \varepsilon) =0$. $\cqfd$
\end{deo}

Maintenant, le probl\`eme est de savoir \`a quoi ressemble le
noyau $K_J=\ker (\pi\Box \pi)$ pour pouvoir bien comprendre
l'hypoth\`ese \`a v\'erifier.

\begin{pro}
Dans une cat\'egorie mono\"\i dale ab\'elienne $(\CA, \, \Box,\,
\II)$, soit $J$ un sous-objet de $M$ tel que le conoyau $M/J$
poss\`ede une section. Alors, le noyau $\ker (\pi \Box \pi)$
correspond \`a l'image de $M\Box (M\oplus \underline{J}) \oplus (M
\oplus \underline{J}) \Box M $ dans $M\Box M$ via l'application
$M\Box(M + i_J) + (M + i_J)\Box M$ que nous noterons $K_J=M\Box
(M\oplus \underline{J}) + (M \oplus \underline{J}) \Box M $
\end{pro}

\begin{deo}
Cette proposition est une cons\'equence directe de la remarque qui
suit le lemme~\ref{lemmeconoyau}.$\cqfd$
\end{deo}

Tout ceci permet de donner une autre d\'efinition, \'equivalente
dans le cas scind\'e, de la notion d'id\'eal.

\begin{cor}
\label{caracterisationdideal}
 Un sous-objet $\xymatrix{J\
\ar@{^{(}->}[r]^-{i}& M}$ d'un mono\"\i de $M$, tel qu'il existe
un objet $N$ v\'erifiant $J\oplus N=M$, est un id\'eal de $M$ si
et seulement si
$$\left\{ \begin{array}{l}
\xymatrix@C=40pt{ M\Box (M\oplus \underline{J}) \ar[r]^-{M\Box (M+i)}
&M\Box M \ar[r]^-{\mu}& J} \\
\xymatrix@C=40pt{ (M\oplus \underline{J})\Box M \ar[r]^-{(M+i)\Box
M}&M\Box M \ar[r]^-{\mu}& J.}
\end{array} \right.$$
\end{cor}

\textsc{Remarque :} Il est \'equivalent de dire que $J$ est un
bimodule multilin\'eaire sur $M$ pour la re\-pr\'e\-sen\-ta\-tion
r\'eguli\`ere.\\

Ceci ressemble plus \`a la d\'efinition classique d'un id\'eal. En
effet, dans le cas d'un produit mono\"\i dal biadditif, cela
revient \`a exiger que
$$\left\{ \begin{array}{l}
\xymatrix{M\Box J \ar[r]^-{\mu} & J  } \\
\xymatrix{J\Box M \ar[r]^-{\mu}& J.  }
\end{array}\right. $$
On reconnait bien la notion d'id\'eal pour un anneau ou une
alg\`ebre.\\

Dans le cas des op\'erades (cat\'egorie des $\Sy$-modules munie du
produit $\circ$), la d\'efinition devient
$$\left\{ \begin{array}{l}
\xymatrix{M\otimes_{\Sy_n} \underbrace{M(i_1)\otimes \cdots
\otimes J(i_k) \otimes
\cdots \otimes M(i_n)}_{\textrm{au moins un J}}  \ar[r]^-{\mu} & J  } \\
\xymatrix{J\circ M \ar[r]^-{\mu}& J.  }
\end{array}\right. $$
On retrouve la d\'efinition donn\'ee dans \cite{GK} et par M.
Markl dans \cite{Markl1}.\\

\subsection{Id\'eal engendr\'e}
Supposons maintenant que la cat\'egorie $\CA$ soit petite et
compl\`ete, pour pouvoir d\'efinir l'intersection d'un certain
ensemble d'objets (\emph{cf.} \cite{MacLane1}).

\begin{dei}[Id\'eal engendr\'e]
\index{id\'eal engendr\'e}
 Soit $R$ un sous-objet d'un mono\"\i de $M$. On appelle
\emph{id\'eal engendr\'e par $R$}, le plus petit id\'eal de $M$
contenant $R$. Ce dernier existe et est donn\'e par l'intersection
$\bigcap_{J} J$ pour $J$ id\'eal de $M$ contenant $R$ ($R
\hookrightarrow J \hookrightarrow M$). On le note $(R)_M$ voire
$(R)$ lorsqu'il n'y a pas d'amibigu\"\i t\'e.
\end{dei}

Pour $R$ un sous-objet de $M$, on consid\`ere la partie
multilin\'eaire en $R$ de $M\Box (M\oplus R)\Box M$ c'est-\`a-dire
$$M\Box (M\oplus \underline{R})\Box M = \coker (M \
\Box M\Box M \to M\Box (M\oplus R) \Box M).$$

Le description de l'id\'eal engendr\'e sur $R$ \`a l'aide d'une
intersection n'\'etant pas tr\`es explicite, on en donne une autre
forme. Comme la d\'efinition d'id\'eal repose sur la notion de
partie multilin\'eaire, l'id\'eal libre se construit aussi avec
cette notion.

\begin{pro}
Soit $R$ un sous-objet d'un mono\"\i de  $M$. Alors l'id\'eal
engendr\'e par $R$ correspond \`a l'image du morphisme
$$ \xymatrix@C=80pt{ M\Box (M\oplus \underline{R}) \Box M
\ar[r]^-{\mu^2\circ (M\Box (M+i) \Box M)} & M}$$ que l'on notera $(R)$.
\end{pro}

\begin{deo}
A l'aide de la deuxi\`eme caract\'erisation d'un id\'eal
(corollaire~\ref{caracterisationdideal}), on voit que $(R)$ est un
id\'eal de $M$. Puis, pour tout id\'eal $J$ de $M$ contenant $R$,
on a que
$$ \xymatrix{M\Box (M\oplus \underline{R})\Box M \ar[r]^-{\mu^2}& J.}$$
Donc, (R) est inclus dans tout $J$ et ainsi $(R)=\bigcap_J J$.
$\cqfd$
\end{deo}

\textsc{Remarque :} L'id\'eal libre $(R)$ correspond au bimodule
multilin\'eaire libre engendr\'e par $R$.\\

On retrouve le cas des anneaux et des alg\`ebres
$$ (R) = \mu^2(A\otimes R \otimes A),$$
ainsi que celui des op\'erades. Par exemple, V. Ginburg et M. M.
Kapranov dans \cite{GK} d\'ecrivent l'id\'eal engendr\'e par un
$\Sy$-module $R$ avec les arbres dont les sommets sont indic\'es
par des \'el\'ements de $M$, en imposant qu'au moins un sommet
soit indic\'e par un \'el\'ement de $R$.

\section{Foncteurs polynomiaux et analytiques}

La notion de foncteur analytique, qui est une colimite de foncteurs
polynomiaux, est essentielle dans le reste de ce travail.

Dans la suite de cette th\`ese, nous introduirons un nouveau
produit mono\"\i dal $\boxtimes_c$ que nous \'e\-tu\-di\-erons en
d\'etails. Une propri\'et\'e fondamentale de ce produit est que
les foncteurs de multiplication induits sont des foncteurs
analytiques scind\'es. Comme les foncteurs analytiques scind\'es
pr\'eservent les co\'egalisateurs r\'eflexifs, on pourra appliquer
la construction du mono\"\i de libre \`a cette cat\'egorie.

En outre, le fait de reconnaitre sur les produits $A\boxtimes_c B$
une structure analytique permet d'introduire une graduation
suppl\'ementaire (voire une bigraduation) sur de tels objets.
C'est cette id\'ee qui nous permettra de d\'emontrer les lemmes
homologiques
sur lesquels repose toute cette th\`ese.\\

On se place dans une cat\'egorie mono\"\i dale sym\'etrique
ab\'elienne $(\CA,\, \otimes,\, k)$. Soit $\Delta_n$ le foncteur
diagonal $\CA \to \CA^{\times n}$.

\begin{dei}[Foncteurs polynomiaux homog\`enes]
\index{foncteurs polynomiaux homog\`enes}
On appelle
\emph{foncteur polynomial homog\`ene de degr\'e $n$} tout foncteur
$f \, :\, \CA \to \CA$ qui s'\'ecrit sous la forme
$f_{(n)}=f_n\circ \Delta_n$ avec $f_n$ un foncteur de $\CA^{\times
n} \to \CA$ additif en chacune de ses entr\'ees.
\end{dei}

\begin{dei}[Foncteurs polynomiaux scind\'es]
\index{foncteurs polynomiaux scind\'es} Un foncteur $f\,  : \, \CA
\to \CA$ est dit \emph{polynomial scind\'e} s'il se d\'ecompose en
somme directe de foncteurs polynomiaux homog\`enes
$f=\bigoplus_{n=0}^N f_{(n)}$.
\end{dei}

Les foncteurs que nous rencontrerons par la suite ne s'expriment
pas tous \`a l'aide de sommes finies.

\begin{dei}[Foncteurs analytiques scind\'es]
\index{foncteurs analytiques scind\'es}
 On appelle \emph{foncteur
analytique scind\'e}, tout foncteur $f\, : \, \CA \to \CA$ qui
s'\'ecrit sous la forme $f=\bigoplus_{n=0}^\infty f_{(n)}$ o\`u
$f_{(n)}$ est un foncteur polynomial homog\`ene de degr\'e $n$.
\end{dei}

\textsc{Exemple :} Le foncteur de Schur $\mathcal{S}_\Po$
associ\'e \`a une op\'erade $\Po$
$$\mathcal{S}_\Po (V) =\bigoplus_{n=0}^\infty \Po(n)\otimes_{\Sy_n}
V^{\otimes n}$$
est un foncteur analytique scind\'e.\\

Dans la suite, nous utiliserons la graduation naturelle fourni par
de tels foncteurs. Nous la noterons toujours entre parenth\`eses
$(n)$. Et, par abus de langage, on utilisera dans la suite le
terme de foncteur analytique pour parler de foncteurs analytiques
scind\'es.

\begin{pro}
\label{foncteuranalytiquereflexif} Tout foncteur analytique
scind\'e pr\'eserve les co\'egalisateurs r\'eflexifs.
\end{pro}

\begin{deo}
Soient $\xymatrix{X_1 \ar@<-1ex>[r]_{d_0} \ar[r]^{d_1} &
\ar@/_1pc/[l]_{s_0}X_0 \ar@{>>}[r]^{\pi}& X}$ un co\'egalisateur
r\'eflexif et $f=\bigoplus_{n=0}^\infty f_{(n)}$ un foncteur
analytique scind\'e. Posons $f_{(n)}=f_n\circ \Delta_n$ o\`u
$f_n\, :\, \CA^{\times n} \to \CA$ est un foncteur $n$-additif. Le
r\'esultat vient de l'\'egalit\'e
$$\sum_{i=1}^n f_n(X_0,\ldots,\,
\underbrace{(d_0-d_1)(X_1)}_{i^{\textrm{\`eme}}\ \textrm{place}},
\,\ldots ,\ X_0)=\left( f_n(d_0,\ldots,\ d_0)-f_n(d_1,\ldots,\,
d_1)\right)\circ \Delta_n(X_1).$$ L'inclusion $\supset$ est
toujours vraie et vient de la formule
$$f_n(d_0,\ldots,\
d_0)-f_n(d_1,\ldots,\, d_1)=\sum_{i=1}^n f_n(d_0,\ldots ,\, d_0
,\, \underbrace{d_0-d_1}_{i^{\textrm{\`eme}}\ \textrm{place}},\,
d_1,\ldots ,\,d_1).$$

L'inclusion inverse $\subset$ repose sur le rel\`evement $s_0$ et
vient de
\begin{eqnarray*}
&& f_n(X_0,\ldots,\,X_0,\,  (d_0-d_1)(X_1),\, X_0 ,\,\ldots ,\
X_0)\\
&=&  f_n(X_0,\ldots,\,  d_0(X_1) ,\,\ldots ,\
X_0)-f_n(X_0,\ldots,\,
d_1(X_1) ,\,\ldots ,\ X_0)\\
&=& f_n(d_0s_0(X_0),\ldots,\,  d_0(X_1) ,\,\ldots ,\
d_0s_0(X_0))-f_n(d_1s_0(X_0),\ldots,\, d_1(X_1) ,\,\ldots ,\
d_1s_0(X_0)).
\end{eqnarray*} $\cqfd$
\end{deo}

\chapter{Prop\'erades et PROPs}

\thispagestyle{empty}

On poursuit ici la m\^eme d\'emarche qui a men\'e \`a
l'introduction des op\'erades. Les op\'erades ont \'et\'e
d\'efinies pour mod\'eliser les op\'erations \`a $n$ entr\'ees et
une sortie $\vcenter{\xymatrix@M=0pt@R=6pt@C=6pt{\ar@{-}[dr] &
\ar@{-}[d] &\ar@{-}[dl]  \\
&\ar@{-}[d] & \\  & &}}$ sur les diff\'erents types d'alg\`ebres.
Pour re\-pr\'e\-sen\-ter alg\'ebriquement l'ensemble des
op\'erations agissant sur un type alg\`ebres $A$, on utilise des
$\Sy_n$-modules : $\Po(n)\otimes_{\Sy_n} A^{\otimes n} \to A$.

Dans certains cas, comme ceux des big\`ebres et des big\`ebres de
Lie, on veut pouvoir repr\'esenter des op\'erations \`a plusieurs
entr\'ees et plusieurs sorties
$\vcenter{\xymatrix@M=0pt@R=6pt@C=6pt{\ar@{-}[dr]
& \ar@{-}[d] &\ar@{-}[dl]  \\
&\ar@{-}[dl] \ar@{-}[dr] & \\  & &}}$ agissant sur un module $A$.
Pour cela, on introduit ici la notion de $(\Sy_m,\,
\Sy_n)$-bimodule : $\Po(m,\, n)\otimes_{\Sy_n} A^{\otimes n} \to
A^{\otimes m}$. On d\'efinit ensuite un produit $\boxtimes$ dans
la cat\'egorie des $\Sy$-bimodules qui repr\'esente
alg\'ebriquement les compositions d'op\'erations. Ce produit peut
aussi s'\'ecrire \`a l'aide de graphes dirig\'es.

En partant de l'observation que les diff\'erents types de
g\`ebres, que l'on consid\`ere en pratique, sont d\'efinies par
des g\'en\'erateurs et des relations bas\'es sur des graphes
connexes, il suffit de prendre en compte les compositions
\'ecrites \`a l'aide de graphes connexes pour obtenir toute
l'information escompt\'ee. On d\'efinit ainsi le produit
$\boxtimes_c$ en se restreignant aux graphes connexes. A la
diff\'erence du produit $\boxtimes$, le produit $\boxtimes_c$ est
un produit mono\"\i dal dans la cat\'egorie des $\Sy$-bimodules.
R\'eciproquement, on retrouve le produit $\boxtimes$ \`a partir du
produit $\boxtimes_c$ par concat\'enation.

Nous d\'efinissons une \emph{prop\'erade} comme un mono\"\i de
 dans la cat\'egorie mono\"\i dale des
$\Sy$-bi\-mo\-du\-les munie du produit connexe $\boxtimes_c$. Un
PROP correspond \`a un ``mono\"\i de'' pour le produit $\boxtimes$
avec en plus un morphisme de concat\'enation des op\'erations.
Nous montrons que ces deux notions sont reli\'ees par une paire de
foncteurs adjoints.

\section{La cat\'egorie des $\Sy$-bimodules}

Afin de repr\'esenter les op\'erations \`a $n$ entr\'ees et $m$
sorties, on introduit une cat\'egorie dont les objets sont des
$\mathbb{S}$-bimodules.

\subsection{$\Sy$-bimodules}

\begin{dei}[$\Sy$-bimodule]
\index{$\Sy$-bimodule}
 On appelle $\Sy$\emph{-bimodule}, une
collection de $\left( \Po(m,\, n)\right)_{m,\, n \in \mathbb{N}}$,
o\`u chaque $k$-module $\Po(m,\, n)$ est muni d'une action de
$\Sy_m$ \`a gauche et d'une action de $\Sy_n$ \`a droite, telles
que ces deux actions commutent entre elles.

Un morphisme entre deux $\Sy$-bimodules $\Po$, $\Qo$ est une
collection d'applications lin\'eaires $f_{m,\, n} \, :\, \Po(m,\,
n) \to \Qo(m,\, n)$ \'equivariantes \`a gauche par $\mathbb{S}_m$
et \`a droite par $\mathbb{S}_n$.

Les $\Sy$-bimodules et leurs morphismes forment une cat\'egorie
que l'on note $\Sy$-biMod.
\end{dei}

\begin{dei}[$\Sy$-bimodule r\'eduit]
\index{$\Sy$-bimodule r\'eduit}
Lorsqu'un $\Sy$-bimodule $\Po$
v\'erifie $\Po(0,\, n)=0$ et $\Po(m,\, 0)=0$ pour tous les entiers
$n$ et $m$, on dit qu'il est \emph{r\'eduit}.
\end{dei}

Les $\Sy$-bimodules servent \`a coder les op\'erations sur un
certain type de g\`ebre. Il faut maintenant expliquer comment on
repr\'esente alg\'ebriquement les compositions entre ces
op\'erations.

\subsection{Permutations connexes}

A la diff\'erence du cas g\'en\'eral des PROPs (\emph{cf.} section
$3$), on ne consid\`ere ici que les compositions d'op\'erations
bas\'ees sur des graphes connexes.

On cherche \`a \'ecrire ces compositions \`a l'aide d'un produit
mono\"\i dal. Et l'\'ecriture alg\'ebrique de ce produit repose
sur une sous-classe de permutations de $\Sy_N$.

\begin{dei}[Permutations connexes]
\index{permutations connexes} Soit $N$ un nombre entier. Soient
$\ok=(k_1,\, \ldots,\, k_b)$ un $b$-uplet et $\oj=(j_1,\, \ldots
,\, j_a)$ un $a$-uplet tels que $|\ok| = k_1+\cdots+k_b
=|\oj|=j_1+\cdots+ j_a= N$.

On d\'efinit les \emph{permutations $(\ok,\, \oj)$-connexes} de
$\Sy_N$
 comme l'ensemble des permutations de
$\Sy_N$ dont le graphe est connexe, si on
 relie les entr\'ees indic\'ees par $j_1+\cdots+ j_i +1, \,
 \ldots,\, j_1+\cdots +j_i + j_{i+1}$, pour $0\leq i \leq a-1$, et les
 sorties indic\'ees par $k_1+\cdots+ k_i +1, \,
 \ldots,\, k_1+\cdots +k_i + k_{i+1}$, pour $0\leq i \leq b-1$.

 On note cet ensemble $\Sc_\kj$.
\end{dei}

\textsc{Exemple :} Dans $\Sy_4$, consid\`erons la permutation
$(1324)$ et sa repr\'esentation g\'eom\`etrique suivante :
$$\xymatrix@C=15pt{1 \ar@{{*}-{*}}[d] &2 \ar@{{*}-{*}}[dr]&
3 \ar@{{*}-{*}}[dl] | \hole &4 \ar@{{*}-{*}}[d]\\
1 &2 & 3&4.}$$ Si on prend $\ok=(2,\, 2)$ et $\oj=(2,\, 2)$, on
relie les entr\'ees $1$, $2$ et $3$, $4$ ainsi que les sorties
$1$, $2$ et $3$, $4$. Ce qui donne le graphe connexe
$$
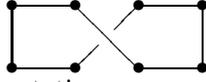
\xymatrix{ *=0{\bullet} \ar@{-}[r] \ar@{-}[d] &*=0{\bullet}\ar@{-}[dr] &*=0{\bullet} \ar@{-}[dl] |\hole \ar@{-}[r]
&*=0{\bullet} \ar@{-}[d] \\
 *=0{\bullet} \ar@{-}[r] &*=0{\bullet} &*=0{\bullet}\ar@{-}[r] &*=0{\bullet} } $$
 Ainsi, la permutation $(1324)$ est une permutation connexe
pour $(2,\, 2)$ et $(2,\, 2)$, $(1324)\in \Sc_{(2,\,2),\, (2,\,
2)}$.\\

\textsc{Contre-exemple :} On prend toujours la permutation
$(1324)$ de $\Sy_4$ mais maintenant $\ok=(1,\, 1,\, 2)$ et
$\oj=(2,\, 1,\, 1)$ cette fois-ci. Ceci donne le graphe non
connexe suivant :
$$
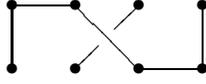
\xymatrix{ *=0{\bullet} \ar@{-}[r] \ar@{-}[d] &*=0{\bullet}\ar@{-}[dr] &*=0{\bullet} \ar@{-}[dl] |\hole
&*=0{\bullet} \ar@{-}[d] \\
 *=0{\bullet}  &*=0{\bullet} &*=0{\bullet}\ar@{-}[r] &*=0{\bullet} }
 $$\\

Citons enfin une proposition \'evidente qui permettra notamment de
faire le lien avec les op\'erades.

\begin{pro}
\label{permconnexesarbre} Lorsque $\ok$ est r\'eduit \`a $(N)$,
$\ok=(N)$, toutes les permutations de $\Sy_N$ sont $\left( (N),
\oj \right)$-connexes, c'est-\`a-dire $\Sc_{\left( (N), \oj
\right)}=\Sy_N$.

Et si $\ok$ est diff\'erent de $(N)$, on a $\Sc_{\ok,\, (1,\,
\ldots,\, 1)}=\emptyset$.
\end{pro}

\subsection{Composition verticale connexe des $\Sy$-bimodules}

On peut maintenant  d\'efinir, sur la cat\'egorie des
$\Sy$-bimodules, la
 structure mono\"\i dale qui nous int\'eresse.\\

Pour deux $a$-uplets $\oj$ et $\oi$, on utilise les conventions
d'\'ecriture suivantes. On note $\Po(\oj,\, \oi)$ le produit
tensoriel $\Po(j_1,\, i_1)\otimes_k \cdots \otimes_k \Po(j_a,\,
i_a)$ et $\Sy_\oj$ l'image du produit direct des groupes
$\Sy_{j_1}\times\cdots \times \Sy_{j_b}$ dans $\Sy_{|\oj|}$.

\begin{dei}[Produit mono\"\i dal $\boxtimes_c$]

\index{produit mono\"\i dal $\boxtimes_c$} Soient $\Qo$, $\Po$
deux $\Sy$-bimodules. Le \emph{produit mono\"\i dal connexe} de
$\Qo$ et $\Po$ est le $\Sy$-bimodule $\Qo \boxtimes_c \Po$ donn\'e
par

$$\mathcal{Q}\boxtimes_c \mathcal{P}(m,\, n) =
 \bigoplus_{N\in \mathbb{N}} \left( \bigoplus_{\ol,\, \ok,\, \oj,\, \oi} k[\mathbb{S}_m]
 \otimes_{\mathbb{S}_\ol}
\mathcal{Q}(\ol,\, \ok)\otimes_{\mathbb{S}_{\ok}} k[\Sc_\kj]
\otimes_{\mathbb{S}_\oj} \mathcal{P}(\oj,\, \oi)
\otimes_{\mathbb{S}_\oi} k[\mathbb{S}_n] \right) \Bigg/_\sim \ ,$$

o\`u la somme directe court sur les $b$-uplets $\ol$, $\ok$ et les
$b$-uplets $\oj$, $\oi$ tels que $|\ol|=m$, $|\ok|=|\oj|=N$,
$|\oi|=n$ et o\`u la relation d'\'equivalence $\sim$ est donn\'ee
par
\begin{eqnarray*}
&&\theta \otimes q_1\otimes \cdots \otimes q_b \otimes \sigma
\otimes p_1 \otimes \cdots \otimes p_a \otimes\omega
  \sim \\
&&\theta \,\tau^{-1}_\ol \otimes q_{\tau^{-1}(1)}\otimes \cdots
\otimes q_{\tau^{-1}(b)} \otimes
 \tau_\ok\,
 \sigma \, \nu_\oj \otimes p_{\nu(1)} \otimes \cdots \otimes p_{\nu(a)}
 \otimes \nu^{-1}_{\ol} \,
 \omega,
\end{eqnarray*}
pour $\theta \in \Sy_m$, $\omega \in \Sy_n$, $\sigma \in \Sc_\kj$
et pour $\tau \in \mathbb{S}_b$ avec $\tau_{k_1,\ldots ,\, k_b}$
la permutation par blocs correspondante (\emph{cf.} conventions),
$\nu \in \mathbb{S}_a$ et $\nu_{j_1,\ldots ,\, j_a}$ la
permutation par blocs correspondante.
\end{dei}

\textsc{Remarque :} A cause de la lourdeur de l'\'ecriture, nous
omettrons souvent dans la suite les repr\'esentations induites.
(C'est souvent le cas pour les op\'erades).

\begin{pro}
L'objet $\Qo\boxtimes_c\Po$ est bien d\'efini.
\end{pro}

\begin{deo}
Pour $\ol$ et $\ok$ deux $b$-uplets tels que $|\ol|=m$ et
$|\ok|=N$, on pose $\ol'=(l_{\tau^{-1}(1)}$, $\ldots$,
$l_{\tau^{-1}(b)})$ et $\ok'=(k_{\tau^{-1}(1)},\, \ldots , \,
k_{\tau^{-1}(b)})$ qui v\'erifient aussi $|\ol'|=m$ et $|\ok'|=N
$. De plus, toute permutation $(\ok,\, \oj)$-connexe $\sigma$
 donne par composition \`a gauche
avec une permutation par blocs du type $\tau_\ok$ et \`a droite
par une permutation par blocs du type
 $\nu_\oj$ une permutation $(\ok',\, \oj')$-connexe. En effet, les
 repr\'esentations g\'eom\'etriques
de $\sigma\in\Sy_N$ et de $\tau_\ok \, \sigma\, \nu_\oj$ sont
hom\'eomorphes et on relie exactement les m\^emes points. En
r\'esum\'e, si $\sigma\in\Sc_\kj$ on a encore  $\tau_\ok \, \sigma
\, \nu_\oj \in \Sc_{\ok',\ \ol'}$. $\cqfd$
\end{deo}

Ce produit mono\"\i dal est d\'efini ainsi pour correspondre \`a
la composition verticale de graphes connexes.

\begin{dei}[Graphes dirig\'es]
\index{graphes dirig\'es} Ce que l'on appelle ici par
\emph{graphes dirig\'es} sont des graphes non planaires, dirig\'es
par un flot, dont les ar\^etes entrant et sortant d'un noeud (ou
sommet) sont indic\'ees par des entiers $\{1,\, \ldots,\, n \}$,
tout comme le sont les entr\'ees et sorties du graphe. On suppose
de plus que chaque noeud admette au moins une entr\'ee et une
sortie.
\end{dei}

\begin{dei}[Graphes connexes]
\index{graphes connexes}
On dit qu'un graphe est \emph{connexe}
s'il est connexe en tant qu'espace topologique.
\end{dei}

\begin{figure}[h]
$$ \xymatrix{
 & \ar[dr]_(0.7){1} 1&2 \ar[d]_(0.6){3} &3\ar[dl]^(0.7){2} &4\ar[dr]_(0.7){1} & &5\ar[dl]^(0.7){2} \\
*+[o][F-]{1} & \ar@{--}[r]& *+[F-,]{\nu_1} \ar@{-}[dl]_(0.3){1}
\ar@{-}[drr]^(0.3){2} \ar@{--}[rrr]& & & *+[F-,]{\nu_2}
\ar@{-}[dr]^(0.3){1}   \ar@{-}[dll]_(0.3){2} |(0.75) \hole  \ar@{--}[r]& \\
 & *=0{} \ar[dr]_(0.7){1} & &*=0{}\ar[dl]^(0.7){2} &*=0{}\ar[dr]_(0.7){1} & &*=0{} \ar[dl]^(0.7){2}\\
*+[o][F-]{2}& \ar@{--}[r] & *+[F-,]{\nu_3} \ar[d]_(0.3){1}
\ar@{--}[rrr]& & & *+[F-,]{\nu_4} \ar[dl]_(0.3){3} \ar[d]_(0.3){2}
\ar[dr]^(0.3){1}
\ar@{--}[r] & \\
& &4 & &1 &2 &3 } $$ \caption{Exemple de graphe connexe \`a deux
niveaux.} \label{figrapheniveau}
\end{figure}
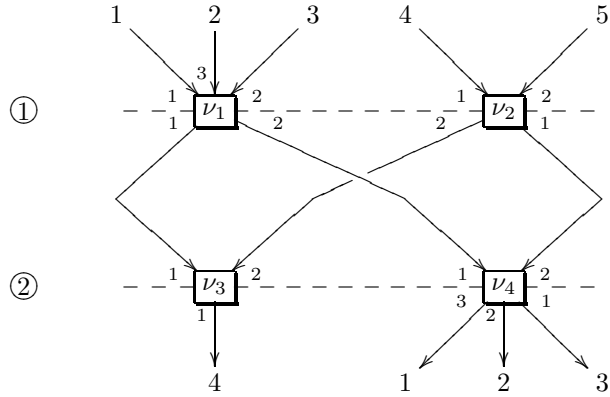

Pour d\'ecrire le produit mono\"\i dal $\boxtimes_c$, on se sert
des \emph{graphes \`a niveaux} connexes \index{graphes \`a
niveaux}, c'est-\`a-dire des graphes dont les noeuds se
r\'epartissent sur des niveaux. La figure~\ref{figrapheniveau}
repr\'esente un graphe \`a deux niveaux o\`u les noeuds sont
indic\'es par $\nu_i$.

A l'aide des graphes connexes \`a deux niveaux $\mathcal{G}_2^c$,
on construit le produit mono\"\i dal $\boxtimes^\mathcal{G}_c$ sur
les $\Sy$-bimodules. Pour un graphe $g$, on appelle $\N_1$
l'ensemble des noeuds appartenant au premier niveau (en fonction
de la direction donn\'ee par le flot global) et $\N_2$ l'ensemble
des noeuds appartenant au second niveau. Pour un noeud $\nu$, on
consid\`ere les deux ensembles $In(\nu)$ et $Out(\nu)$ compos\'es
des ar\^etes entrant et sortant du noeud. On note $|Out(\nu)|$ et
$|In(\nu)|$ les cardinaux de ces ensembles.

\begin{dei}[Produit mono\"\i dal $\boxtimes^\G_c$]
\index{produit mono\"\i dal $\boxtimes^\G_c$}
A $\Qo$ et $\Po$
deux $\Sy$-bimodules, on associe le $\Sy$-bimodule $\Qo
\boxtimes^\G_c \Po$ donn\'e par la formule

$$\Qo\boxtimes^\G_c \Po = \left( \bigoplus_{g\in\G_2^c}
\bigotimes_{\nu \in \N_2} \Qo(|Out(\nu)|,\, |In(\nu)|) \otimes
\bigotimes_{\nu \in \N_1} \Po(|Out(\nu)|,\, |In(\nu)|)\right)
\Bigg/_\approx \ ,$$ o\`u la relation d'\'equivalence $\approx$
est engendr\'ee par
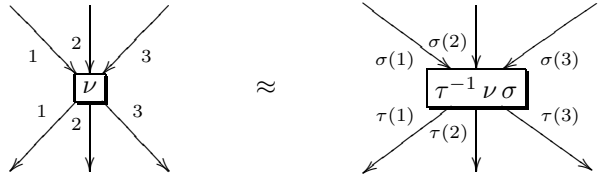
$$\xymatrix{\ar[dr]_(0.5){1} &\ar[d]_(0.5){2} & \ar[dl]^(0.5){3}
& & \ar[dr]_(0.5){\sigma(1)}& \ar[d]_(0.5){\sigma(2)}& \ar[dl]^(0.5)
{\sigma(3)} \\
 & *+[F-,]{\nu} \ar[dl]_(0.4){1}\ar[d]_(0.4){2}\ar[dr]^(0.4){3}&
 &\approx & &
 *+[F-,]{\tau^{-1}\, \nu\,\sigma} \ar[dl]_(0.5){\tau(1)}\ar[d]_(0.5){\tau(2)}
 \ar[dr]^(0.5){\tau(3)}&\\
 & & & & & &\ .}$$
\end{dei}

Ceci qui revient \`a indicer les sommets des graphes \`a deux
niveaux par des \'el\'ements de $\Qo$ et de $\Po$, et \`a relier
les indices des sorties (resp. entr\'ees) d'un sommet \`a l'action
de $\Sy_m$ \`a gauche (resp. $\Sy_n$ \`a droite) sur l'\'el\'ement
correspondant.\\

\textsc{Remarque :} Par souci de concision, on omettra souvent
dans la suite d'\'ecrire la relation d'\'e\-qui\-va\-len\-ce.
Ainsi, lorsque l'on parlera de graphes indic\'es par des
$\Sy$-bimodules, on quotientera implicitement par cette relation
d'\'equivalence.

\begin{pro}
Pour tout couple $(\Qo,\, \Po)$ de $\Sy$-bimodules, les deux
produits $\Qo \boxtimes_c \Po $ et $\Qo \boxtimes^\G_c \Po$ sont
naturellement isomorphes.
\end{pro}

\begin{deo}
A tout \'el\'ement $ \theta \otimes q_1\otimes \cdots \otimes q_b
\otimes \sigma \otimes p_1 \otimes \cdots \otimes p_a
\otimes\omega$ de $k[\mathbb{S}_m] \otimes \mathcal{Q}(\ol,\,
\ok)\otimes k[\Sc_\kj] \otimes \mathcal{P}(\oj,\, \oi) \otimes
k[\mathbb{S}_n] $, on associe un graphe $\Psi(\theta \otimes
q_1\otimes \cdots \otimes q_b \otimes \sigma \otimes p_1 \otimes
\cdots \otimes p_a \otimes\omega)$ dont les noeuds sont indic\'es
par les $q_\beta$ et les $p_\alpha$. Pour cela, on consid\`ere une
repr\'esentation g\'eom\'etrique de $\sigma$. On regroupe et
r\'eindice les sorties en fonction de $\ok$ et les entr\'ees en
fonction de $\oj$. On indice les noeuds ainsi cr\'e\'es \`a l'aide
des $q_\beta$ et $p_\alpha$. Puis pour chaque $q_\beta$ on
construit $l_\beta$ ar\^etes sortant que l'on indice par $1,\,
\ldots ,\, l_\beta$. On proc\`ede de m\^eme avec $p_\alpha$ et les
entr\'ees du graphe. Enfin, on num\'erote les sorties du graphe
avec $\theta$ et les entr\'ees avec $\omega$. Par exemple,
l'\'el\'ement $(4123)\otimes q_1\otimes q_2\otimes (1324) \otimes
p_1 \otimes p_2 \otimes (12435)$ donne le graphe repr\'esent\'e
\`a la figure~\ref{figraphecompo}.

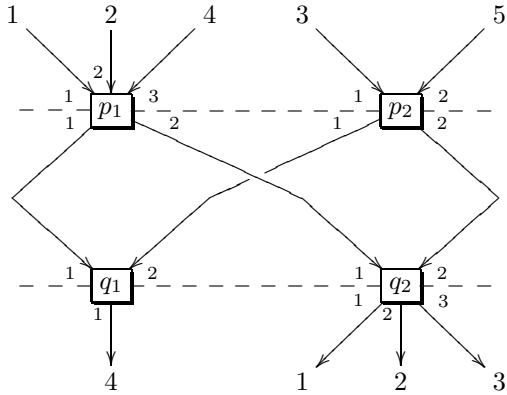
\begin{figure}[h]
$$ \xymatrix{
\ar[dr]_(0.7){1} 1&2 \ar[d]_(0.6){2} &4\ar[dl]^(0.7){3} &3\ar[dr]_(0.7){1} & &5\ar[dl]^(0.7){2} \\
\ar@{--}[r]& *+[F-,]{p_1} \ar@{-}[dl]_(0.3){1}
\ar@{-}[drr]^(0.3){2} \ar@{--}[rrr]& & & *+[F-,]{p_2}
\ar@{-}[dr]^(0.3){2}   \ar@{-}[dll]_(0.3){1} |(0.75) \hole  \ar@{--}[r]& \\
 *=0{} \ar[dr]_(0.7){1} & &*=0{}\ar[dl]^(0.7){2} &*=0{}\ar[dr]_(0.7){1} & &*=0{} \ar[dl]^(0.7){2}\\
\ar@{--}[r] & *+[F-,]{q_1} \ar[d]_(0.3){1}  \ar@{--}[rrr]& & &
*+[F-,]{q_2} \ar[dl]_(0.3){1} \ar[d]_(0.3){2} \ar[dr]^(0.3){3}
\ar@{--}[r] & \\
&4 & &1 &2 &3 } $$ \caption{Image via $\Psi$ de $(4123)\otimes
q_1\otimes q_2\otimes (1324) \otimes p_1 \otimes p_2 \otimes
(12435)$.} \label{figraphecompo}
\end{figure}

Soit $\bar{\Psi}(\theta \otimes q_1\otimes \cdots \otimes q_b
\otimes \sigma \otimes p_1 \otimes \cdots \otimes p_a
\otimes\omega)$ la classe d'\'equivalence de $\Psi(\theta \otimes
q_1\otimes \cdots \otimes q_b \otimes \sigma \otimes p_1 \otimes
\cdots \otimes p_a \otimes\omega)$ pour la relation
 $\approx$. Gr\^ace \`a cette relation
d'\'equivalence, l'application $\bar{\Psi}$ passe naturellement au
quotient $k[\mathbb{S}_m]
 \otimes_{\mathbb{S}_\ol}
\mathcal{Q}(\ol,\, \ok)\otimes_{\mathbb{S}_\ok} k[\Sc_\kj]
\otimes_{\mathbb{S}_\oj} \mathcal{P}(\oj,\, \oi)
\otimes_{\mathbb{S}_\oi} k[\mathbb{S}_n]$. Quant \`a la relation
d'\'equivalence $\sim$, elle correspond \`a un r\'earrangement
(hom\'eomorphe) dans l'espace du graphe engendr\'e. Le graphe
ainsi cr\'e\'e \'etant non-planaire, il est invariant sur les
classes d'\'equivalences pour la relation $\sim$. On a finalement
une application $\widetilde{\Psi}$ de $\Qo\boxtimes_c\Po$ vers
$\Qo\boxtimes^\G_c \Po$. Comme tout graphe admet une
repr\'esentation de la forme $ \theta \otimes q_1\otimes \cdots
\otimes q_a \otimes \sigma \otimes p_1 \otimes \cdots \otimes p_b
\otimes\omega$, on peut exhiber une r\'eciproque \`a
$\widetilde{\Psi}$. Ainsi, les deux produits $\Qo\boxtimes_c\Po$
et $\Qo\boxtimes^\G_c \Po$ sont naturellement isomorphes et
repr\'esentent la m\^eme information. $\cqfd$
\end{deo}

Il reste \`a d\'efinir l'objet qui jouera le r\^ole d'unit\'e dans
cette cat\'egorie. On pose
$$I=\left\{
\begin{array}{l}
I(1,\, 1)=k,    \\
I(m,\,n)=0 \quad \textrm{sinon}.
\end{array} \right.$$

\begin{pro}
La cat\'egorie $(\mathbb{S}\textrm{-biMod}, \, \boxtimes_c, \, I)$
est une cat\'egorie mono\"\i dale.
\end{pro}

\begin{deo}
Pour montrer la relation d'unit\'e sur $\Po\boxtimes_c I(m,\, n)$,
il faut \'etudier $\Sc_\kj$ pour $\oj=(1,\, \ldots,\, 1)$. Ce
dernier est vide sauf si on r\'eunit toutes les sorties, soit pour
$\ok=(n)$. Et dans ce cas $\Sc_{(n),\, (1,\, \ldots ,\, 1)}$ vaut
$\Sy_n$ (\emph{cf.} Proposition~\ref{permconnexesarbre}). Ainsi,
on a
$$ \mathcal{P}\boxtimes_c I (m,\, n) =
\mathcal{P}(m,\, n) \otimes_{\mathbb{S}_n} k\lbrack \mathbb{S}_n
\rbrack \otimes_{\mathbb{S}_1^{\times n}} k^{\otimes n} =
\mathcal{P}(m,\, n) \otimes k .id_n \otimes k^{\otimes n} =
\mathcal{P}(m,\, n).$$ (Le cas $I\boxtimes_c \Po=\Po$ est
parfaitement sym\'etrique.)

Pour montrer la relation d'associativit\'e $\mathcal{R}\boxtimes_c
(\mathcal{Q} \boxtimes_c \mathcal{P}) = (\mathcal{R}\boxtimes_c
\mathcal{Q}) \boxtimes_c \mathcal{P}$ du produit $\boxtimes_c$, on
utilise le produit $\boxtimes_c^\G$. En effet, il suffit de voir
que les $\Sy$-bimodules $\mathcal{R}\boxtimes_c^\G (\mathcal{Q}
\boxtimes_c^\G \mathcal{P}) $ et $(\mathcal{R}\boxtimes_c^\G
\mathcal{Q}) \boxtimes_c^\G \mathcal{P}$ correspondent aux graphes
\`a $3$ niveaux indic\'es par des \'el\'ements de $\mathcal{R}$,
$\Qo$ et $\Po$, soit
$$\left( \bigoplus_{g\in\G_3^c} \bigotimes_{\nu \in \N_3}
\mathcal{R}(|Out(\nu)|,\, |In(\nu)|) \otimes \bigotimes_{\nu \in
\N_2} \Qo(|Out(\nu)|,\, |In(\nu)|) \otimes \bigotimes_{\nu \in
\N_1} \Po(|Out(\nu)|,\, |In(\nu)|)\right) \Bigg/_\approx \ .$$ (La
composition verticale de graphes connexes donne encore un graphe
connexe). $\cqfd$
\end{deo}

\subsection{Les sous-cat\'egories mono\"\i dales $(k\textrm{-Mod},\ \otimes,\, k)$
et $(\Sy\textrm{-Mod},\, \circ,\, I)$}

Les deux ca\-t\'e\-go\-ries mono\"\i dales $k$-Mod et $\Sy$-Mod
apparaissent comme des sous-cat\'egories mono\"\i dales pleines de
la cat\'egories des $\Sy$-bimodules avec le produit $\boxtimes_c$
d\'efini
pr\'ec\'edemment. \\

En ce qui concerne la cat\'egorie des modules sur $k$ munie du
produit tensoriel classique, il suffit d'associer \`a tout module
$V$ le $\mathbb{S}$-bimodule suivant :
$$\left\{ \begin{array}{l}
V(1,\, 1) =V, \\
V(j,\, i) = 0 \quad \textrm{sinon}.
\end{array} \right. $$

On retrouve alors le produit tensoriel sur $k$ car pour $V$ et $W$
deux modules, on a $V\boxtimes_c W (1,\, 1)= V\otimes_k W$ et
comme il n'existe pas de permutations connexes associ\'ees \`a des
$a,\, b$-uplets de la forme $(1,\, \ldots ,\, 1)$ (\emph{cf.}
proposition~\ref{permconnexesarbre}), on a $V\boxtimes W (j,\, i)=
0$ pour $(j,\, i) \ne (1,\, 1)$. Et les morphismes entre deux
$\Sy$-bimodules compos\'es uniquement d'un $(\Sy_1,\,
\Sy_1)$-module
correspondent aux morphismes de $k$-modules.\\

La cat\'egorie des $\Sy$-bimodules contient aussi les
$\mathbb{S}$-modules li\'es aux op\'erades. De la m\^eme
mani\`ere, \`a partir d'un $\mathbb{S}$-module $\mathcal{P}$, on
d\'efinit un $\mathbb{S}$-bimodule :
$$\left\{ \begin{array}{l}
\mathcal{P}(1,\, n) =\mathcal{P}(n) \quad \textrm{pour} \ n\in \mathbb{N}^*, \\
\mathcal{P}(j,\, i) = 0 \quad \textrm{si} \ j\ne 1.
\end{array} \right. $$

Nous avons vu \`a la proposition~\ref{permconnexesarbre} que
$\Sc_{(N),\, (1,\, \ldots,\, 1)})=\Sy_N$. Ce qui se traduit ici
par
$$ \mathcal{Q}\boxtimes_c \mathcal{P} (1,\, n)=
\bigoplus_{N\leq n} \left( \bigoplus_{i_1+\cdots +i_N =n}
\mathcal{Q}(1,\, N) \otimes_{\mathbb{S}_N} k[\mathbb{S}_N]
\otimes_{\mathbb{S}_1^{\times N}} \mathcal{P}(\overline{1},\, \oi)
\otimes_{\mathbb{S}_\oi} k[\mathbb{S}_n] \right) \Bigg/_\sim \ ,$$
o\`u la relation d'\'equivalence $\sim$ s'\'ecrit ici
$$q\otimes \sigma \nu
\otimes p_1 \otimes \cdots \otimes p_N \sim q\otimes \sigma
\otimes p_{\nu(1)}\otimes \cdots \otimes p_{\nu(N)}.$$ Cette
relation revient \`a prendre les coinvariants pour l'action de
$\Sy_N$ dans l'expression $ \mathcal{Q}(1,\, N) \otimes_k
\mathcal{P}(1,\, i_1) \otimes_k \cdots \otimes_k \mathcal{P}(1,\,
i_N)$. On retombe bien sur le produit mono\"\i dal des
$\mathbb{S}$-modules (\emph{cf.} \cite{GK} et J.-P. May
\cite{May}), que l'on connait plus sous la forme alg\'ebrique
suivante (en omettant les induites)
$$  \mathcal{Q}\boxtimes \mathcal{P} (1,\, n)=
\bigoplus_{N\leq n} \left( \bigoplus_{i_1+\cdots +i_N =n}
\mathcal{Q}(1,\, N) \otimes \mathcal{P}(1,\, i_1) \otimes \cdots
\otimes \mathcal{P}(1,\, i_N) \right)_{\mathbb{S}_N} =
\mathcal{Q}\circ\mathcal{P}(n).$$ Remarquons qu'avec le produit
$\boxtimes_c^\G$, on retrouve la composition des $\Sy$-modules
\'ecrite \`a l'aide des arbres (\emph{cf.} \cite{GK} et
\cite{Loday3}).

Quant \`a $\mathcal{Q}\boxtimes_c\mathcal{P}(j,\, i)$, pour $j\ne
1$, cette composition correspond \`a la juxtaposition d'arbres
(for\^et) et repose donc sur des graphes non connexes. Ainsi,
$\mathcal{Q}\boxtimes_c\mathcal{P}(j,\, i)=0$, pour $j\ne 1$.

\subsection{Composition horizontale et verticale des $\Sy$-bimodules}

Afin de repr\'esenter la concat\'enation de deux op\'erations, on
introduit un produit mono\"\i dal $\otimes$ de la mani\`ere
suivante.

\begin{dei}[Produit de concat\'enation $\otimes$]
\index{produit de concat\'enation $\otimes$}
 Soient $\Po$, $\Qo$
deux $\Sy$-bimodules. On d\'efinit le \emph{produit de
concat\'enation $\otimes$} par la formule suivante :
$$ \Po\otimes \Qo (m,\,n) = \bigoplus_{m'+m''=m \atop n'+n''=n}
\Po(m',\, n')\otimes \Qo(m'',\, n''),$$
$$ \textrm{o\`u}  \ \
\Po(m,\ n)\otimes \Qo(m',\,
n')=k[\Sy_{m+m'}]\otimes_{\Sy_m\times\Sy_{m'}} \Po(m,\,n)
\otimes_k \Qo(m',\, n') \otimes_{\Sy_n\times
\Sy_{n'}}k[\Sy_{n+n'}].$$
\end{dei}

\begin{figure}[h]
$$ \xymatrix@C=10pt{ 1\ar[dr]_(0.6){1} & 2 \ar[dr]^(0.85){1} & 3
\ar[dl]_(0.85){2} | \hole &4 \ar[dl]^(0.6){2} \\
\ar@{..}[r]& *+[F-,]{\quad p\quad} \ar[d]^(0.35){1} \ar[dl]_(0.5){2}
\ar@{..}[r]& *+[F-,]{\quad q\quad}\ar[d]^(0.35){1} \ar@{..}[r] & \\
1& 3& 2& }$$ \caption{Composition horizontale de deux
op\'erations.} \label{concatenation}
\end{figure}
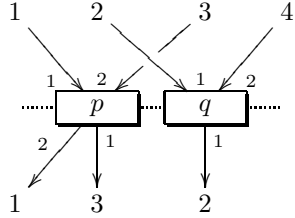

Le produit $\otimes$ correspond \`a la notion intuitive de
concat\'enation d'op\'erations r\'epr\'esent\'ee \`a la
figure~\ref{concatenation}. Pour ce produit on parle de
\emph{composition horizontale} \index{composition horizontale}.

\begin{pro}
Le produit $\otimes$ d\'efinit un produit mono\"\i dal
sym\'etrique dans la cat\'egorie des $\Sy$-bimodules. L'unit\'e
est donn\'ee par le $\Sy$-bimodule $k$ d\'efini par
$$ \left\{
\begin{array}{ll}
k(0,\, 0)=k, & \\
k(m,\, n)=0 & \textrm{sinon}.
\end{array}
\right.$$
\end{pro}

\begin{deo}
La relation d'unit\'e vient de l'\'egalit\'e
$$ \Po\otimes k (m,\, n) = k[\Sy_m]\otimes_{\Sy_m} \Po(m,\, n) \otimes_k k
\otimes_{\Sy_n} k[\Sy_n]=\Po(m,\, n). $$

Et l'associativit\'e vient de la relation
\begin{eqnarray*}
&&\big(\Po(m,\,n)\otimes \Qo(m',\, n')\big)\otimes R(m'',\, n'') =
\Po(m,\,n)\otimes \big(\Qo(m',\, n')\otimes R(m'',\, n'')\big) = \\
&& k[\Sy_{m+m'+m''}]\otimes_{\Sy_m\times \Sy_{m'}\times \Sy_{m''}}
\Po(m,\,n)\otimes_k \Qo(m',\, n')\otimes_k R(m'',\, n'')
\otimes_{\Sy_n\times \Sy_{n'}\times \Sy_{n''}} k[\Sy_{n+n'+n''}].
\end{eqnarray*}

L'isomorphisme de sym\'etrie est donn\'e par
$$ \Po(m,\, n) \otimes \Qo(m',\, n') \xrightarrow{\tau} (1,\,2)_{m,\, m'}
\big( \Po(m,\, n) \otimes \Qo(m',\, n') \big) (1,\, 2)_{n,\, n'} =
\Qo(m',\, n') \otimes \Po(m,\, n).$$
\end{deo}

\textsc{Remarque :} Ce produit mono\"\i dal est bilin\'eaire.\\

En mimant ce que l'on a fait dans les parties pr\'ec\'edentes, on
peut d\'efinir un produit qui repr\'esente les compositions
verticales d'op\'erations mais suivant des sch\'emas non
n\'ecessairement connexes.

\begin{dei}[Produit de composition $\boxtimes$]
\index{produit de composition $\boxtimes$}
 Soient $\Qo$, $\Po$
deux $\Sy$-bimodules. On d\'efinit le \emph{produit de composition
$\Qo \boxtimes \Po$} par

$$  \Qo\boxtimes \Po (m, \, n) = \bigoplus_{N \in \mathbb{N}}
\left( \bigoplus_{\ol,\, \ok,\, \oj,\, \oi} k[\Sy_m]
\otimes_{\Sy_\ol} \Qo(\ol,\, \ok) \otimes_{S_\ok}  k\lbrack \Sy_N
\rbrack \otimes_{S_\oj} \mathcal{P}(\oj,\, \oi) \otimes_{\Sy_\oi}
k[\Sy_n]\right) \Bigg/_\sim \ ,$$ o\`u la relation d'\'equivalence
$\sim$ est la m\^eme que dans le cas connexe.
\end{dei}

A la diff\'erence du produit $\otimes$, le produit $\boxtimes$
repr\'esente les \emph{compositions verticales} \index{composition
verticale} d'op\'erations.\\

Comme dans le cas connexe, on a une autre \'ecriture de ce produit
en terme de graphes dirig\'es, \`a la diff\'erence pr\`es qu'ici
les graphes ne seront pas suppos\'es connexes.

On consid\`ere l'ensemble $\mathcal{G}_2$ des graphes dirig\'es
 \`a deux niveaux non n\'ecessairement
connexes. On d\'efinit alors le produit bas\'e sur de tels
graphes.

\begin{dei}[Produit de composition $\boxtimes^\mathcal{G}$]
\index{produit de composition $\boxtimes^\mathcal{G}$}
 Soient $\Qo$, $\Po$ deux $\Sy$-bimodules. On d\'efinit le
produit $\Qo \boxtimes^\mathcal{G} \Po$ par
$$  \Qo\boxtimes^\mathcal{G}\Po  =
\left( \bigoplus_{g\in \mathcal{G}^2} \bigotimes_{\nu \in \N_2}
\Qo(|Out(\nu)|,\, |In(\nu)|) \otimes  \bigotimes_{\nu \in \N_1}
\Po(|Out(\nu)|,\, |In(\nu)|) \right) \Bigg/_\approx \ ,$$ o\`u la
relation d'\'equivalence $\approx$ est la m\^eme que dans le cas
connexe.
\end{dei}

Comme dans le cas connexe, ces deux d\'efinitions sont
\'equivalentes.

\begin{pro}
Pour tout couple $(\Qo,\, \Po)$ de $\Sy$-bimodules, les deux
produits de composition $\Qo \boxtimes \Po$ et $\Qo
\boxtimes^\mathcal{G} \Po$ sont naturellement isomorphes.
\end{pro}

\begin{deo}
La d\'emonstration est identique. On introduit le m\^eme type
d'isomorphisme $\widetilde{\Psi}$ entre $\Qo \boxtimes \Po$ et
$\Qo \boxtimes^\mathcal{G} \Po $. $\cqfd$
\end{deo}

Par la suite, nous aurons besoin de consid\'erer l'ensemble des
concat\'enations possibles d'op\'erations de $\Po$.

\begin{dei}[Les $\Sy$-bimodules $T_\otimes(\Po)$ et $S_\otimes(\Po)$]
\index{alg\`ebre sym\'etrique libre $S_\otimes(\Po)$} Comme le
produit mono\"\i dal de $\Sy$-bimodules $\otimes$ est
bilin\'eaire, le mono\"\i de libre sur $\Po$ pour la
concat\'enation est donn\'e par
$$T_\otimes(\Po)=\bigoplus_{n\in \mathbb{N}} \Po^{\otimes n}.$$
Comme $\Po\otimes \Qo$ est isomorphe \`a $\Qo \otimes \Po$ et que
l'on on aura \`a consid\'erer un morphisme commutatif sur $\Po
\otimes \Qo$ par la suite, on introduit l'\emph{alg\`ebre
sym\'etrique libre tronqu\'ee} sur $\Po$.
$$S_\otimes(\Po)= \bigoplus_{n\in \mathbb{N}^*}
(\Po^{\otimes n})_{\Sy_n}.$$
\end{dei}

\textsc{Remarque :} On a exclut l'op\'eration scalaire $k(0,\, 0)$
de la d\'efinition de l'alg\`ebre sym\'etrique consid\'er\'ee.\\

Le foncteur $S_\otimes$ permet de relier les deux produits de
composition $\boxtimes_c$ et $\boxtimes$.

\begin{pro}
\label{lienentreproduits}
 Soient $\Qo$ et $\Po$ deux
$\Sy$-bimodules. On a l'\'egalit\'e
$$S_\otimes (\Qo \boxtimes_c \Po) = \Qo \boxtimes \Po. $$
\end{pro}

\begin{deo}
Le r\'esultat repose sur le fait que tout graphe de
$\mathcal{G}_2$ peut se d\'ecomposer en produit de graphes
connexes appartenant \`a $\mathcal{G}^c_2$ et que l'on ne tient
pas compte ici de l'ordre suivant lequel sont concat\'en\'es ces
graphes connexes. $\cqfd$
\end{deo}

Comme le produit mono\"\i dal $\otimes$ est bien connu, par
exemple du point de vue homologique, pour \'etudier le produit
$\boxtimes$, il suffira de faire l'\'etude sur $\boxtimes_c$ et de
passer \`a la concat\'enation \`a la fin. Ce sera par exemple le
cas des diff\'erents complexes de cha\^\i nes introduits dans la
suite de cette th\`ese (bar constructions, complexes de Koszul)
dont les diff\'erentielles se font composante connexe par
composante connexe.\\

\textsc{Remarque :} Le produit de composition $\boxtimes$ n'est
pas un produit mono\"\i dal. L'associativit\'e ne fait aucun
doute. Par contre, la relation d'unit\'e fait d\'efaut. On a
$$ \Po\boxtimes I= S_\otimes(\Po).$$
Or, en g\'en\'eral, $\Po$ est diff\'erent de
$S_\otimes(\Po)$.

\subsection{Repr\'esentation des \'el\'ements de
$\Qo \boxtimes_c \Po$ et de $\Qo \boxtimes \Po$} Les produits $\Qo
\boxtimes_c \Po$ et $\Qo \boxtimes \Po$ sont d\'efinis comme des
quotients par la relation d'\'equivalence $\sim$ qui permute
l'ordre des \'el\'ements de $\Qo$ et de $\Po$. Afin notamment de
pouvoir \'etudier le comportement homologique de ces produits, on
cherche des repr\'esentants naturels des classes d'\'equivalence
pour la relation $\sim$.\\

Pour cela on consid\`ere les couples de partitions ordonn\'ees de
$([m],\, [n])$.

\begin{dei}[Partition ordonn\'ee de $\lbrack n \rbrack$]
\index{partition ordonn\'ee de $[n]$}
 Une \emph{partition ordonn\'ee} de $[n]$ est une suite
$(\Pi_1,\ldots,\, \Pi_k)$ d'ensembles qui forment une partition,
au sens usuel du terme, de $[n]=\{1,\, \ldots ,\, n\}$.
\end{dei}

Parmi les partitions ordonn\'ees de $[n]$, on choisit celles qui
sont \emph{croissantes}, c'est-\`a-dire qui v\'erifient
$\min(\Pi_1)<\min(\Pi_2)<\cdots <min(\Pi_k)$.

Soit $\Theta(m,\, n)$ l'ensemble des couples
$\big((\Pi_1',\ldots,\, \Pi_b'),\, (\Pi_1,\ldots,\, \Pi_a)
\big)$de partitions ordonn\'ees croissantes de $\left( [m],\,
[n]\right)$.

Lorsqu'il n'y a pas d'ambigu\"\i t\'e sur les entiers $m$ et $n$,
on note cet ensemble $\Theta$.

\begin{pro}
\label{Representation} Le produit  $\Qo\boxtimes_c\Po(m,\,n)$ est
isomorphe, en tant que $(\Sy_m,\Sy_n)$-bimodule, \`a
$$\bigoplus_{\left((\Pi_1',\ldots,\, \Pi_b'),\, (\Pi_1,\ldots,\, \Pi_a)
\right)\in \Theta(m,\, n) \atop |\ok|=|\oj|}
k[\Sy_m]\otimes_{\Sy_\ol} \mathcal{Q}(\ol,\,
\ok)\otimes_{\mathbb{S}_{\ok}} k[ \mathbb{S}_c^\kj]
\otimes_{\mathbb{S}_{\oj}} \mathcal{P}(\oj,\,
\oi)\otimes_{\Sy_\oi} k[\Sy_n],$$ o\`u $l_\beta$ est \'egal au
cardinal de l'ensemble $\Pi_\beta'$ et $i_\alpha$ \`a celui de
$\Pi_\alpha$.
\end{pro}

Dit autrement, tout \'el\'ement de $\mathcal{Q}\boxtimes_c
\mathcal{P}$ peut se repr\'esenter sous la forme
$$ q_1^{\Pi_1'}\otimes \cdots \otimes q_b^{\Pi_b'} \otimes
\sigma \otimes p_1^{\Pi_1} \otimes \cdots \otimes p_a^{\Pi_a} =
(q_1^{\Pi_1'} ,\, \ldots ,\,  q_b^{\Pi_b'}) \, \sigma
(p_1^{\Pi_1},\, \ldots ,\, p_a^{\Pi_a}), $$ avec
$\left((\Pi_1',\ldots,\, \Pi_b'),\, (\Pi_1,\ldots,\, \Pi_a)
\right)\in \Theta$.\\

L'action de $\Sy_m$ \`a gauche sur ce module revient \`a permuter
les entiers d'une partition $\left(\Pi_1',\ldots ,\, \Pi_b'
\right)$. Et s'il faut permuter des $\Pi_\beta'$ pour retomber sur
une partition ordonn\'ee, on permute de la m\^eme mani\`ere les
\'el\'ements $q_\beta$ correspondant.

\begin{deo}
On utilise la description du produit $\boxtimes_c$ en termes de
graphes dont les entr\'ees et les sorties sont indic\'ees par des
entiers de $\{1,\ldots,\, n\}$ et $\{1,\ldots ,\,m\}$
respectivement. Les partitions $(\Pi_1,\, \ldots ,\ \Pi_a)$
repr\'esentent les entr\'ees des op\'erations $(p_1,\ldots ,\,
p_a)$ et les partitions $(\Pi_1',\, \ldots ,\ \Pi_b')$
repr\'esentent les sorties des op\'erations $(q_1,\ldots ,\,
q_b)$. $\cqfd$
\end{deo}

Cette \'ecriture nous permettra de montrer que le
produit mono\"\i dal $\boxtimes_c$ v\'erifie certaines
propri\'et\'es homologiques (\emph{cf.} chapitre $3$ section $1.2$).\\

\textsc{Remarque :} La premi\`ere op\'eration $p_1$ sur la ligne
des $p$ est celle qui re\c coit l'entr\'ee indic\'ee par $1$. Il
en va de m\^eme pour $q_1$ avec la sortie indic\'ee par $1$. Cette
propri\'et\'e permet notamment de diff\'erencier sur chaque ligne
une op\'eration des autres. Ceci nous permettra de construire des
homotopies sur les bar et cobar constructions augment\'ees
(\emph{cf.} chapitre $4$ section $3$).\\

On a le m\^eme r\'esultat pour le produit $\boxtimes$.

\begin{pro}
\label{Representationboxtimes} Le produit  $\Qo\boxtimes\Po(m,\,
n)$ est isomorphe, en tant que $(\Sy_m,\, \Sy_n)$-bimodule, \`a
$$\bigoplus_{\left((\Pi_1',\ldots,\, \Pi_b'),\, (\Pi_1,\ldots,\, \Pi_a)
\right)\in \Theta(m,\, n) \atop |\ok|=|\oj|} k[\Sy_m]
\otimes_{\Sy_\ol} \mathcal{Q}(\ol,\,
\ok)\otimes_{\mathbb{S}_{\ok}} k[\Sy_{|\ok|} ]
\otimes_{\mathbb{S}_{\oj}} \mathcal{P}(\oj,\, \oi)
\otimes_{\Sy_\oi} k[\Sy_n].$$
\end{pro}

\section{Bifoncteurs de Schur}

A tout $\Sy$-bimodule $\Po$, on peut associer
un bifoncteur dit de Schur. La composition de tels bifoncteurs est
li\'ee au produit mono\"\i dal $\boxtimes_c$. Les bifoncteurs
introduits ici g\'en\'eralisent la notion de foncteur de Schur des
op\'erades.

\subsection{D\'efinition}

Soit $\Po$ un $\Sy$-bimodule.
\begin{dei}[Bifoncteur de Schur complet]
\index{bifoncteur de Schur complet}
On appelle \emph{bifoncteur de
Schur complet}, associ\'e au $\Sy$-bimodule $\Po$, le foncteur
d\'efini par
\begin{eqnarray*}
\Sy\textrm{-biMod}\times\Sy\textrm{-biMod} &\to& \Sy\textrm{-biMod}\\
(W,\, V) &\mapsto& W\boxtimes_c \Po \boxtimes_c V.
\end{eqnarray*}
Et on le note $\mathcal{S}_\Po$.
\end{dei}

\begin{dei}[Bifoncteur de Schur r\'eduit]
\index{bifoncteur de Schur r\'eduit}
On appelle \emph{bifoncteur
de Schur r\'eduit} le foncteur suivant
\begin{eqnarray*}
\Sy\textrm{-biMod}\times\Sy\textrm{-biMod} &\to& \textrm{bigr-}k\textrm{-Mod}\\
(W,\, V) &\mapsto& \bigoplus_{m,\, n}
k\otimes_{\Sy_m}((W\boxtimes_c\Po \boxtimes_c
V)(m,\, n))\otimes_{\Sy_n}k \\
&=& \bigoplus_{m,\, n}\bigoplus_N \bigoplus_{\left((\ol,\, \ok),\,
( \oj,\, \oi )\right )\in \, \Theta} W\boxtimes_c\Po(\ol,\,
\ok)\otimes_{\mathbb{S}_\ok} k\lbrack \mathbb{S}_c^{\ok, \, \oj}
\rbrack \otimes_{\mathbb{S}_{\oj}} V(\oj,\,
\oi) \\
&=& \bigoplus_{m,\, n}\bigoplus_N \bigoplus_{\left((\ol,\, \ok),\,
(\oj,\, \oi)\right )\in \, \Theta} W(\ol,\,
\ok)\otimes_{\mathbb{S}_{\ok}} k\lbrack \mathbb{S}_c^{\ok, \, \oj}
\rbrack \otimes_{\mathbb{S}_\oj} \Po\boxtimes_c V(\oj,\, \oi),
\end{eqnarray*}
que l'on note $\overline{\Sh}_\Po$.
\end{dei}

Ces deux bifoncteurs de Schur induisent un foncteur $\Sh$
(respectivement $\overline{\Sh}$) entre la cat\'egorie des
$\Sy$-bimodules et celle des bifoncteurs sur
$\Sy\textrm{-biMod}\times\Sy\textrm{-biMod}$
($\textrm{biFonct}_{\Sy\textrm{-biMod}}$).

\begin{dei}[Foncteur de Schur]
\index{foncteur de Schur}
On appelle \emph{foncteur de Schur},
not\'e $\Sh$, le foncteur
\begin{eqnarray*}
\Sy\textrm{-biMod} &\to& \textrm{biFonct}_{\Sy\textrm{-biMod}} \\
\Po &\mapsto& \Sh_\Po.
\end{eqnarray*}
\end{dei}

\subsection{Op\'erations sur les bifoncteurs}

Nous allons \'etudier la comportement de ce foncteur $\Sh$
vis-\`a-vis des
op\'erations respectives des deux cat\'egories $\Sy$-biMod et biFonct.\\

La cat\'egorie des bifoncteurs de
$\Sy\textrm{-biMod}\times\Sy\textrm{-biMod}$ vers $\Sy$-biMod (ou
bigr-$k$-Mod) est munie de coproduits $\oplus$. En effet, on
d\'efinit $(\mathcal{F}\oplus\mathcal{F}')(W,\, V)$ par
$\mathcal{F}(W,\, V) \oplus \mathcal{F}'(W,\, V)$.
Malheureusement, le foncteur $\Sh$ ne pr\'eserve pas toujours les
coproduits. La d\'ependance de $\Sh$ en $\Po$ est pas lin\'eaire.  \\

De la m\^eme mani\`ere, si on d\'efinit le produit tensoriel
 de deux $\Sy$-bimodules $\Qo$, $\Po$ par
$$(\Qo\bar{\otimes}\Po) (m,\, n) = \Qo(m,\, n) \otimes_k \Po(m,\,
n),$$ o\`u les actions de $\Sy_m$ et $\Sy_n$ sont les actions
diagonales (\`a ne pas confondre avec la concat\'enation
$\otimes$), alors la cat\'egorie des bifoncteurs est munie d'un
produit tensoriel sym\'etrique
$(\mathcal{F}\bar{\otimes}\mathcal{F}')(W,\, V)= \mathcal{F}(W,\,
V) \bar{\otimes} \mathcal{F}'(W,\, V)$. Cette fois-ci, le foncteur
de Schur $\Sh$ est un foncteur mono\"\i dal pour le produit
tensoriel $\bar{\otimes}$. On pose $(k)_{m,\, n}$ comme \'etant le
$\Sy$-bimodule ayant $k$ comme module en tout degr\'e. Cette objet
est l'unit\'e dans la cat\'egorie mono\"\i dale
$(\Sy\textrm{-biMod},\, \bar{\otimes})$. De m\^eme, on pose
$(k)_{m,\, n}$ pour le bifoncteur constant d'image $(k)_{m,\, n}$
qui est l'unit\'e pour la cat\'egorie mono\"\i dale
$(\textrm{biFonct}_{\Sy\textrm{-biMod}},\, \bar{\otimes)}$.

\begin{pro}
Le foncteur de Schur $\Sh$ est un foncteur mono\"\i dal entre la
cat\'egorie $(\Sy\textrm{-biMod}$, $\bar{\otimes}$, $(k)_{m,\, n})
$ et $(\textrm{biFonct}_{\Sy\textrm{-biMod}},\, \bar{\otimes},\,
(k)_{m,\, n}) $.
\end{pro}

 En ce qui concerne la composition des bifoncteurs, on a les
relations suivantes.

\begin{lem}
Pour tout  $W,\, V$ et $\Qo,\, \Po$ dans $\Sy$-biMod, les
bifoncteurs de Schur v\'erifient
$$\Sh_\Qo(W,\, \Sh_\Po (I,\, V))= \Sh_\Qo(\Sh_\Po (W,\, I),\, V)=\Sh_{\Qo\boxtimes_c\Po}(W,\, V).$$
\end{lem}

\begin{deo}
Par d\'efinition du bifoncteur de Schur on a
$$\Sh_\Qo(W,\, \Sh_\Po (I,\, V))= \Sh_\Qo(\Sh_\Po (W,\, I),\, V)=\Sh_{\Qo\boxtimes_c\Po}(W,\, V)=
W\boxtimes_c \Qo \boxtimes_c \Po \boxtimes_c V .$$ $\cqfd$
\end{deo}

\begin{dei}[Foncteur de Schur \`a droite]
\index{foncteur de Schur \`a droite}
On appelle \emph{foncteur de
Schur \`a droite} associ\'e \`a un $\Sy$-bimodule $\Po$ la
restriction suivante du bifoncteur $\Sh_\Po$
\begin{eqnarray*}
\Sh'_\Po \ : \ \Sy\textrm{-bimod} &\to& \Sy\textrm{-bimod} \\
V &\mapsto & \Sh_\Po (I,\, V).
\end{eqnarray*}
\end{dei}
Avec cette d\'efinition, le lemme pr\'ec\'edent s'\'ecrit :

\begin{pro}
\label{foncteurdeSchur} Le foncteur $\Sh'$ est un foncteur
mono\"\i dal entre les cat\'egories $(\Sy\textrm{-biMod},\,
\boxtimes_c ,\, I)$ et $(Fonct_{\Sy\textrm{-biMod}},\, \circ,\,
id)$.
\end{pro}

\begin{deo}
Dit autrement, on a $\Sh'_\Qo\circ\Sh'_\Po
(V)=\Sh'_{\Qo\boxtimes_c \Po}(V) $ et $\Sh'_I(V)=V$.$\cqfd$
\end{deo}

\textsc{Remarque :} Cette \'ecriture relie l'associativit\'e du
produit $\boxtimes_c$  \`a celle de la composition des foncteurs.

\section{D\'efinitions des prop\'erades et des PROPs}

On peut maintenant d\'efinir les notions de \emph{prop\'erade} et
de \emph{PROP} ainsi que celle de \emph{g\`ebre sur une
prop\'erade}. Une prop\'erade nous servira \`a coder les
op\'erations \`a plusieurs entr\'ees et plusieurs sorties qui
r\'egissent les diff\'erents types de g\`ebres (alg\`ebres,
cog\`ebres, big\`ebres, etc ...).

\subsection{D\'efinition et premiers exemples de prop\'erades}

\begin{dei}[Prop\'erade]
\index{prop\'erade} Une \emph{prop\'erade} $(\Po,\, \mu,\, \eta)$
est un mono\"\i de dans la cat\'egorie
$(\mathbb{S}\textrm{-biMod}$, $\boxtimes_c$, $I)$, o\`u le produit
mono\"\i dal $\boxtimes_c$ est celui d\'efinit \`a la section
$1.3$.
\end{dei}

Toutes les prop\'erades que nous consid\'ererons ici seront
\emph{r\'eduites}, c'est-\`a-dire qu'elles verifient $\Po(m,\,
n)=0$ si $m=0$ ou $n=0$.\index{prop\'erade r\'eduite}\\

\textsc{Exemples :}
\begin{itemize}

\item Les premiers exemples viennent des sous-cat\'egories
mono\"\i dales pleines $k$-Mod et $\Sy$-Mod. Ainsi, les alg\`ebres
et les op\'erades sont des prop\'erades particuli\`eres.

\item Les autres exemples de prop\'erades que nous traitons ici
viennent de PROPs (\emph{cf.} section~\ref{exemples}).
\end{itemize}

\begin{dei}[$\Sy$-bimodule gradu\'e par un poids]
\index{$\Sy$-bimodule gradu\'e par un poids} On appelle
\emph{$\Sy$-bimodule gradu\'e par un poids} toute somme directe,
indic\'ee par $\rho\in \mathbb{N}$, de $\Sy$-bimodules,
c'est-\`a-dire $M=\bigoplus_{\rho \in \mathbb{N}} M^{(\rho)}$.

Les morphismes de $\Sy$-bimodules gradu\'es par un poids sont les
morphismes de $\Sy$-bimodules qui pr\'eservent cette
d\'ecomposition. L'ensemble des $\Sy$-bimodules
gradu\'es par un poids muni des morphismes correspondant forme
une cat\'egorie que l'on
note gr-$\Sy$-biMod.
\end{dei}

\textsc{Remarque :} On d\'efinit le produit tensoriel $\otimes_k$
dans la cat\'egorie des $\Sy$-bimodules gradu\'es par un poids par
$$\left(\Po(m,\, n)\otimes_k \Qo(m',\, n')\right)^{(\rho)} =
\bigoplus_{s+t=\rho} \Po(m,\, n)^{(s)}\otimes_k\Qo(m',\,
n')^{(t)}.$$ A l'aide de cette g\'en\'eralisation, on  peut
\'etendre les produits $\otimes$, $\boxtimes_c$ et $\boxtimes$ \`a
la cat\'egorie gr-$\Sy$-biMod.

\begin{dei}[Prop\'erade gradu\'ee par un poids]
\index{prop\'erade gradu\'ee par un poids}
 Une \emph{prop\'erade
gradu\'ee par un poids} $(\Po,\, \mu,\, \eta)$ est un mono\"\i de
dans la cat\'egorie $(\textrm{gr-}\Sy\textrm{-biMod},\,
\boxtimes,\, I)$. Lorsque les \'el\'ements de degr\'e $0$, pour
cette graduation, correspondent \`a $\Po^{(0)}=I$, on parle de
\emph{prop\'erade connexe}.\index{prop\'erade connexe}
\end{dei}

Dualement, on d\'efinit la notion de coprop\'erade.

\begin{dei}[Coprop\'erade]
\index{coprop\'erade}
Une \emph{coprop\'erade} est un comono\"\i
de dans la cat\'egorie $(\Sy\textrm{-bimod},\, \boxtimes_c,\ I)$.
\end{dei}

De mani\`ere explicite, cela signifie que $(\mathcal{C},\
\Delta,\, \varepsilon)$ est une coprop\'erade si et seulement si
$\Delta\, : \, \mathcal{C} \to \mathcal{C}\boxtimes_c \mathcal{C}$
est un morphisme de $\Sy$-bimodules coassociatif
$$\xymatrix{ \mathcal{C} \ar[d]^{\Delta} \ar[r]^{\Delta} & \mathcal{C}
\boxtimes_c\mathcal{C} \ar[d]^{\Delta \boxtimes_c \mathcal{C}}\\
\mathcal{C}\boxtimes_c\mathcal{C}\ar[r]^(0.4){\mathcal{C}
\boxtimes_c \Delta}
&\mathcal{C}\boxtimes_c\mathcal{C}\boxtimes_c\mathcal{C}}$$ et
$\varepsilon \, :\, \mathcal{C}\to I$ une counit\'ee
$$\xymatrix{ I \boxtimes_c \mathcal{C}& \mathcal{C} \boxtimes_c \mathcal{C}
\ar[l]_{\varepsilon \boxtimes_c \mathcal{C}}
\ar[r]^{\mathcal{C} \boxtimes_c \varepsilon}& \mathcal{C} \boxtimes_c I \\
 &\mathcal{C} \ar[u]^{\Delta} \ .  \ar[ur]_{\rho_\mathcal{C}^{-1}}
\ar[ul]^{\lambda_\mathcal{C}^{-1}}& } $$

\subsection{D\'efinition de PROP}

Dans le cas des prop\'erades, on ne consid\`ere qu'une
composition, la composition verticale $\mu$. Dans le cas des
PROPs, il faut en plus prendre en compte une autre composition, la
composition horizontale.

\begin{dei}[PROP]
\index{PROP}
Une structure de \emph{PROP} sur un $\Sy$-bimodule
correspond aux donn\'ees suivantes :
\begin{itemize}
\item une composition verticale associative $\Po\boxtimes \Po
\xrightarrow{\mu} \Po$,

\item une composition horizontale associative et commutative $\Po\otimes \Po
\xrightarrow{conc} \Po$,

\item une unit\'e pour la composition verticale $S_\otimes(I)
\xrightarrow{\eta} \Po.$
\end{itemize}

En outre, on impose que les deux compositions commutent,
c'est-\`a-dire qu'elles v\'erifient la relation
d'\emph{Interchange Law} \index{interchange law}

$$\xymatrix{ (\Po\boxtimes_c \Po) \otimes (\Po\boxtimes_c \Po)
\ar[dd]^{\mu \otimes \mu} \ar[r]& \left(\bigoplus_{m,\, n,\, N}
S_\otimes(\Po)(m,\, N) \otimes_{\Sy_N} k[\Sy_N] \otimes_{\Sy_N}
S_\otimes(\Po)(N,\, m)\right) \Big/_\sim
  \ar[d]^{conc\, \boxtimes_c conc}  \\
 &  \Po \boxtimes_c \Po \ar[d]^{\mu}\\
\Po \otimes \Po \ar[r]^{conc}  &  \Po.} $$
\end{dei}

\textsc{Remarque :} Le $\Sy$-bimodules $S_\otimes(I)$ est
isomorphe au $\Sy$-bimodule $\bigoplus_{n \in\mathbb{N}^*}
k[\Sy_n]$.

\begin{pro}
Cette d\'efinition de PROP est \'equivalente \`a la d\'efinition
classique donn\'ee par Lawvere et Mac Lane (\emph{cf.}
\cite{Lawvere} et \cite{MacLane2}).
\end{pro}

\begin{deo}
La seule diff\'erence entre les deux d\'efinitions vient de la composition
verticale. Dans la d\'efiniton de Mac Lane, la composition verticale est la
donn\'ee d'un morphisme de $\Sy$-bimodules associatif de la forme
$$\circ \ : \ \Po(m,\,N) \otimes_{\Sy_N} \Po(N,\, n) \to \Po(m,\, n).$$

Si on se donne un morphisme $\Po\boxtimes \Po \xrightarrow{\mu} \Po$, on
construit une composition $\circ$ par restriction
$$ \circ \ : \ \Po(m,\, N) \otimes_{\Sy_N} \Po(N,\, n)=
\Po(m,\, N) \otimes_{\Sy_N} k[\Sy_N] \otimes_{\Sy_N} \Po(N,\, n)
\to \Po \boxtimes \Po \xrightarrow{\mu} \Po. $$

R\'eciproquement, si on a une composition de la forme $\circ$, on d\'efinit
$\mu$ par
\begin{eqnarray*}
 \mu \ &:& \ \Po\boxtimes \Po(m,\,n)= \bigoplus_{N\in\mathbb{N}}
S_\otimes
(\Po)(m,\, N) \otimes_{\Sy_N} S\otimes (\Po)(N,\, n)\\
&& \xrightarrow{\bigoplus_N
conc\otimes_{\Sy_N}conc} \bigoplus_N \Po(m,\, N) \otimes_{\Sy_N} \Po(N,\, n)
\xrightarrow{\circ} \Po(m,\, n).
\end{eqnarray*}
Ces deux constructions sont inverses l'une de l'autre gr\^ace \`a
l'interchange law.$\cqfd$
\end{deo}

\textsc{Remarque :} La terminologie de PROP a \'et\'e introduite
par S. Mac Lane. Elle vient de ``PROduits et Permutations''.\\

A partir d'un PROP $\Po$, si on oublie la composition horizontale
et que l'on ne parle que de compositions verticales connexes,
alors on retrouve une structure de prop\'erade.

\begin{dei}[Le foncteur oubli $U_c$]
\index{foncteur oubli $U_c$}
 On appelle \emph{foncteur oubli}
$U_c$, le foncteur
$$U_c \ : \ \textrm{PROPs} \to \textrm{prop\'erades}$$
qui d\'efinit une prop\'erade \`a partir d'un PROP.
\end{dei}

\begin{pro}
\label{concatenationlibre} Le foncteur oubli entre les PROPs et
les prop\'erades admet un adjoint \`a gauche donn\'e par la
construction $S_\otimes$

$$\xymatrix{\textrm{PROPs \ } \ar@<0.5ex>@^{>}[r]^{U_c}   &
\ar@<0.5ex>@^{>}[l]^{S_\otimes}{\ \textrm{prop\'erades}}.} $$
\end{pro}

\begin{deo}
Soit $(\Po,\, \mu ,\, \eta)$ une prop\'erade, le morphisme de
concat\'enation $S_\otimes(\Po)\otimes S_\otimes(\Po) \to \Po$ est
celui donn\'ee par la construction de l'alg\`ebre sym\'etrique
libre. Quant \`a la composition verticale $\widetilde{\mu}$, on la
d\'efinit, gr\^ace \`a la proposition~\ref{lienentreproduits},
par
$$\widetilde{\mu} \ :\ S_\otimes(\Po)\boxtimes S_\otimes(\Po) =
S_\otimes(\Po\boxtimes_c \Po) \xrightarrow{S_\otimes(\mu)}
S_\otimes(\Po). $$ Et l'unit\'e de la prop\'erade permet d'obtenir
l'inclusion $S_\otimes(I) \xrightarrow{S_\otimes(\eta)}
S_\otimes(\Po)$. Ainsi d\'efini $(S_\otimes(\Po)$,
$\widetilde{\mu}$, $conc$, $S_\otimes(\eta))$ forme un PROP. Et on
v\'erife facilement la relation d'adjonction. $\cqfd$
\end{deo}

\subsection{L'exemple fondamental $End(V)$}

\begin{dei}[$End(V)$]
\index{$End(V)$}
 On pose $End(V) = \lbrace \textrm{applications
lin\'eaires} \ : \ V^{\otimes n} \rightarrow V^{\otimes m}
\rbrace$ o\`u les actions de $\mathbb{S}_m$ et de $\Sy_n$ se font
\`a la source $V^{\otimes n}$ et \`a l'arriv\'ee $V^{\otimes m}$
par permutation des variables. La composition
$\xymatrix{End(V)\boxtimes End(V) \ar[r]^-{\chi} & End(V) }$
revient \`a composer les applications multilin\'eaires suivant le
sch\'ema donn\'e par le produit mono\"\i dal $\boxtimes$. La
composition horizontale $End(V)\otimes End(V) \xrightarrow{conc}
End(V)$ est exactement la concat\'enation des applications
lin\'eaires. Et, l'inclusion $S_\otimes (I) \to End(V)$ correspond
aux applications lin\'eaires de $V^{\otimes n} \to V^{\otimes n}$
obtenues par permutations des variables.
\end{dei}

Par exemple, pour $p_\alpha \in Hom_{k\textrm{-Mod}}(V^{\otimes
i_\alpha},\, V^{\otimes j_\alpha})$ et $q_\beta\in
Hom_{k\textrm{-Mod}}(V^{\otimes k_\beta},\, V^{\otimes l_\beta})$,
l'application $\mu(\theta\otimes q_1 \otimes \cdots \otimes q_b
\otimes \sigma \otimes p_1\otimes\cdots \otimes p_a \otimes
\omega)\in Hom_{k\textrm{-Mod}}(V^{\otimes n},\, V^{\otimes m})$
correspond \`a la composition suivante
$$ \xymatrix@C=35pt{V^{\otimes n} \ar[r]^-{\omega}&V^{\otimes n}
\ar[r]^-{p_1\otimes \cdots \otimes p_a }&V^{\otimes N}
\ar[r]^-{\sigma}&V^{\otimes N}\ar[r]^-{q_1\otimes \cdots \otimes
q_b}& V^{\otimes m} \ar[r]^-{\theta}  &V^{\otimes m},}$$ o\`u
 l'application $p_1\otimes\cdots \otimes p_a\, : \, V^{\otimes n}=
 V^{\otimes i_1} \otimes \cdots \otimes V^{\otimes i_b} \to
 V^{\otimes j_1} \otimes \cdots \otimes V^{\otimes j_b}=
 V^{\otimes N}$ est la concat\'enation des
applications $p_\alpha$.\\

Ainsi d\'efini, $(End(V),\, \chi,\, conc,\, \eta)$ est un PROP. Et
si l'on se restreint aux compositions verticales connexes $
End(V)\boxtimes_c End(V) \xrightarrow{\chi_c} End(V) $, on obtient
une prop\'erade.

\subsection{Lien avec les foncteurs de Schur}
On a la proposition suivante
\begin{pro}
Si $\Po$ est une prop\'erade, alors le foncteur de Schur $\Sh_\Po
'$ est une monade, c'est-\`a-dire un mono\"\i de dans la
cat\'egorie $(Fonct_{\Sy\textrm{-bimod}},\, \circ,\, id)$.
\end{pro}
\begin{deo}
On se sert de la proposition~\ref{foncteurdeSchur}. $\Box$
\end{deo}

\subsection{$\Po$-g\`ebre}

La notion de g\`ebre sur une prop\'erade (ou un PROP) est la
g\'en\'eralisation naturelle de celle d'alg\`ebre sur une
op\'erade.

\begin{dei}[$\Po$-g\`ebre]
\index{$\Po$-g\`ebre} \index{g\`ebre sur une prop\'erade}

Soit $\mathcal{P}$ une prop\'erade (respectivement un PROP). Une
structure de $\mathcal{P}$\emph{-g\`ebre} sur le $k$-module $V$
est la donn\'ee d'un morphisme de prop\'erades (respectivement de
PROPs) $\mathcal{P} \rightarrow End(V)$.
\end{dei}

\textsc{Remarque :} Une $\Po$-g\`ebre est une ``alg\`ebre'' sur
une prop\'erade $\Po$. Mais le terme d'alg\`ebre est ici \`a
prendre dans un sens large. En effet, les $\Po$-g\`ebres peuvent
aussi \^etre munies de coproduits. C'est le cas des big\`ebres et des
big\`ebrse de Lie, par exemple.\\

On consid\`ere le $\Sy$-bimodule $T_*(V) (m,\,n)=V^{\otimes m}$,
o\`u l'action de $\Sy_n$ est l'action triviale et celle de $\Sy_m$
correspond \`a la permutation des variables. Un morphisme de
$\Sy$-bimodules $\mu_V\, : \, \Po\boxtimes_c V \to T_*(V)$
s'\'etend en un unique morphisme $T(\mu_V)$ de $\Po\boxtimes
T_*(V) \to T_*(V)$ par concat\'enation. En effet, tout \'el\'ement
de $\Po\boxtimes T_*(V)$ peut se voir comme la concat\'enation
d'\'el\'ements de $\Po\boxtimes_c V$. Par exemple, l'\'el\'ement
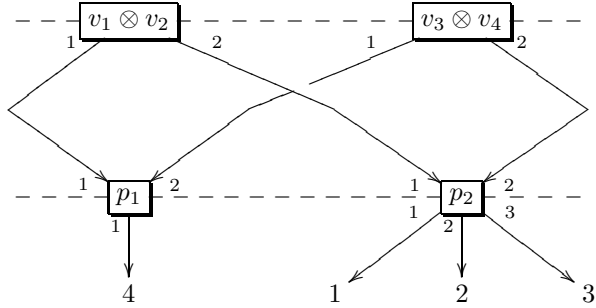
$$ \xymatrix{\ar@{--}[r]& *+[F-,]{v_1\otimes v_2} \ar@{-}[dl]_(0.4){1}
\ar@{-}[drr]^(0.4){2} \ar@{--}[rrr]& & & *+[F-,]{v_3\otimes v_4}
\ar@{-}[dr]^(0.4){2}   \ar@{-}[dll]_(0.4){1} |(0.75) \hole  \ar@{--}[r]& \\
 *=0{} \ar[dr]_(0.7){1} & &*=0{}\ar[dl]^(0.7){2} &*=0{}\ar[dr]_(0.7){1}
 & &*=0{} \ar[dl]^(0.7){2}\\
\ar@{--}[r] & *+[F-,]{p_1} \ar[d]_(0.3){1}  \ar@{--}[rrr]& & &
*+[F-,]{p_2} \ar[dl]_(0.3){1} \ar[d]_(0.3){2} \ar[dr]^(0.3){3}
\ar@{--}[r] & \\
&4 & &1 &2 &3 } $$ s'\'ecrit aussi
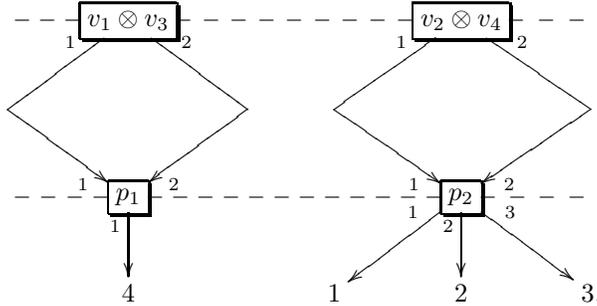
$$ \xymatrix{\ar@{--}[r]& *+[F-,]{v_1\otimes v_3} \ar@{-}[dl]_(0.4){1}
\ar@{-}[dr]^(0.4){2} \ar@{--}[rrr]& & & *+[F-,]{v_ 2\otimes v_4}
\ar@{-}[dr]^(0.4){2}   \ar@{-}[dl]_(0.4){1}  \ar@{--}[r]& \\
 *=0{} \ar[dr]_(0.7){1} & &*=0{}\ar[dl]^(0.7){2} &*=0{}\ar[dr]_(0.7){1}
 & &*=0{} \ar[dl]^(0.7){2}\\
\ar@{--}[r] & *+[F-,]{p_1} \ar[d]_(0.3){1}  \ar@{--}[rrr]& & &
*+[F-,]{p_2} \ar[dl]_(0.3){1} \ar[d]_(0.3){2} \ar[dr]^(0.3){3}
\ar@{--}[r] & \\
&4 & &1 &2 &3. } $$ \\

Par exemple, pour $\Po=End(V)$, on a naturellement une application
$\xi\, : \, End(V)\boxtimes_c V \to T_*(V)$ qui correspond \`a
l'\'evaluation d'une application de
$Hom_{k\textrm{-Mod}}(V^{\otimes n},\, V^{\otimes m})$ par un
\'el\'ement de $V^{\otimes n}$. Par d\'efinition de $\chi$, on a
imm\'ediatement le lemme suivant.

\begin{lem}
Le diagramme ci-dessous est commutatif
$$\xymatrix@C=40pt{End(V)\boxtimes_c End(V)\boxtimes_c V
\ar[r]^-{End(V)\boxtimes_c \xi}\ar[d]^-{\chi \boxtimes_c V}
& End(V) \boxtimes_c T_*(V)  \ar[d]^-{T(\xi)}\\
End(V) \boxtimes_c V \ar[r]^-{\xi}& T_*(V).} $$
\end{lem}
Ce lemme se g\'en\'eralise de la mani\`ere suivante.

\begin{pro}
Un $k$-module $V$ est une g\`ebre sur $\Po$ si et seulement si il
existe un morphisme de $\Sy$-bimodules $\mu_V\, : \,
\Po\boxtimes_c V \to T_*(V)$ tel que le diagramme suivant commute
$$\xymatrix{ \Po \boxtimes_c \Po \boxtimes_c V  \ar[r]^-{\Po
\boxtimes_c \mu_V} \ar[d]^-{\mu \boxtimes_c V} &
\Po \boxtimes_c T_ *(V)\ar[d]^-{T(\mu_V)} \\
\Po \boxtimes_c V \ar[r]^-{\mu_V} & T_*(V). } $$
\end{pro}

\textsc{Remarque :} Lorsque $\Po$ est une op\'erade, on retrouve
la notion classique d'alg\`ebre sur l'op\'erade $\Po$.

\section{Prop\'erade et PROP libres, quadratiques}

\subsection{Prop\'erade et PROP libres}

Le but de cette partie est de d\'ecrire la prop\'erade libre sur
un $\Sy$-bimodule $V$. Pour cela, on utilise le travail plus
g\'en\'eral, effectu\'e sur le mono\"\i de libre, au chapitre $1$
section~\ref{monolibre}. Puis en utilisant la
proposition~\ref{concatenationlibre}, on obtient le PROP libre sur $V$.\\
\index{prop\'erade libre}

\begin{lem}
 Pour tout couple $(A,\, B)$ de
$\Sy$-bimodules, le foncteur
$$\Phi_{A,\, B}\, :\, X \mapsto
A\boxtimes_c X \boxtimes_c B$$
 est un foncteur analytique
(\emph{cf.} Chapitre $1$ section $8$).
\end{lem}

\begin{deo}
Le $\Sy$-bimodule $A\boxtimes_c X \boxtimes_c B$ est donn\'e par
la somme directe sur les graphes \`a $3$ niveaux $\mathcal{G}_3^c$
dont les sommets du premier niveau sont indic\'es par des
\'el\'ements de $B$, ceux du deuxi\`eme niveau par des
\'el\'ements de $X$ et ceux du troisi\`eme niveau par des
\'el\'ements de $A$. Si on pose $\mathcal{G}_3^n$ l'ensemble des
graphes \`a trois niveaux ayant $n$ sommets sur le deuxi\`eme
niveau, alors le foncteur $\Phi_{A,\, B}$ s'\'ecrit
\begin{eqnarray*}
\Phi_{A,\, B}(X) &=& A\boxtimes_c X \boxtimes_c B \\
&=& \bigoplus_{n\in \mathbb{N}} \Big( \bigoplus_{g \in
\mathcal{G}^n_3} \bigotimes_{\nu \in \mathcal{N}_1}
A(|Out(\nu)|,\, |In(\nu)|) \otimes \bigotimes_{i=1}^n
X(|Out(\nu_i)|,\, |In(\nu_i)|) \otimes\\
&& \bigotimes_{\nu \in \mathcal{N}_3}
B(|Out(\nu)|,\, |In(\nu)|) \Big) \Big/_\approx  \\
&=& \bigoplus_{n\in \mathbb{N}} {\Phi_n}(X),
\end{eqnarray*}
o\`u $\Phi_n$ est un foncteur polynomial homog\`ene de degr\'e
$n$.$\cqfd$
\end{deo}

\begin{pro}
La partie multilin\'eaire en $Y$ not\'ee $ A\boxtimes_c (X\oplus
\underline{Y}) \boxtimes_c B$ correspond au sous-$\Sy$-bimodule de
$A\boxtimes_c (X\oplus Y) \boxtimes_c B $ d\'efini par la somme
directe sur les graphes connexes \`a $3$ niveaux dont les sommets
du deuxi\`eme niveau sont indic\'es par des \'el\'ements de $X$ et
de $Y$ mais avec au moins un \'el\'ement de $Y$.
\end{pro}

\begin{lem}
\label{SbiModMultilineaire} La cat\'egorie $(\Sy\textrm{-biMod},\,
\boxtimes_c,\, I)$ est une cat\'egorie mono\"\i dale ab\'elienne
qui pr\'eserve les co\'egalisateurs r\'eflexifs ainsi que les
colimites s\'equentielles.
\end{lem}

\begin{deo}
Pout tout $\Sy$-bimodule $A$, les foncteurs de multiplication \`a
gauche et \`a droite $L_A$ et $R_A$ par $A$ sont des foncteurs
analytiques par le lemme pr\'ec\'edent. Ils pr\'eservent donc les
colimites s\'equentielles et les co\'egalisateurs r\'eflexifs par
la proposition~\ref{foncteuranalytiquereflexif}. $\cqfd$
\end{deo}

Cette proposition permet d'appliquer les r\'esultats sur le
mono\"\i de libre de la partie~\ref{monolibre} du chapitre $1$.
Donnons les interpr\'etations en termes de $\Sy$-bimodules de
cette partie. Soit $V$ un $\Sy$-bimodule que l'on augmente en
posant $V_+=I\oplus V$. Ensuite, $V_n=(V_+)^{\boxtimes_c n}$
correspond aux graphes connexes \`a $n$ \'etages dont les noeuds
(sommets) sont indic\'es par des \'el\'ements de $V$ et $I$.

Le $\Sy$-bimodule $\widetilde{V_n}=\coker \left( \bigoplus_i
R_{V_i,\, V_{n-i-2}}\to V_n \right)$ correspond aux graphes \`a
$n$ niveaux quotient\'es par la relation engendr\'ee par
$V\boxtimes_c I \sim I \boxtimes_c V$, ce qui revient \`a oublier
les niveaux.

\begin{thm}
\label{proplibre}
 La prop\'erade libre sur $V$, not\'ee $\F(V)$, est la
somme directe sur l'ensemble des graphes (sans niveau) connexes
$\mathcal{G}^c$ o\`u l'on indice les sommets par $V$, soit
$$\F(V)=\left( \bigoplus_{g\in \mathcal{G}^c} \bigotimes_{\nu \in
\mathcal{N}} V(|Out(\nu)|,\, |In(\nu)|)\right) \Bigg/_\approx \
.$$ Quant \`a la composition $\mu$, elle vient de la composition
des graphes dirig\'es.
\end{thm}

\textsc{Remarque : } Dans le cas o\`u $V$ est un $\Sy$-module, on
obtient la construction donn\'ee dans \cite{GK} \`a l'aide des
arbres.\\

Comme annonc\'e dans l'introduction, en concat\'enant ensuite toutes
les op\'erations de $\F(V)$, on trouve le PROP libre.

\begin{cor}
Le $\Sy$-bimodule $S_\otimes(\left( \F(V)\right)$ est le PROP
libre sur $V$.
\end{cor}

\begin{deo}
On utilise la proposition~\ref{concatenationlibre} pour montrer,
par composition, que le foncteur $S_\otimes(\F)$ est un adjoint
\`a gauche de $U\circ U_c$
$$\xymatrix{\textrm{PROPs \ } \ar@<0.5ex>@^{>}[r]^{U_c}   &
\ar@<0.5ex>@^{>}[l]^{S_\otimes}{\ \textrm{prop\'erades}}
\ar@<0.5ex>@^{>}[r]^{U}  &
\ar@<0.5ex>@^{>}[l]^{\F}\Sy\textrm{-biMod}.}$$ $\cqfd$
\end{deo}

\textsc{Remarque :} On retrouve bien la m\^eme construction que
celle du PROP libre exprim\'ee avec des graphes (non
n\'ecessairement connexes, sans niveau) de B. Enriquez et P.
Etingof
(\emph{cf.} \cite{EE}).\index{PROP libre}\\

On peut d\'ecomposer l'ensemble des graphes en fonction du nombre
total de sommets. Soit $\mathcal{G}^n$ l'ensemble des graphes \`a
$n$ sommets. Cette d\'efinition permet d'affiner l'\'ecriture de
la prop\'erade libre (respectivement du PROP libre). En effet, on
a
$$ \F(V)=\left( \bigoplus_{g\in \mathcal{G}_c}
\bigotimes_{\nu \in \mathcal{N}} V(|Out(\nu)|,\, |In(\nu)|)
\right) \Bigg/_\approx= \bigoplus_{n\in \mathbb{N}} \underbrace{
\left( \bigoplus_{g\in \mathcal{G}_c^{(n)}} \bigotimes_{i=1}^n
V(|Out(\nu_i)|,\, |In(\nu_i)|)\right)
\Bigg/_\approx}_{\F_{(n)}(V)}.$$

\begin{pro}
\label{Fanalytique} Le foncteur $\F \, : \, V \mapsto \F(V)$ est
un foncteur analytique en $V$, dont la partie de degr\'e $n$ est
not\'ee $\F_{(n)}(V)$. Cette graduation est stable pour la
composition de la prop\'erade libre $\mu$ (respectivement pour les
compositions $\widetilde{\mu}$ et $\textit{conc}$ du PROP libre). Ainsi,
toute prop\'erade libre (respectivement PROP libre) est gradu\'ee
par un poids.
\end{pro}

\subsection{Coprop\'erade colibre connexe}

\index{coprop\'erade colibre connexe}
Le $\Sy$-bimodule sur lequel
on d\'efinit une structure de coprop\'erade colibre connexe est le
m\^eme que pour la prop\'erade libre. On pose
$$\F^c(V)= \F(V)=\left( \bigoplus_{g\in \mathcal{G}^c}
\bigotimes_{\nu \in \mathcal{N}} V(|Out(\nu)|,\, |In(\nu)|)\right)
\Bigg/_\approx.$$

La projection sur la composante engendr\'ee par le graphe
trivial $\ \vcenter{\xymatrix@M=0pt@R=10pt@C=6pt{ {}\ar@{-}[d] \\
{}}} \ $ fournit la counit\'e $\varepsilon \, : \, \F^c(V) \to I$.

Et la comultiplication repose sur l'ensemble des d\'ecoupages en
deux des graphes. Sur un \'el\'ement $g(V)$ qui repr\'esente un
graphe $g$ indic\'e par des op\'erations de $V$, on d\'efinit le
morphisme $\Delta$ par
$$\Delta(g(V))=\sum_{\left (g_1(V),\, g_2(V) \right)} g_1(V) \boxtimes_c g_2(V),$$
o\`u la somme court sur les couples $(g_1(V),\, g_2(V))$ tels que
$\mu(g_1(V) \boxtimes_c g_2(V))=g(V)$. Il faut faire attention ici
que $g_1(V)$ (et $g_2(V)$) repr\'esente une famille de graphes
indic\'es.

\begin{pro}
Pour tout $\Sy$-bimodule $V$, $(\F^c(V),\, \Delta,\, \varepsilon)$
est une coprop\'erade gradu\'ee par un poids. Cette coprop\'erade
est colibre dans le cadre des coprop\'erades connexes.
\end{pro}

\begin{deo}
La relation de counit\'e vient de
\begin{eqnarray*}
(\varepsilon \boxtimes_c id)\circ \Delta \left( g(V) \right) &=&
(\varepsilon \boxtimes_c id)\circ \Delta \left(\sum_{\left
(g_1(V),\, g_2(V) \right)} g_1(V) \boxtimes_c g_2(V)\right) \\
&=& \sum_{\left (g_1(V),\, g_2(V) \right)} \varepsilon(g_1(V))
\boxtimes_c g_2(V) \\
&=& I \boxtimes_c g(V)\\
 &=& g(V).
\end{eqnarray*}
La relation de coassociativit\'e vient de
\begin{eqnarray*}
&& (\Delta\boxtimes_c id)\circ \Delta (g(V)) =(id \boxtimes_c
\Delta)\circ \Delta (g(V)) = \\
&& \sum_{(g_1(V),\ g_2(V),\, g_3(V))} g_1(V)\boxtimes_c
g_2(V)\boxtimes_c g_3(V),
\end{eqnarray*}
o\`u la somme court sur les triplets $(g_1(V),\ g_2(V),\, g_3(V))$
tels que $\mu(g_1(V) \boxtimes_c g_2(V) \boxtimes_c g_3(V))=g(V)$.

Soient $C$ une coprop\'erade connexe et $f\, : \, C \to V$ un
morphisme de $\Sy$-bimodules. La d\'ecomposition analytique de
$\F^c(V)$, en fonction du nombre de sommets des graphes
utilis\'es, permet par r\'e\-cu\-rren\-ce de d\'efinir un
morphisme de coprop\'erades connexes $\bar{f} \, : \, C \to
\F^c(V)$. Ce morphisme est enti\`erement d\'etermin\'e par $f$ et
est l'unique morphimes de coprop\'erades connexes \`a faire
commuter le diagramme
$$\xymatrix{V & \ar[l] \F^c(V) \\
& \ar[ul]^{f} \quad C.  \ar[u]_{\bar{f}}}$$ $\cqfd$
\end{deo}

\textsc{Remarque :} Cette construction g\'en\'eralise la
construction de la cog\`ebre colibre connexe donn\'ee par T. Fox
dans \cite{F}.

\subsection{Prop\'erades d\'efinies par g\'en\'erateurs et relations, prop\'erades quadratiques}
\label{exemples}
 Continuons l'interpr\'etation des notions
mono\"\i dales dans le cadre des $\Sy$-bimodules. Soit $R$ un
sous-$\Sy$-bimodule de $\F(V)$. La
proposition~\ref{SbiModMultilineaire} permet d'appliquer les
r\'esultats montr\'es au chapitre $1$ qui affirment
 que l'id\'eal engendr\'e par $R$ correspond \`a la somme
sur l'ensemble des graphes, o\`u l'on indice les sommets avec des
\'el\'ements de $V$, et qui admet au moins un sous-graphe
appartenant \`a $R$.

\begin{pro}
La prop\'erade quotient $\F(V)/(R)$ correspond \`a la somme
directe des graphes indic\'es par des \'el\'ements de $V$ mais
dont aucun sous-graphe n'appartient \`a $R$.
\end{pro}

\begin{deo}
Il suffit de voir que $R$, comme sous-objet de $\F(V)$, se repr\'esente
avec des sommes de graphes connexes et que la substitution d'un sous-graphe
connexe par un autre sous-graphe connexe dans un graphe connexe donne toujours
un graphe connexe. $\Box$
\end{deo}

\textsc{Remarque :} Dans le cas o\`u $V$ est un $\Sy$-module, les
seuls graphes qui interviennent sont des arbres. On retrouve alors
les constructions op\'eradiques donn\'ees dans \cite{GK}.\\

Comme nous l'avons vu \`a la proposition~\ref{Fanalytique}, le
foncteur $\F(V)$ est analytique. Ceci permet de distinguer une
classe particuli\`erement int\'eressante de prop\'erades
d\'efinies par g\'en\'erateurs et relations.

\begin{dei}[Prop\'erades quadratiques]
\index{prop\'erade quadratique}
 On appelle \emph{prop\'erade
quadratique} une prop\'erade de la forme $\F(V)/(R)$, o\`u $R$ est
un sous-$\Sy$-bimodule de $\F_{(2)}(V)$ (partie de degr\'e $2$ de
la prop\'erade libre sur $V$).
\end{dei}

\textsc{Exemples :} Les premiers exemples viennent encore une fois
des sous-cat\'egories pleines $k$-Mod et $\Sy$-Mod. Les alg\`ebres
quadratiques, comme les alg\`ebres sym\'etriques $S(V)$ et
ext\'erieures $\Lambda(V)$, sont des prop\'erades quadratiques. Il
en va de m\^eme pour les op\'erades quadratiques, citons les
op\'erades des alg\`ebres associatives $\A$, commutatives $\C$,
des alg\`ebres de Lie $\Li$, etc ... \\

\textsc{Exemples :} Parmi les exemples nouveaux que cette
th\'eorie permet de traiter, citons
\begin{itemize}
\item \label{BiLie} \index{big\`ebres de Lie}La prop\'erade $\BLi$
codant les big\`ebres de Lie. Cette prop\'erade est engendr\'ee
par le $\Sy$-bimodule
$$V=\left\{ \begin{array}{l}
V(1,\, 2)=\lambda.k\otimes sgn_{\Sy_2}, \\
V(2,\, 1)=\Delta.sgn_{\Sy_2}\otimes k, \\
V(m,\, n) = 0 \quad \textrm{sinon,}
\end{array}
\right.= \vcenter{\xymatrix@M=0pt@R=6pt@C=4pt{{\scriptstyle 1} & &
{\scriptstyle 2}
\\\ar@{-}[dr] &
&\ar@{-}[dl]  \\  &\ar@{-}[d] & \\
& & \\ & & }}\otimes sgn_{\Sy_2} \oplus
\vcenter{\xymatrix@M=0pt@R=6pt@C=4pt{  & \ar@{-}[d]&  \\
&\ar@{-}[dl] \ar@{-}[dr]& \\ & & \\{\scriptstyle 1}&
&{\scriptstyle 2} }} \otimes sgn_{\Sy_2},$$
o\`u $sgn_{\Sy_2}$ est
la repr\'esentation signature de $\Sy_2$. Elle est soumise aux
relations quadratiques (c'est-\`a-dire s'\'ecrivant avec deux
(co)op\'erations)
\begin{eqnarray*}
R&=&\left\{
\begin{array}{l}
\lambda (\lambda,\, 1) \big( (123)+(231)+(312)  \big) \\
\oplus \ \big((123)+(231)+(312) \big)(\Delta,\, 1)\Delta  \\
\oplus \ \Delta \otimes \lambda - (\lambda,\, 1)\otimes(213)\otimes(1,\, \Delta) \\
\quad -(1,\, \lambda)(\Delta,\, 1) - (1,\,\lambda
)(\Delta,\,1)-(1,\, \lambda )\otimes(132)\otimes(\Delta,1)
\end{array}
\right. \\
&=&\left\{ \begin{array}{l}
\vcenter{\xymatrix@M=0pt@R=6pt@C=4pt{{\scriptstyle 1} & &
{\scriptstyle 2}& & {\scriptstyle 3}\\
\ar@{-}[dr] &
&\ar@{-}[dl] & & \ar@{-}[dl]  \\
& \ar@{-}[dr] & &\ar@{-}[dl]  & \\
& &\ar@{-}[d] & & \\
& & \\ & & }} + \vcenter{\xymatrix@M=0pt@R=6pt@C=4pt{{\scriptstyle
2} & &
{\scriptstyle 3}& & {\scriptstyle 1}\\
\ar@{-}[dr] &
&\ar@{-}[dl] & & \ar@{-}[dl]  \\
& \ar@{-}[dr] & &\ar@{-}[dl]  & \\
& &\ar@{-}[d] & & \\
& & \\ & & }} +\vcenter{\xymatrix@M=0pt@R=6pt@C=4pt{{\scriptstyle
3} & &
{\scriptstyle 1}& & {\scriptstyle 2}\\
\ar@{-}[dr] &
&\ar@{-}[dl] & & \ar@{-}[dl]  \\
& \ar@{-}[dr] & &\ar@{-}[dl]  & \\
& &\ar@{-}[d] & & \\
& & \\ & & }} \\
\oplus \vcenter{\xymatrix@M=0pt@R=6pt@C=4pt{ & & & & \\
 & &\ar@{-}[d] & & \\
& & \ar@{-}[dl]  \ar@{-}[dr] & & \\
& \ar@{-}[dl]  \ar@{-}[dr] & &\ar@{-}[dr] & \\
& & & &\\ {\scriptstyle 1} & & {\scriptstyle 2}& &{\scriptstyle 3}
}} +
\vcenter{\xymatrix@M=0pt@R=6pt@C=4pt{ & & & & \\
 & &\ar@{-}[d] & & \\
& & \ar@{-}[dl]  \ar@{-}[dr] & & \\
& \ar@{-}[dl]  \ar@{-}[dr] & &\ar@{-}[dr] & \\
& & & &\\ {\scriptstyle 2} & & {\scriptstyle 3}& &{\scriptstyle 1}
}} + \vcenter{\xymatrix@M=0pt@R=6pt@C=4pt{ & & & & \\
 & &\ar@{-}[d] & & \\
& & \ar@{-}[dl]  \ar@{-}[dr] & & \\
& \ar@{-}[dl]  \ar@{-}[dr] & &\ar@{-}[dr] & \\
& & & &\\ {\scriptstyle 3} & & {\scriptstyle 1}& &{\scriptstyle
2}}} \\
\ \\
\oplus \vcenter{\xymatrix@M=0pt@R=6pt@C=4pt{{\scriptstyle 1} &
& {\scriptstyle 2}
\\\ar@{-}[dr] &  &\ar@{-}[dl]  \\  &\ar@{-}[d] & \\
& \ar@{-}[dl]\ar@{-}[dr]& \\ & & \\ {\scriptstyle 1} & &
{\scriptstyle 2}}} -
 \vcenter{\xymatrix@M=0pt@R=6pt@C=2pt{ & {\scriptstyle 1} & &{\scriptstyle 2}
 \\ & \ar@{-}[d] &  &\ar@{-}[d]  \\  & \ar@{-}[dl]
 \ar@{-}[dr] &  & \ar@{-}[dl]\\
\ar@{-}[d]&  & \ar@{-}[d] & \\ & & & \\ {\scriptstyle 1} & &
{\scriptstyle 2} &}} +
 \vcenter{\xymatrix@M=0pt@R=6pt@C=2pt{ & {\scriptstyle 2} & & {\scriptstyle 1}
  \\ & \ar@{-}[d] &  &\ar@{-}[d]  \\  & \ar@{-}[dl]
 \ar@{-}[dr] &  & \ar@{-}[dl]\\
\ar@{-}[d]&  & \ar@{-}[d] & \\ & & & \\ {\scriptstyle 1} & &
{\scriptstyle 2} &}} -
\vcenter{\xymatrix@M=0pt@R=6pt@C=2pt{{\scriptstyle 1} & & {\scriptstyle 2} & \\
\ar@{-}[d] & &\ar@{-}[d] &  \\ \ar@{-}[dr] &
& \ar@{-}[dl] \ar@{-}[dr]& \\
& \ar@{-}[d]&  & \ar@{-}[d] \\ & & & \\ &{\scriptstyle 1} &
&{\scriptstyle 2}}}+
\vcenter{\xymatrix@M=0pt@R=6pt@C=2pt{{\scriptstyle 2} & &
{\scriptstyle 1} &
\\ \ar@{-}[d] & &\ar@{-}[d] &  \\ \ar@{-}[dr] &
& \ar@{-}[dl] \ar@{-}[dr]& \\
& \ar@{-}[d]&  & \ar@{-}[d] \\ & & & \\ & {\scriptstyle 1} & &
{\scriptstyle 2}}}\ .
\end{array}
\right. \end{eqnarray*}

Les $\BLi$-g\`ebres correspondent aux big\`ebres de Lie de
Drinfeld (\emph{cf.} \cite{Drinfeld}).

\item \label{InfBi} \index{big\`ebres de Hopf infinit\'esimales}
La
prop\'erade $\IBi$ codant les big\`ebres de Hopf
infinit\'esimales. Cette prop\'erade est engendr\'ee par
$$V=\left\{ \begin{array}{l}
V(1,\, 2)=m.k\otimes k[\Sy_2], \\
V(2,\, 1)=\Delta.k[\Sy_2]\otimes k, \\
V(m,\, n) = 0 \quad \textrm{sinon,}
\end{array}
\right.  = \vcenter{\xymatrix@M=0pt@R=6pt@C=6pt{\ar@{-}[dr] &
&\ar@{-}[dl]  \\  &\ar@{-}[d] & \\
& & \\ & & }} \oplus
\vcenter{\xymatrix@M=0pt@R=6pt@C=6pt{  & \ar@{-}[d]&  \\
&\ar@{-}[dl] \ar@{-}[dr]& \\ & &\\&&}}\ ,$$. Elle est soumise aux
relations
\begin{eqnarray*}
R&=&\left\{ \begin{array}{l}
m(m,\, 1)-m(1,\, m) \\
\oplus \ (\Delta,\,  1)\Delta - (1,\Delta)\Delta  \\
\oplus \ \Delta \otimes m - (m,\, 1)(1,\, \Delta) - (1,\,
m)(\Delta,\, 1).
\end{array}
\right. \\
&=& \left\{ \begin{array}{l} \vcenter{\xymatrix@M=0pt@R=6pt@C=6pt{
& & & & \\ \ar@{-}[ddrr] &
&\ar@{-}[dl] & & \ar@{-}[ddll]  \\
&  & &  & \\
& &\ar@{-}[d] & & \\
& & }} - \vcenter{\xymatrix@M=0pt@R=6pt@C=6pt{& & & & \\
\ar@{-}[dr] &
&\ar@{-}[dr] & & \ar@{-}[dl]  \\
& \ar@{-}[dr] & &\ar@{-}[dl]  & \\
& &\ar@{-}[d] & & \\
& & }} \\
\oplus \vcenter{\xymatrix@M=0pt@R=6pt@C=6pt{ & & & & \\
 & &\ar@{-}[d] & & \\
& & \ar@{-}[dl]  \ar@{-}[dr] & & \\
& \ar@{-}[dl]  \ar@{-}[dr] & &\ar@{-}[dr] & \\
& & & &}} -
\vcenter{\xymatrix@M=0pt@R=6pt@C=6pt{ & & & & \\
 & &\ar@{-}[d] & & \\
& & \ar@{-}[dl]  \ar@{-}[dr] & & \\
& \ar@{-}[dl]   & &\ar@{-}[dl] \ar@{-}[dr] & \\
& & & & }}  \\
\ \\
\oplus \vcenter{\xymatrix@M=0pt@R=6pt@C=6pt{\ar@{-}[dr] &  &\ar@{-}[dl]  \\  &\ar@{-}[d] & \\
& \ar@{-}[dl]\ar@{-}[dr]& \\ & & }} -
 \vcenter{\xymatrix@M=0pt@R=6pt@C=6pt{  & \ar@{-}[d] &  &\ar@{-}[d]  \\  & \ar@{-}[dl]
 \ar@{-}[dr] &  & \ar@{-}[dl]\\
\ar@{-}[d]&  & \ar@{-}[d] & \\ & & & }} -
\vcenter{\xymatrix@M=0pt@R=6pt@C=6pt{ \ar@{-}[d] & &\ar@{-}[d] &
\\ \ar@{-}[dr] &
& \ar@{-}[dl] \ar@{-}[dr]& \\
& \ar@{-}[d]&  & \ar@{-}[d] \\ & & & }}\ .

\end{array} \right. \end{eqnarray*}
 Les $\IBi$-g\`ebres sont des analogues
associatifs des big\`ebres de Lie. Elles correspondent \`a ce que
M. Aguiar appellent les big\`ebres de Hopf infinit\'esimales
(\emph{cf.} \cite{Ag1}, \cite{Ag2} et \cite{Ag3}).

\item La prop\'erade $\frac{1}{2}\Bi$, analogue d\'eg\'en\'er\'e
de la prop\'erade des big\`ebres. L'espace g\'en\'erateur $V$ est
le m\^eme que celui de $\IBi $, c'est-\`a-dire constitu\'e d'une
op\'eration et d'une coop\'eration. Quant \`a celui des relations
$R$, il vaut
\begin{eqnarray*}
R&=&\left\{ \begin{array}{l}
m(m,\, 1)-m(1,\, m) \\
\oplus \ (\Delta,\,  1)\Delta - (1,\Delta)\Delta  \\
\oplus \ \Delta \otimes m \qquad \qquad \qquad (\kappa).
\end{array}
\right. \\
&=& \left\{ \begin{array}{l} \vcenter{\xymatrix@M=0pt@R=6pt@C=6pt{
& & & & \\ \ar@{-}[ddrr] &
&\ar@{-}[dl] & & \ar@{-}[ddll]  \\
&  & &  & \\
& &\ar@{-}[d] & & \\
& & }} - \vcenter{\xymatrix@M=0pt@R=6pt@C=6pt{& & & & \\
\ar@{-}[dr] &
&\ar@{-}[dr] & & \ar@{-}[dl]  \\
& \ar@{-}[dr] & &\ar@{-}[dl]  & \\
& &\ar@{-}[d] & & \\
& & }} \\
\oplus \vcenter{\xymatrix@M=0pt@R=6pt@C=6pt{ & & & & \\
 & &\ar@{-}[d] & & \\
& & \ar@{-}[dl]  \ar@{-}[dr] & & \\
& \ar@{-}[dl]  \ar@{-}[dr] & &\ar@{-}[dr] & \\
& & & &}} -
\vcenter{\xymatrix@M=0pt@R=6pt@C=6pt{ & & & & \\
 & &\ar@{-}[d] & & \\
& & \ar@{-}[dl]  \ar@{-}[dr] & & \\
& \ar@{-}[dl]   & &\ar@{-}[dl] \ar@{-}[dr] & \\
& & & & }}  \\
\ \\
\oplus \vcenter{\xymatrix@M=0pt@R=6pt@C=6pt{\ar@{-}[dr] &  &\ar@{-}[dl]  \\  &\ar@{-}[d] & \\
& \ar@{-}[dl]\ar@{-}[dr]& \\ & & }} \qquad \qquad \qquad
(\kappa).\end{array} \right.
\end{eqnarray*}
Cette exemple
a \'et\'e introduit par M. Markl dans \cite{Markl3} afin de
trouver le mod\`ele minimal pour le PROP des big\`ebres.
\end{itemize}

\textsc{Contre-exemples :}
$ $
\begin{itemize}

\item \index{big\`ebres} La prop\'erade pr\'ec\'edente est un cas
quadratique issu de la prop\'erade $\Bi$ des big\`ebres. Cette
derni\`ere est bien d\'efinie par g\'en\'erateurs et relations
mais n'est pas quadratique. En effet, la d\'efinition de la
prop\'erade $\Bi$ est la m\^eme que celle de $\frac{1}{2}\Bi$ sauf
que la relation $\kappa$ est remplac\'ee par
\begin{eqnarray*}
\kappa'\, &:& \, \Delta \otimes m - (m,\,m)\otimes (1324)\otimes
(\Delta,\,\Delta)\\ &=&
\vcenter{\xymatrix@M=0pt@R=6pt@C=6pt{\ar@{-}[dr] &  &\ar@{-}[dl]  \\
&\ar@{-}[d] & \\
& \ar@{-}[dl]\ar@{-}[dr]& \\ & & }} -
\vcenter{\xymatrix@R=6pt@C=2pt{ &\ar@{-}[d] & &\ar@{-}[d] & \\
*=0{}& *=0{} \ar@{-}[dl] \ar@{-}[ddrr] |\hole & *=0{}&*=0{}
\ar@{-}[ddll]\ar@{-}[dr] &*=0{} \\
*=0{}\ar@{-}[dr] & *=0{}&*=0{} &*=0{} &*=0{}\ar@{-}[dl] \\
*=0{}&*=0{} \ar@{-}[d]&*=0{} &*=0{} \ar@{-}[d]&*=0{} \\
& & & & \\}}
\end{eqnarray*}
 Cette prop\'erade n'est ni quadratique (on
utilise $4$ op\'erations pour \'ecrire le second \'el\'ement), ni
homog\`ene.

\item \index{big\`ebres de Hopf infinit\'esimales} Un autre
exemple qui ne rentre pas stricto sensu dans la th\'eorie des
prop\'erades est donn\'ee par les big\`ebres de Hopf
infinit\'esimales \emph{unitaires} (\emph{cf.} \cite{Loday4}).
Leur d\'efinition est la m\^eme que celle des big\`ebres de Hopf
infinit\'esimales, \`a la seule diff\'erence pr\`es que la
troisi\`eme relation est
$$\Delta
\otimes m - (m,\, 1)(1,\, \Delta) - (1,\, m)(\Delta,\, 1)- (1,\,
1)= \vcenter{\xymatrix@M=0pt@R=6pt@C=6pt{\ar@{-}[dr] &  &\ar@{-}[dl]  \\  &\ar@{-}[d] & \\
& \ar@{-}[dl]\ar@{-}[dr]& \\ & & }} -
 \vcenter{\xymatrix@M=0pt@R=6pt@C=6pt{  & \ar@{-}[d] &  &\ar@{-}[d]  \\  & \ar@{-}[dl]
 \ar@{-}[dr] &  & \ar@{-}[dl]\\
\ar@{-}[d]&  & \ar@{-}[d] & \\ & & & }} -
\vcenter{\xymatrix@M=0pt@R=6pt@C=6pt{ \ar@{-}[d] & &\ar@{-}[d] &
\\ \ar@{-}[dr] &
& \ar@{-}[dl] \ar@{-}[dr]& \\ & \ar@{-}[d]&  & \ar@{-}[d] \\ & &
& }} -
\vcenter{\xymatrix@M=0pt@R=6pt@C=6pt{ \ar@{-}[d] \\
\ar@{-}[d] \\ \ar@{-}[d] \\{} }}
\ \  \vcenter{\xymatrix@M=0pt@R=6pt@C=6pt{ \ar@{-}[d] \\
 \ar@{-}[d] \\ \ar@{-}[d] \\ {} }} \ .$$ Cette relation
n'\'etant pas connexe, on a l\`a un exemple de PROP qui n'est pas
une prop\'erade.
\end{itemize}

Le fait que les relations soient homog\`enes, c'est-\`a-dire,
toutes de m\^eme degr\'e, permet de conserver la graduation de la
prop\'erade libre par passage au quotient.

\begin{pro}
Soit $\F(V)/(R)$ une prop\'erade d\'efinie par g\'en\'erateurs et
relations. Si $R \subset \F_{(n)}(V)$ pour un certain $n$, alors
la prop\'erade $\F(V)/(R)$ est gradu\'ee par un poids donn\'e par
le nombre de sommets.
\end{pro}

\begin{cor}
Toute prop\'erade quadratique est gradu\'ee en
fonction du nombre d'op\'erations.
\end{cor}

\subsection{PROPs d\'efinis par g\'en\'erateurs et relations connexes}

Au niveau des PROPs, on consid\`ere le PROP libre $S_\otimes(\F(V))$
sur un $\Sy$-bimodule $V$ soumis \`a des relations $R$. Si les relations
ne font pas intervenir la concat\'enation des op\'erations, c'est-\`a-dire
si $R$ appartient \`a $\F(V)$, le PROP quotient
$S_\otimes(\F(V))/(R)_{S_\otimes(\F(V))}$ correspond
\`a l'alg\`ebre sym\'etrique libre
pour le produit $\otimes$ sur la prop\'erade quotient $\F(V)/(R)_{\F(V)}$

\begin{pro}
\label{Lambdapropéradequotient}
Soit $V$ un $\Sy$-bimodule et soit $R$ un sous-$\Sy$-bimodule de $\F(V)$.
Le PROP quotient $S_\otimes(\F(V))/(R)_{S_\otimes(\F(V))}$ est
naturellement
isomorphe au PROP  $S_\otimes(\F(V)/(R)_{\F(V)})$ construit \`a partir
de la prop\'erade quotient $\F(V)/(R)_{\F(V)}$.
\end{pro}

\begin{deo}
On part du morphisme naturel
$$\phi \ :  \ S_\otimes(\F(V)) \twoheadrightarrow
S_\otimes(\F(V)/(R)_{\F(V)}).$$ Il admet pour noyau l'id\'eal
libre sur $R$ dans le PROP libre $S_\otimes(\F(V))$. Cette
application induit donc un isomorphisme $\bar{\phi}$ de
$\Sy$-bimodules. On conclut en remarquant que $\phi$ est un
morphisme de PROP et qu'il en est de m\^eme pour $\bar{\phi}$.
$\cqfd$
\end{deo}

Cette proposition permet de g\'en\'eraliser aux PROPs quotients tous les
r\'esultats obtenus au niveau des prop\'erades.

\begin{pro}
Soit $S_\otimes(\F(V))/(R)$ un PROP d\'efini par g\'en\'erateurs et
relations,
o\`u $R$ est un sous-$\Sy$-bimodules de $\F_{(n)}(V)$.
Alors le PROP $\ S_\otimes(\F(V))/(R)$ est gradu\'e par le nombre de
sommets
des graphes.
\end{pro}

\begin{dei}[PROPs quadratiques]
On appelle \emph{PROP quadratique} tout PROP de la forme
$S_\otimes(\F(V))/(R)$ avec $R \subset \F_{(2)}(V)$.
\end{dei}

\begin{cor}
Tout PROP quadratique est gradu\'e par un poids.
\end{cor}

\chapter{Prop\'erades et PROPs diff\'erentiels}

\thispagestyle{empty}

Nous avons, pour l'instant, d\'efini les prop\'erades et les PROPS
en partant de la cat\'egorie mono\"\i dale sym\'etrique des
$k$-modules, puis en consid\'erant les $k$-modules qui admettent
une action de $\Sy_m$ \`a gauche et de $\Sy_n$ \`a droite. Cette
d\'emarche est g\'en\'eralisable \`a toute cat\'egorie mono\"\i
dale sym\'etrique. On peut donc d\'efinir des prop\'erades
ensemblistes, des prop\'erades topologiques, etc ... (Cette
d\'efinition est une version prop\'eradique de celle donn\'ee par
J.-P. May dans \cite{May} pour les op\'erades). Le but du
pr\'esent chapitre est de d\'efinir des prop\'erades et les PROPS
non plus lin\'eaires mais diff\'erentiels, c'est-\`a-dire \`a
valeurs dans la cat\'egorie des complexes de cha\^\i nes (modules
diff\'erentiels gradu\'es), et d'\'etudier les propri\'et\'es
homologiques de tels objets.

\section{La cat\'egorie des $\Sy$-bimodules diff\'erentiels gradu\'es}

\subsection{Les $\Sy$-bimodules diff\'erentiels gradu\'es et
les produits $\boxtimes_c$ et $\boxtimes$}

On se place maintenant dans la cat\'egorie des $k$-modules
diff\'erentiels gradu\'es, not\'ee dg-Mod. Rappelons que l'on
munit cette cat\'egorie d'une structure mono\"\i dale en
d\'efinissant le produit $V\otimes_k W$ par
$$(V\otimes_k W)_d=\bigoplus_{i+j=d}V_i\otimes_k W_j, $$
et en posant pour diff\'erentielle sur $V\otimes_k W$
$$\delta(v\otimes w)=\delta(v)\otimes w  + (-1)^{|v|}v\otimes \delta(w),$$
o\`u $v\otimes w$ est un tenseur \'el\'ementaire. Pour d\'ecrire
les isomorphismes de sym\'etrie, on utilise les r\`egles de signes
de Koszul-Quillen, \`a savoir
$$\tau_{v,\, w}\ : \ v\otimes w \mapsto (-1)^{|v||w|}w\otimes v.$$
Plus g\'en\'eralement, la commutation de deux \'el\'ements
(morphismes, diff\'erentielles, objets, etc ...) de degr\'e
homologique $d$ et $e$ entra\^\i ne un signe $(-1)^{de}$.
\index{r\`egles de signes de Koszul-Quillen}

\begin{dei}[dg-$\Sy$-bimodule]
\index{dg-$\Sy$-bimodule} Un \emph{dg-$\Sy$-bimodule} $\Po$ est
une collection de modules diff\'erentiels gradu\'es $(\Po(m,\,
n))_{m,\, n \in \mathbb{N}}$ munis d'une action de $\Sy_m$ \`a
gauche et de $\Sy_n$ \`a droite, ces deux actions commutent entre
elles et agissent comme des morphismes de dg-modules. On note
$\Po_d(m,\, n)$, le sous-$(\Sy_m,\, \Sy_n)$-bimodule compos\'e des
\'el\'ements de degr\'e homologique $d$ du complexe $\Po(m,\, n)$.
\end{dei}

Le produit mono\"\i dal $\boxtimes_c$ sur les $\Sy$-bimodules
s'\'etend au cas diff\'erentiel. La d\'efinition reste la m\^eme,
\`a la seule diff\'erence que la relation $\sim$, bas\'ee sur les
isomorphismes de sym\'etrie, fait maintenant intervenir des signes
(r\`egle de Koszul-Quillen). Par exemple, on a
\begin{eqnarray*}
\theta \otimes q_1\otimes \cdots \otimes q_i \otimes q_{i+1}
\otimes \cdots \otimes q_b \otimes \sigma \otimes p_1 \otimes
\cdots \otimes p_a \otimes \omega
\sim \\
(-1)^{|q_i||q_{i+1}|} \theta \otimes q_1\otimes \cdots \otimes
q_{i+1} \otimes q_i \otimes \cdots \otimes q_b \otimes \sigma'
\otimes p_1 \otimes \cdots \otimes p_a\otimes \omega,
\end{eqnarray*}

o\`u $\sigma'=(1\ldots i+1\ i\ldots\, a)_\ok\ \sigma$.

\textsc{Remarque :} Comme les induites $\theta$ et $\omega$
n'influent pas sur les r\`egles de signes, nous omettrons de les
\'ecrire par la suite pour all\'eger les notations.\\

Pour un \'el\'ement $q_1\otimes \cdots \otimes q_b \otimes \sigma
\otimes p_1 \otimes \cdots \otimes p_a$ de $\Qo(\ol,\,
\ok)\otimes_{\Sy_\ok} k[\Sc_\kj] \otimes_{\Sy_\ol} \Po(\oj,\,
\oi)$, on d\'efinit la diff\'erentielle $\delta$ par
\begin{eqnarray*}
&& \delta(q_1\otimes \cdots \otimes q_b \otimes \sigma \otimes p_1
\otimes \cdots \otimes p_a)= \\
&&\sum_{\beta=1}^b (-1)^{|q_1|+\cdots +|q_{\beta-1}|}q_1\otimes
\cdots \otimes \delta(q_\beta)\otimes \cdots \otimes q_b \otimes
\sigma \otimes p_1 \otimes \cdots \otimes p_a +\\
&&\sum_{\alpha=1}^a (-1)^{|q_1|+\cdots
+|q_b|+|p_1|+\cdots+|p_{\alpha-1}|}  q_1 \otimes \cdots \otimes
q_b \otimes \sigma \otimes p_1 \otimes \cdots \otimes
\delta(p_\alpha) \otimes \cdots \otimes p_a.
\end{eqnarray*}

\begin{lem}
\label{differentielle}
 La diff\'erentielle $\delta$ est constante
sur les classes d'\'equivalence pour la relation $\sim$.
\end{lem}

\begin{deo}
Comme les transpositions de la forme $(1 \ldots i+1 \ i \ldots b)$
engendrent $\Sy_b$, il suffit de faire les v\'erifications
suivantes
\begin{eqnarray*}
&&\delta((-1)^{|q_i||q_{i+1}|}(q_1,\ldots,\,q_{i+1},\,
q_i,\ldots,\, q_b)\, \sigma'\, (p_1,\ldots,\,p_a))=\\
&&\sum_{\beta=1}^{i-1}(-1)^{|q_i||q_{i+1}|+|q_1|+\cdots+|q_{\beta-1}|}(q_1,\ldots,\,
\delta(q_\beta),\ldots,\, q_{i+1},\, q_i,\ldots,\,
q_b)\, \sigma'\, (p_1,\ldots,\, p_a) +\\
&&(-1)^{|q_i||q_{i+1}|+|q_1|+\cdots+|q_{i-1}|}(q_1,\ldots,\,
\delta(q_{i+1}),\, q_i,\ldots,\, q_b)\, \sigma'\, (p_1,\ldots,\,
p_a) + \\
&&(-1)^{|q_i||q_{i+1}|+|q_1|+\cdots+|q_{i-1}|+|q_{i+1}|}(q_1,\ldots,\,
q_{i+1},\, \delta(q_i),\ldots,\, q_b)\, \sigma'\, (p_1,\ldots,\,
p_a) + \\
&&\sum_{\beta=i+2}^b(-1)^{|q_i||q_{i+1}|+|q_1|+\cdots+|q_{\beta-1}|}(q_1,\ldots,\,
q_{i+1},\, q_i,\ldots,\, \delta(q_\beta),\ldots,\, q_b)\,
\sigma'\, (p_1,\ldots,\, p_a) +\\
&&\sum_{\alpha=1}^a
(-1)^{|q_i||q_{i+1}|+|q_1|+\cdots+|q_b|+|p_1|+\ldots+|p_{\alpha-1}|}(q_1,\cdots,\,q_{i+1},\,
q_i,\ldots,\,q_b)\,\sigma' \, (p_1,\ldots,\,
\delta(p_\alpha),\ldots,\, p_a) \sim \\
&&\sum_{\beta=1}^{i-1}(-1)^{|q_1|+\cdots+|q_{\beta-1}|}(q_1,\ldots,\,
\delta(q_\beta),\ldots,\, q_{i},\, q_{i+1},\ldots,\,
q_b)\, \sigma\, (p_1,\ldots,\, p_a) +\\
&&(-1)^{|q_1|+\cdots+|q_{i-1}|+|q_i|}(q_1,\ldots,\, q_i
,\,\delta(q_{i+1}),\ldots,\, q_b)\, \sigma\, (p_1,\ldots,\,
p_a) + \\
&&(-1)^{|q_1|+\cdots+|q_{i-1}|}(q_1,\ldots,\,
\delta(q_i),\,q_{i+1} ,\ldots,\, q_b)\, \sigma\, (p_1,\ldots,\,
p_a) + \\
&&\sum_{\beta=i+2}^b(-1)^{|q_1|+\cdots+|q_{\beta-1}|}(q_1,\ldots,\,
q_i,\, q_{i+1},\ldots,\, \delta(q_\beta),\ldots,\, q_b)\,
\sigma\, (p_1,\ldots,\, p_a) +\\
&&\sum_{\alpha=1}^a
(-1)^{|q_1|+\cdots+|q_b|+|p_1|+\cdots+|p_{\alpha-1}|}(q_1,\ldots,\,q_b)\,\sigma'
\, (p_1,\ldots,\, \delta(p_\alpha),\ldots,\, p_a) \\
&& =\delta\big((q_1,\ldots,\, q_b)\, \sigma\, (p_1,\ldots,\,
p_a)\big),
\end{eqnarray*}
o\`u $\sigma'=(1\ldots i+1\ i\ldots\, b)_\ok\ \sigma$. Le calcul
faisant intervenir des permutations sur les $p_\beta$ est
similaire. $\cqfd$
\end{deo}

\begin{pro}
La cat\'egorie $(\textrm{dg-}\Sy\textrm{-biMod},\,  \boxtimes_c,\,
I)$ est une cat\'egorie mono\"\i dale qui admet la cat\'egorie
$(\Sy\textrm{-biMod},\, \boxtimes_c,\, I)$ comme sous-cat\'egorie
mono\"\i dale pleine.
\end{pro}

\begin{deo}
Tout $(\Sy_m,\, \Sy_n)$-bimodule $\Po(m,\, n)$ peut \^etre vu
comme un $(\Sy_m,\, \Sy_n)$-bimodule diff\'erentiel gradu\'e. Il
suffit pour cela de consid\'erer tous ses \'el\'ements comme
\'etant de degr\'e $0$ et $\Po(m,\,n)$ muni d'une diff\'erentielle
nulle. $\cqfd$
\end{deo}

On peut \'etendre sans difficult\'e la
proposition~\ref{Representation} au cadre diff\'erentiel.

\begin{pro}
\label{representationboxtimes}
 Le produit $\Qo\boxtimes_c \Po$ est
isomorphe, en tant que dg-$\Sy$-bimodule, \`a
$$\bigoplus _{\left((\Pi_1',\ldots,\, \Pi_b'),\, (\Pi_1,\ldots,\, \Pi_a)
\right)\in \Theta \atop |\ok|=|\oj|} k[\Sy_m] \otimes_{\Sy_\ol}
\Qo(\ol,\, \ok)\otimes_{\Sy_\ok} k[\Sc_\kj] \otimes_{\Sy_\oj}
\Po(\oj,\, \oi)\otimes_{\Sy_\oi} k[\Sy_n].$$
\end{pro}

Sur les parties connexes, on vient de pr\'eciser les r\`egles de
signes, les relations d'\'equivalences \`a consid\'erer et les
diff\'erentielles. On d\'efinit les m\^emes notions pour le
produit $\boxtimes$ en g\'en\'eralisant les r\`egles de
Koszul-Quillen du produit tensoriel $\otimes_k$ au produit
$\otimes$.

Soient $\Po$ et $\Qo$ deux $\Sy$-bimodules diff\'erentiels gradu\'es. Pour
$p\otimes q$ appartenant \`a $\Po_d(m,\,n)\otimes \Qo_e(m',\, n')$, on pose
$$\delta(p\otimes q)=\delta(p)\otimes q + (-1)^{d}p\otimes \delta(q).$$
Les isomorphismes de sym\`etrie s'\'ecrivent ici
$$ \tau_{p,\, q} \ : \ p\otimes q \mapsto (-1)^{de} q\otimes p.$$

Avec ces r\`egles de signes, on d\'efinit une diff\'erentielle naturelle sur
le produit $\Qo \boxtimes \Po$.

\begin{pro}
Le bifoncteur $\boxtimes$ de la cat\'egorie
des $\Sy$-bimodules s'\'etend \`a la cat\'egorie des $\Sy$-bimodules diff\'erentiels
gradu\'es.
\end{pro}

Comme dans le cas connexe, on peut \'etendre aux
dg-$\Sy$-bimodules la proposition~\ref{Representationboxtimes}.

\begin{pro}
\label{Representationboxtimesdiff}
 Le produit $\Qo\boxtimes \Po$ est
isomorphe, en tant que dg-$\Sy$-bimodule, \`a
$$\bigoplus_{\left((\Pi_1',\ldots,\, \Pi_a'),\, (\Pi_1,\ldots,\, \Pi_b)
\right)\in \Theta' \atop |\ok|=|\oj|} k[\Sy_m]
\otimes_{\Sy_\ol}\Qo(\ol,\, \ok)\otimes_{\Sy_\ok} k[\Sy_{|\ok|}]
\otimes_{\Sy_\oj} \Po(\oj,\, \oi) \otimes_{\Sy_\oi} k[\Sy_n].$$
\end{pro}

\subsection{Structure diff\'erentielle des produits de dg-$\Sy$-bimodules}
\label{structurediffproduit} On \'etudie ici l'homologie des
produits de deux $\Sy$-bimodules diff\'erentiels.\\

Remarquons d'abord que le complexe $(\Qo\boxtimes_c \Po,\, \delta)$
correspond au complexe total d'un bicomplexe. En effet, on
d\'efinit le bidegr\'e d'un \'el\'ement
$$(q_1,\, \ldots,\, q_b)\,
\sigma\, (p_1,\ldots,\, p_a) \quad  \textrm{de} \quad \Qo(\ol,\,
\ok)\otimes_{\Sy_{\ok}} k[\Sc_\kj] \otimes_{\Sy_{\oj}} \Po (\oj,\,
\oi)$$ par $(|q_1|+\cdots+|q_b|,\, |p_1|+\cdots+|p_a|)$. Et, la
diff\'erentielle horizontale $\delta_h \, :\, \Qo\boxtimes_c \Po
\to \Qo\boxtimes_c \Po$ est induite par celle de $\Qo$. Quant \`a
la diff\'erentielle verticale $\delta_v \, :\, \Qo\boxtimes_c \Po
\to \Qo\boxtimes_c \Po$, elle est induite par celle de $\Po$. De
mani\`ere explicite, on a
\begin{eqnarray*}
&& \delta_h\big((q_1,\, \ldots,\, q_b)\, \sigma\, (p_1,\ldots,\,
p_a) \big)= \\
&&\sum_{\beta=1}^b (-1)^{|q_1|+\cdots+|q_{\beta-1}|}(q_1,\ldots,\,
\delta(q_\beta),\ldots,\, q_b)\, \sigma\, (p_1,\ldots,\, p_a)\\
\textrm{et} && \delta_v\big((q_1,\, \ldots,\, q_b)\, \sigma\,
(p_1,\ldots,\, p_a) \big)=\\
&& \sum_{\alpha=1}^a
(-1)^{|q_1|+\cdots+|q_b|+|p_1|+\cdots+|p_{\alpha-1}|}(q_1,\ldots,\,
q_b)\, \sigma\, (p_1,\ldots,\,\delta(p_\alpha),\ldots,\, p_a).
\end{eqnarray*}
On voit clairement que $\delta=\delta_h+\delta_v$ et
$\delta_h\circ \delta_v =-\delta_v\circ \delta_h$.

Comme tout bicomplexe, celui-ci donne naissance \`a deux suites
spectrales $I^r(\Qo\boxtimes_c \Po)$ et $II^r(\Qo\boxtimes_c \Po)$
qui convergent vers l'homologie totale $H_*(\Qo\boxtimes_c \Po, \,
\delta)$ (\emph{cf.} Conventions). Rappelons que ces derni\`eres
v\'erifient
\begin{eqnarray*}
 && I^1(\Qo\boxtimes_c \Po)=H_*(\Qo\boxtimes_c \Po, \, \delta_v), \qquad
I^2(\Qo\boxtimes_c \Po)=H_*(H_*(\Qo\boxtimes_c \Po, \,
\delta_v),\,\delta_h) \\
&et& II^1(\Qo\boxtimes_c \Po)=H_*(\Qo\boxtimes_c \Po, \,
\delta_h), \qquad II^2(\Qo\boxtimes_c \Po)=H_*(H_*(\Qo\boxtimes_c
\Po, \, \delta_h),\, \delta_v).
\end{eqnarray*}
Gr\^ace \`a la proposition~\ref{representationboxtimes}, on sait
que produit mono\"\i dal $\Qo\boxtimes_c \Po$ s'\'ecrit \`a l'aide
de la somme directe
$$ \bigoplus_{\Theta}
\Qo(\ol,\, \ok)\otimes_{\Sy_\ok} k[\Sc_\kj] \otimes_{\Sy_\oj}
\Po(\oj,\, \oi).$$ Cette d\'ecomposition est compatible avec la
structure de bicomplexe. Ainsi, les suites spectrales se
d\'ecomposent de la m\^eme mani\`ere
\begin{eqnarray*}
&& I^r(\Qo\boxtimes_c\Po)=\bigoplus_{\Theta} I^r\left(\Qo(\ol,\,
\ok)\otimes_{\Sy_\ok} k[\Sc_\kj] \otimes_{\Sy_\oj} \Po(\oj,\, \oi)
\right) \\
&\textrm{et}&II^r(\Qo\boxtimes_c\Po)=\bigoplus_{\Theta}
II^r\left(\Qo(\ol,\, \ok)\otimes_{\Sy_\ok} k[\Sc_\kj]
\otimes_{\Sy_\oj} \Po(\oj,\, \oi) \right)
\end{eqnarray*}
On a imm\'ediatement les deux expressions suivantes.

\begin{pro}
\label{I1}
 Si $\Qo(\ol,\, \ok)\otimes_{\Sy_\ok} k[\Sc_\kj]$ est un
$k[\Sy_\oj]$-module projectif (\`a droite), alors on a
$$I^1\left(\Qo(\ol,\,
\ok)\otimes_{\Sy_\ok} k[\Sc_\kj] \otimes_{\Sy_\oj} \Po(\oj,\, \oi)
\right)=\Qo(\ol,\,
\ok)\otimes_{\Sy_\ok}k[\Sc_\kj]\otimes_{\Sy_\oj}
H_*\left(\Po(\oj,\, \oi)\right).$$ Et si,
$k[\Sc_\kj]\otimes_{\Sy_\oj}\Po(\oj,\,\oi)$ est un
$k[\Sy_\ok]$-module projectif (\`a gauche), alors on a
$$II^1\left(\Qo(\ol,\,
\ok)\otimes_{\Sy_\ok} k[\Sc\kj] \otimes_{\Sy_\oj} \Po(\oj,\, \oi)
\right)=H_*\left(\Qo(\ol,\,
\ok)\right)\otimes_{\Sy_\ok}k[\Sc_\kj]\otimes_{\Sy_\oj} \Po(\oj,\,
\oi).$$
\end{pro}

\begin{cor}
Si $k$ est un corps de caract\'eristique nulle, on a toujours
$$I^1\left(\Qo(\ol,\,
\ok)\otimes_{\Sy_\ok} k[\Sc_\kj] \otimes_{\Sy_\oj} \Po(\oj,\, \oi)
\right)=\Qo(\ol,\,
\ok)\otimes_{\Sy_\ok}k[\Sc_\kj]\otimes_{\Sy_\oj}
H_*\left(\Po(\oj,\, \oi)\right)\ \textrm{et}$$
$$II^1\left(\Qo(\ol,\,
\ok)\otimes_{\Sy_\ok} k[\Sc_\kj] \otimes_{\Sy_\oj} \Po(\oj,\, \oi)
\right)=H_*\left(\Qo(\ol,\,
\ok)\right)\otimes_{\Sy_\ok}k[\Sc_\kj]\otimes_{\Sy_\oj} \Po(\oj,\,
\oi).$$
\end{cor}

\begin{deo}
Par le th\'eor\`eme de Maschke, on sait que tous les anneaux
$k[\Sy_\oj]$ sont semisimples. Ainsi, tout $k[\Sy_\oj]$-module est
projectif. $\cqfd$
\end{deo}

\begin{pro}
\label{I2II2}
Si $k$ est de caract\'eristique nulle, alors on a
\begin{eqnarray*}
&& I^2\left(\Qo(\ol,\, \ok)\otimes_{\Sy_\ok} k[\Sc_\kj]
\otimes_{\Sy_\oj} \Po(\oj,\, \oi) \right)= II^2\left(\Qo(\ol,\,
\ok)\otimes_{\Sy_\ok} k[\Sc_\kj] \otimes_{\Sy_\oj} \Po(\oj,\, \oi)
\right)=\\ && H_*\Qo(\ol,\,
\ok)\otimes_{\Sy_\ok}k[\Sc_\kj]\otimes_{\Sy_\oj} H_*\Po(\oj,\,
\oi).
\end{eqnarray*}
\end{pro}

\begin{deo}
Comme tout $k[\Sy_\ok]$-module est projectif, on a
imm\'editatement
$$I^2\left(\Qo(\ol,\, \ok)\otimes_{\Sy_\ok} k[\Sc_\kj]
\otimes_{\Sy_\oj} \Po(\oj,\, \oi) \right)= H_*\left(\Qo(\ol,\,
\ok)\right)\otimes_{\Sy_\ok}k[\Sc_\kj]\otimes_{\Sy_\oj}
H_*\left(\Po(\oj,\, \oi)\right).$$ Et on conclut avec la formule
de K\"unneth
\begin{eqnarray*}
H_*\left(\Qo(\ol,\, \ok)\right)&=& H_*\left(\Qo(l_1,\, k_1)\otimes
\cdots \otimes \Qo(l_a,\, k_a) \right) = \\ && H_*\Qo(l_1,\,
k_1)\otimes \cdots \otimes H_*\Qo(l_a,\, k_a) = H_*\Qo(\ol,\,
\ok).
\end{eqnarray*}$\cqfd$
\end{deo}

\textsc{Remarque :} Soient $\Phi \, : \, M \to M'$ et $\Psi \, :
\, N \to N'$ deux morphismes de dg-$\Sy$-bimodules. Alors le
morphisme $\Phi\boxtimes_c \Psi \, :\, M\boxtimes_c N \to M'
\boxtimes N'$ est un morphisme de bicomplexes. Ainsi, il induit
des morphismes de suites spectrales $I^r(\Phi\boxtimes_c \Psi)\,
:\, I^r(M\boxtimes_c N) \to I^r(M'\boxtimes_c N')$ et $
II^r(\Phi\boxtimes_c \Psi)\, :\, II^r(M\boxtimes_c N) \to
II^r(M'\boxtimes_c N')$.\\

De mani\`ere g\'en\'erale, on a la proposition suivante.

\begin{pro}
\label{homologieproduit}
 Lorsque le corps $k$ est de
caract\'eristique nulle, on a toujours
$$H_*(\Qo \boxtimes_c \Po)(m,\, n)= \bigoplus_{\Theta}
H_*\Qo(\ol,\, \ok)\otimes_{\Sy_\ok} k[\Sc_\kj] \otimes_{\Sy_\oj}
H_*\Po(\oj,\, \oi).$$ C'est-\`a-dire, $H_*(\Qo \boxtimes_c \Po)
=(H_*\Qo) \boxtimes_c (H_*\Po)$.
\end{pro}

\begin{deo}
C'est une application directe des th\'eor\`emes de Mashke (pour
les coinvariants) et de K\"unneth (pour les produits tensoriels).
$\cqfd$
\end{deo}

On peut reprendre les m\^emes raisonnements dans le cas du produit
$\boxtimes$.\\

Le complexe $(\Qo \boxtimes \Po,\ \delta)$ est encore le complexe
total d'un bicomplexe. Les \'el\'ements de $\Qo$ fournissent la
graduation et la diff\'erentielle horizontales et ceux de $\Po$ la
graduation et la diff\'erentielle verticales. En ce qui concerne
les signes, il faut faire attention aux composantes connexes. Par
exemple, pour un objet
$$(q_1,\, \ldots,\, q_b)\,
\sigma\, (p_1,\ldots,\, p_a) \otimes (q'_1,\, \ldots,\, q'_{b'})\,
\sigma\, (p'_1,\ldots,\, p'_{a'})$$ de $\Qo \boxtimes_c \Po
\otimes\Qo \boxtimes_c \Po$, on a

\begin{eqnarray*}
&&\delta_h\left(  (q_1,\, \ldots,\, q_b)\, \sigma\, (p_1,\ldots,\,
p_a) \otimes (q'_1,\, \ldots,\, q'_{b'})\, \sigma\,
(p'_1,\ldots,\, p'_{a'})\right) = \\
&&\sum_{\beta=1}^b (-1)^{|q_1|+\cdots+|q_{\beta-1}|}(q_1,\ldots,\,
\delta(q_\beta),\ldots,\, q_b)\, \sigma\, (p_1,\ldots,\, p_a)
\\
&&\quad \quad \quad  \quad \quad \quad \quad \quad \quad \otimes
(q'_1,\, \ldots,\, q'_{b'})\, \sigma\, (p'_1,\ldots,\,
p'_{a'})+\\
&&\sum_{\beta=1}^{b'}
(-1)^{|q_1|+\cdots+|q_b|+|p_1|+\cdots+|p_a|+|q'_1|+\cdots+|q'_{\beta-1}|}
(q_1,\ldots,\, q_b)\, \sigma\, (p_1,\ldots,\, p_a)\\
&& \quad \quad \quad  \quad \quad \quad \quad \quad \quad \otimes
(q'_1,\, \ldots,\, \delta(q'_\beta),\ldots,\, q'_{b'})\, \sigma\,
(p'_1,\ldots,\, p'_{a'})
\end{eqnarray*}
et
\begin{eqnarray*}
&& \delta_v\left( (q_1,\, \ldots,\, q_b)\, \sigma\, (p_1,\ldots,\,
p_a) \otimes (q'_1,\, \ldots,\, q'_{b'})\, \sigma\,
(p'_1,\ldots,\, p'_{a'})\right) = \\
&&\sum_{\alpha=1}^a
(-1)^{|q_1|+\cdots+|q_b|+|p_1|+\cdots+|p_{\alpha-1}|}(q_1,\ldots,\,
 q_b)\, \sigma\, (p_1,\ldots,\,\delta(p_\alpha),\, \ldots ,\,
p_a)\\
&& \quad \quad \quad  \quad \quad \quad \quad \quad \quad \otimes
(q'_1,\, \ldots,\, q'_{b'})\, \sigma\, (p'_1,\ldots,\,
p'_{a'})+\\
&&\sum_{\beta=1}^{b'}
(-1)^{|q_1|+\cdots+|q_b|+|p_1|+\cdots+|p_a|+|q'_1|+\cdots+|q'_{b'}|
+|p'_1|+\cdots+|p'_{\alpha-1}|} \\
&& (q_1,\ldots,\, q_b)\, \sigma\, (p_1,\ldots,\, p_a)\otimes
(q'_1,\, \ldots ,\,  q'_{b'})\, \sigma\,
(p'_1,\ldots,\,\delta(p'_\alpha),\ldots,\,  p'_{a'}).
\end{eqnarray*}

On a alors les m\^emes r\'esultats que dans le cas connexe.

\begin{pro}
Si $k$ est un corps de caract\'eristique nulle, on a les
\'egalit\'es
$$I^1\left(\Qo(\ol,\,
\ok)\otimes_{\Sy_\ok} k[\Sy_N] \otimes_{\Sy_\oj} \Po(\oj,\, \oi)
\right)=\Qo(\ol,\, \ok)\otimes_{\Sy_\ok}k[\Sy_N]\otimes_{\Sy_\oj}
H_*\left(\Po(\oj,\, \oi)\right)\ \textrm{et}$$
$$II^1\left(\Qo(\ol,\,
\ok)\otimes_{\Sy_\ok} k[\Sy_N] \otimes_{\Sy_\oj} \Po(\oj,\, \oi)
\right)=H_*\left(\Qo(\ol,\,
\ok)\right)\otimes_{\Sy_\ok}k[\Sy_N]\otimes_{\Sy_\oj} \Po(\oj,\,
\oi).$$
\end{pro}

\begin{pro}
Si $k$ est de caract\'eristique nulle, alors on a
\begin{eqnarray*}
&& I^2\left(\Qo(\ol,\, \ok)\otimes_{\Sy_\ok} k[\Sy_N]
\otimes_{\Sy_\oj} \Po(\oj,\, \oi) \right)= II^2\left(\Qo(\ol,\,
\ok)\otimes_{\Sy_\ok} k[\Sy_N] \otimes_{\Sy_\oj} \Po(\oj,\, \oi)
\right)=\\ && H_*\Qo(\ol,\,
\ok)\otimes_{\Sy_\ok}k[\Sy_N]\otimes_{\Sy_\oj} H_*\Po(\oj,\, \oi).
\end{eqnarray*}
\end{pro}

\begin{pro}
\label{homologieproduitboxtimes}
 Lorsque le corps $k$ est de
caract\'eristique nulle, on a
$$H_*(\Qo \boxtimes \Po) =(H_*\Qo)
\boxtimes (H_*\Po).$$
\end{pro}

\begin{cor}
Soient $\Qo$ et $\Po$ deux  dg-$\Sy$-bimodules. On a la relation
$$H_*\left( \Qo\boxtimes \Po \right) = S_\otimes\left( H_*(\Qo \boxtimes_c \Po )\right).$$
\end{cor}

\textsc{Remarque :} Toutes ces propositions sont des
g\'en\'eralisations des r\'esultats de B. Fresse sur le produit de
composition $\circ$ des op\'erades aux produits $\boxtimes_c$ et
$\boxtimes$ des prop\'erades et des PROPs (\emph{cf.}
\cite{Fresse} 2.3.).

\subsection{Prop\'erades et PROPs diff\'erentiels}

On peut donner une version diff\'erentielle \`a la notion de
prop\'erade et de PROP.

\begin{dei}[Prop\'erades diff\'erentielles]
\index{prop\'erade diff\'erentielle}
On appelle \emph{prop\'erade
diff\'erentielle} un mono\"\i de $(\Po,\, \mu,\, \eta)$ dans la
cat\'egorie $(\textrm{dg-}\Sy\textrm{-biMod},\, \boxtimes_c, \,
I)$.
\end{dei}

Cela signifie qu'en plus de la donn\'ee d'une prop\'erade
lin\'eaire, le morphisme de composition $\mu$ est un morphisme de
dg-modules de degr\'e $0$ pr\'eservant les actions de $\Sy_m$ et
$\Sy_n$. Ainsi, on a la relation du type d\'erivation :
\begin{eqnarray*}
&& \delta \left( \mu((p_1,\ldots,\, p_b)\, \sigma \,
(p'_1,\ldots,\, p'_a) ) \right) = \\
&& \sum_{\beta=1}^b (-1)^{|p_1|+\cdots +|p_{\beta-1}|} \mu((
p_1,\ldots ,\, \delta(p_\beta),\ldots ,\,
p_b) \, \sigma \, (p'_1,\ldots ,\,p'_a)) +\\
&& \sum_{\alpha=1}^a (-1)^{|p_1|+\cdots
+|p_b|+|p'_1|+\cdots+|p'_{\alpha-1}|}  \mu((p_1,\ldots,\, p_b) \,
\sigma \, (p'_1,\ldots,\, \delta(p'_\alpha),\ldots,\,p'_a)).
\end{eqnarray*}

La proposition~\ref{homologieproduit} donne la propri\'et\'e
suivante.

\begin{pro}
Lorsque le corps $k$ est de caract\'eristique nulle, pour toute
prop\'erade diff\'erentielle $(\Po,\, \mu,\, \eta)$, l'homologie
de cette prop\'erade $(H_*(\Po),\ H_*(\mu),\, H_*(\eta))$ est une
prop\'erade gradu\'ee.
\end{pro}

\begin{dei}[$\Sy$-bimodules diff\'erentiels gradu\'es par un poids]
\index{$\Sy$-bimodule diff\'erentiel gradu\'e par un poids}

On appelle \emph{$\Sy$-bimodule diff\'erentiel gradu\'e par un
poids} toute somme directe sur $\rho \in \mathbb{N}$ de
$\Sy$-bimodules diff\'erentiels $M=\bigoplus_{\rho\in \mathbb{N}}
M^{(\rho)}$.

Les morphismes de $\Sy$-bimodules diff\'erentiels gradu\'es par un
poids sont les morphismes de $\Sy$-bimodules diff\'erentiels qui
pr\'eservent cette d\'ecomposition. L'ensemble des $\Sy$-bimodules
diff\'erentiels gradu\'es par un poids, muni des morphismes
correspondant, forme une cat\'egorie que l'on note
gr-dg-$\Sy$-biMod.
\end{dei}

Remarquons que la graduation est une information suppl\'ementaire
ind\'ependante du degr\'e homologique.

\begin{pro}
\label{bigraduation}
 Soit $f=\bigoplus_{n=0}^\infty f_{(n)}$ un
foncteur analytique dans la cat\'egorie des $\Sy$-bimodules. Soit
$M=\bigoplus_{\rho\in \mathbb{N}} M^{(\rho)}$ un $\Sy$-bimodule
gradu\'e par un poids. Alors le $\Sy$-bimodule $f(M)$ est
bigradu\'e, d'une part avec la graduation venant du foncteur
analytique, d'autre part avec la graduation totale venant du
poids.
\end{pro}

\begin{deo}
Posons $f_{(n)}=f_n\circ \Delta_n$ o\`u $f_n$ est une application
$n$-lin\'eaire. La premi\`ere graduation s'\'ecrit
$f(M)_{(n)}=f_n(M,\ldots,\, M)$. La deuxi\`eme correspond \`a
$$f(M)^{(\rho)}=\sum_{i_1+\cdots+i_n=\rho} f_n(M^{(i_1)},\ldots,\, M^{(i_n)}).$$
$\cqfd$
\end{deo}

\begin{cor}
Toute prop\'erade libre (respectivement PROP libre) sur un
$\Sy$-bimodule $M$ gradu\'e par un poids est bigradu\'e.
\end{cor}

\begin{deo}
Le foncteur $\F$ (respectivement $S_\otimes (\F)$) est un foncteur
analytique dont les degr\'es sont donn\'es par le nombre de
sommets des graphes. La premi\`ere graduation est donc donn\'ee
par le nombre d'op\'erations de $M$ qui servent \`a repr\'esenter
un \'el\'ement de $\F(M)$. La deuxi\`eme graduation est \'egale
\`a la somme des poids de ces op\'erations. $\cqfd$
\end{deo}

\begin{dei}[Prop\'erades diff\'erentielle gradu\'ees par un poids]
\index{prop\'erade diff\'erentielle gradu\'ee par un poids}

Une \emph{prop\'erade diff\'erentielle gradu\'ee par un poids}
$(\Po,\, \mu,\, \eta)$ est un mono\"\i de dans la cat\'egorie
$(\textrm{gr-dg-}\Sy\textrm{-biMod},\, \boxtimes_c,\, I)$.

Une telle prop\'erade sera dite \emph{connexe} si de plus
$\Po^{(0)}=I$. \index{prop\'erade connexe}
\end{dei}

\textsc{Remarque :} Duallement, on a la notion de coprop\'erade
diff\'erentielle \index{coprop\'erade diff\'erentielle} (et
diff\'erentielle gradu\'ee par un poids). \index{coprop\'erade
diff\'erentielle gradu\'ee par un poids}Une coprop\'erade
diff\'erentielle correspond \`a une coprop\'erade dont le
coproduit est un morphisme de dg-$\Sy$-bimodules,
c'est-\`a-dire v\'erifie une relation du type cod\'erivation.\\

De la m\^eme mani\`ere, on peut consid\'erer des PROPs
diff\'erentiels.

\begin{dei}[PROPs diff\'erentiels]
\index{PROP diff\'erentiel}

On appelle \emph{PROP diff\'erentiel} la donn\'ee d'une structure
de PROP $(\Po,\, \widetilde{\mu},\, conc,\ \eta)$ dans la
cat\'egorie $(\textrm{dg-}\Sy\textrm{-biMod},\, \boxtimes,\,
\otimes)$.
\end{dei}

On exige donc en plus que les morphismes de compositions
$\widetilde{\mu}$ et $conc$ soient des morphismes de degr\'e $0$
de $\Sy$-bimodules diff\'erentiels, c'est-\`a-dire qu'ils
v\'erifient des relations du type d\'erivation. Par exemple, pour
la composition horizontale, on a
$$\delta(conc(p\otimes q))=conc(\delta(p)\otimes q) +(-1)^{|p|}conc(p\otimes \delta(q)).$$

On a des versions PROPiques des r\'esultats pr\'ec\'edents. Par
exemple, avec la proposition~\ref{homologieproduitboxtimes}, on
montre la proposition suivante.

\begin{pro}
Lorsque le corps $k$ est de caract\'eristique nulle, pour tout
PROP diff\'erentiel $(\Po,\, \widetilde{\mu},\, conc)$,
l'homologie de ce PROP $(H_*(\Po),\ H_*(\widetilde{\mu)},\,
H_*(conc))$ est un PROP gradu\'e.
\end{pro}

On peut aussi d\'efinir la notion de PROP diff\'erentiel gradu\'e
par un poids.

\begin{dei}[PROP diff\'erentiel gradu\'e par un poids]
\index{PROP diff\'erentiel gradu\'e par un poids}

Un \emph{PROP diff\'erentiel gradu\'e par un poids} est la donn\'e
d'un PROP dans la cat\'egorie des $\Sy$-bimodules diff\'erentiels
gradu\'es par un poids $(\textrm{gr-dg-}\Sy\textrm{-biMod},\,
\boxtimes ,\, \otimes)$.

Un tel PROP est dit \emph{connexe} si
$\Po^{(0)}=S_\otimes(I)$.\index{PROP connexe}
\end{dei}

\section{$\Po$-modules libres et quasi-libres}

Soit $\Po$ une prop\'erade diff\'erentielle. La notion de
$\Po$-module s'\'etend naturellement au cas diff\'erentiel. On
s'int\'eresse ici plus particuli\`erement aux $\Po$-modules libres
et \`a une version l\'eg\`erement diff\'erente, celle des
\emph{$\Po$-modules quasi-libres}. On ne traitera ici que les
$\Po$-modules \`a droite, le cas des $\Po$-modules \`a gauche
\'etant parfaitement sym\'etrique.

\subsection{$\Po$-modules diff\'erentiels libres}

\index{$\Po$-module diff\'erentiel libre} Les r\'esultats
g\'en\'eraux sur les cat\'egories mono\"\i dales donn\'es au
chapitre $1$, montrent que le $\Po$-module diff\'erentiel libre
\`a droite sur un dg-$\Sy$-module $M$ correspond \`a $L=M\boxtimes
\Po$. La diff\'erentielle $\delta$ sur $L$ est la diff\'erentielle
canonique sur le produit
$M\boxtimes \Po$ explicit\'ee pr\'ec\'edemment. \\

De la m\^eme mani\`ere, si $(\Po,\, \mu,\, \eta,\, \varepsilon)$
est une prop\'erade diff\'erentielle augment\'ee, c'est-\`a-dire
que l'on a en plus un morphisme mono\"\i dal de dg-$\Sy$-bimodules
$\varepsilon \, : \, \Po \to I$, on peut appliquer la
proposition~\ref{quotientindecomposable}. Dans ce cas, le quotient
ind\'ecomposable du $\Po$-module diff\'erentiel libre
$L=M\boxtimes \Po$ est isomorphe, en tant que dg-$\Sy$-bimodule,
 \`a $M$.

\subsection{$\Po$-modules quasi-libres}

Par la suite, nous traiterons des exemples qui ne rentrent pas
strictement dans le cadre des $\Po$-modules libres. En effet, dans
ces cas l\`a, la diff\'erentielle $\delta_\theta$ correspond \`a
la somme de la diff\'erentielle canonique $\delta$ avec un
morphisme homog\`ene de degr\'e $-1$.

\begin{dei}[$\Po$-module quasi-libre \`a droite]
\index{$\Po$-module quasi-libre}

Un \emph{$\Po$-module quasi-libre \`a droite} $L$ est un
dg-$\Sy$-bimodule de la forme $M\boxtimes \Po$ mais o\`u la
diff\'erentielle $\delta_\theta$ est la somme $\delta+d_\theta$ de
la diff\'erentielle canonique $\delta$ sur le produit $M \boxtimes
\Po$ avec un morphisme  homog\`ene de degr\'e $-1$, $d_\theta \, :
\, M\boxtimes \Po \to M\boxtimes \Po$ tel que
\begin{eqnarray*}
&& d_\theta\left((m_1,\ldots,\, m_b)\, \sigma \, (p_1,\ldots,\,
p_a) \right) =\\
&& \sum_{\beta=1}^b (-1)^{|m_1|+\cdots+|m_{\beta-1}|}
(m_1,\ldots,\, d_\theta(m_\beta),\ldots,\, m_b)\, \sigma \,
(p_1,\ldots,\, p_a).
\end{eqnarray*}
\end{dei}

Ainsi, $\delta_\theta \, : \, M\boxtimes \Po \to M\boxtimes \Po$
v\'erifie la relation de d\'erivation suivante
\begin{eqnarray*}
&& \delta_\theta\left((m_1,\ldots,\, m_b)\, \sigma \,
(p_1,\ldots,\,
p_a) \right) =\\
&& \sum_{\beta=1}^b (-1)^{|m_1|+\cdots+|m_{\beta-1}|}
(m_1,\ldots,\, \delta_\theta(m_\beta),\ldots,\, m_b)\, \sigma \,
(p_1,\ldots,\, p_a) +\\
&& \sum_{\alpha=1}^a
(-1)^{|m_1|+\cdots+|m_b|+|p_1|+\cdots+|p_{\alpha-1}|}
(m_1,\ldots,\, m_b)\, \sigma \, (p_1,\ldots,\, \delta(p_\alpha)
,\ldots,\, p_a).
\end{eqnarray*}

Comme $\delta^2=0$, l'identit\'e ${\delta_\theta}^2=0$ est
\'equivalente \`a $\delta d_\theta + d_\theta \delta +
{d_\theta}^2=0$. \\

Si on \'ecrit $d_\theta(m_\beta)=(m'_1,\ldots,\, m'_{b'})\,
\sigma'\,(p'_1,\ldots,\, p'_{a'})\in M\boxtimes \Po$, alors
$$(m_1,\ldots,\, d_\theta(m_\beta),\ldots,\, m_b)\, \sigma \,
(p_1,\ldots,\, p_a)$$ correspond en fait \`a
$$(m_1,\ldots,\, m'_1,\ldots,\, m'_{b'},\ldots ,\ m_b)\, \sigma' \,
\mu\left( (p'_1,\ldots,\, p'_{a'}) \, \sigma \, (p_1,\ldots,\,
p_a) \right) \in M\boxtimes \Po.$$

On peut remarquer que le morphisme $d_\theta$ est d\'etermin\'e
par sa restriction sur $\xymatrix{M \ar[r]^-{M\boxtimes \eta}&
M\boxtimes \Po}$.

\begin{pro}
Le morphisme $\xymatrix{M \ar[r]^-{M\boxtimes \eta}& M\boxtimes
\Po}$ est un morphisme de dg-$\Sy$-bimodules si et seulement si
$d_\theta=0$.
\end{pro}

Cependant si $\Po$ est une prop\'erade diff\'erentielle
augment\'ee, la condition pour que le morphisme $M\boxtimes
\varepsilon$ soit un morphisme de dg-module est moins restrictive.

\begin{pro}
Soit $(\Po,\, \mu,\, \eta,\, \varepsilon)$ une prop\'erade
augment\'ee. Le morphisme $\xymatrix{M\boxtimes\Po
\ar[r]^-{M\boxtimes \varepsilon}& M}$ est un morphisme de
dg-$\Sy$-bimodules si et seulement si
$$\xymatrix{ M\boxtimes \Po \ar[r]^-{d_\theta}& M \boxtimes
\Po \ar[r]^-{M\boxtimes \varepsilon}& M }=0.$$
\end{pro}

\begin{cor}
Dans le cas o\`u la composition $\xymatrix{ M\boxtimes \Po
\ar[r]^-{d_\theta}& M \boxtimes \Po \ar[r]^-{M\boxtimes
\varepsilon}& M }$ est nulle, le quotient ind\'ecomposable de $L$
est isomorphe \`a $M$ en tant que dg-$\Sy$-bimodule.
\end{cor}

Dans la suite, les exemples de $\Po$-modules quasi-libres que nous
aurons \`a traiter auront souvent une forme particuli\`ere.

\begin{dei}[$\Po$-module quasi-libre analytique]
\index{$\Po$-module quasi-libre analytique}

Soit $\Po$ une prop\'erade diff\'erentielle augment\'ee. Un
\emph{$\Po$-module quasi-libre analytique \`a droite} est un
$\Po$-module quasi-libre de la forme $L=\Upsilon(V)\boxtimes \Po$,
o\`u $V$ est un dg-$\Sy$-bimodule et $\Upsilon(V)$ un foncteur
analytique en $V$. On pose $\Upsilon=\bigoplus_n \Upsilon_{(n)}$
o\`u $\Upsilon_{(n)}$ est un foncteur polynomial homog\`ene de
degr\'e $n$ tel que $\Upsilon_{(0)}=I$. En outre, on suppose de
$V$ est gradu\'e par un poids et que $V^{(0)}=0$. Comme nous
l'avons vu \`a la proposition~\ref{bigraduation}, $\Upsilon(V)$
est bigradu\'e, d'abord par les degr\'es du foncteur analytique
puis par la graduation totale venant du poids de $V$.

On impose de plus que le morphisme $d_\theta$ v\'erifie
$$d_\theta \ : \ \Upsilon(V)^{(n)} \to \underbrace{\Upsilon(V)}_{(<n)}\boxtimes
(I\oplus \underbrace{\oPo}_1),$$ o\`u la graduation $(n)$ est ici
la graduation totale venant du poids de $V$.

De plus, si $\Po$ est gradu\'ee par un poids on exige que le
morphisme $d_\theta$ pr\'eserve globalement la graduation totale
venant des poids de $\Po$ et de $V$.
\end{dei}

On a imm\'ediatement la proposition suivante.

\begin{pro}
Soit $L=\Upsilon(V)\boxtimes \Po$ un $\Po$-module quasi-libre
analytique. Le morphisme $d_\theta$ v\'erifie
$$d_\theta \
:\ \underbrace{\Upsilon(V)}_{(n)}\boxtimes \, \Po \to
\underbrace{\Upsilon(V)}_{(<n)}\boxtimes \, \Po.$$ Et,
$$\xymatrix{ M\boxtimes \Po \ar[r]^-{d_\theta}& M \boxtimes \Po
\ar[r]^-{M\boxtimes \varepsilon}& M }=0.$$ On peut donc identifier
le quotient ind\'ecomposable de $L$ \`a $\Upsilon(V)$.
\end{pro}

\textsc{Exemples :} Les exemples de $\Po$-modules quasi-libres
analytiques que nous aurons \`a traiter se divisent en deux
cat\'egories.
\begin{itemize}
\item Dans le cas o\`u $V=\overline{\Po}$, on trouve la bar
construction $\Upsilon(V)=\bar{\B}(\Po)=\F^c(\Sigma\oPo)$
(\emph{cf.} chapitre $4$) ainsi que la bar construction
simpliciale
$\Upsilon(V)=\bar{\mathcal{C}}(\Po)=\mathcal{FS}(\oPo)$
(\emph{cf.} chapitre $6$). Alors la bar construction augment\'ee
$\bar{\B}(\Po)\boxtimes \Po$, et la bar construction simpliciale
augment\'ee $\bar{\mathcal{C}}(\Po)\boxtimes \Po$ sont deux
exemples de $\Po$-modules quasi-libres analytiques.

\item Si $\Po$ est une prop\'erade quadratique engendr\'ee par un
dg-$\Sy$-bimodule $V$. La prop\'erade duale de Koszul $\Po^{\ac}$
est un foncteur analytique en $V$ (\emph{cf.} chapitre $7$). Le
complexe de Koszul $\Po^{\ac}\boxtimes \Po$ est alors un exemple
de $\Po$-module quasi-libre analytique.
\end{itemize}

\begin{pro}
\label{décompositionLrho}
 Lorsque $\Po$ est une prop\'erade
diff\'erentielle augment\'ee gradu\'ee par un poids, tout
$\Po$-module quasi-libre analytique $L=\Upsilon(V)\boxtimes \Po$
se d\'ecompose en somme directe de sous-complexes
$L=\bigoplus_{\rho \in \mathbb{N}} L^{(\rho)}$, o\`u
$L^{(\rho)}=(\Upsilon(V)\boxtimes \Po)^{(\rho)}$ repr\'esente le
sous-module de $L$ compos\'e des \'el\'ements de graduation totale
$\rho$.
\end{pro}

\begin{deo}
Les diff\'erentielles $\delta_{\Upsilon(V)}$ et $\delta_\Po$
pr\'eservent cette graduation sur $\Upsilon(V)$ et $\Po$
respectivement. Et, par d\'efinition, le morphisme $d_\theta$
pr\'eserve cette graduation globalement.$\cqfd$
\end{deo}

\begin{dei}[Morphisme de $\Po$-modules quasi-libres analytiques]
Soient $L=\Upsilon(V)\boxtimes \Po$ et $L'=\Upsilon'(V')\boxtimes
\Po$ deux $\Po$-modules quasi-libres analytiques. Et soit $\Phi \,
:\, L \to L'$ un morphisme de $\Po$-modules diff\'erentiels. On
dit que $\Phi$ est un \emph{morphisme de $\Po$-modules
quasi-libres analytiques} si, de plus, $\Phi$ pr\'eserve la
graduation en $V$ :
$$\Phi \ :\ \underbrace{\Upsilon(V)}_n \boxtimes \Po \to
\underbrace{\Upsilon'(V')}_n\boxtimes \Po.$$

Si, en outre, la prop\'erade diff\'erentielle $\Po$ est gradu\'ee
par un poids, comme $V$ et $V'$ sont aussi gradu\'es par un poids,
on impose en plus que $\Phi$ pr\'eserve la graduation naturelle
globale de $\Upsilon(V)\boxtimes \Po$ :
$$\Phi \ : \  L^{(\rho)}=(\Upsilon(V)\boxtimes \Po)^{(\rho)} \to
{L'}^{(\rho)}=(\Upsilon'(V')\boxtimes \Po)^{(\rho)}.$$ On note
$\Phi^{(\rho)}$ la restriction de $\Phi$ au sous-complexe
$L^{(\rho)}$.
\end{dei}

\textsc{Remarque :} On a les notions duales de
$\mathcal{C}$-comodule diff\'erentiel colibre sur un
dg-$\Sy$-bimodule $M$ (donn\'e par $M\boxtimes \mathcal{C}$), de
$\mathcal{C}$-comodule quasi-colibre et de $\mathcal{C}$-comodule
 quasi-colibre analytique. La cobar construction aumgent\'ee \`a droite $\bar{\B}^c(\mathcal{C})\boxtimes \mathcal{C}$
 est un exemple de $\mathcal{C}$-comodule
 quasi-colibre analytique. \index{$\mathcal{C}$-comodule diff\'erentiel
 colibre} \index{$\mathcal{C}$-comodule
 quasi-colibre analytique}

\section{Prop\'erades et PROPs quasi-libres}

Pour expliciter la construction de la prop\'erade et du PROP
libres dans le cas diff\'erentiel, il reste \`a comprendre
d\'efinir la diff\'erentielle. En outre, comme pour les
$\Po$-modules, les notions de prop\'erade et de PROP
 diff\'erentiels libres n'englobent pas tous les cas que nous aurons
\`a traiter. C'est pourquoi, on consid\`ere les
\emph{prop\'erades} et les \emph{PROPs quasi-libres}, qui
apparaissent comme des prop\'erades et des PROPs libres mais dont
la diff\'erentielle est la somme de la diff\'erentielle de la
prop\'erade libre avec une \emph{d\'erivation}.

\subsection{Prop\'erades diff\'erentielles libres}
\index{prop\'erade diff\'erentielle libre}

Soit $(V,\, \delta)$ un dg-$\Sy$-bimodule. Une application directe
des r\'esultats de la section $6$ du chapitre $1$, sur la
construction du mono\"\i de libre, montre que le $\Sy$-bimodule
gradu\'e $\F(V)$ est donn\'e par la somme directe sur les graphes
connexes (sans niveau) dont les sommets sont indic\'es par des
\'el\'ements de $V$ (\emph{cf.} th\'eor\`eme~\ref{proplibre}),
soit
$$ \F(V)=\left( \bigoplus_{g\in\mathcal{G}} \bigotimes_{\nu\in
\mathcal{N}}V(|Out(\nu)|,\, |In(\nu)|)\right) \Bigg/_\approx   \ .
$$ Cette construction revient \`a consid\'erer l'objet simplicial
libre $\mathcal{FS}(V)=\bigoplus_{n\in\mathbb{N}} (V_+)^{\boxtimes
n}$ quotient\'e par la relation engendr\'ee par $I\boxtimes_c V
\equiv V \boxtimes_c I$. Or, dans la cadre diff\'erentiel, pour
pouvoir d\'efinir la diff\'erentielle d'un tel objet, il faut
\^etre plus fin au niveau des signes. On quotiente alors
$\mathcal{FS}(V)$ par la relation engendr\'ee par $(\nu_1,\,
\nu,\, \nu_2)\, \sigma \, (\nu'_1,\, 1,\, \nu'_2)\equiv
(-1)^{|\nu_2|+|\nu'_1|}(\nu_1,\, 1,\, 1,\, \nu_2)\, \sigma' \,
(\nu'_1,\, \nu,\, \nu'_2)  $, par exemple, o\`u $\nu$ appartient
\`a $V(2,\, 1)$ ici.

\begin{pro}
La diff\'erentielle canonique $\delta$ d\'efinie sur
$\mathcal{FS}(V)$ \`a partir de celle de $V$, passe au quotient
pour la relation d'\'equivalence $\equiv$.
\end{pro}

\begin{deo}
La v\'erification des signes est imm\'ediate et le calcul est
semblable \`a celui effectu\'e pour d\'emontrer le
lemme~\ref{differentielle}.$\cqfd$
\end{deo}

\begin{thm}
La prop\'erade diff\'erentielle $\F(V)$ munie de la
diff\'erentielle $\delta$ est la prop\'erade diff\'erentielle
libre sur $V$.
\end{thm}

\textsc{Remarque :} La d\'ecomposition analytique du foncteur
$\F(V)$ donn\'ee par la proposition~\ref{Fanalytique} est stable
par la diff\'erentielle $\delta$. Cette d\'ecomposition reste donc
toujours vraie dans le cas diff\'erentiel. De plus, la projection
$\F(V)\to \F_{(1)}(V)=V$ permet d'identifier le
dg-$\Sy$-bimodule $V$ avec le quotient ind\'ecomposable de $\F(V)$.\\

Nous allons voir que le foncteur $\F$ poss\`ede de bonnes
propri\'et\'es homologiques.

\begin{pro}
\label{Ffoncteurexact} Lorsque $k$ est un corps de
caract\'eristique nulle, le foncteur $\F \, : \,
\textrm{dg-}\Sy\textrm{-biMod} \to \textrm{dg-prop\'erades}$ est
un foncteur exact, c'est-\`a-dire $H_*(\F(V))=\F(H_*(V))$.
\end{pro}

\begin{deo}
La d\'emonstation de cette proposition est essentiellement la
m\^eme que celle du th\'eor\`eme $4$ de l'article de M. Markl et
A. A. Voronov \cite{MV} et repose encore une fois sur le
th\'eor\`eme de Mashke. $\cqfd$
\end{deo}

\textsc{Remarque :} On fait exactement le m\^eme raisonnement pour
les coprop\'erades colibres puisque le dg-$\Sy$-bimodule
sous-jacent $\F(V)=\F^c(V)$ est le m\^eme. \\

On poursuit la m\^eme d\'emarche pour d\'ecrire le PROP libre.
Pour cela, on consid\`ere les signes dus aux r\`egles de
Koszul-Quillen dans les isomorphismes de sym\'etrie du produit
mono\"\i dal $\otimes$ des $\Sy$-bimodules. Avec ces r\`egles de
signes, on d\'efinit imm\'ediatement une diff\'erentielle $\delta$
sur $S_\otimes(\F(V))$.

\begin{thm}
Le PROP diff\'erentiel libre sur un dg-$\Sy$-bimodule $V$ est
donn\'e par la construction $S_\otimes(\F(V))$ munie de la
diff\'erentielle $\delta$.
\end{thm}

En utilisant le th\'eor\`eme de Mashke, on montre que le foncteur
$S_\otimes$ est un foncteur exact. Et par composition de
foncteurs exacts, le foncteur PROP libre est un foncteur exact.

\begin{pro}
Lorsque $k$ est un corps de caract\'eristique nulle, le foncteur
$$S_\otimes \circ \F \ : \  \textrm{dg-}\Sy\textrm{-biMod}
\to \textrm{dg-PROPs}$$ est un foncteur exact, c'est-\`a-dire
$$H_*( S_\otimes(\F(V)))=S_\otimes(\F(H_*(V))) .$$
\end{pro}

\subsection{Prop\'erades quasi-libres}

Une prop\'erade quasi-libre est une prop\'erade libre mais dont la
diff\'erentielle est compos\'ee de la somme de deux termes : la
diff\'erentielle canonique $\delta$ plus une d\'erivation.

\begin{dei}[D\'erivation]
\index{d\'erivation}

Soit $\Po$ une prop\'erade diff\'erentielle. Un morphisme $d\, :
\, \Po \to \Po$ homog\`ene de degr\'e $|d|$ de $\Sy$-bimodules est
appel\'e \emph{d\'erivation} s'il v\'erifie l'identit\'e suivante
:
\begin{eqnarray*}
&& \delta \left( \mu((p_1,\ldots,\, p_b)\, \sigma \,
(p'_1,\ldots,\, p'_a) ) \right) = \\
&& \sum_{\beta=1}^b (-1)^{(|p_1|+\cdots +|p_{\beta-1}|)|d|} \mu((
p_1,\ldots ,\, \delta(p_\beta),\ldots ,\,
p_b) \, \sigma \, (p'_1,\ldots ,\,p'_a)) +\\
&& \sum_{\alpha=1}^a (-1)^{(|p_1|+\cdots
+|p_b|+|p'_1|+\cdots+|p'_{\alpha-1}|)|d|}  \mu((p_1,\ldots,\, p_b)
\, \sigma \, (p'_1,\ldots,\, \delta(p'_\alpha),\ldots,\,p'_a)).
\end{eqnarray*}
\end{dei}

Par exemple, la diff\'erentielle $\delta$ d'une prop\'erade
diff\'erentielle est une d\'erivation homog\`ene de degr\'e $-1$
telle que $\delta^2=0$. Remarquons que si $d$ est une d\'erivation
homog\`ene de degr\'e $-1$ alors la somme $\delta + d$ est encore
une d\'erivation de m\^eme degr\'e. Et dans ce cas, comme
$\delta^2=0$, l'\'equation $(\delta + d)^2=0$ est \'equivalente
\`a $\delta d + d \delta + d^2=0$.\\

\textsc{Remarque :} Dualement, pour une coprop\'erade
diff\'erentielle $\mathcal{C}$, on d\'efinit la notion de
\emph{cod\'erivation}.\index{cod\'erivation} \\

On s'int\'eresse maintenant \`a la forme des d\'erivations sur la
prop\'erade diff\'erentielle libre $\F(V)$.

\begin{lem}
\label{derivation}
 Soit $\F(V)$ la prop\'erade libre sur le
dg-$\Sy$-bimodule $V$. Pour tout morphisme homog\`ene $\theta \, :
\, V \to \F(V)$, il existe une unique d\'erivation homog\`ene de
m\^eme degr\'e $d_\theta\, : \, \F(V) \to \F(V)$ telle que sa
restriction \`a $V\subset \F(V) $ corresponde \`a $\theta$. Cette
correspondance est bijective. De plus, si $\theta\,  : \, V \to
\F_{(r)}(V)$ alors, on a $d_\theta\left( \F_{(s)}(V)\right)\subset
\F_{(r+s-1)}(V)$.
\end{lem}

\begin{deo}
Soit $g$ un graphe \`a $n$ sommets. On d\'efinit $d_\theta$ sur un
repr\'esentant d'une classe d'\'equivalence pour la relation
$\equiv$ (d\'eplacement vertical des sommets au signe pr\`es)
associ\'e au graphe $g$. Cela revient \`a choisir un ordre pour
\'ecrire les \'el\'ements de $V$ qui indicent les sommets de $g$.
Ecrivons $\bigotimes_{i=1}^n \nu_i$ ce repr\'esentant, avec
$\nu_i\in V$. On pose alors
$$d_\theta(\bigotimes_{i=1}^n \nu_i)=\sum_{i=1}^n (-1)^{|\nu_1|
+\cdots+|\nu_{i-1}|}\left( \bigotimes_{j=1}^{i-1}\nu_j
\right)\otimes \theta(\nu_i) \otimes \left( \bigotimes_{j=i+1}^{n}
\nu_j \right).$$ Un calcul rapide montre que cette application est
bien constante sur la classe d'\'equivalence pour la relation
$\equiv$ (la v\'erification est du m\^eme que celle du
lemme~\ref{differentielle} et fait appel aux signes d\'efinis pour
la relation $\equiv$ dans le cadre diff\'erentiel). En fait,
l'application $d_\theta$ revient \`a effectuer $\theta$ \`a chaque
\'el\'ement de $V$ indi\c cant un sommet du graphe $g$. (Comme
$\theta$ est un morphisme de $\Sy$-bimodules, entre autre,
l'application $d_\theta$ passe bien au quotient pour la relation
$\approx$.)

La surjectivit\'e du produit $\mu$ sur la prop\'erade libre
$\F(V)$ montre l'unicit\'e de la d\'erivation $d_\theta$ ainsi que
le fait que toute d\'erivation soit de cette forme.

Enfin, si l'application $\theta$ associe \`a un \'el\'ement de
$V$, des \'el\'ements de $\F(V)$ qui s'\'ecrivent avec des graphes
\`a $r$ sommets ($\theta$ cr\'ee $r-1$ sommets), alors on voit de
la forme de $d_\theta$ donn\'ee pr\'ec\'edemment que
$d_\theta\left( \F_{(s)}(V)\right)\subset \F_{(r+s-1)}(V)$.
$\cqfd$
\end{deo}

On a le lemme dual dans le cadre des coprop\'erades.

\begin{lem}
\label{lemmecodérivation}
 Soit $\F^c(V)$ la coprop\'erade colibre connexe
sur le dg-$\Sy$-bimodule $V$. Pour tout morphisme homog\`ene
$\theta \, : \, \F^c(V) \to V$, il existe une unique
cod\'erivation homog\`ene de m\^eme degr\'e $d_\theta\, : \,
\F^c(V) \to \F^c(V)$ telle que sa projection sur $V$ corresponde
\`a $\theta$. Cette correspondance est bijective. De plus, si
$\theta\,  : \, \F^c(V) \to V$ est nulle sur toutes les
composantes $\F^c_{(s)}(V)\subset \F^c(V)$, pour $s\ne r$, alors,
on a $d_\theta\left( \F_{(s+r-1)}(V)\right)\subset \F_{(s)}(V)$,
pour tout $s>0$.
\end{lem}

Graphiquement, effectuer la cod\'erivation $d_\theta$ sur un
\'el\'ement de $\F^c(V)$ repr\'esent\'e par un graphe $g$ revient
\`a appliquer $\theta$ \`a tous les sous-graphes possibles de $g$.\\

On peut maintenant donner la d\'efinition d'une prop\'erade
quasi-libre.

\begin{dei}(Prop\'erade quasi-libre)
\index{prop\'erade quasi-libre}

Une \emph{prop\'erade quasi-libre} $\Po$ est une prop\'erade
diff\'erentielle dont le $\Sy$-bimodule gradu\'e sous-jacent est
de la forme $\Po=\F(V)$ mais dont la diff\'erentielle
$\delta_\theta \, : \, \F(V) \to \F(V)$ est \'egale \`a la somme
de la diff\'erentielle canonique $\delta$ sur la prop\'erade libre
avec une d\'erivation $d_\theta$, soit $\delta_\theta = \delta +
d_\theta$.
\end{dei}

Gr\^ace au lemme pr\'ec\'edent, on sait que toute d\'erivation
$d_\theta$ est d\'etermin\'ee par sa restriction sur $V$, $\theta
\, : \, V \to \F(V)$. Ainsi, l'inclusion $V \to \F(V)$ est un
morphisme de dg-$\Sy$-bimodules si et seulement si $\theta$ est
nul.  Par contre, la condition pour que la projection $\F(V) \to
V$ soit un morphisme de dg-$\Sy$-bimodules est moins restrictive.

\begin{pro}
La projection $\F(V) \to V$ est un morphisme de dg-$\Sy$-bimodules
si et seulement si $\theta(V) \subset \bigoplus_{r\ge
2}\F_{(r)}(V)$.
\end{pro}

Dans ce cas, on a $d_\theta\left( \F(V)\right)\subset
\bigoplus_{r\ge 2} \F_{(r)}(V)$, et alors le quotient
ind\'ecomposable de la prop\'erade quasi-libre $\Po=\F(V)$ est
isomorphe, en tant que dg-$\Sy$-bimodule, \`a $V$.

En dualisant, on a la notion de coprop\'erade quasi-colibre.

\begin{dei}[Coprop\'erade quasi-colibre]
\index{coprop\'erade quasi-colibre}

Une \emph{coprop\'erade quasi-colibre} $\mathcal{C}$ est une
coprop\'erade diff\'erentielle dont le $\Sy$-bimodule gradu\'e
sous-jacent est de la forme $\mathcal{C}=\F^c(V)$ mais dont la
diff\'erentielle $\delta_\theta \, : \, \F^c(V) \to \F^c(V)$ est
\'egale \`a la somme de la diff\'erentielle canonique $\delta$ sur
la coprop\'erade libre avec une cod\'erivation $d_\theta$, soit
$\delta_\theta = \delta + d_\theta$.
\end{dei}

Le lemme pr\'ec\'edent sur les cod\'erivations des coprop\'erades
colibres montrent que la projection $\F^c(V)\to V$ est un
morphisme de dg-$\Sy$-bimodules si et seulement si $\theta$ est
nulle. Alors que l'inclusion $V \to \F^c (V)$ est un morphisme de
dg-$\Sy$-bimodules si et seulement si $\theta$ est nulle sur
$\F^c_{(1)}(V)=V$.\\

\textsc{Exemples :} Les deux exemples fondamentaux de
coprop\'erades et prop\'erades quasi-libres sont donn\'es par la
bar construction $\F^c\left( \Sigma \oPo \right)$ et la cobar
construction $\F\left( \Sigma^{-1} \overline{\mathcal{C}} \right)$
r\'eduites. Ces deux constructions sont \'etudi\'ees en d\'etail
dans le chapitre suivant.

\subsection{PROPs quasi-libres}
On fait la m\^eme \'etude pour les PROPs quasi-libres.

\begin{dei}
Soit $\Po$ un PROP diff\'erentiel. Un morphisme $d \, :\, \Po \to
\Po$ homog\`ene de degr\'e $|d|$ de $\Sy$-bimodules est appel\'e
\emph{d\'erivation} si $d$ est une d\'erivation pour la
prop\'erade $U_c(\Po)$ et pour l'alg\`ebre $(\Po, \otimes,\,
conc)$. Cette derni\`ere condition s'\'ecrit
$$\delta(conc(p\otimes q)) =conc(\delta(p)\otimes q) +(-1)^{|p||d|}conc(p\otimes \delta(q)) .$$
\end{dei}

On a le m\^eme type de lemme que dans le cas des prop\'erades.

\begin{lem}
Soit $S_\otimes(\F(V))$ le PROP libre sur le
dg-$\Sy$-bimodule $V$. Pour tout morphisme homog\`ene $\theta \, :
\, V \to S_\otimes(\F(V))$, il existe une unique
d\'erivation homog\`ene de m\^eme degr\'e $d_\theta\, : \,
S_\otimes(\F(V)) \to S_\otimes(\F(V))$ telle que sa
restriction \`a $V\subset S_\otimes(\F(V)) $ corresponde \`a
$\theta$. Cette correspondance est bijective. De plus, si
$\theta\,  : \, V \to S_\otimes(\F(V))_{(r)}$ alors, on a
$d_\theta\left( S_\otimes(\F(V))_{(s)}\right)\subset
S_\otimes(\F(V))_{(r+s-1)}$.
\end{lem}

La d\'erivation $d_\theta$ appliqu\'ee \`a un \'el\'ement de
$S_\otimes(\F(V))$ repr\'esent\'e par un graphe $g$
s'\'ecrit avec une sommme indic\'ee par l'ensemble des sommets
$\nu$ du graphe $g$ de graphes o\`u on remplace l'op\'eration
plac\'ee au sommet $\nu$ par son image via $\theta$.

\begin{deo}
La d\'emonstration est la m\^eme que dans le cas des prop\'erades.
$\Box$
\end{deo}

Dualement, on la notion de \emph{cod\'erivation} sur les coPROPs
et le lemme suivant.

\begin{lem}
\label{lemmecodérivationPROP}
 Soit $S_\otimes(\F^c(V))$ le coPROP colibre
sur le dg-$\Sy$-bimodule $V$. Pour tout morphisme homog\`ene
$\theta \, : \, S_\otimes(\F^c(V)) \to V$, il existe une
unique cod\'erivation homog\`ene de m\^eme degr\'e $d_\theta\, :
\, S_\otimes(\F^c(V)) \to S_\otimes(\F^c(V))$ telle
que sa projection sur $V$ corresponde \`a $\theta$. Cette
correspondance est bijective. De plus, si $\theta\,  : \,
S_\otimes(\F^c(V)) \to V$ est nulle sur toutes les
composantes $S_\otimes(\F^c(V))_{(s)}\subset
S_\otimes(\F^c(V))$, pour $s\ne r$, alors, on a
$d_\theta\left( S_\otimes(\F(V))_{(s+r-1)}\right)\subset
S_\otimes(\F(V))_{(s)}$, pour tout $s>0$.
\end{lem}

Graphiquement, effectuer la cod\'erivation $d_\theta$ sur un
\'el\'ement de $S_\otimes(\F^c(V))$ repr\'esent\'e par un
graphe $g$ revient
\`a appliquer $\theta$ \`a tous les sous-graphes possibles de $g$.\\

La d\'efinition de PROP quasi-libre est semblable \`a celle de
prop\'erade quasi-libre.

\begin{dei}[PROP quasi-libre]
\index{PROP quasi-libre}

On appelle \emph{PROP quasi-libre}, un PROP diff\'erentiel $\Po$
dont le $\Sy$-bimodule gradu\'e sous-jacent est un PROP libre
$S_\otimes(\F(V))$ et dont la diff\'erentielle $\delta_\theta$ est
la somme de la diff\'erentielle canonique $\delta$ et d'une
d\'erivation $d_\theta$.
\end{dei}

Dualement, on a la notion de coPROP quasi-colibre.

\begin{dei}[coPROP quasi-colibre]
\index{coPROP quasi-colibre}

Un \emph{coPROP quasi-colibre} $\mathcal{C}$ est un coPROP
diff\'erentiel dont le $\Sy$-bimodule gradu\'e sous-jacent est de
la forme $S_\otimes(F^c(V))$ mais dont la diff\'erentielle
$\delta_\theta \, : \, S_\otimes(\F^c(V)) \to S_\otimes(\F^c(V))$
est \'egale \`a la somme de la diff\'erentielle canonique $\delta$
sur le coPROP colibre avec une cod\'erivation $d_\theta$, soit
$\delta_\theta = \delta + d_\theta$.
\end{dei}

\textsc{Exemples :} Les deux principaux exemples de PROPs et de
coPROPs quasi-libres sont donn\'es par la bar et la cobar
constructions r\'eduites qui font l'objet du chapitre suivant.

\chapter{Bar et cobar constructions}

\thispagestyle{empty}

On g\'en\'eralise ici les bar et cobar constructions des
alg\`ebres et des op\'erades aux prop\'erades et aux PROPs. Pour
cela, on \'etudie d'abord les propri\'et\'es des produits et
coproduits partiels. Puis, on d\'efinit la bar et la cobar
construction en commen\c cant par leurs versions r\'eduites pour
passer ensuite \`a la version \`a coefficients. Ces bar
constructions poss\`edent de bonnes propri\'et\'es homologiques,
elles permettront d'obtenir des r\'esolutions pour les
prop\'erades et les PROPs diff\'erentiels (\emph{cf.}
th\'eor\`eme~\ref{barcobarresolution} et le chapitre 7). Afin
d'\'etablir ces r\'esolutions, nous commen\c cons par montrer dans
ce chapitre que les bar et les cobar constructions augment\'ees
 sont acycliques.

\section{Produit et coproduit de composition partiel}

Avant de donner la d\'efinition des deux bar constructions, nous
rappelons ce qu'est la \emph{suspension} et la
\emph{d\'esuspension} d'un dg-$\Sy$-bimodule ainsi que les
propri\'et\'es v\'erifi\'ees par les \emph{produits} et les
{coproduits de composition partiels}.

\subsection{Suspension d'un dg-$\Sy$-bimodule}
Soit $\Sigma$ le $\Sy$-bimodule gradu\'e d\'efini par
$$\left\{ \begin{array}{ll}
\Sigma(0,\,0) = k.s &\textrm{o\`u $s$ est un \'el\'ement de
degr\'e +1,} \\
\Sigma(m,\, n)=0 & \textrm{sinon.}
\end{array}\right.$$

\begin{dei}[Suspension $\Sigma V$]
\index{suspension}

On appelle \emph{suspension} du dg-$\Sy$-bimodule $V$, le
dg-$\Sy$-bimodule $\Sigma V = \Sigma \otimes V$.
\end{dei}

Cette op\'eration revient \`a tensoriser par $s$ de tous les
\'el\'ements de $V$. A $v\in V_{d-1}$ on associe donc $s\otimes
v\in (\Sigma V)_d$ que l'on note souvent $\Sigma v$. Ainsi,
$(\Sigma V)_d$ est naturellement isomorphe \`a $V_{d-1}$. Selon
les principes \'enonc\'es au chapitre pr\'ec\'edent, la
diff\'erentielle sur $\Sigma V$ est donn\'ee par la formule
$\delta (\Sigma v)= -\Sigma \delta (v)$, pour tout $v$ dans $V$.
La suspension $\Sigma V$ correspond donc \`a l'introduction d'un
\'el\'ement de degr\'e $+1$ qu'il faut prendre en compte lors des
permutations faisant intervenir des signes (r\`egles de
Koszul-Quillen).\\

De la m\^eme mani\`ere, on d\'efinit $\Sigma^{-1}$ par
$$\left\{ \begin{array}{ll}
\Sigma^{-1}(0,\,0) = k.s^{-1} & \textrm{o\`u $s$ est un
\'el\'ement de
degr\'e -1,} \\
\Sigma^{-1}(m,\, n)=0 & \textrm{sinon.}
\end{array}\right.$$

\begin{dei}[D\'esuspension $\Sigma^{-1} V$]
\index{d\'esuspension}

On appelle \emph{d\'esuspension} du dg-$\Sy$-bi\-mo\-du\-le $V$,
le dg-$\Sy$-bi\-mo\-du\-le $\Sigma^{-1} V = \Sigma^{-1} \otimes
V$.
\end{dei}

\subsection{Produit de composition partiel}

Soit $(\Po,\, \mu,\, \eta,\, \varepsilon)$ une prop\'erade
augment\'ee. Comme les graphes \`a deux niveaux sont munis d'une
bigraduation en fonction du nombre de sommets sur chaque niveau,
on peut d\'ecomposer le produit de composition $\mu$ de la
mani\`ere suivante $\mu=\bigoplus_{r,\, s \, \in
 \mathbb{N}} \mu_{(r,\,s)}$, o\`u
 $$\mu_{(r,\,s)} \ : \ (I\oplus \underbrace{\oPo}_{r}) \boxtimes_c
 (I \oplus \underbrace{\oPo}_{s}) \to \oPo.$$

\begin{dei}[Produit de composition partiel]
\index{produit de composition partiel}

On appelle \emph{produit de composition partiel} la restriction du
produit de composition $\mu$ \`a $ (I\oplus
\underbrace{\oPo}_{1})\boxtimes_c (I\oplus \underbrace{\oPo}_{1})$.
Le produit de composition partiel correspond donc au $\mu_{(1,\,
1)}$ pr\'ec\'edent.
\end{dei}

On peut remarquer que le produit de composition partiel correspond
au produit ne faisant intervenir que deux \'el\'ements non
triviaux de $\Po$.\\

\textsc{Remarque :} Dans le cadre des op\'erades \cite{GK}, ce
produit de composition partiel n'est autre que le produit partiel
not\'e $\circ_i$. Une particularit\'e fondamentale des op\'erades
est que ce produit partiel engendre le produit global $\circ$
(toute composition peut s'\'ecrire avec un nombre fini de
compositions partielles successives).

Dans le cas des diop\'erades (\emph{cf.} \cite{Gan}), l'auteur ne
consid\`ere que les produits de composition partiels ne faisant
intervenir qu'une seule branche entre chaque sommet. Et dans le
cas des $\frac{1}{2}$-PROPs (\emph{cf.} \cite{MV}), les auteurs
restreignent encore les produits partiels aux couples dont au
moins une op\'eration n'admet qu'une seule entr\'ee ou sortie.
Deux ces deux cas, le produit mono\"\i dal \'etudi\'e est celui
engendr\'e par les produits partiels.

Notons que le produit de composition des prop\'erades est aussi
engendr\'e par les produits partiels (\emph{cf.} chapitre $6$).

\begin{lem}
\label{thetainduitproduit}
 Si $(\Po,\, \mu,\, \eta,\,
\varepsilon)$ est une prop\'erade diff\'erentielle augment\'ee,
alors le produit de composition partiel $\mu_{(1,\, 1)}$ induit un
morphisme homog\`ene de degr\'e $-1$
$$\theta \ : \ \F_{(2)}^c(\Sigma \oPo) \to \Sigma \oPo.$$
\end{lem}

Rappelons que si $\Po$ est un $\Sy$-bimodule diff\'erentiel, la
prop\'erade libre $\F(\Sigma \oPo)$ sur $\Sigma \oPo$ est
bigradu\'ee. La premi\`ere graduation vient du nombre
d'op\'erations de $\Sigma \oPo$ utilis\'ees pour repr\'esenter un
\'el\'ement de $\F(\Sigma \oPo)$. On note cette graduation
$\F_{(n)}(\Sigma \oPo)$. La deuxi\`eme graduation, ici le degr\'e
homologique, est obtenu en effectuant la somme des degr\'es des
op\'erations de $\Sigma \oPo$.

\begin{deo}
Soit $(1,\ldots,\, 1,\, \Sigma q,\, 1,\ldots ,\,
1)\sigma(1,\ldots,\, 1,\, \Sigma p,\, 1,\ldots ,\, 1)$ un
\'el\'ement de degr\'e homologique $|q|+|p|+2$ de
$\F_{(2)}^c(\Sigma \oPo)$. On d\'efinit $\theta$ par
\begin{eqnarray*}
&& \theta\left((1,\ldots,\, 1,\, \Sigma q,\, 1,\ldots ,\,
1)\sigma(1,\ldots,\, 1,\, \Sigma p,\, 1,\ldots ,\, 1)\right)=\\
&& (-1)^{|q|}\Sigma\mu_{(1,\, 1)}(1,\ldots,\, 1,\,  q,\, 1,\ldots
,\, 1)\sigma(1,\ldots,\, 1,\, p,\, 1,\ldots ,\, 1).
\end{eqnarray*}
Et comme $\Po$ est une prop\'erade gradu\'ee, ce dernier
\'el\'ement est de degr\'e $|q|+|p|+1$ soit un de moins que son
ant\'ec\'edent. $\cqfd$
\end{deo}

Au morphisme homog\`ene $\theta$, on peut associer une
cod\'erivation $d_\theta \, : \,\F^c(\Sigma \oPo) \to \F^c(\Sigma
\oPo) $ gr\^ace au lemme~\ref{lemmecodérivation}.

\begin{pro}
\label{propriétés-codérivation}
 La cod\'erivation $d_\theta$
v\'erifie les propri\'et\'es suivantes :
\begin{enumerate}
\item L'\'equation $\delta d_\theta +d_\theta \delta =0 $ est
vraie si et seulement si le produit de composition partiel
$\mu_{(1,\, 1)}$ est un morphisme de dg-$\Sy$-bimodules. \item On
a toujours ${d_\theta}^2=0$.
\end{enumerate}
\end{pro}

\begin{deo}
$\ $
\begin{enumerate}
\item D'apr\`es le lemme~\ref{lemmecodérivation}, appliquer
$d_\theta$ \`a un \'el\'ement de $\F^c(\Sigma\oPo)$ revient \`a
effectuer $\theta$ \`a tous les sous-graphes possibles du graphe
repr\'esentant cet \'el\'ement. Ici, il suffit donc d'appliquer
$\theta$ \`a tous les couples de sommets reli\'es par au moins une
branche et n'admettant aucun sommet interm\'ediaire. Or pour un
couple de sommets indic\'es par $\Sigma q$ et $\Sigma p$, on a
d'une part
\begin{eqnarray*}
&& \delta \circ d_\theta \left( (1,\ldots,\, 1,\, \Sigma q,\,
1,\ldots ,\, 1)\sigma(1,\ldots,\, 1,\, \Sigma p,\, 1,\ldots ,\, 1)
\right)  = \\
&&(-1)^{|q|}\delta \left( \Sigma \mu_{(1,\ 1)}((1,\ldots,\, 1,\,
q,\, 1,\ldots ,\, 1)\sigma(1,\ldots,\, 1,\, p,\, 1,\ldots ,\, 1))
\right)  = \\
&& (-1)^{|q|+1}\Sigma \, \delta \left( \mu_{(1,\ 1)}((1,\ldots,\,
1,\, q,\, 1,\ldots ,\, 1)\sigma(1,\ldots,\, 1,\, p,\, 1,\ldots ,\,
1)) \right)
\end{eqnarray*}
et d'autre part
\begin{eqnarray*}
&& d_\theta \circ \delta \left( (1,\ldots,\, 1,\, \Sigma q,\,
1,\ldots ,\, 1)\sigma(1,\ldots,\, 1,\, \Sigma p,\, 1,\ldots ,\, 1)
\right)  = \\
&& d_\theta \big( -(1,\ldots,\, 1,\, \Sigma \delta(q),\, 1,\ldots
,\, 1)\sigma(1,\ldots,\, 1,\, \Sigma p,\, 1,\ldots ,\, 1) +
 \\
 &&(-1)^{|q|} (1,\ldots,\, 1,\, \Sigma q,\, 1,\ldots ,\,
1)\sigma(1,\ldots,\, 1,\, \Sigma \delta(p),\, 1,\ldots ,\, 1)
\big)= \\
&& (-1)^{|q|}\Sigma \mu_{(1,\, 1)}((1,\ldots,\, 1,\, \delta(q),\,
1,\ldots ,\, 1)\sigma(1,\ldots,\, 1,\, p,\, 1,\ldots ,\, 1))+ \\
&& \Sigma \mu_{(1,\, 1)} ((1,\ldots,\, 1,\, q,\, 1,\ldots ,\,
1)\sigma(1,\ldots,\, 1,\, \delta(p),\, 1,\ldots ,\, 1)).
\end{eqnarray*}
Ainsi, $\delta d_\theta +d_\theta \delta =0$ implique
\begin{eqnarray*}
&&\delta \left( \mu_{(1,\, 1)} ((1,\ldots,\, 1,\,  q,\, 1,\ldots
,\, 1)\sigma(1,\ldots,\, 1,\,  p,\, 1,\ldots ,\, 1))   \right) = \\
&&\mu_{(1,\, 1)}( (1,\ldots,\, 1,\,  \delta(q),\, 1,\ldots ,\,
1)\sigma(1,\ldots,\, 1,\, p,\, 1,\ldots ,\, 1)) + \\
&& (-1)^{|q|}\mu_{(1,\, 1)}( (1,\ldots,\, 1,\,  q,\, 1,\ldots ,\,
1)\sigma(1,\ldots,\, 1,\, \delta(p),\, 1,\ldots ,\, 1)),
\end{eqnarray*}
c'est-\`a-dire que le produit de composition partiel
$\mu_{(1,\,1)}$ est un morphisme de dg-$\Sy$-bimodules. Dans
l'autre sens, on conclut en remarquant qu'effectuer $\delta
d_\theta +d_\theta \delta$ revient \`a faire une somme
d'expressions du type pr\'ec\'edent.\\

\item Pour montrer que ${d_\theta}^2=0$, il faut s'int\'eresser
aux paires de couples distincts de sommets d'un graphe. Deux cas
de figure sont possibles : soit on a affaire \`a deux couples de
sommets dont les quatres sommets sont disctincts (a), soit les
deux couples ont un sommet en commun (b).
\begin{enumerate}
\item Sur un \'el\'ement de la forme
$$X\otimes \Sigma q_1 \otimes
\Sigma p_1 \otimes Y \otimes \Sigma q_2 \otimes \Sigma p_2 \otimes
Z, $$ on effectue $\theta$ deux fois en commen\c cant par un
couple diff\'erent \`a chaque fois, ce qui donne
\begin{eqnarray*}
&&\big( (-1)^{|X|+|q_1|}(-1)^{|X|+|q_1|+|p_1|+1+|Y|+|q_2|} +
(-1)^{|X|+|q_1|+|p_1|+|Y|+|q_2|} (-1)^{|X|+|q_1|} \big) \\
&&X\otimes \Sigma \mu_{(1,\, 1)}(q_1\otimes p_1) \otimes Y \otimes
\Sigma \mu_{(1,\, 1)}(q_2\otimes p_2)\otimes Z =0.
\end{eqnarray*}

\item Dans ce cas, deux configurations sont possibles
\begin{itemize}
\item Le sommet commun peut \^etre entre les deux autres sommets
(dans le sens du graphe) (\emph{cf.} figure~\ref{entresommets} ).

\begin{figure}[h]
$$ \xymatrix{ & & \\
 & *+[F-,]{\Sigma p}  \ar@{.}[u] \ar@{.}[ul]\ar@{.}[ur] \ar@/_/@{-}[d]\ar@/^/@{-}[d]& \\
 & *+[F-,]{\Sigma q} \ar@{-}[d]& \\
 & *+[F-,]{\Sigma r} \ar@{.}[ul]\ar@{.}[ur]\ar@{.}[dl]\ar@{.}[dr]& \\
 & &} $$ \caption{}  \label{entresommets}
\end{figure}

Alors, sur un \'el\'ement de la forme
$$X\otimes \Sigma r \otimes \Sigma q \otimes \Sigma p \otimes Y,$$
si on applique $\theta$ deux fois en commen\c cant par un couple
diff\'erent \`a chaque fois, on trouve
\begin{eqnarray*}
&& (-1)^{|X|+|r|}(-1)^{|X|+|r|+|q|} X\otimes \Sigma \mu_{(1,\,
1)}(\mu_{(1,\, 1)}(r\otimes q)\otimes p)\otimes Y + \\
&& (-1)^{|X|+|r|+1+|q|}(-1)^{|X|+|r|} X\otimes \Sigma
\mu_{(1,\,1)} (r\otimes \mu_{(1,\,1)}(q\otimes p))\otimes Y=0,
\end{eqnarray*}
par associativit\'e du produit partiel $\mu_{(1,\,1)}$ (qui vient
de celle de $\mu$).

\item Sinon, le sommet commun est au-dessus (\emph{cf.}
figure~\ref{audessussommets}) ou en dessous des deux autres.

\begin{figure}[h]
$$ \xymatrix{& & & \\
& &*+[F-,]{\Sigma p} \ar@{.}[u] \ar@{.}[ul] \ar@{.}[ur] \ar@{.}[d] & \\
& *+[F-,]{\Sigma r}
\ar@{-}[ur]\ar@{.}[u]\ar@{.}[ul]\ar@{.}[d]\ar@{.}[dl]  & *=0{}
\ar@{.}[d]&
*+[F-,]{\Sigma q} \ar@{-}@/_/[ul]\ar@{-}@/^/[ul]\ar@{.}[u] \ar@{.}[d] \\
 & & & } $$ \caption{}  \label{audessussommets}
\end{figure}
On fait le m\^eme calcul que pr\'ec\'edemment pour un \'el\'ement
de cette forme pour obtenir
\begin{eqnarray*}
&&(-1)^{|X|+|r|+1+|q|}(-1)^{|X|+|r|}X\otimes \Sigma
\mu_{(1,\,1)}(r\otimes
\mu_{(1,\,1)}(q\otimes p))\otimes Y+ \\
&&(-1)^{(|r|)+1(|q|+1)+|X|+|q|+1+|r|}(-1)^{|X|+|q|}X\otimes
\Sigma\mu_{(1,\, 1)}(q\otimes \mu_{(1,\, 1)} (r\otimes p))\otimes
Y=0.
\end{eqnarray*}
En effet, l'associativit\'e de $\mu_{(1,\, 1)}$ donne
$$\mu_{(1,\, 1)}(q\otimes \mu_{(1,\, 1)} (r\otimes p))= (-1)^{|r||q|}
\mu_{(1,\,1)}(r\otimes \mu_{(1,\,1)}(q\otimes p)). $$
\end{itemize}
\end{enumerate}
On conclut en remarquant que le r\'esultat de l'application
${d_\theta}^2$ est une somme d'expressions des deux formes
pr\'ec\'edentes.$\cqfd$
\end{enumerate}
\end{deo}

\textsc{Remarque :} Ce sont les r\`egles naturelles sur les signes
qui permettent d'avoir cette proposition. Dans les articles de
\cite{GK} et \cite{Gan}, ces signes sont exprim\'es sous la forme
de l'op\'erade (resp. diop\'erade) $Det$.

\begin{dei}[Sommets adjacents]
\index{sommets adjacents}

On appelle \emph{sommets adjacents} tout couple de sommets d'un graphe
reli\'es par au moins une branche et n'admettant aucun sommet interm\'ediaire.
Ils correspondent \`a des sous-graphes de la forme $\F_{(2)}(V)$ et sont
les couples composables par la cod\'erivation $d_\theta$.
\end{dei}

\textsc{Remarque :} La cod\'erivation $d_\theta$ consiste \`a
composer les couples de sommets adjacents. Pour d\'efinir
l'homologie des graphes, M. Kontsevich avait introduit dans
\cite{Ko} la notion d'\emph{edge contraction}. Cette notion
revient \`a composer deux sommets reli\'es par une seule branche.
Dans le cas des op\'erades (arbres) et des diop\'erades (graphes
de genre $0$), les sommets adjacents sont reli\'es entre eux par
une seule branche, on retrouve alors cette notion (\emph{cf.}
\cite{GK} et \cite{Gan}). La cod\'erivation $d_\theta$ est donc la
g\'en\'eralisation naturelle de cette
notion d'edge contraction. \index{edge contraction} \\

\subsection{Coproduit partiel}

On peut dualiser les r\'esultats de la partie pr\'ec\'edente.

\begin{dei}[Coproduit partiel]
\index{coproduit partiel}

Soit $(\mathcal{C},\, \Delta,\, \varepsilon,\, \eta)$ une
coprop\'erade coaugment\'ee. On appelle \emph{coproduit partiel}
la composition d'applications
$$\overline{\mathcal{C}}  \xrightarrow{\Delta } \mathcal{C}\boxtimes_c
 \mathcal{C} \twoheadrightarrow
(I\oplus \underbrace{\overline{\mathcal{C}}}_1) \boxtimes_c
(I\oplus \underbrace{\overline{\mathcal{C}}}_1),$$ not\'ee
$\Delta_{(1,\,1)}$.
\end{dei}

\begin{lem}
\label{thetainduitcoproduit} Pour toute coprop\'erade
diff\'erentielle coaugment\'ee $(\mathcal{C},\, \Delta ,\,
\varepsilon,\, \eta)$, le coproduit partiel induit un morphisme
homog\`ene de degr\'e $-1$
$$\theta' \ : \ \Sigma^{-1} \oC \to \F_{(2)}(\Sigma^{-1} \oC).$$
\end{lem}

\begin{deo}
Soit $c$ un \'el\'elment de $\oC$. En s'inspirant des notations de
Sweedler, posons
$$\Delta_{(1,\,1)}(c)=\sum_{(c',\, c'')}
(1,\ldots,\, 1,\,  c',\, 1,\ldots ,\, 1)\sigma(1,\ldots,\, 1,\,
 c'',\, 1,\ldots ,\, 1).$$
Alors $\theta'(\Sigma^{-1}c)$ d\'efini par
$$\theta'(\Sigma^{-1}c)=-\sum_{(c',\, c'')} (-1)^{|c'|}
(1,\ldots,\, 1,\,  \Sigma^{-1}\, c',\, 1,\ldots ,\,
1)\sigma(1,\ldots,\, 1,\, \Sigma^{-1}\,c'',\, 1,\ldots ,\, 1),$$
convient. $\cqfd$
\end{deo}

\textsc{Remarque :} Le signe $-$ apparaissant devant le symbole $\Sigma$ dans
 la derni\`ere expression n'est pas essentiel ici. Il deviendra par
contre fondamental dans la d\'emonstration de la r\'esolution
bar-cobar(\emph{cf.} th\'eor\`eme~\ref{barcobarresolution}).\\

Gr\^ace au lemme~\ref{derivation}, au morphisme $\theta'$ on peut
associer une d\'erivation $d_{\theta'} \, :\, \F(\Sigma^{-1} \oC)\to
 \F(\Sigma^{-1} \oC)$ de degr\'e $-1$ v\'erifiant la proposition suivante :

\begin{pro}
$ \ $
\begin{enumerate}
\item L'\'equation $\delta d_{\theta'} +d_{\theta'} \delta =0 $ est
vraie si et seulement si le coproduit partiel $\Delta_{(1,\, 1)}$
est un morphisme de dg-$\Sy$-bimodules. \item On a toujours
${d_{\theta'}}^2=0$.
\end{enumerate}
\end{pro}

\begin{deo}
La d\'emonstration du lemme~\ref{derivation} montre qu'effectuer
$d_{\theta'}$ sur un \'el\'ement de $\F(\Sigma^{-1} \oC)$ revient
\`a appliquer $\theta'$ \`a chaque sommet du graphe repr\'esentant
le dit \'el\'ement. Ainsi, les arguments pour montrer cette
proposition sont du m\^eme type que pour la proposition
pr\'ec\'edente. Notamment le deuxi\`eme point  vient de la
coassociativit\'e de $\Delta_{(1,\, 1)}$ qui vient de celle de
$\Delta$ et des r\`egles de signes. $\cqfd$
\end{deo}

\textsc{Remarque :} La d\'erivation $d_{\theta'}$ revient \`a
effectuer le coproduit partiel sur toutes les op\'erations indi\c
cant un sommet. Cette d\'erivation est la g\'en\'eralisation
naturelle de la notion de \emph{vertex expansion} introduite par
M. Kontsevich \cite{Ko} dans le cadre de la cohomologie des
graphes. Dans le cas des alg\`ebres et des op\'erades, la
d\'erivation $d_{\theta'}$ correspond \`a la "vertex expansion"
des arbres. \index{vertex expansion}

\section{Bar et cobar constructions} A l'aide
des deux propositions pr\'ec\'edentes, on peut maintenant
d\'efinir la bar et la cobar construction des prop\'erades. Puis nous
g\'en\'eralisons ces deux constructions au cadre des PROPs.

\subsection{Bar et cobar constructions r\'eduites}

\begin{dei}[Bar construction r\'eduite]
\index{bar construction r\'eduite}

A tout prop\'erade diff\'erentielle augment\'ee $(\Po,\, \mu,\,
\eta,\, \varepsilon)$, on peut associer la coprop\'erade
quasi-colibre $\bar{\B}(\Po) = \F^c(\Sigma \oPo)$ munie de la
diff\'erentielle $\delta_\theta=\delta +d_\theta$, o\`u $\theta$
est le morphisme induit par le produit partiel sur $\Po$
(\emph{cf.} lemme~\ref{thetainduitproduit}). Cette coprop\'erade
quasi-colibre est appel\'ee \emph{bar construction r\'eduite de
$\Po$ }.
\end{dei}

L'\'egalit\'e ${\delta_\theta}^2=0$ vient de $\delta d_\theta
+d_\theta \delta =0$ et de ${d_\theta}^2=0$.

Si on pose $\bar{\B}_{(s)}(\Po) = \F^c_{(s)}(\Sigma \oPo)$, alors
$d_\theta$ d\'efinit le complexe
$$\xymatrix{\cdots \ar[r]^(0.4){d_\theta} & \bar{\B}_{(s)}(\Po)
\ar[r]^(0.5){d_\theta} & \bar{\B}_{(s-1)}(\Po)
\ar[r]^(0.6){d_\theta} & \cdots \ar[r]^(0.4){d_\theta} &
\bar{\B}_{(1)}(\Po) \ar[r]^(0.5){d_\theta}&
\bar{\B}_{(0)}(\Po).}$$

En raisonnant de mani\`ere duale, on obtient la d\'efinition
suivante :

\begin{dei}[Cobar construction r\'eduite]
\index{cobar construction r\'eduite}

A tout coprop\'erade diff\'erentielle coaugment\'ee
$(\mathcal{C}$, $\Delta$, $\varepsilon$, $\eta)$, on peut associer
la prop\'erade quasi-libre $\bar{\B}^c(\mathcal{C}) =
\F(\Sigma^{-1} \oC)$ munie de la diff\'erentielle
$\delta_{\theta'}=\delta +d_{\theta'}$, o\`u $\theta'$ est le
morphisme induit par le coproduit partiel sur $\mathcal{C}$
(\emph{cf.} lemme~\ref{thetainduitcoproduit}). Cette prop\'erade
quasi-libre est appel\'ee \emph{cobar construction r\'eduite de
$\mathcal{C}$ }.
\end{dei}

A nouveau $d_{\theta'}$ d\'efinit un complexe
$$\xymatrix{\bar{\B}^c_{(0)}(\mathcal{C}) \ar[r]^{d_{\theta'}} &
\bar{\B}^c_{(1)}(\mathcal{C}) \ar[r]^(0.6){d_{\theta'}} & \cdots
\ar[r]^(0.45){d_{\theta'}}&
\bar{\B}^c_{(s)}(\mathcal{C})\ar[r]^(0.45){d_{\theta'}}&
\bar{\B}^c_{(s+1)}(\mathcal{C}) \ar[r]^(0.6){d_{\theta'}}& \cdots
.}$$

\textsc{Remarque :} Lorsque $\Po$ (resp. $\mathcal{C}$) est une
alg\`ebre ou une op\'erade (resp. cog\`ebre ou une coop\'erade),
on retombe sur les d\'efinitions classiques (\emph{cf.}
S\'eminaire H. Cartan \cite{Cartan} et \cite{GK}).

\subsection{Bar construction \`a coefficients}

Soit $(\Po, \, \mu ,\, \eta,\, \varepsilon)$ une prop\'erade
diff\'erentielle augment\'ee. Sur $\bar{\B}(\Po)$, on d\'efinit
deux morphismes homog\`enes $\theta_r\, :\,
\bar{\B}(\Po)\to\bar{\B}(\Po)\boxtimes_c \Po $ et $\theta_l\, :\,
\bar{\B}(\Po)\to \Po \boxtimes_c \bar{\B}(\Po)$ de degr\'e $-1$ par
$$\theta_r \ : \ \bar{\B}(\Po)=\F^c(\Sigma
\oPo)\xrightarrow{\Delta}  \F^c(\Sigma \oPo)\boxtimes_c
\F^c(\Sigma \oPo) \twoheadrightarrow \F^c(\Sigma \oPo)\boxtimes_c
(I \oplus \underbrace{\oPo}_1),$$ et par
$$\theta_l \ : \ \bar{\B}(\Po)=\F^c(\Sigma
\oPo)\xrightarrow{\Delta}  \F^c(\Sigma \oPo)\boxtimes_c
\F^c(\Sigma \oPo) \twoheadrightarrow (I \oplus
\underbrace{\oPo}_1)\boxtimes_c \F^c(\Sigma \oPo).$$

Remarquons que ces deux morphismes reviennent \`a extraire, un par
un, les sommets extr\'emaux (en haut et en bas) d'un graphe
repr\'esentant un \'el\'ement de $\F^c(\Sigma \Po)$.

\begin{lem}
\label{lemme-bar-construction-augmentée} Le morphisme $\theta_r$
induit un morphisme homog\`ene $d_{\theta_r}$ de degr\'e $-1$
$$d_{\theta_r}\ : \ \bar{\B}(\Po)\boxtimes_c \Po \to
\bar{\B}(\Po)\boxtimes_c \Po .$$ Et le dg-$\Sy$-bimodule
$\bar{\B}(\Po)\boxtimes_c \Po$ muni de la somme de la
diff\'erentielle canonique $\delta$ avec la co\-d\'e\-ri\-va\-tion
$d_\theta$ d\'efinie sur $\bar{\B}(\Po)$ et du morphisme
$d_{\theta_r}$ est un $\Po$-module quasi-libre analytique (\`a
droite).
\end{lem}

\begin{deo}
Le d\'emonstration du lemme d\'ecoule principalement des
d\'efinitions du chapitre $3$. Seul point difficile, l'\'equation
$(\delta +d_\theta +d_{\theta_r})^2=0$. On sait d\'ej\`a, gr\^ace
aux propositions de la section pr\'ec\'edente que
$(\delta+d_\theta)^2=0$. Quant aux \'equations
$$\left\{  \begin{array}{l}
(\delta+d_\theta)d_{\theta_r}+d_{\theta_r}(\delta+d_\theta)=0 \quad \textrm{et}\\
{d_{\theta_r}}^2=0,
\end{array} \right. $$
 elles se d\'emontrent de la m\^eme mani\`ere que la
proposition~\ref{propriétés-codérivation} en \'ecrivant
soigneusement les r\`egles sur les signes. $\cqfd$
\end{deo}

On peut remarquer que le morphisme $d_{\theta_r}$ revient \`a
\'ecr\^eter les sommets se situant en haut des graphes
repr\'esentant des \'el\'ements de $\bar{\B}(\Po)$ puis \`a
composer les op\'erations de $\Po$ ainsi obtenues avec celle de
$\Po$ issues de la premi\`ere ligne de $\bar{\B}(\Po) \boxtimes
\Po$.

\begin{dei}[Bar construction augment\'ee \`a droite]
\index{bar construction augment\'ee}

Le $\Po$-module quasi-libre analytique $\bar{\B}(\Po)\boxtimes_c
\Po$ est appel\'e \emph{bar construction augment\'ee \`a droite}.
\end{dei}

On peut faire exactement le m\^eme raisonnement pour $\theta_l$,
$d_{\theta_l}$ et la bar construction augment\'ee \`a gauche $\Po
\boxtimes_c \bar{\B}(\Po)$.

Nous d\'emontrons dans la prochaine section que, comme dans le cas
des alg\`ebres et des op\'erades,  les deux  bar constructions
augment\'ees (\`a gauche et \`a droite) sont toujours acycliques.\\

Posons $\B(\Po,\, \Po,\, \Po)=\Po\boxtimes_c \bar{\B}(\Po) \boxtimes_c
\Po$. De la m\^eme mani\`ere que pr\'ec\'edemment, les morphismes
$\theta_r$ et $\theta_l$ induisent sur $\B(\Po,\, \Po,\, \Po)$
deux morphismes homog\`enes de degr\'e $-1$ :
$$d_{\theta_r} ,\, d_{\theta_l}\ : \\B(\Po,\, \Po,\, \Po) \to \B(\Po,\, \Po,\, \Po) .$$
Ainsi, on d\'efinit sur $\B(\Po,\, \Po,\, \Po)$ une
diff\'erentielle $d$ par la somme de plusieurs termes:
\vspace{6pt}

\begin{tabular}{l}
$\bullet \quad$ la diff\'erentielle canonique $\delta$ induite par celle de $\Po$,\\
$\bullet\quad$ la cod\'erivation $d_\theta$ d\'efinie sur $\bar{\B}(\Po)$,\\
$\bullet\quad$ du morphisme homog\`ene $d_{\theta_r}$ de degr\'e $-1$,\\
$\bullet\quad$ du morphisme homog\`ene $d_{\theta_l}$ de degr\'e
$-1$.
\end{tabular}

\begin{lem}
Le morphisme $d$ ainsi d\'efini v\'erifie l'\'equation $d^2=0$.
\end{lem}

\begin{deo}
Encore une fois, la d\'emonstration repose sur les r\`egles de
signes. Les calcus sont du m\^eme style que ceux de la
proposition~\ref{propriétés-codérivation}. $\cqfd$
\end{deo}

\begin{dei}[Bar construction \`a coefficients]
\index{bar construction \`a coefficients}

Soient $L$ un $\Po$-module diff\'erentiel \`a droite et $R$ un
$\Po$-module diff\'erentiel \`a gauche. A toute prop\'erade
diff\'erentielle augment\'ee $(\Po,\, \mu,\ \eta,\, \varepsilon)$,
on peut associer la \emph{bar construction \`a coefficients dans
$L$ et $R$} d\'efinie par le dg-$\Sy$-bimodule $\B(L,\, \Po,\,
R):= L{\boxtimes_c}_\Po \B(\Po,\, \Po,\, \Po) {\boxtimes_c}_\Po R$ muni de
la diff\'erentielle induite par $d$ et par celles de $L$ et $R$
not\'ees $\delta_L$ et $\delta_R$.
\end{dei}

\begin{pro}
La bar construction $\B(L,\, \Po,\, R)$ \`a coefficients dans les
modules $L$ et $R$ est isomorphe, en tant que dg-$\Sy$-bimodule,
\`a $L\boxtimes_c \bar{\B}(\Po) \boxtimes_c R$ muni de la
diff\'erentielle $d$, d\'efinie par la somme des termes suivants
\vspace{6pt}

\begin{tabular}{l}
$\bullet \quad$ la diff\'erentielle canonique $\delta$ induite par celle de $\Po$ sur $\bar{\B}(\Po)$,\\
$\bullet \quad$ la diff\'erentielle canonique $\delta_L$ induite par celle de $L$,\\
$\bullet \quad$ la diff\'erentielle canonique $\delta_R$ induite par celle de $R$,\\
$\bullet\quad$ la cod\'erivation $d_\theta$ d\'efinie sur $\bar{\B}(\Po)$,\\
$\bullet\quad$ d'un morphisme homog\`ene $d_{\theta_R}$ de degr\'e $-1$,\\
$\bullet\quad$ d'un morphisme homog\`ene $d_{\theta_L}$ de degr\'e
$-1$.
\end{tabular}
\end{pro}

\begin{deo}
Il s'agit de bien comprendre comment le morphisme $d_{\theta_r}$
donne, apr\`es passage au produit mono\"\i dal relatif, un
morphisme homog\`ene $d_{\theta_R}$. En fait, le morphisme
$d_{\theta_R}$ correspond \`a la composition
$$\xymatrix@C=60pt{ L\boxtimes_c \bar{\B}(\Po) \boxtimes_c R
\ar[r]^{\widetilde{\theta_r}}& L \boxtimes_c \bar{\B}(\Po)
\boxtimes_c \Po \boxtimes_c R \ar[r]^(0.55){L\boxtimes_c
\bar{\B}(\Po)\boxtimes_c r} & L\boxtimes_c \bar{\B}(\Po)
\boxtimes_c R,}$$ o\`u le morphisme $\widetilde{\theta_r}$ est
induit par $\theta_r$ de la mani\`ere suivante : sur un
\'el\'ement de $ L\boxtimes_c \bar{\B}(\Po) \boxtimes_c R$, que
l'on repr\'esente par
$$l_1 \otimes \cdots\otimes l_b \otimes b_1\otimes \cdots \otimes b_s
\otimes r_1 \otimes \cdots\otimes r_a,$$ le morphisme
$\widetilde{\theta_r}$ vaut
\begin{eqnarray*}
&&\sum_{i=1}^s (-1)^{|b_i|(|b_{i+1}|+\cdots+|b_s|)}
(-1)^{|l_1|+\cdots+|l_b|+|b_1|+\cdots+|b_{i-1}|+|b_{i+1}|+\cdots+|b_s|} \\
&&l_1 \otimes \cdots\otimes l_b \otimes b_1\otimes \cdots \otimes
b_{i-1} \otimes b_{i+1} \otimes \cdots \otimes b_s \otimes
\theta_r(b_i)\otimes r_1 \otimes \cdots\otimes r_b.
\end{eqnarray*}
Il en va de m\^eme pour $d_{\theta_L}$. $\cqfd$
\end{deo}

De cette expression de la bar construction \`a coefficients
d\'ecoule imm\'ediatement le corollaire suivant.

\begin{cor}
$\ $
\begin{itemize} \item La bar construction r\'eduite
$\bar{\B}(\Po)$ correspond \`a la bar construction avec des
coefficients trivaux $L=R=I$, c'est-\`a-dire
$\bar{\B}(\Po)=\B(I,\, \Po,\, I)$.

\item Le complexe de cha\^\i nes $\B(L,\, \Po,\, \Po)=L\boxtimes_c
\bar{\B}(\Po)\boxtimes_c \Po$ (resp. $\B(\Po,\, \Po\, R)=\Po
\boxtimes_c \bar{\B} (\Po) \boxtimes_c R$) est un $\Po$-module
quasi-libre analytique \`a droite (resp. \`a gauche)
\end{itemize}
\end{cor}

\subsection{Cobar construction \`a coefficients}

En dualisant les arguments pr\'ec\'edents, on obtient la cobar
construction \`a coefficients.\\

Soit $(\mathcal{C},\,  \Delta,\ \varepsilon,\, \eta)$ une
coprop\'erade diff\'erentielle coaugment\'ee et soient $(L,\ l)$
et $(R,\, r)$ deux $\mathcal{C}$-comodules diff\'erentiels
respectivement \`a droite et \`a gauche.

Ces deux $\mathcal{C}$-comodules induisent les applications
$\theta'_L$ et $\theta'_R$ suivantes :
$$\theta'_R\ : \ R \xrightarrow{r} \mathcal{C}\boxtimes_c R
\twoheadrightarrow (I\oplus \underbrace{\oC}_1)\boxtimes_c R,$$
$$\theta'_L \ : \ L \xrightarrow{r}  L\boxtimes_c \mathcal{C}
\twoheadrightarrow L\boxtimes_c (I\oplus \underbrace{\oC}_1).$$

\begin{lem}
Les deux applications $\theta'_R$ et $\theta'_L$ induisent chacune
un morphisme homog\`ene de degr\'e $-1$ de la forme suivante sur
$L\boxtimes_c \bar{\B}^c(\mathcal{C})\boxtimes_c R$ :
\begin{eqnarray*}
d_{\theta'_R} \ &:& \ \xymatrix{ \ L\boxtimes_c
\bar{\B}^c(\mathcal{C})\boxtimes_c
R\ar[r]^(0.4){\widetilde{\theta'_R}} & \ {\begin{array}{l} \\ L
\boxtimes_c \bar{\B}^c(\mathcal{C}) \boxtimes_c (I\oplus
\underbrace{\oC}_1) \boxtimes_c R \end{array}} \ar[r]& }\\
&&\xymatrix@C=80pt{{\begin{array}{l} \\ L \boxtimes_c
\bar{\B}^c(\mathcal{C}) \boxtimes_c (I\oplus
\underbrace{\Sigma^{-1} \oC}_1) \boxtimes_c R \end{array}}
\ar[r]^(0.6){L \boxtimes_c \mu_{\bar{\B}^c(\mathcal{C})}
\boxtimes_c R} & L\boxtimes_c \bar{\B}^c(\mathcal{C})\boxtimes_c
R},
\end{eqnarray*}

\begin{eqnarray*}
\theta'_L \ &:& \ \xymatrix{d_{\theta'_L} \ : \ L\boxtimes_c
\bar{\B}^c(\mathcal{C})\boxtimes_c
R\ar[r]^(0.4){\widetilde{\theta'_L}} & {\begin{array}{l} \\ L
\boxtimes_c (I\oplus \underbrace{\oC}_1) \boxtimes_c
\bar{\B}^c(\mathcal{C}) \boxtimes_c R \end{array}}
\ar[r]& }\\
&&\xymatrix@C=80pt{{\begin{array}{l} \\ L \boxtimes_c (I\oplus
\underbrace{\Sigma^{-1} \oC}_1) \boxtimes_c
\bar{\B}^c(\mathcal{C}) \boxtimes_c R \end{array}}\ar[r]^(0.6){L
\boxtimes_c \mu_{\bar{\B}^c(\mathcal{C})} \boxtimes_c R} &
L\boxtimes_c \bar{\B}^c(\mathcal{C})\boxtimes_c R}.
\end{eqnarray*}
\end{lem}

\begin{pro}
Sur le $\Sy$-bimodule gradu\'e $L\boxtimes_c
\bar{\B}^c(\mathcal{C})\boxtimes_c R$, on a une diff\'erentielle
$d$ donn\'ee par la somme des termes suivants :
\vspace{6pt}

\begin{tabular}{l}
$\bullet \quad$ la diff\'erentielle canonique $\delta$ induite par celle de $\mathcal{C}$ sur $\bar{\B}^c(\mathcal{C})$,\\
$\bullet \quad$ la diff\'erentielle canonique $\delta_L$ induite par celle de $L$,\\
$\bullet \quad$ la diff\'erentielle canonique $\delta_R$ induite par celle de $R$,\\
$\bullet\quad$ la d\'erivation $d_\theta$ d\'efinie sur $\bar{\B}^c(\mathcal{C})$,\\
$\bullet\quad$ du morphisme homog\`ene $d_{\theta'_R}$ de degr\'e $-1$,\\
$\bullet\quad$ du morphisme homog\`ene $d_{\theta'_L}$ de degr\'e
$-1$.
\end{tabular}
\end{pro}

\begin{dei}[Cobar construction \`a coefficients]
\index{cobar construction \`a coefficients}

Le dg-$\Sy$-bimodule $L\boxtimes_c \bar{\B}^c(\mathcal{C})\boxtimes_c
R$ muni de la diff\'erentielle $d$ est appel\'e \emph{cobar
construction \`a coefficients dans $L$ et $R$} et est not\'e
$\B^c(L,\, \mathcal{C},\ R)$.
\end{dei}

\begin{cor}
$ \ $
\begin{itemize}
\item La cobar construction r\'eduite $\bar{\B}^c(\mathcal{C})$
correspond \`a la cobar construction \`a coefficients dans le
$\mathcal{C}$-comodule trivial $I$.

\item Le complexe de cha\^\i nes $\B^c(L,\, \mathcal{C},\,
\mathcal{C})=L\boxtimes_c \bar{\B}^c(\mathcal{C})\boxtimes_c
\mathcal{C}$ (resp. $\B(\mathcal{C},\, \mathcal{C},\,
R)=\mathcal{C} \boxtimes_c \bar{\B}^c (\mathcal{C}) \boxtimes_c R$)
est un $\mathcal{C}$-comodule quasi-colibre analytique \`a droite
(resp. \`a gauche)

\end{itemize}
\end{cor}

\subsection{Bar et cobar constructions des PROPs}

Par concat\'enation, on \'etend naturellement les d\'efinitions
pr\'ec\'edentes au cadre des PROPs.

\begin{dei}[Bar construction r\'eduite d'un PROP]
Soit $(\Po,\, \mu,\, conc)$ un PROP diff\'erentiel augment\'e. On appelle
\emph{bar construction r\'eduite du PROP $\Po$} le coPROP quasi-colibre
d\'efini par le $\Sy$-bimodule $S_\otimes (\F^c(\Sigma \oPo))$ muni
de la diff\'erentielle $\delta_\theta$, somme de la diff\'erentielle
canonique $\delta$ avec l'unique cod\'erivation $d_\theta$ qui prolonge le
morphisme $\theta$ induit par le produit partiel
$$\mu_{(1,\,1)} \, : \,
(I \oplus \underbrace{\oPo}_1) \boxtimes_c (I \oplus
\underbrace{\oPo}_1)
  \to \Po.$$
On la note aussi $\bar{\B}(\Po)$.
\end{dei}

Etant donn\'e que la cod\'erivation se fait composante connexe par
composante connexe, la bar construction PROPique est un complexe
obtenu par concat\'enation de la bar construction
pro\-p\'e\-ra\-di\-que.

\begin{pro}
\label{Lambdabar}
On a l'isomorphisme de coPROPs diff\'erentiels suivant
$$\bar{\B}(\Po)=S_\otimes(\bar{\B}(U_c(\Po))). $$
\end{pro}

Dualement, on d\'efinit la cobar construction sur un coPROP.

\begin{dei}[Cobar construction r\'eduite d'un coPROP]
Soit $(\mathcal{C},\, \Delta,\, deconc)$ un coPROP
diff\'e\-ren\-ti\-el coaugment\'e . On appelle \emph{cobar
construction r\'eduite du coPROP $\mathcal{C}$} le PROP
quasi-libre d\'efini par le $\Sy$-bimodule $S_\otimes
(\F(\Sigma^{-1} \bar{\mathcal{C}}))$ muni de la diff\'erentielle
$\delta_{\theta'}$, somme de la diff\'erentielle canonique
$\delta$ avec l'unique d\'erivation $d_{\theta'}$ qui prolonge le
morphisme $\theta'$ induit par le coproduit partiel
$$\Delta_{(1,\,1)} \, : \,
\bar{\mathcal{C}} \xrightarrow{\Delta} \mathcal{C} \boxtimes \mathcal{C}
\twoheadrightarrow (I \oplus \underbrace{\bar{\mathcal{C}}}_1)
\boxtimes_c (I \oplus \underbrace{\bar{\mathcal{C}}}_1 ).$$
On la note aussi $\bar{\B}^c(\mathcal{C})$.
\end{dei}

Comme la d\'erivation $d_{\theta'}$ pr\'eserve les composantes connexes, la
cobar
construction d'un coPROP est l'alg\`ebre sym\'etrique libre sur
 la cobar construction
de la coprop\'erade associ\'ee.

\begin{pro}
On a un isomorphisme de PROPs diff\'erentiels
$$\bar{\B}^c(\mathcal{C})=S_\otimes(\bar{\B}^c(U_c(\mathcal{C}))).  $$
\end{pro}

De la m\^eme mani\`ere que pour les prop\'erades, on d\'efinit la bar et
la cobar constructions \`a coefficients dans
un $\Po$-module (comodule) \`a droite $L$ et un $\Po$-module (comodule)
\`a gauche $R$ ainsi que les bar et cobar constructions augment\'ees.

\begin{pro}
Soit $(\Po,\, \mu, conc)$ un PROP diff\'erentiel augment\'e.
Soient $L$ un $\Po$-module diff\'e\-ren\-ti\-el \`a droite et $R$
un $\Po$-module diff\'erentiel \`a gauche. On a l'isomorphisme de
$\Sy$-bimodules diff\'erentiels
$$L\boxtimes \bar{\B}(\Po)\boxtimes R = S_\otimes \left(
L\boxtimes_c \bar{\B}(Uc(\Po)) \boxtimes_c R \right) .$$
\end{pro}

\section{Acyclicit\'e des bar et cobar constructions augment\'ees}

Lorsque $\Po$ est une alg\`ebre, c'est-\`a-dire que le
dg-$\Sy$-bimodule $\Po$ est nul en dehors de $\Po(1,\,1)=A$, on
sait que la bar construction augment\'ee (\`a gauche comme \`a
droite) sur l'alg\`ebre unitaire $A$ est acyclique. Pour
d\'emontrer cela, on introduit directement une homotopie
contractante (\emph{cf.} S\'eminaire H. Cartan \cite{Cartan}).
Dans le cas des op\'erades, B. Fresse montre que la bar
construction augment\'ee \`a gauche $\Po\circ \bar{\B}(\Po)$ est
acyclique, \`a nouveau, en exhibant une homotopie contractante.
Cette homotopie repose sur le fait que le produit mono\"\i dal
$\circ$ est lin\'eaire \`a gauche. Et l'acyclicit\'e de la bar
construction augment\'ee \`a droite d\'ecoule ensuite du
r\'esultat pr\'ec\'edent et des lemmes de comparaisons sur les
modules quasi-libres (\emph{cf.} \cite{Fresse} section $4.6$).
Dans le cadre g\'en\'eral des prop\'erades, comme le produit
mono\"\i dal $\boxtimes_c$ n'est lin\'eaire ni \`a gauche, ni \`a
droite, il faut affiner ces arguments. On commence ainsi par
d\'efinir une filtration sur le complexe $\left(\Po \boxtimes_c
\bar{\B}(\Po),\, d \right)$. Cette filtration induit une suite
spectrale $E^*_{p,\, q}$ convergente. Enfin on montre que les
complexes $\left( E^0_{p,\, *},\, d^0 \right)$ sont acycliques
pour $p>0$ en introduisant une homotopie contractante du m\^eme
type que dans le cas des alg\`ebres et des op\'erades. On
proc\`ede de la m\^eme mani\`ere pour montrer que les deux cobar
constructions coaugment\'ees \`a droite et \`a droite sur une
coprop\'erade gradu\'ee par un poids sont acycliques. Et comme le
foncteur $S_\otimes$ est un foncteur exact, on en conclut
l'acyclicit\'e des bar et cobar constructions augment\'ees dans le
cas PROPique.

\subsection{Acyclicit\'e de la bar construction augment\'ee}

Le dg-$\Sy$-bimodule $\Po\boxtimes_c\bar{\B}(\Po)$ est l'image du
foncteur
$$\oPo \mapsto (I\oplus \oPo)\boxtimes_c\F^c(\Sigma \oPo).$$
Ce foncteur est analytique scind\'e, c'est-\`a-dire qu'il est la
somme directe de foncteurs polynomiaux en $\oPo$ scind\'es
(\emph{cf.} chapitre $1$ section $8$). Pour voir cela, il suffit
de consid\'erer le nombre de sommets non r\'eduits (indic\'es par
des \'el\'ements de $\oPo$) qui composent un repr\'esentant d'un
\'el\'ement de $(I\oplus \oPo)\boxtimes_c\F^c(\Sigma \oPo)$. On
note cette d\'ecomposition :
$$\Po\boxtimes_c\bar{\B}(\Po)=\bigoplus_{s\in\mathbb{N}} \left(
\Po\boxtimes_c\bar{\B}(\Po)\right)_{(s)}.$$ On consid\`ere alors la
filtration d\'efinie par
$$ F_i=F_i\left(\Po\boxtimes\bar{\B}(\Po)\right)=\bigoplus_{s\leq i} \left(
\Po\boxtimes_c\bar{\B}(\Po) \right)_{(s)}.$$ Ainsi, $F_i$ est
compos\'ee des \'el\'ements de $\Po\boxtimes\bar{\B}(\Po)$
repr\'esentables par des graphes admettant au plus $i$ sommets non
r\'eduits.

\begin{lem}
La filtration $F_i$ du dg-$\Sy$-bimodule
$\Po\boxtimes_c\bar{\B}(\Po)$ est stable par la diff\'erentielle $d$
de la bar construction augment\'ee.
\end{lem}

\begin{deo}
Le lemme~\ref{lemme-bar-construction-augmentée} donne la forme de
la diff\'erentielle $d$. Celle-ci est la somme de trois termes :
\begin{enumerate}
\item La diff\'erentielle $\delta$ issue de celle de $\Po$. Cette
derni\`ere laisse invariant le nombre d'\'el\'ements de $\oPo$.
Ainsi, on a $\delta(F_i)\subset F_i$.

\item La cod\'erivation $d_\theta$ de la bar construction
r\'eduite $\bar{\B}(\Po)$. Cette application consiste \`a composer
des pairs de sommets. Elle v\'erifie donc $d_\theta(F_i)\subset
F_{i-1}$.

\item Le morphisme $d_{\theta_l}$. Ce morphisme a pour effet
d'\'ecr\^eter une op\'eration $\Sigma \oPo$ de $\bar{\B}(\Po)$ par
le bas, pour la composer ensuite avec des \'el\'ements de $\Po$
$$ \Po\boxtimes_c\bar{\B}(\Po) \to  \Po \boxtimes_c(I\oplus \oPo)
\boxtimes_c\bar{\B}(\Po) \to  \Po\boxtimes_c\bar{\B}(\Po).$$ Le
nombre global de sommets en $\oPo$ est donc d\'ecroissant. Ce qui
s'\'ecrit $d_{\theta_l}(F_i)\subset F_i$.$\cqfd$
\end{enumerate}
\end{deo}

De ce lemme, on obtient que la filtration $F_i$ induit une suite
spectrale not\'ee $E^*_{p,\, q}$, dont le premier terme vaut
$$E^0_{p,\, q}=F_p\left( (\Po\boxtimes_c\bar{\B}(\Po))_{p+q}\right)/
F_{p-1}\left( (\Po\boxtimes_c\bar{\B}(\Po))_{p+q}\right),$$ o\`u
$p+q$ repr\'esente le degr\'e homologique. Ainsi, le module
$E^0_{p,\, q}$ est donn\'e par les graphes \`a $p$ sommets non
r\'eduits $E^0_{p,\, q}=\left(
 \left(\Po\boxtimes_c\bar{\B}(\Po)\right)_{(p)}\right)_{p+q}$ et la
 diff\'erentielle $d^0$ est la somme de deux termes $d^0=\delta
 +d'_{\theta_l}$. Cette  derni\`ere application $d'_{\theta_l}$
 est \'equivalente \`a l'application $d_{\theta_l}$ lorsqu'elle ne
 diminue par le nombre global de sommets. De mani\`ere explicite, le
 morphisme $d'_{\theta_l}$ consiste \`a \'ecr\^eter une
 op\'eration $\Sigma \oPo$ de $\bar{\B}(\Po)$ par le
bas et \`a l'ins\'erer dans la ligne des \'el\'ements de $\Po$,
sans composition avec des op\'erations de $\oPo$ . (La composition
revient \`a faire $I\boxtimes_c \Po \to \Po$). Dans tous les autres
cas o\`u $d_{\theta_l}$ exige une composition non triviale avec
des \'el\'ements de $\Po$, $d'_{\theta_l}$ est nulle (\emph{cf.}
figure~\ref{d'}).

\begin{figure}[h]
$$ \xymatrix{ \ar@{-}[dr]&\ar@{-}[d] & &\ar@{-}[dd] &\ar@{-}[d] &\ar@{-}[dl] \\
 &*+[F-,]{\Sigma p_4} \ar@{-}[ddd] \ar@{-}[ddr] \ar@{-}[ur]& & & *+[F-,]{\Sigma p_5} \ar@{-}[dl] \ar@{-}[dd]&\\
 & & & *+[F-,]{\Sigma p_3}\ar@{-}[dl]\ar@{-}[dr] & &\\
 & & *+[F-,]{\Sigma p_1}\ar@{-}[dd] \ar@/^1pc/@{-->}[d]^{d'_{\theta_l}}& & *+[F-,]{\Sigma p_2} \ar@{-}@/_/[d]\ar@{-}@/^/[d]&\\
 &*+[F-,]{p'_1}\ar@{-}[d]\ar@{-}[dl]  \ar@{.}[l] \ar@{.}[rrr]& & & *+[F-,]{p'_2}\ar@{-}[d]\ar@{-}[dl]\ar@{-}[dr] \ar@{.}[r]&\\
& & & & & } $$ \caption{}  \label{d'}
\end{figure}
Remarquons que ce morphisme $d'_{\theta_l}$ est homog\`ene de
degr\'e $-1$ ($\Sigma \oPo \to \Po$).\\

La suite spectrale $E^*_{p,\, q}$ se situe dans le demi-plan
$p\ge0$.

On calcule l'homologie des complexes de cha\^\i nes $\left(
E^0_{p,\, *},\, d^0 \right)$ pour montrer que la suite spectrale
d\'eg\'en\`ere au rang $E^1_{p,\, q}$.

\begin{lem}
Au rang $E^1_{p,\, q}$ on a
$$ E^1_{p,\, q}=\left\{ \begin{array}{l}
I \quad \textrm{si} \quad p=q=0, \\
0 \quad \textrm{sinon}.
\end{array} \right.$$
\end{lem}

\begin{deo}
Lorsque $p=0$ on a
$$ E^0_{0,\, q}=\left\{ \begin{array}{l}
I \quad \textrm{si} \quad q=0, \\
0 \quad \textrm{sinon}.
\end{array} \right.$$
et la diff\'erentielle $d^0$ est nulle sur les modules $E^0_{0,\,
q}$. L'homologie de ces modules vaut donc
$$ E^1_{0,\, q}=\left\{ \begin{array}{l}
I \quad \textrm{si} \quad q=0, \\
0 \quad \textrm{sinon}.
\end{array} \right.$$\\

Lorsque $p>0$, pour chaque complexe de cha\^\i nes
$\left(E^0_{p,*},\ d^0 \right)$, on va exhiber une homotopie
contractante $h$.

Gr\^ace \`a la proposition~\ref{representationboxtimes}, on peut
repr\'esenter un \'el\'ement de $\Po \boxtimes \bar{\B}(\Po)$ par
$(p_1,\, \ldots ,\, p_r)\sigma (b_1,\, \ldots ,\, b_s)$. Si
$p_1\in I$, on pose
$$h\left((p_1,\, \ldots ,\, p_r)\sigma (b_1,\, \ldots ,\, b_s)
\right)=0,$$ sinon on d\'efinit $h$ par
\begin{eqnarray*}
h\left( (p_1,\, \ldots ,\, p_r)\sigma (b_1,\, \ldots ,\,
b_s)\right)= (-1)^{(|p_1|+1)(|p_2|+\cdots+|p_r|)} (1,\,\ldots,\
1,\,  p_2,\, \ldots,\, p_r)\sigma'(b',\,b_{i+1},\ \ldots ,\, b_s),
\end{eqnarray*}
o\`u $b'=\mu_{\F(\Sigma \oPo)} \left(\Sigma p_1\otimes
(b_1,\,\ldots,\,b_i)\right)$ avec $b_1,\, \ldots,\, b_i$ les
\'el\'ements de $\bar{\B}(\Po)$ reli\'es au sommet indic\'e par
$p_1$ dans la repr\'esentation graphique de $(p_1,\, \ldots ,\,
p_r)\sigma (b_1,\, \ldots ,\, b_s)$. Cette application $h$ ne
change pas le nombre d'op\'erations $\oPo$ et est de degr\'e
homologique $+1$ ($\oPo \to \Sigma \oPo$). Ainsi, on a $h\,:\,
E^0_{p,\, q} \to E^0_{p,\, q+1}$. Intuitivement, l'application $h$
revient \`a prendre une op\'eration (non triviale) parmi la ligne
de $\Po$, \`a la suspendre et \`a la remonter d'un cran pour
l'inclure dans celles de $\bar{\B}(\Po)$. Cette d\'emarche est la
d\'emarche inverse de celle de $d'_{\theta_l}$. V\'erifions
maintenant que $h$ est bien une homotopie contractante,
c'est-\`a-dire que $hd^0+d^0h=i_d$.
\begin{enumerate}
\item L'application $h$ anticommute avec $\delta$, $h\delta+\delta
h =0$. Le calcul est similaire \`a ceux effectu\'es
pr\'ec\'edemment. Encore une fois, le r\'esultat vient des
r\`egles de signes et de la suspension $\Sigma p_1$.

\item On a $h d'_{\theta_l}+d'_{\theta_l} h=i_d$. D\'ej\`a, on
calcule $d'_{\theta_l}$ :
\begin{eqnarray*}
&& d'_{\theta_l}\left( (p_1,\, \ldots ,\, p_r)\sigma (b_1,\,
\ldots
,\, b_s)\right)= \\
&&\qquad  \sum
(-1)^{|p'|(|p_{j+1}|+\cdots+|p_r|+|b_1|+\cdots+|b_k|)}\\
&& \qquad (p_1,\ldots,\, p',\, p_{j+1} \ldots,\
p_r)\widetilde{\sigma}(b_1,\ldots,\, b_k,\, b'_{k+1},\ldots ,\,
b_s),
\end{eqnarray*}
o\`u $b_{k+1}=\mu_{\F(\Sigma \oPo)}(\Sigma p'\otimes b'_{k+1})$.
Ensuite, on obtient
\begin{eqnarray*}
&&hd'_{\theta_l}\left( (p_1,\, \ldots ,\, p_r)\sigma (b_1,\,
\ldots ,\, b_s)
\right) =\\
&&\qquad  \sum (-1)^{|p'|(|p_{j+1}|+\cdots+|p_r|+|b_1|+\cdots+|b_k|)+(|p_1|+1)(|p_2|+\cdots+|p'|+\cdots+|p_r|)} \\
&&\qquad (1,\dots,\, 1,\ p_2,\ldots,\, p',\, p_{j+1} \ldots,\
p_r)\widetilde{\sigma}'(b',\, b_{i+1},\ldots,\, b_k,\,
b'_{k+1},\ldots ,\, b_s).
\end{eqnarray*}
Alors que le calcul de $d'_{\theta_l} h$ donne
\begin{eqnarray*}
&& d'_{\theta_l} h\left( (p_1,\, \ldots ,\, p_r)\sigma (b_1,\,
\ldots
,\, b_s)\right)=\\
&& \qquad (p_1,\, \ldots ,\, p_r)\sigma (b_1,\, \ldots ,\, b_s) - \\
&& \qquad \sum (-1)^{|p'|(|p_{j+1}|+\cdots+|p_r|+|b_1|+\cdots+|b_k|)+(|p_1|+1)(|p_2|+\cdots+|p'|+\cdots+|p_r|)} \\
&&\qquad (1,\dots,\, 1,\ p_2,\ldots,\, p',\, p_{j+1} \ldots,\
p_r)\widetilde{\sigma}'(b',\, b_{i+1},\ldots,\, b_k,\,
b'_{k+1},\ldots ,\, b_s),
\end{eqnarray*}
le signe $-1$ vient ici de la commutation de $p'$ avec $\Sigma
p_1$.
\end{enumerate}
L'existence de cette homotopie permet d'affirmer que les complexes
$\left( E^0_{p,\, *},\, d^0\right)$ sont acycliques pour $p>0$ et
donc que
$$  E^1_{p,\, q}=0 \quad \textrm{pour tout} \ q \ \textrm{si} \ p>0.$$
$\cqfd$
\end{deo}

On peut maintenant conclure la d\'emonstration de l'acyclicit\'e
de la bar construction augment\'ee \`a gauche.

\begin{thm}
\label{acyclicitébaraugmentée}
L'homologie du complexe $\left( \Po
\boxtimes_c \bar{\B}(\Po),\, d\right)$ est la suivante :
$$ \left\{ \begin{array}{l}
H_0\left( \Po \boxtimes_c \bar{\B}(\Po)\right)=I , \\
H_n\left( \Po \boxtimes_c \bar{\B}(\Po)\right)=0 \quad
\textrm{sinon}.
\end{array} \right.$$
\end{thm}

\begin{deo}
Comme la filtration $F_i$ est exhaustive $\Po \boxtimes_c
\bar{\B}(\Po)=\bigcup_i F_i$ et born\'ee in\-f\'e\-rieu\-re\-ment
$F_{-1}\left(\Po \boxtimes_c \bar{\B}(\Po)\right)=0$, par les
th\'eor\`emes classiques de convergence des suites spectrales
(\emph{cf.} \cite{Weibel} $5.5.1$), on sait que la suite spectrale
$E^*_{p,\, q}$ converge vers l'homologie de $\Po \boxtimes_c
\bar{\B}(\Po)$
$$E^1_{p,\, q} \Longrightarrow H_{p+q}\left(\Po \boxtimes_c \bar{\B}(\Po),\,
d\right) .$$ Et la forme d\'eg\'en\'er\'ee de $E^1_{p,\,q}$ permet
de conclure. $\cqfd$
\end{deo}

\begin{cor}
Le morphisme d'augmentation
$$\xymatrix@C=50pt{\Po\boxtimes_c \bar{\B}(\Po) \ar[r]^{\varepsilon_\Po
\boxtimes_c \varepsilon_{\F^c(\Sigma \oPo)}} & I\boxtimes_c I = I} $$
est un quasi-isomorphisme.
\end{cor}

Ce dernier r\'esultat se g\'en\'eralise de la mani\`ere suivante :

On d\'efinit le \emph{morphisme d'augmentation}
$\varepsilon(\Po,\, \Po,\, \Po)$ par la composition
$$\B(\Po,\, \Po,\,
\Po)=\Po\boxtimes_c \bar{\B}(\Po) \boxtimes_c \Po
\xrightarrow{\Po\boxtimes_c
\varepsilon_{\F^c(\Sigma \oPo)}\boxtimes_c \Po}  \Po\boxtimes_c
I\boxtimes_c \Po=\Po\boxtimes_c \Po \twoheadrightarrow \Po
\boxtimes_\Po \Po= \Po.$$

\begin{thm}
Soit $R$ un module diff\'erentiel \`a gauche sur $\Po$.
Le morphisme d'augmentation $\varepsilon(\Po,\, \Po,\, R)$
$$ \xymatrix@C=75pt{ \B(\Po,\Po,\,R)=\Po\boxtimes_\Po
\B(\Po,\, \Po,\, \Po) \boxtimes_\Po R \ar[r]^(0.6){\Po
\boxtimes_\Po \varepsilon(\Po,\, \Po,\, \Po) \boxtimes_\Po R }
 & \Po\boxtimes_\Po
 \Po \boxtimes_\Po R = R} $$ est un
quasi-isomorphime.
\end{thm}

\begin{deo}
On introduit la m\^eme filtration que pr\'ec\'edemment en comptant le
nombre de sommets indic\'es par des op\'erations de $\oPo$
des rep\'esentants des \'el\'ements de
$\Po\boxtimes_c \bar{\B}(\Po)\boxtimes_c R$. Alors, en utilisant
la m\^eme homotopie contractante on montre
que la suite spectrale d\'eg\'en\'ere au rang $E^1_{p,\, q}$ sous la forme
$$ E^1_{p,\, q}=\left\{ \begin{array}{ll}
H_q(R) & \textrm{si} \quad p=0, \\
0 & \textrm{si} \quad p>0.
\end{array} \right.$$
On conclut ensuite de la m\^eme mani\`ere, en utilisant la
convergence de la suite spectrale.$\cqfd$
\end{deo}

Pour d\'emontrer l'acyclicit\'e de la bar
construction augment\'ee \`a droite $\bar{\B}(\Po)\boxtimes_c \Po $, on
introduit la m\^eme filtration et
on d\'efinit l'homotopie contractante $h$ par
$$h((b_1,\, \ldots,\, b_r)\, \sigma \, (p_1,\,\ldots,\,p_s))=0 ,$$
lorsque $p_1\in I$, et dans le cas contraire par
\begin{eqnarray*}
&& h((b_1,\, \ldots,\, b_r)\, \sigma \, (p_1,\,\ldots,\,p_s))=\\
&&(-1)^{|b_1|+\cdots+|b_i|+|p_1|(|b_{i+1}|+\cdots+|b_r|)}(b',\,
b_{i+1}, \,\ldots ,\,b_r) \, \sigma'
(1,\,\ldots,\,1,\,p_2,\,\ldots,\,p_s),
\end{eqnarray*}
o\`u $b'=\mu_{\F(\Sigma \oPo)}((b_1,\, \ldots ,\, b_i)\otimes\Sigma p_1)$.\\

On a alors le m\^eme type de r\'esultat pour $L$ un module \`a
droite sur $\Po$.

\begin{thm}
Soit $L$ un module diff\'erentiel \`a droite sur $\Po$.
Le morphisme d'augmentation $\varepsilon(L,\, \Po,\, \Po)$
$$ \xymatrix@C=75pt{ \B(L,\Po,\,\Po)=L\boxtimes_\Po
\B(\Po,\, \Po,\, \Po) \boxtimes_\Po \Po \ar[r]^(0.6){L
\boxtimes_\Po \varepsilon(\Po,\, \Po,\, \Po) \boxtimes_\Po \Po }
 & L \boxtimes_\Po
 \Po \boxtimes_\Po \Po = L} $$ est un
quasi-isomorphime.
\end{thm}

Ces r\'esultats peuvent se reformuler dans le cadre des PROPs.

\begin{cor}
Soit $\Po$ un PROP diff\'erentiel augment\'e. Les morphismes
d'augmentation suivant sont des quasi-isomorphismes :

$$\Po\boxtimes \bar{\B}(\Po) \xrightarrow{\varepsilon_\Po
\boxtimes \varepsilon_{S_\otimes(\F^c(\Sigma \oPo))}}
I\boxtimes S_\otimes(I) = S_\otimes(I) $$
$$\bar{\B}(\Po)\boxtimes \Po \xrightarrow{\varepsilon_{
S_\otimes(\F^c(\Sigma \oPo))}
\boxtimes \varepsilon_\Po}  S_\otimes(I)\boxtimes I =
S_\otimes(I). $$
\end{cor}

\begin{deo}
Il suffit de voir que
$$\Po\boxtimes \bar{\B}(\Po)=S_\otimes(\Po \boxtimes_c \bar{\B}(U_c(\Po)))
$$
et que le foncteur $S_\otimes$ est un foncteur exact pour pouvoir
appliquer les th\'eor\`emes pr\'ec\'edents. $\cqfd$
\end{deo}

\subsection{Acyclicit\'e de la cobar construction coaugment\'ee}

De la m\^eme mani\`ere, on peut montrer l'acyclicit\'e de la cobar
construction augment\'ee \`a droite
$\bar{\B}^c(\mathcal{C})\boxtimes_c \mathcal{C} $ en utilisant des
arguments duaux. Par contre ici, pour des probl\`emes de
convergence de suites spectrales, on se place dans le cas o\`u la
coprop\'erade $\mathcal{C}$ est gradu\'ee par un poids.

\begin{thm}
\label{acyclicitécobarcoaugmentée}
Pour toute coprop\'erade coaugment\'ee gradu\'ee par un poids
$\mathcal{C}$, l'homologie du complexe
$\left(\bar{\B}^c(\mathcal{C})\boxtimes_c \mathcal{C} ,\, d\right)$
est la suivante :
$$ \left\{ \begin{array}{l}
H_0\left( \bar{\B}^c(\mathcal{C})\boxtimes_c \mathcal{C}\right)=I , \\
H_n\left( \bar{\B}^c(\mathcal{C})\boxtimes_c \mathcal{C}\right)=0
\quad \textrm{sinon}.
\end{array} \right.$$
\end{thm}

\begin{deo}
On d\'efinit une filtration $F_{i}$ par
$$F_{i}=\bigoplus_{s\ge -i} \left( \bar{\B}^c(\mathcal{C})\boxtimes_c
 \mathcal{C} \right)_{(s)} .$$
On v\'erifie que cette filtration est stable sous l'action de la
diff\'erentielle $d$ :
\begin{enumerate}
\item  La diff\'erentielle $\delta$ laisse invariant le nombre
d'\'el\'ements de $\oC$. Ainsi, on a $\delta(F_i)\subset F_i$.

\item La d\'erivation $d_{\theta'}$ de la cobar construction
r\'eduite $\bar{\B}^c(\mathcal{C})$ consiste \`a d\'ecomposer en
deux des op\'erations indi\c cant des sommets. La nombre de
sommets augmente donc de $1$. La  d\'erivation $d_{\theta'}$
v\'erifie donc $d_{\theta'}(F_i)\subset F_{i-1}$.

\item Le morphisme $d_{\theta'_r}$ a pour effet de d\'ecomposer un
\'el\'ement de $\mathcal{C}$ en deux pour inclure une des deux
parties dans $\bar{\B}^c(\mathcal{C})$. Le nombre global de
sommets en $\oPo$ est donc croissant. Ce qui s'\'ecrit
$d_{\theta'_r}(F_i)\subset F_i$.
\end{enumerate}

Cette filtration induit donc une suite spectrale $E^*_{p,\, q}$
dont le premier terme est donn\'e par
$$E^0_{p,\, q}=F_p\left((\bar{\B}^c(\mathcal{C})\boxtimes_c
\mathcal{C})_{p+q}\right)/
F_{p-1}\left((\bar{\B}^c(\mathcal{C})\boxtimes_c
\mathcal{C})_{p+q}\right)=\left((\bar{\B}^c(\mathcal{C})\boxtimes_c
\mathcal{C})_{(-p)} \right)_{p+q}.$$

La diff\'erentielle $d^0$ est alors la somme de deux termes :
$d^0=\delta+d'_{\theta'_r}$ o\`u $d'_{\theta'_r}$ correspond \`a
$d_{\theta'_r}$ lorsque celle-ci n'augmente pas le nombre de
sommets indic\'es par $\oC$. Le morphisme $d'_{\theta'_r}$ revient
donc juste \`a prendre un \'el\'ement $\oC$ de la premi\`ere ligne
(sans le d\'ecomposer), \`a le d\'esuspendre et \`a l'inclure dans
$\bar{\B}^c(\mathcal{C})$.

On montre ensuite de la m\^eme mani\`ere que
$$ E^1_{p,\, q}=\left\{ \begin{array}{l}
I \quad \textrm{si} \quad p=q=0, \\
0 \quad \textrm{sinon}.
\end{array} \right.$$
Ici l'homotopie contractante $h$ est d\'efinie comme l'op\'eration
inverse de $d'_{\theta'_r}$. Pour un \'el\'ement de
$\bar{\B}^c(\mathcal{C})\boxtimes \mathcal{C}$ repr\'esent\'e par
$$(b_1,\ldots ,\,b_r)\, \sigma \,(c_1,\ldots,\, c_s),$$
on d\'efinit $h$ en d\'ecomposant $b_1\in
\bar{\B}^c(\mathcal{C})_{(n)}$ :
$$\xymatrix{b_1 \ar@{|->}[r]&\sum b'_1\otimes \Sigma^{-1}
c },$$ o\`u $b'_1 \in \bar{\B}^c(\mathcal{C})_{(n-1)}$ et $c\in
\oC$. Si l'op\'eration $\Sigma^{-1}c$ n'est reli\'e par le haut
qu'\`a des \'el\'ements de $I$ alors on la suspend et on l'inclut
dans la ligne des \'el\'ements de $\mathcal{C}$. Un calcul du
m\^eme type que pr\'ec\'edemment montre que l'on a bien
$hd^0+d^0h=i_d$. \\

Comme la coprop\'erade $\mathcal{C}$ est gradu\'ee par un poids, on
peut d\'ecomposer le complexe $\bar{\B}^c(\mathcal{C})\boxtimes_c \mathcal{C}$
en somme directe de sous-complexes \`a l'aide de la graduation totale
induite, ce qui donne
$$\bar{\B}^c(\mathcal{C})\boxtimes_c \mathcal{C}=
\bigoplus_{\rho \in \mathbb{N}} (\bar{\B}^c(\mathcal{C})\boxtimes
\mathcal{C})^{(\rho)}. $$
Cette d\'ecomposition est compatible avec la filtration pr\'ec\'edente, aisni
qu'avec l'homotopie contractante $h$. On a donc
$$ E^1_{p,\, q}\big( (\bar{\B}^c(\mathcal{C})\boxtimes_c \mathcal{C})^{(\rho)}
\big)=0\quad \textrm{pour tout}\  p,\,q \ \textrm{d\`es que}\ \rho>0
\ \textrm{et}$$
$$E^1_{p,\, q}\big( (\bar{\B}^c(\mathcal{C})\boxtimes_c \mathcal{C})^{(0)}
\big)=\left\{ \begin{array}{l}
I \quad \textrm{si} \quad p=q=0, \\
0 \quad \textrm{sinon}.
\end{array} \right.$$

Enfin, comme la filtration $F_i$ sur le sous-complexe
$(\bar{\B}^c(\mathcal{C})\boxtimes_c \mathcal{C})^{(\rho)}$ est born\'ee
$$F_0\big ((\bar{\B}^c(\mathcal{C})\boxtimes_c \mathcal{C})^{(\rho)}\big) =
(\bar{\B}^c(\mathcal{C})\boxtimes_c \mathcal{C})^{(\rho)} \quad \textrm{et}
\quad
F_{-\rho-1}\big ((\bar{\B}^c(\mathcal{C})\boxtimes_c
\mathcal{C})^{(\rho)}\big) =0,$$
par le th\'eor\`eme classique de convergence des suites spectrales
(\emph{cf.} \cite{Weibel} $5.5.1$), on a que la suite spectrale
$E^*_{p,\, q}\big( (\bar{\B}^c(\mathcal{C})\boxtimes_c
\mathcal{C})^{(\rho)}
\big)$ converge vers l'homologie de $(\bar{\B}^c(\mathcal{C})\boxtimes_c
 \mathcal{C})^{(\rho)}$. On obtient alors le r\'esultat en faisant la
somme directe sur $\rho$. $\cqfd$
\end{deo}

\textsc{Remarque :} En introduisant une homotopie du m\^eme type, on montre
que la cobar construction coaugment\'ee \`a gauche est aussi
acyclique.

Et de la m\^eme mani\`ere, on a le th\'eor\`eme suivant

\begin{thm}
Soit $L$ un comodule diff\'erentiel \`a droite sur $\mathcal{C}$.
Si $\mathcal{C}$ est une coprop\'erade diff\'e\-ren\-ti\-elle
gradu\'ee par un poids, alors l'homologie de la cobar construction
\`a coefficients dans $L$ et $\mathcal{C}$ vaut
$$H_n\big( B^c(L,\,\mathcal{C} ,\,\mathcal{C}),\, d\big)=H_n(L,\,
\delta_L).$$
\end{thm}

On a bien sur le m\^eme r\'esultat pour les comodules \`a gauche sur
$\mathcal{C}$.\\

Comme corollaire, on a l'acyclicit\'e de la cobar construction coaugment\'ee
sur un coPROP gradu\'e par un poids.

\begin{cor}
Pour tout coPROP  $\mathcal{C}$coaugment\'e gradu\'e par un poids, l'homologie
du complexe $\left(\bar{\B}^c(\mathcal{C})\boxtimes \mathcal{C} ,\, d\right)$
est la suivante :
$$ \left\{ \begin{array}{l}
H_0\left( \bar{\B}^c(\mathcal{C})\boxtimes \mathcal{C}\right)=S_\otimes(I)
, \\
H_n\left( \bar{\B}^c(\mathcal{C})\boxtimes \mathcal{C}\right)=0
\quad \textrm{sinon}.
\end{array} \right.$$
\end{cor}

\begin{deo}
La cobar construction PROPique est l'image par le foncteur exacte
$S_\otimes$ de la cobar construction prop\'eradique sur la coprop\'erade
induite
$$\bar{\B}^c(\mathcal{C})\boxtimes \mathcal{C} = S_\otimes
\left( \bar{\B}^c(U_c(\mathcal{C}))\boxtimes_c \mathcal{C}\right).
$$
$\cqfd$
\end{deo}

\chapter{Lemmes de comparaison}

\index{lemmes de comparaison} 
\thispagestyle{empty}

Nous produisons ici des lemmes techniques de comparaison entre les
$\Po$-modules quasi-libres et entre les prop\'erades quasi-libres.
Ces lemmes sont les g\'en\'eralisations pour le produit mono\"\i
dal $\boxtimes_c$ de lemmes classiques dans le cadre du produit
tensoriel $\otimes_k$. Remarquons que B. Fresse avait d\'ej\`a
g\'en\'eralis\'e ces lemmes pour le produit mono\"\i dal $\circ$
des op\'erades dans \cite{Fresse}. Ce dernier, tout comme W. L.
Gan, dans \cite{Gan} fonde ses d\'emonstrations sur les
propri\'et\'es des arbres (respectivement des graphes de genre
$0$). Par exemple, le fait qu'un arbre \`a $n$ feuilles n'admette
qu'au plus $n-1$ noeuds non triviaux s'av\`ere crucial dans la
convergence des suites spectrales en jeu. Comme les graphes de
genre quelconque ne poss\`ede pas ces propri\'et\'es, nous avons
\'et\'e oblig\'e d'utiliser un autre outil, la graduation
naturelle par un poids de certains dg-$\Sy$-bimodules ainsi que
celles des $\Po$-modules quasi-libres analytiques et des
prop\'erades diff\'erentielles quasi-libres qu'ils engendrent. Ces
lemmes techniques nous permettrons, entre autre, de montrer que la
construction bar-cobar sur une prop\'erade (respectivement un
PROP) gradu\'ee par un poids fournit une r\'esolution de la
prop\'erade (respectivement du PROP) de d\'epart.

\section{Au niveau des $\mathcal{P}$-modules quasi-libres}

Gr\^ace \`a une filtration bien choisie, nous d\'emontrons les
lemmes de comparaison entre les $\Po$-modules quasi-libres
analytiques. Une autre filtration permet de montrer le m\^eme type
de lemme pour les $\mathcal{C}$-comodules quasi-colibres
analytiques. Enfin, ce dernier r\'esultat permet de montrer que la
construction bar-cobar sur une prop\'erade (respectivement un
PROP) gradu\'ee par un poids en fournit une r\'esolution.

\subsection{Filtration et suite spectrale associ\'ees \`a
un $\Po$-module quasi-libre analytique}

Soit $\Po$ une prop\'erade diff\'erentielle et soit
$L=M\boxtimes_c \Po$ un $\Po$-module quasi-libre analytique \`a
droite, o\`u $M=\Upsilon(V)$. Posons $M_{d}^{(\alpha)}$ le
sous-module de $M$ compos\'e des \'el\'ements de degr\'e
homologique $d$ et de graduation $\alpha$ (d\'ecomposition
polynomiale du foncteur $\Upsilon$ ou graduation totale si $V$ est
gradu\'e par un poids).

Sur $L$, on d\'efinit la filtration suivante :
$$F_s(L)=\bigoplus_\Xi \bigoplus_{|\bar{d}|+|\bar{\alpha}|\le s}
M_{\bar{d}}^{(\bar{\alpha})} (\ol,\, \ok) \otimes_{\Sy_\ok}
k[\Sc_\kj ] \otimes_{\Sy_\oj} \Po_{\bar{e}}(\oj,\, \oi),$$ o\`u la
somme directe $(\Xi)$ porte sur les entiers $m,\,n$ et $N$ et sur
les uplets $\ol,\, \ok,\, \oj,\, \oi$.

\begin{lem}
La filtration $F_s$ est stable par la diff\'erentielle $d$ de $L$.
\end{lem}

\begin{deo}
La diff\'erentielle $d$ est compos\'ee de trois termes :
\begin{itemize}
\item La diff\'erentielle $\delta_M$ induite par celle de $M$.
Elle a pour effet de diminuer de $1$ le degr\'e homologique
$|\bar{d}|\to |\bar{d}|-1$. Ainsi, on a $\delta_M(F_s)\subset
F_{s-1}$. \item La diff\'erentielle $\delta_\Po$ induite par celle
de $\Po$. Elle diminue le degr\'e homologique $|\bar{e}|$ de $1$,
soit$|\bar{e}| \to |\bar{e}|-1$. Ainsi, on a
$\delta_\Po(F_s)\subset F_s$. \item Par d\'efinition d'un
$\Po$-module quasi-libre analytique, le morphisme $d_\theta$  fait
baisser la graduation $|\bar{\alpha}|$ de $1$ et le degr\'e
homologique $|\bar{d}|$ d'au moins $1$ d'o\`u
$d_\theta(F_s)\subset F_{s-2}$.$\cqfd$
\end{itemize}
\end{deo}

Cette filtration induit donc une suite spectrale $E^*_{s,\, t}$,
dont le premier terme $E^0_{s,\, t}$ correspond au sous dg-$\Sy$-bimodule
 de $L$ suivant :
$$E^0_{s,\, t}=\bigoplus_\Xi \bigoplus_{r=0}^s \bigoplus_{|\bar{d}|=s-r
\atop |\bar{e}|=t+r}  M_{\bar{d}}^{(\bar{\alpha})} (\ol,\, \ok)
\otimes_{\Sy_\ok} k[\Sc_\kj ] \otimes_{\Sy_\oj}
\Po_{\bar{e}}(\oj,\, \oi),$$ o\`u la premi\`ere somme $(\Xi)$
porte sur $m,\,n,\, N$ et $\ol,\, \ok,\, \oj,\, \oi$. On peut
remarquer que la filtration $F_s$ ainsi que sa suite spectrale se
d\'ecomposent toujours gr\^ace \`a cette somme directe. En
reprenant les notations de la section~\ref{structurediffproduit}
du chapitre $3$, on a
$$E^0_{s,\,t}=\bigoplus_{r=0}^s I^0_{s-r,\, t+r}\left(\bigoplus_{|\bar{\alpha}|=r}
 M^{(\bar{\alpha})}\boxtimes_c \Po\right).$$

\textsc{Remarque :} Lorsque $\Po$ est gradu\'ee par un poids, on
exige que $V$ le soit aussi, ainsi le complexe $L$ se d\'ecompose
\`a l'aide de cette graduation sous la forme $L=\bigoplus_{\rho}
L^{(\rho)}$ (\emph{cf.} proposition~\ref{décompositionLrho}). De
plus, cette d\'ecomposition
 est stable par la filtration
$F_s$. Ainsi, la suite spectrale $E^*_{s,_, t}$ se d\'ecompose en
somme directe suivant la graduation par le poids de $L$ :
$$E^*_{s,\,t}(L)=\bigoplus_{\rho \in \mathbb{N}} E^*_{s,\, t}(L^{(\rho)}).$$
Par d\'efinition de $\Po$-module quasi-libre analytique, les
\'el\'ements de $V$ sont de degr\'e au moins $1$, on a alors
$$E^0_{s,\,t}(L^{(\rho)})=\bigoplus_{r=0}^{\textrm{min}(s,\, \rho)}
I^0_{s-r,\, t+r}\left(\big( \bigoplus_{|\bar{\alpha}|=r}
M^{(\bar{\alpha})} \boxtimes_c \Po \big) \bigcap L^{(\rho)}\right)
.$$

\begin{pro}
\label{SuiteSpectrale}
 Soit $L$ un $\Po$-module quasi-libre
analytique. Alors la suite spectrale $E^*_{s,\,t}$ con\-ver\-ge
vers l'homologie de $L$
$$E^*_{s,\, t}\Longrightarrow H_{s+t}(L,\,d) .$$
De plus, la diff\'erentielle $d^0$ sur $E^0_{s,\, t}$ est donn\'ee par
$\delta_\Po$ et la diff\'erentielle $d^1$ sur $E^1_{s,\, t}$ est donn\'ee par
$\delta_M$. D'o\`u, $E^2_{s,\, t}$ est reli\'e \` a $I^2$ par la formule
$$ E^2_{s,\,t}=\bigoplus_{r=0}^s I^2_{s-r,\, t+r}\left(\bigoplus_{|\bar{\alpha}|=r}
M^{(\bar{\alpha})}\boxtimes_c \Po\right)  .$$ Ce qui s'\'ecrit
$$E^2_{s,\,t}=\bigoplus_{r=0}^s \bigoplus_{|\bar{\alpha}|=r}
\bigoplus_{ |\bar{d}|=s-r \atop |\bar{e}|=t-r} (
H_*(M))^{(\bar{\alpha})}_{\bar{d}} \boxtimes_c
(H_*(\Po))_{\bar{e}} .$$
\end{pro}

\begin{deo}
Comme la filtration $F_s$ est exhaustive $\bigcup_s F_s = L$ et born\'ee
inf\'erieurement $F_{-1}=0$, le th\'eor\`eme classique de convergence des
suites spectrales assure que la suite spectrale $E^*_{s,\,t}$ converge vers
l'homologie de L (\emph{cf.} \cite{Weibel} 5.5.1).

La forme des diff\'erentielles $d^0$ et $d^1$ vient de l'\'etude
pr\'ec\'edente de l'action de la diff\'erentielle $d=\delta_M+\delta_\Po+d_\theta$
sur la filtration $F_s$.

Enfin, la formule de $E^2_{s,\,t}$ vient de la
proposition~\ref{I2II2}. $\cqfd$
\end{deo}

\subsection{Lemme de comparaison des $\Po$-modules quasi-libres analytiques}

Pour pouvoir d\'emontrer le lemme de comparaison proprement dit,
nous avons besoin du lemme technique suivant.

\begin{lem}
Soit $\Psi\, :\, \Po \to \Po'$ un morphisme de prop\'erades
diff\'erentielles augment\'ees et soient $(L,\, \lambda)$ et
$(L',\, \lambda')$ deux modules quasi-libres analytiques sur $\Po$
et $\Po'$. Posons $\bar{L}=M$ et $\bar{L'}=M'$, les quotients
ind\'ecomposables respectifs. Soit $\Phi \, :\, L \to L'$ un
morphisme de $\Po$-modules analytiques, o\`u la structure de
$\Po$-module sur $L'$ est celle donn\'ee par le foncteur de
restriction $\Psi^!$ (\emph{cf.} chapitre $1$ section $4.2$).
\begin{enumerate}

\item Un tel morphisme pr\'eserve la filtration $F_s$ et donc donne naissance
\`a un morphisme de suites spectrales
$$E^*(\Phi)\ :\ E^*(L) \to E^*(L').$$

\item Soit $\bar{\Phi} \, :\, M\to M'$ le morphisme de
dg-$\Sy$-bimodules induit par $\Phi$. Si de plus, les deux
prop\'erades $\Po$ et $\Po'$ sont gradu\'ees par un poids, on a
alors que le morphisme $E^0(\Phi)\, : \, M\boxtimes_c \Po \to
M'\boxtimes_c \Po'$ vaut $E^0=\bar{\Phi}\boxtimes_c \Psi$.
\end{enumerate}
\end{lem}

\begin{deo}
$\ $
\begin{enumerate}

\item Soit $(m_1,\ldots,\, m_b)\, \sigma \,(p_1,\ldots,\, p_a)$ un
\'el\'ement de $F_s(M\boxtimes_c \Po)$, c'est-\`a-dire que
$|\bar{d}|+|\bar{\alpha}|\le s$. Comme $\Phi$ est un morphisme de
$\Po$-modules, on a
\begin{eqnarray*}
&& \Phi\big((m_1,\ldots,\, m_b)\, \sigma \,(p_1,\ldots,\,
p_a)\big) =
\Phi\circ \lambda\big((m_1,\ldots,\, m_b)\, \sigma \,(p_1,\ldots,\, p_a)\big) =\\
&& \lambda' \big( ( \Phi(m_1),\ldots,\, \Phi(m_b))\, \sigma
\,(\Psi(p_1),\ldots,\, \Psi(p_a))\big).
\end{eqnarray*}
Si on pose $$\Phi(m_i)=(m_1^i,\ldots,\, m_{b_i}^i)\, \sigma^i\,
(p_1^i,\ldots,\, p_{a_i}^i),$$ comme $\Phi$ est un morphisme de
dg-$\Sy$-bimodules, on a $d_i=|\bar{d^i}|+|\bar{e^i}|$, d'o\`u
$|\bar{d^i}|\le d_i$ et donc $\sum_i |\bar{d^i}| \le |\bar{d}|$.
Et comme $\Phi$ est un morphisme de modules quasi-libres
anlaytiques, il pr\'eserve, par d\'efinition, la graduation en
$(\bar{\alpha})$, d'o\`u $\sum_i |\bar{\alpha^i}|\le
|\bar{\alpha}|$. Ainsi, on a
$$\Phi\big((m_1,\ldots,\, m_b)\, \sigma \,(p_1,\ldots,\, p_a)\big)\in F_s(L').$$

\item  L'application $E^0_{s,\,t}(\Phi)$ correspond au passage au
quotient suivant
$$E^0_{s,\,t} \ : \ F_s(L_{s+t})/F_{s-1}(L_{s+t}) \to F_s(L'_{s+t})/F_{s-1}(L'_{s+t}).$$
Soit $(m_1,\ldots,\, m_b)\, \sigma \,(p_1,\ldots,\, p_a)$ un
\'el\'ement de $E^0_{s,\,t}$, c'est-\`a-dire que
$|\bar{\alpha}|=r\le s$, $|\bar{d}|=s-r$ et $|\bar{e}|=t+r$. On a
donc
$$E^0_{s,\, t}(\Phi)\big( (m_1,\ldots,\, m_b)\, \sigma
\,(p_1,\ldots,\, p_a)\big)=0$$
 si et seulement si $\Phi\big(
(m_1,\ldots,\, m_b)\, \sigma \,(p_1,\ldots,\, p_a)\big)\in
F_{s-1}(L')$. Or, nous avons vu pr\'e\-c\'e\-demment que
$$ \Phi\big((m_1,\ldots,\, m_b)\, \sigma \,(p_1,\ldots,\, p_a)\big) =
\lambda' \big( ( \Phi(m_1),\ldots,\, \Phi(m_b))\, \sigma
\,(\Psi(p_1),\ldots,\, \Psi(p_a))\big).$$ Ainsi,
$\Phi\big((m_1,\ldots,\, m_b)\, \sigma \,(p_1,\ldots,\, p_a)\big)$
appartient \`a $F_s(L'_{s+t})\backslash F_{s-1}(L'_{s+t})$ si et
seulement si $\Phi(m_i)=m_i=\bar{\Phi}(m_i)$. En effet, si on a
$\Phi(m_i)=(m_1^i,\ldots,\, m_{b_i}^i)\, \sigma^i\,
(p_1^i,\ldots,\, p_{a_i}^i)$ avec au moins un $p_j^i$
n'appartenant pas \`a $I$. Alors, par la d\'efinition de
prop\'erade gradu\'ee par un poids, le degr\'e pour la graduation
par le poids de $p_j^i$ est au moins de $1$ et par conservation
globale de cette graduation par $\Phi$, on a
$|\bar{\alpha^i}|<\alpha^i$. Ainsi, le morphisme $E^0(\Phi)$
correspond bien \`a $\bar{\Phi}\boxtimes_c \Po$.$\cqfd$
\end{enumerate}
\end{deo}

\begin{thm}[Lemme de comparaison des modules quasi-libres analytiques]
Soit $\Psi\, :\, \Po \to \Po'$ un morphisme de prop\'erades
diff\'erentielles gradu\'ees par un poids augment\'ees et soient
$(L,\, \lambda)$ et $(L',\, \lambda')$ deux modules quasi-libres
analytiques sur $\Po$ et $\Po'$. Posons $\bar{L}=M$ et
$\bar{L'}=M'$ les quotients ind\'ecomposables respectifs,
c'est-\`a-dire $L=M\boxtimes_c \Po$ et $L'=M'\boxtimes_c \Po'$.
Soit $\Phi \, :\, L \to L'$ un morphisme de $\Po$-modules
analytiques, o\`u la structure de $\Po$-module sur $L'$ est celle
donn\'ee par le foncteur de restriction $\Psi^!$. Posons
$\bar{\Phi} \, :\, M\to
M'$ le morphisme de dg-$\Sy$-bimodules induit par $\Phi$.\\

Si deux des trois morphismes suivants $\left\{ \begin{array}{l}
\Psi \ : \ \Po \to \Po', \\
\bar{\Phi} \ : \ M \to M', \\
\Phi \ : \ L \to L',
\end{array}\right.$
sont des quasi-isomorphismes, alors le troisi\`eme est aussi un
quasi-isomorphisme.
\end{thm}

\textsc{Remarque :} Si les prop\'erades ne sont pas gradu\'ees par
un poids et si $\Phi=\bar{\Phi}\boxtimes_c \Psi$, lorsque
$\bar{\Phi}$ et $\Psi$ sont des quasi-isomorphismes, $\Phi$ est
aussi un quasi-isomorphisme. La d\'emonstration est rapide. On a
imm\'editatement $E^0(\Phi)=\bar{\Phi}\boxtimes_c \Psi$. Comme
$\bar{\Phi}$ et $\Psi$ sont des quasi-isomorphismes, on obtient
que $E^2(\Phi)$ est un isomorphisme par la
proposition~\ref{SuiteSpectrale} ($d^0=\delta_\Po$ et
$d^1=\delta_M$). Et toujours gr\^ace \`a cette proposition, la
convergence naturelle de la suite spectrale $E^*$ donne que $\Phi$
est un quasi-isomorphisme entre $L$ et $L'$.

\begin{deo}
$ $
\begin{enumerate}

\item  On suppose ici que $\Psi  \, : \, \Po \to \Po'$ et
$\bar{\Phi} \, : \, M \to M'$ sont des quasi-isomorphismes. En
appliquant le lemme pr\'ec\'edent, on sait que $\Phi$ induit un
morphisme de suites spectrales $E^*(\Phi) \, :\, E^*(L)\to
E^*(L')$ et surtout que $E^0(\Phi)=\bar{\Phi}\boxtimes_c \Psi$,
 on en conclut alors le r\'esultat de la m\^eme mani\`ere.\\

\item On suppose maintenant que $\Psi  \, : \, \Po \to \Po'$ et
$\Phi \, : \, L \to L'$ sont des
quasi-isomorphismes. Nous avons vu pr\'ec\'edemment que la suite spectrale
$E^*_{s,\, t}$ pr\'eservait la d\'ecomposition $L^{(\rho)}$. De
plus, on a
$$E^2_{s,\,t}(L^{(\rho)})=\bigoplus_\Xi \bigoplus_{r=\textrm{max}
(0,\, -t)}^{\textrm{min}(s,\, \rho)}I^2_{s-r,\, t+r}\left(
\bigoplus_{\chi}  M_{\bar{d}}^{(\bar{\alpha})} (\ol,\, \ok)
\otimes_{\Sy_\ok} k[\Sc_\kj ] \otimes_{\Sy_\oj}
\Po_{\bar{e}}^{(\bar{\beta})}(\oj,\, \oi)  \right), $$ o\`u la
deuxi\`eme somme directe $(\chi)$ porte sur $|\bar{\alpha}|=r$,
$|\bar{\beta}|=\rho-r$, $|\bar{d}|=s-r$ et $|\bar{e}|=t+r$. Or,
lorsque $s\ge\rho$, on a pour $r=\rho$
\begin{eqnarray*}
 I^2_{s-\rho,\, t+\rho}\left( \bigoplus_{\chi}
M_{\bar{d}}^{(\bar{\alpha})} (\ol,\, \ok) \otimes_{\Sy_\ok}
k[\Sc_\kj ] \otimes_{\Sy_\oj} \Po_{\bar{e}}^{(\bar{\beta})}(\oj,\,
\oi)\right) &=& I^2_{s-\rho,\, t+\rho} \left(
\bigoplus_{|\bar{\alpha}|=\rho
\atop |\bar{d}|=s-\rho} M_{\bar{d}}^{(\bar{\alpha})} \right)\\
&=& \left\{ \begin{array}{ll} H_{s-\rho}(M^{(\rho)}) &
\textrm{si}\
 t=-\rho, \\
0 & \textrm{sinon}.
\end{array} \right.
\end{eqnarray*}

En r\'esum\'e, on obtient
\begin{eqnarray*}
\left\{ \begin{array}{ll}
 E^2_{s,\,t}(L^{(\rho)})=0 &
\textrm{si}\
 t<-\rho, \\
E^2_{s,\, -\rho}(L^{(\rho)})=H_{s-\rho}(M^{(\rho)}) & \textrm{si}
\ s\ge \rho.
\end{array} \right.
\end{eqnarray*}

On montre ensuite par r\'ecurrence sur l'entier $\rho$ que
$\bar{\Phi}^{(\rho)} \, :\, M^{(\rho)} \to M'^{(\rho)}$ induit un
isomorphisme en homologie $H_*(\bar{\Phi}^{(\rho)}) \, :\,
H_*(M^{(\rho)}) \to H_*(M'^{(\rho)})$.

Pour $\rho=0$, comme $L^{(0)}=M^{(0)}$, on a
$\bar{\Phi}^{(0)}=\Phi^{(0)}$, qui est un quasi-isomorphisme.

Supposons maintenant que le r\'esultat soit vrai pour $r<\rho$ et
tous les indices $*$ ainsi que pour $r=\rho$ et pour les indices
$*<d$, et montrons le pour $*=d$.(Comme tous les complexes de
cha\^\i nes sont ici nuls en degr\'e strictement n\'egatif, on a
toujours que $H_s(\bar{\Phi}^{(\rho)})$ est un isomorphisme pour
$s<0$).

Par le lemme pr\'ec\'edent, on a
$$E^2_{s,\,t}(\Phi^{(\rho)})=\bigoplus_{r=0}^{\textrm{min}(s,\, \rho)}
I^2_{s-r,\, t+r}\left( \bigoplus_{|\bar{\alpha}|=r}
\bar{\Phi}^{{(\bar{\alpha})}}\boxtimes_c \Psi \right).$$ Ce qui
donne, avec l'hypoth\`ese de r\'ecurrence, que
$E^2_{s,\,t}(\Phi^{(\rho)})\, :\, E^2_{s,\,t}(L^{(\rho)})\to
E^2_{s,\,t}(L'^{(\rho)})$ est un isomorphisme pour $s<d+\rho$.
Montrons que $$E^2_{d+\rho,\,-\rho}(\Phi^{(\rho)})\, :\,
E^2_{d+\rho,\,-\rho}(L^{(\rho)})\to
E^2_{d+\rho,\,-\rho}(L'^{(\rho)})
$$ est encore un isomorphisme. Pour cela, on introduit le c\^one
associ\'e au morphisme $\Phi^{(\rho)}$ :
$\cone(\Phi^{(\rho)})=\Sigma^{-1}L^{(\rho)}\oplus L'^{(\rho)}$.
Sur ce c\^one, on d\'efinit la filtration
$$F_s(\cone(\Phi^{(\rho)}))=F_{s-1}(\Sigma^{-1}L^{(\rho)})\oplus F_s(L'^{(\rho)}) .$$
Cette filtration induit une suite spectrale qui v\'erifie
$$E^1_{*,\, t}(\cone(\Phi^{(\rho)}))=\cone(E^1_{*,\,t}(\Phi^{(\rho)})).$$
Cependant, le c\^one de $E^1_{*,\,t}(\Phi^{(\rho)})$ induit une
suite exacte longue
\begin{eqnarray*} && \xymatrix@C=55pt{\cdots
\ar[r]& H_{s+1}\big(\cone(E^1_{*,\, t}(\Phi^{(\rho)}))\big)
\ar[r]& H_s\big(E^1_{*,\,t}(L^{(\rho)})\big)
\ar[r]^(0.7){H_s\big(E^1_{*,\,t}(\Phi^{(\rho)}) \big)}& }\\
&& \xymatrix@C=35pt{H_s\big(E^1_{*,\,t}(L'^{(\rho)})\big) \ar[r]&
H_s\big(\cone(E^1_{*,\, t}(\Phi^{(\rho)}))\big)\ar[r] &
H_{s-1}\big(E^1_{*,\,t}(L^{(\rho)})\big) \ar[r]& \cdots,}
\end{eqnarray*}
ce qui donne finalement le suite exacte longue $(\xi)$ suivante :
\begin{eqnarray*} && \xymatrix@C=35pt{\cdots
\ar[r]& E^2_{s+1,\, t}(\cone(\Phi^{(\rho)})) \ar[r]&
E^2_{s,\,t}(L^{(\rho)}) \ar[r]^{E^2_{s,\,t}(\Phi^{(\rho)})} & E^2_{s,\,t}(L'^{(\rho)}) }\\
&& \xymatrix@C=35pt{\quad  \ar[r]&E^2_{s,\,
t}(\cone(\Phi^{(\rho)})) \ar[r] & E^2_{s-1,\,t}(L^{(\rho)})
\ar[r]& \cdots}
\end{eqnarray*}
Nous avons vu pr\'ec\'edemment que pour tout $t<-\rho$,
$$E^2_{s,\,t}(L^{(\rho)})=E^2_{s,\,t}(L'^{(\rho)})=0.$$
 La suite
exacte longue $(\xi)$ montre alors que
$E^2_{s,\,t}(\cone(\Phi^{(\rho)}))=0$ pour tout $s$, lorsque
$t<-\rho$.

De la m\^eme mani\`ere, nous avons vu que
$E^2_{s,\,t}(\Phi^{(\rho)})$ \'etait un isomorphisme pour
$s<d+\rho$. Gr\^ace \`a la suite exacte longue $(\xi)$, on obtient
que $E^2_{s,\,t}(\cone(\Phi^{(\rho)}))=0$ pour tout $t$, lorsque
$s<d+\rho$.

Ces deux r\'esultats permettent de conclure que
$$
\left\{
\begin{array}{l}
E^2_{d+\rho,\, -\rho}(\cone(\Phi^{(\rho)}))=E^\infty_{d+\rho,\,
-\rho}(\cone(\Phi^{(\rho)})), \\
E^2_{d+\rho+1,\,
-\rho}(\cone(\Phi^{(\rho)}))=E^\infty_{d+\rho+1,\,
-\rho}(\cone(\Phi^{(\rho)})).
\end{array} \right.$$

Comme la filtration $F_s$ du c\^one de $\Phi^{(\rho)}$ est
born\'ee inf\'erieurement et exhaustive, on sait que la suite
spectrale $E^*_{s,\,t}(\cone(\Phi^{(\rho)}))$ converge vers
l'homologie de ce c\^one. Comme $\Phi^{(\rho)}$ est un
quasi-isomorphisme, cette homologie est nulle, d'o\`u les
\'egalit\'es suivantes :
$$\left\{ \begin{array}{l}
E^2_{d+\rho,\,
-\rho}(\cone(\Phi^{(\rho)}))=E^\infty_{d+\rho,\,
-\rho}(\cone(\Phi^{(\rho)}))=0, \\
E^2_{d+\rho+1,\,
-\rho}(\cone(\Phi^{(\rho)}))=E^\infty_{d+\rho+1,\,
-\rho}(\cone(\Phi^{(\rho)}))=0.
\end{array} \right. $$
En r\'einjectant ces \'egalit\'es dans la suite exacte longue
$(\xi)$, on a que $E^2_{d+\rho,\, -\rho}(\Phi^{(\rho)})$ est un
isomorphisme.

On conclut en utilisant le fait que $E^2_{d+\rho,\,
-\rho}(L^{(\rho)})=H_d(M^{(\rho)})$ et que $E^2_{d+\rho,\,
-\rho}(L'^{(\rho)})=H_d(M'^{(\rho)})$. Ainsi, $E^2_{d+\rho,\,
-\rho}(\Phi^{(\rho)})$ correspond \`a $H_d\big(
\bar{\Phi}^{(\rho)}\big)\, : \, H_d(M^{(\rho)}) \to
H_d(M'^{(\rho)})$ qui est donc un isomorphisme. Ce qui conclut
la r\'ecurrence et montre que
$\bar{\Phi}$ est un quasi-isomorphisme.\\

\item Supposons que $\Phi \, : \, L \to L'$ et $\bar{\Phi} \, :\,
M \to M'$ sont des quasi-isomorphismes. Nous utilisons essentiellement
les m\^emes id\'ees que pr\'ec\'edemment.
Comme $M^{(0)}=I$, on tire de la relation
$$E^2_{s,\,t}(L^{(\rho)})=\bigoplus_\Xi
\bigoplus_{r=\textrm{max}(0,\, -t)}^{\textrm{min}(s,\, \rho)}
I^2_{s-r,\, t+r}\left( \bigoplus_\chi
M_{\bar{d}}^{(\bar{\alpha})} (\ol,\, \ok) \otimes_{\Sy_\ok}
k[\Sc_\kj ] \otimes_{\Sy_\oj} \Po_{\bar{e}}^{(\bar{\beta})}(\oj,\,
\oi)  \right), $$ o\`u la deuxi\`eme somme directe $\chi$ porte
sur $|\bar{\alpha}|=r$, $|\bar{\beta}|=\rho-r$, $|\bar{d}|=s-r$ et
$|\bar{e}|=t+r$, que pour $s=0$
$$E^2_{0,\,t}(L^{(\rho)})= I^2_{0,\, t}(\Po_t^{(\rho)})=H_t(\Po^{(\rho)}). $$
On a imm\'ediatement que
$$E^2_{s,\,t}(L^{(\rho)})=0,$$
lorsque $s<0$.

On montre ensuite, par r\'ecurrence sur l'entier $\rho$, que les
morphismes $\Psi^{(\rho)} \, : \, \Po^{(\rho)}\to \Po'^{(\rho)}$
sont des quasi-isomorphismes.

On fonde la r\'ecurrence en remarquant que la d\'efinition de
$\Po$, $\Po'$ gradu\'ees par un poids impose
$\Po^{(0)}=\Po'^{(0)}=I$, d'o\`u $\Psi^{(0)}=id_I$ e st un
quasi-isomorphisme.

Supposons maintenant que le r\'esultat soit vrai pour tout $r<\rho$.
Et on montre par r\'ecurrence sur l'entier $t$ que
$$H_t(\Psi^{(\rho)})\ :\  H_t(\Po^{(\rho)}) \to H_t(\Po'^{(\rho)}) $$
est un isomorphisme. On sait que c'est trivialement vrai pour
$t<0$. On suppose que c'est encore vrai pour $t<e$. Par le lemme
pr\'ec\'edent, on a
$$E^2_{s,\,t}(\Phi^{(\rho)})=\bigoplus_{r=0}^{\textrm{min}(s,\, \rho)}
I^2_{s-r,\, t+r}\left(
\bigoplus_{|\bar{\alpha}|=r \atop |\bar{\beta}|=\rho-r}
\bar{\Phi}^{(\bar{\alpha})}\boxtimes_c \Psi^{(\bar{\beta})}
 \right).$$
Avec les hypoth\`eses de r\'ecurrence, on obtient que
$E^2_{s,\,t}(\Phi^{(\rho)}) \, : \, E^2_{s,\,t}(L^{(\rho)}) \to
E^2_{s,\,t}(L'^{(\rho)})$ est un isomorphisme pour tout $s$,
lorsque $t<e$. En injectant ceci dans la suite exacte longue
$(\xi)$, on montre que $E^2_{s,\,t} (\cone(\Phi^{(\rho)}))=0$ pour
tout $s$, lorsque $t<e$.

La forme de la filtration $F_s(\cone(\Phi^{(\rho)}))=F_{s-1}(\
\Sigma^{-1}L^{(\rho)})\oplus F_s(L'^{(\rho)})$ joint au fait que $F_{-1}(L^{(\rho)})=0$ montrent que $E^2_{s,\,t}(\cone(\Phi^{(\rho)}))=0$ pour $s<0$.

Ces deux r\'esultats, plus la convergence de la suite spectrale
$E^*_{s,\,t}(\cone(\Phi^{(\rho)}))$ vers l'homologie nulle du c\^one de
$\Phi^{(\rho)}$, permettent de conclure que
$$
\left\{
\begin{array}{l}
E^2_{0,\, e}(\cone(\Phi^{(\rho)}))=E^\infty_{0,\,
e}(\cone(\Phi^{(\rho)}))=0 \\
E^2_{1,\,
e}(\cone(\Phi^{(\rho)}))=E^\infty_{1,\,
e}(\cone(\Phi^{(\rho)}))=0.
\end{array} \right.
$$

En r\'einjectant ces deux \'egalit\'es dans la suite exacte longue
$(\xi)$, on obtient
que le morphisme
$$\xymatrix@C=45pt{E^2_{0,\,e}(L^{(\rho)}) \ar[r]^{E^2_{0,\, e}(\Phi^{(\rho)})} &   E^2_{0,\,e}(L'^{(\rho)}) }$$
est un isomorphisme.
On conclut en rappelant que $E^2_{0,\,e}(L^{(\rho)})=H_e(\Po^{(\rho)})$,
 $E^2_{0,\,e}(L'^{(\rho)})=H_e(\Po'^{(\rho)})$ et que
$E^2_{0,\, e}(\Phi^{(\rho)})=H_e(\Psi^{(\rho)}).$ $\cqfd$
\end{enumerate}
\end{deo}

En utilisant les m\^emes arguments, on d\'emontre le m\^eme type
de lemme dans le cadre des PROPs.

\begin{thm}[Lemme de comparaison des modules quasi-libres analytiques sur des PROPs]
Soit $\Psi\, :\, \Po \to \Po'$ un morphisme de PROPs
diff\'erentiels augment\'es gradu\'es par un poids et soient
$(L,\, \lambda)$ et $(L',\, \lambda')$ deux modules quasi-libres
analytiques sur $\Po$ et $\Po'$. Posons $\bar{L}=M$ et
$\bar{L'}=M'$ les quotients ind\'ecomposables respectifs. Soit
$\Phi \, :\, L \to L'$ un morphisme de $\Po$-modules analytiques,
o\`u la structure de $\Po$-module sur $L'$ est celle donn\'ee par
le foncteur de restriction $\Psi^!$. Posons $\bar{\Phi} \, :\,
M\to
M'$ le morphisme de dg-$\Sy$-bimodules induit par $\Phi$.\\

Si deux des trois morphismes suivants $\left\{ \begin{array}{l}
\Psi \ : \ \Po \to \Po', \\
\bar{\Phi} \ : \ M \to M', \\
\Phi \ : \ L \to L',
\end{array} \right. $
sont des quasi-isomorphismes, alors le troisi\`eme est aussi un
quasi-isomorphisme.
\end{thm}

\begin{deo}
On applique les m\^emes arguments que pr\'ec\'edemment \`a la
filtration
$$F_s(L)=\bigoplus_\Xi \bigoplus_{|\bar{d}|+|\bar{\alpha}|\le s}
M_{\bar{d}}^{(\bar{\alpha})} (\ol,\, \ok) \otimes_{\Sy_\ok}
k[\Sy_N] \otimes_{\Sy_\oj} \Po_{\bar{e}}(\oj,\, \oi),$$
 o\`u la
somme directe $(\Xi)$ porte sur les entiers $m,\,n$ et $N$ et sur
les uplets $\ol,\, \ok,\, \oj,\, \oi$. La seule diff\'erence
intervient lorsque l'on consid\`ere des expressions de la forme
$M\boxtimes \Po^{(0)}=M \boxtimes I = S_\otimes(M)$ au lieu de
$M\boxtimes_c I = M$. C'est par exemple le cas dans la partie $2$
de la d\'emonstration pr\'ec\'edente. On conclut de la m\^eme
mani\`ere, en faisant une r\'ecurrence sur le poids de $M$, car
les \'el\'ements de $S_\otimes(M)^{(\rho)}$ sont des sommes de
produits tensoriels d'\'el\'ements de poids inf\'erieur ou \'egal
\`a $\rho$. $\cqfd$
\end{deo}

\subsection{Lemme de comparaison des $\mathcal{C}$-comodules quasi-colibres analytiques}

On peut d\'emontrer le m\^eme type de r\'esultat pour des
comodules quasi-colibres analytiques. Soit $\mathcal{C}$ une
coprop\'erade diff\'erentielle gradu\'ee par un poids et soit
$L=M\boxtimes_c \mathcal{C}$ un $\mathcal{C}$-comodule
quasi-colibre analytique \`a droite, o\`u $M=\Upsilon(V)$. Sur ce
complexe on d\'efinit la graduation
$$F'_s(L)=\bigoplus_{(\Xi)} \bigoplus_{|\bar{e}|+|\bar{\beta}|\le s}
M_{\bar{d}}(\ol,\, \ok)\otimes_{\Sy_\ok} k[\Sc_\kj]
\otimes_{\Sy_\oj}\mathcal{C}_{\bar{e}}^{(\bar{\beta})}(\oj,\,\oi),
$$ o\`u la somme directe $(\Xi)$ porte sur les entiers $m$, $n$ et
$N$ et sur les uplets $\ol$, $\ok$, $\oj$, $\oi$, $\bar{d}$,
$\bar{e}$ et $\bar{\beta}$ tels que $|\bar{e}|+|\bar{\beta}|\le
s$.

En appliquant les m\^emes arguments que pr\'ec\'edemment \`a la
filtration $F'_s$, on obtient le th\'eor\`eme suivant :

\begin{thm}[Lemme de comparaison des comodules quasi-colibres analytiques]
\label{lemmecomparaisoncomodule} Soit $\Psi\, :\, \mathcal{C} \to
\mathcal{C}'$ un morphisme de coprop\'erades diff\'erentielles
coaugment\'ees gradu\'ees par un poids et soient $(L,\, \lambda)$
et $(L',\, \lambda')$ deux comodules quasi-colibres analytiques
sur $\mathcal{C}$ et $\mathcal{C}'$. Posons $\bar{L}=M$ et
$\bar{L'}=M'$ les quotients ind\'ecomposables respectifs. Soit
$\Phi \, :\, L \to L'$ un morphisme de $\mathcal{C}$-comodules
analytiques, o\`u la structure de $\mathcal{C}$-comodule sur $L'$
est celle donn\'ee par le foncteur de restriction $\Psi^!$. Posons
$\bar{\Phi} \, :\, M\to
M'$ le morphisme de dg-$\Sy$-bimodules induit par $\Phi$. \\

Si deux des trois morphismes suivants
$\left\{ \begin{array}{l}
\Psi \ : \ \mathcal{C} \to \mathcal{C}', \\
\bar{\Phi} \ : \ M \to M,' \\
\Phi \ : \ L \to L',
\end{array} \right. $
sont des quasi-isomorphismes, alors le troisi\`eme est aussi un
quasi-isomorphisme.
\end{thm}

De la m\^eme mani\`ere, on a une version PROPique de ce th\'eor\`eme.

\begin{thm}[Lemme de comparaison des comodules quasi-colibres analytiques sur
des coPROPs]
\label{lemmecomparaisoncomodulePROP}

Soit $\Psi\, :\, \mathcal{C} \to \mathcal{C}'$ un morphisme de
coPROPs diff\'erentiels coaugment\'es gradu\'es par un poids et
soient $(L,\, \lambda)$ et $(L',\, \lambda')$ deux comodules
quasi-colibres analytiques sur $\mathcal{C}$ et $\mathcal{C}'$.
Posons $\bar{L}=M$ et $\bar{L'}=M'$ les quotients
ind\'ecomposables respectifs. Soit $\Phi \, :\, L \to L'$ un
morphisme de $\mathcal{C}$-comodules analytiques, o\`u la
structure de $\mathcal{C}$-comodule sur $L'$ est celle donn\'ee
par le foncteur de restriction $\Psi^!$. Posons $\bar{\Phi} \, :\,
M\to
M'$ le morphisme de dg-$\Sy$-bimodules induit par $\Phi$. \\

Si deux des trois morphismes suivants
$\left\{ \begin{array}{l}
\Psi \ : \ \mathcal{C} \to \mathcal{C}', \\
\bar{\Phi} \ : \ M \to M', \\
\Phi \ : \ L \to L',
\end{array} \right. $
sont des quasi-isomorphismes, alors le troisi\`eme est aussi un
quasi-isomorphisme.
\end{thm}

Nous allons utiliser ces deux derniers th\'eor\`emes dans le chapitre suivant
pour \'etablir la r\'esolution bar-cobar.

\subsection{Application : la r\'esolution bar-cobar}

\index{r\'esolution bar-cobar}
Nous g\'en\'eralisons ici aux
pro\-p\'e\-ra\-des et aux PROPs gradu\'es par un poids la
proposition de V. Ginzburg et M.M. Kapranov \cite{GK} qui affirme
que la construction bar-cobar sur une op\'erade $\Po$ fournit une
r\'esolution de $\Po$.\\

Commen\c cons par le cas des prop\'erades. On d\'efinit le
morphisme $\zeta \, :\, \bar{\B}^c(\bar{\B}(\Po))\to \Po$ par la
composition
\begin{eqnarray*}
\xymatrix{\bar{\B}^c(\bar{\B}(\Po))=\F\big(\Sigma^{-1}
\bar{\F}^c(\Sigma \oPo)\big) \ar@{>>}[r] &
\F\big(\Sigma^{-1} \bar{\F}_{(1)}^c(\Sigma
\oPo)\big)=\F(\oPo) \ar[r]^(0.75){c_\Po}
& \Po,}
\end{eqnarray*}
o\`u le morphisme $c_\Po$ correspond \`a la counit\'e dans la
d\'emonstration du mono\"\i de libre (\emph{cf.}
th\'e\-o\-r\`eme~\ref{monoidelibre}). Ce morphisme revient \`a
composer entre elles, via $\mu$, les op\'erations $\Po$ de
$\F(\oPo)$.

\begin{pro}
L'application $\zeta$ ainsi d\'efinie est un morphisme de prop\'erades
diff\'erentielles gradu\'ees par un poids.
\end{pro}

\begin{deo}
La d\'efinition de $\zeta$ repose sur la composition $\mu$ des
op\'erations de $\Po$, ainsi on montre facilement que $\zeta\circ
\mu_{\F(\Sigma^{-1} \bar{\F}^c(\Sigma \Po))} =\mu_\Po\circ
(\zeta\boxtimes_c\zeta)$, c'est-\`a-dire que $\zeta$ est un
morphisme de prop\'erades.

Comme le morphisme $\zeta$ est
soit nul soit correspond \`a des compositions d'op\'erations, actions  qui
pr\'eservent la
graduation par le poids de $\Po$, alors on a que $\zeta$ est un morphisme de
prop\'erades gradu\'ees par un poids.

Reste \`a montrer que $\zeta$ commute avec les diff\'erentielles respectives.

Sur $\bar{\B}^c(\bar{\B}(\Po))=\F\big(\Sigma^{-1}
\bar{\F}^c(\Sigma \oPo)\big)$  la diff\'erentielle $d$ correspond
\` a la somme de la d\'erivation  $d_{\theta'}$ de la cobar
construction induite par le coproduit partiel de $\F^c(\Sigma
\oPo)$ avec la diff\'erentielle canonique
$\delta_{\bar{\B}(\Po)}$. Or, la diff\'erentielle canonique
$\delta_{\bar{\B}(\Po)}$ est la somme de la cod\'erivation
$d_\theta$ de la bar construction $\bar{\B}(\Po)$ induite par le
produit partiel de $\Po$ avec la diff\'erentielle canonique
$\delta_\Po$.

\begin{itemize}

\item Si tous les sommets de $\xi$ sont indic\'es par des
\'el\'ements de $\F^c_{(1)}(\Sigma \oPo)=\Sigma \oPo$, effectuer
la diff\'erentielle $d$ sur $\xi$ consiste \`a n'effectuer que la
diff\'erentielle canonique $\delta_\Po$. Et comme la composition
$\mu$ de la prop\'erade diff\'erentielle $\Po$ commute avec
$\delta_\Po$, on a bien $\zeta\circ d(\xi)= \zeta \circ
\delta_\Po(\xi)=\delta_\Po\circ \zeta(\xi)$.

\item Soit $\xi$ un \'el\'ement de $\F\big(\Sigma^{-1}
\bar{\F}^c(\Sigma \oPo)\big)$ qui se repr\'esente \`a l'aide d'un
graphe dont au moins un sommet est indic\'e par un \'el\'ement de
$\F^c_{(s)}(\Sigma \oPo)$, avec $s\ge 3$. Son image par $\zeta$
est nulle. En outre, l'image de $\xi$ par la diff\'erentielle $d$
est une somme des graphes dont au moins un des sommets est
indic\'e par un \'el\'ement de $\F^c_{(s)} (\Sigma \oPo)$, avec
$s\ge 2$. Ainsi, on a $d\circ \zeta(\xi) + \zeta\circ d(\xi)=0$.

\item Si $\xi$ est repr\'esent\'e par un graphe dont les sommets
sont indic\'es par des \'el\'ements de $\F^c_{(s)}(\Sigma \oPo)$
o\`u $s=1,\, 2$, avec au moins un \'el\'ement de
$\F^c_{(2)}(\Sigma \oPo)$, son image par $\zeta$ est nulle. Reste
\`a montrer qu'il en est de m\^eme pour $\zeta\circ d$. Comme
$\delta_\Po$ pr\'eserve la graduation naturelle de $\F^c$, on a
imm\'ediatement $\zeta\circ \delta_\Po (\xi)=0$. Regardons l'effet
de $d_\theta+d_{\theta'}$ sur un \'el\'ement de $\F^c_{(2)}(\Sigma
\oPo)$. Posons $\xi=X\otimes \Sigma^{-1}(\Sigma p_1\otimes \Sigma
p_2)\otimes Y$ o\`u $(\Sigma p_1\otimes \Sigma p_2)$ appartient
\`a $\F_{(2)}^c(\Sigma \oPo)$ et $X$, $Y$ \`a $\F(\Sigma^-1
\bar{\F}^c(\Sigma \oPo))$. En appliquant $d_{\theta'}$, on obtient
un terme de la forme
$$(-1)^{|X|+1}(-1)^{|p_1|+1}X\otimes \Sigma^{-1}\Sigma p_1\otimes
 \Sigma^{-1}\Sigma p_2\otimes Y, $$
o\`u $\Sigma^{-1}\Sigma p_1 $ et $\Sigma^{-1}\Sigma p_2 $ sont des
\'el\'ements
de $\Sigma^{-1} \F^c_{(1)}(\Sigma \oPo)$.
Et, en appliquant $d_\theta$, on obtient un terme de la forme
$$(-1)^{|X|+1}(-1)^{|p_1|}X\otimes \Sigma^{-1} \Sigma \mu(p_1\otimes p_2)
\otimes
Y, $$
o\`u $\Sigma^{-1} \Sigma \mu(p_1\otimes p_2) $ appartient \`a $\Sigma^{-1}
\F^c_{(1)}(\Sigma \oPo)$.
Ainsi, $\zeta\circ (d_\theta +d_{\theta'}) (\xi)$ est une somme de termes de la
forme
$$\big( (-1)^{|X|+|p_1|} +(-1)^{|X|+|p_1|+1}\big) X\otimes \mu(p_1\otimes p_2)
\otimes Y=0,$$ d'o\`u la conclusion. $\cqfd$
\end{itemize}
\end{deo}

\textsc{Remarque :} On comprend, gr\^ace \`a cette d\'emonstration,
pourquoi on a
introduit un signe $-$ suppl\'ementaire dans la d\'efintion de la d\'erivation
$d_{\theta'}$  de la cobar construction.\\

Si on pose $\mathcal{C}=\bar{\B}(\Po)$, la coprop\'erade
diff\'erentielle coaugment\'ee d\'efinie par la bar construction
sur $\Po$, on peut alors consid\'erer la cobar construction
coaugment\'ee \`a droite $\bar{\B}^c(\mathcal{C})\boxtimes_c
\mathcal{C}=\bar{\B}^c(\bar{\B}(\Po))\boxtimes_c \bar{\B}(\Po)$.
Remarquons que la coprop\'erade $\bar{\B}(\Po)=\F^c(\Sigma \oPo)$
est gradu\'ee par un poids (par la graduation de $\F^c$, en
fonction du nombre de sommets, joint \`a celle de $\Po$) et que la
cobar construction coaugment\'ee
$\bar{\B}^c(\mathcal{C})\boxtimes_c \mathcal{C}$ est un
$\mathcal{C}$-comodules quasi-colibre analytique.

\begin{lem}
Le morphisme $ \zeta\boxtimes_c \bar{\B}(\Po) \ : \
\bar{\B}^c(\bar{\B}(\Po))\boxtimes_c \bar{\B}(\Po) \to
\Po\boxtimes_c \bar{\B}(\Po)$ entre la cobar construction
coaugment\'ee $\bar{\B}^c(\mathcal{C})\boxtimes_c \mathcal{C}$ et
la bar construction augment\'ee $\Po\boxtimes_c \bar{\B}(\Po)$ est
un morphisme de dg-$\Sy$-bimodules. En outre, il s'agit d'un
morphisme de $\bar{\B}(\Po)$-comodules quasi-colibres analytiques.
\end{lem}

\begin{deo}
Les applications $\zeta$ et $id_{\bar{\B}(\Po)}$ \'etant des
morphismes de dg-$\Sy$-bimodules, le morphisme $\zeta\boxtimes_c
id_{\bar{\B}(\Po)}$ pr\'eserve les diff\'erentielles canoniques
$\delta_{\bar{\B}^c(\mathcal{C})}+\delta_{\mathcal{C}}$ et
$\delta_\Po + \delta_{\bar{\B}(\Po)}$. En outre, le morphisme
$d_{\theta_r}$ intervenant dans la d\'efinition de la
diff\'erentielle de la cobar construction
$\bar{\B}^c(\mathcal{C})\boxtimes_c \mathcal{C}$ (\emph{cf.}
chapitre $4$ section $2.2$) correspond bien, via $\zeta
\boxtimes_c \bar{\B}(\Po)$ au morphisme $d_{\theta_l}$ intervenant
dans la d\'efinition de la diff\'erentielle de la bar construction
$\Po\boxtimes_c \bar{\B}(\Po)$. (Le morphisme $d_{\theta_r}$
revient \`a extraire une op\'eration $\Po$ de $\bar{\B}(\Po)$ par
le bas pour la composer, par le haut, \`a
$\bar{\B}^c(\bar{\B}(\Po))$, alors que $d_{\theta_l}$ revient \`a
 extraire
une op\'eration $\Po$ de
$\bar{\B}(\Po)$ par le bas pour la composer, par le haut, \`a
$\Po$.)

Enfin, comme $\zeta$ et $id_{\bar{\B}(\Po)}$ pr\'eservent la
graduation par le poids  venant de celle de $\Po$, on a que
$\zeta\boxtimes_c \bar{\B}(\Po)$ est un morphisme de
$\bar{\B}(\Po)$-comodules quasi-colibres analytiques. $\cqfd$
\end{deo}

On peut maintenant conclure le th\'eor\`eme escompt\'e.

\begin{thm}[R\'esolution bar-cobar]
\label{barcobarresolution}
Pour toute prop\'erade diff\'erentielle augment\'ee gradu\'ee par un poids
$\Po$, le
morphisme
$$\zeta \ :\ \bar{\B}^c(\bar{\B}(\Po))\to \Po$$
est un quasi-isomorphisme de prop\'erades diff\'erentielles
gradu\'ees par un poids.
\end{thm}

\begin{deo}
On sait d'apr\`es le th\'eor\`eme~\ref{acyclicitécobarcoaugmentée}
que le complexe $ \bar{\B}^c(\bar{\B}(\Po))\boxtimes_c
\bar{\B}(\Po)$ est acyclique (cobar construction augment\'ee). De
m\^eme, le th\'eor\`eme~\ref{acyclicitébaraugmentée} affirme que
le complexe $\Po\boxtimes_c \bar{\B}(\Po)$ est lui aussi acyclique
(bar construction augment\'ee). On en conclut que le morphisme de
comodules quasi-colibres analytiques
$\zeta\boxtimes_c\bar{\B}(\Po)$
 est un quasi-isomorphisme.
En posant $\Psi=id_{\bar{\B}(\Po)}$ et $\Phi=\zeta\boxtimes_c
\bar{\B}(\Po)$, on peut appliquer la partie $(2)$ du th\'eor\`eme
de comparaison des comodules quasi-colibres analytiques, ce qui
donne que $\zeta$ est un quasi-isomorphisme. $\cqfd$
\end{deo}

\textsc{Remarque :} Nous utiliserons principalement ce r\'esultat dans le
cas o\`u la prop\'erade $\Po$ est quadratique (donc gradu\'ee par le nombre de
sommets des graphes en jeu).
Gra\^ce au lemme de comparaison entre les prop\'erades
quasi-libres, d\'emontr\'e dans la section suivante, on montrera
au chapitre $7$ que, dans le cas o\`u la prop\'erade quadratique de d\'epart
est de Koszul,  la r\'esolution bar-cobar peut se ``simplifier'' pour
donner le mod\`ele minimal sur $\Po$.\\

On peut faire sensiblement le m\^eme travail dans le cas des
PROPs.\\

Soit $(\Po,\, \mu,\, conc)$ un PROP diff\'erentiel augment\'e. On d\'efinit le
morphisme $\zeta'$ par la composition

$$ \bar{\B}^c(\bar{\B}(\Po) \twoheadrightarrow S_\otimes \F(\Sigma^{-1}
\bar{\F}^c(\Sigma \oPo))=S_\otimes(\F(\oPo)) \xrightarrow{c'_\Po}
\Po,$$ o\`u $c'_\Po =conc \circ S_\otimes(c_\Po)$. Cette
derni\`ere application revient \`a effectuer toutes les
compositions possibles (verticales et horizontales) d'op\'erations
de $\oPo$ suivant le sch\'ema de composition donn\'e par
$S_\otimes (\F)$.

\begin{pro}
Le morphisme $\zeta'$ est un morphisme de PROPs diff\'erentiels gradu\'es par
un poids.
\end{pro}

\begin{deo}
Comme la d\'efinition de $\zeta'$ repose sur les compositions $\mu$ et $conc$
du PROP $\Po$, on voit que ce morphisme est un morphisme de PROPs.

En outre, le PROP $\Po$ est un PROP gradu\'e par un poids, ce qui implique que les
compositions $\mu$ et $conc$ pr\'eserve ce poids. Il en est donc de m\^eme pour
$\zeta'$ qui est un morphisme de PROPs gradu\'e par un poids.

Pour voir que $\zeta'$ pr\'eserve les diff\'erentielles resp\'ectives, on l'\'ecrit
comme compos\'e des morphismes de $\Sy$-bimodules diff\'erentiels. Le diagramme
suivant est commutatif.
\begin{center}
$$\xymatrix{{\bar{\B}^c(\bar{\B}(\Po))=\bar{\B}^c(S_\otimes(
\F^c(\Sigma \oPo)))} \ar[r]^(0.7){\zeta'} \ar@{>>}[d]   & {\Po} \\
{\bar{\B}^c(\F^c(\Sigma \oPo))=S_\otimes(\bar{\F}(\Sigma^{-1}
\bar{\F}^c(\Sigma \oPo)))} \ar[r]^(0.7){S_\otimes(\zeta)}
 & {S_\otimes(\Po).} \ar[u]_{conc} }$$
\end{center}
$\cqfd$
\end{deo}

\begin{lem}
Le morphisme $\zeta'\boxtimes \bar{\B}(\Po) \, : \,
\bar{\B}^c(\bar{\B}(\Po))\to \Po \boxtimes \bar{\B}(\Po)$ est un
morphisme de $\bar{\B}(\Po)$-comodules quasi-colibres analytiques.
\end{lem}

\begin{deo}
Comme le morphisme $\zeta'$ est une version concat\'en\'ee du
morphisme $\zeta$, la d\'e\-mon\-stra\-tion reste la m\^eme.
$\cqfd$
\end{deo}

\begin{thm}[R\'esolution bar-cobar pour les PROPs]
\label{barcobarresolutionPROP} Pour tout PROP diff\'erentiel
augment\'e gradu\'e par un poids, le morphisme
$$\zeta' \ : \ \bar{\B}^c(\bar{\B}(\Po))\to \Po $$
est un quasi-isomorphisme de PROPs diff\'erentiels gradu\'es par un poids.
\end{thm}

\begin{deo}
La d\'emonstration est la m\^eme que dans le cas des prop\'erades.
On utilise ici que la bar et la cobar constructions augment\'ees
sont toujours acycliques dans le cas des PROPs. Et on conclut en
utilisant les versions PROPiques des lemmes de comparaison.
$\cqfd$
\end{deo}

\section{Au niveau des prop\'erades quasi-libres}

On montre ici un lemme de comparaison du m\^eme type que les
pr\'ec\'edents, mais au niveau des prop\'erades quasi-libres.

\begin{thm}[Lemme de comparaison des prop\'erades quasi-libres]
\label{lemmecomparaisonproperades} Soient $M$ et $M'$ deux
dg-$\Sy$-bimodules gradu\'es par un poids et de degr\'e au moins
$1$. Soient $\Po$ et $\Po'$ deux prop\'erades quasi-libres de la
forme $\Po=\F(M)$ et $\Po'=\F(M')$, munies de d\'erivations
$d_\theta$ et $d_{\theta'}$ provenant de morphismes $\theta \, :\,
M \to \bigoplus_{s\ge 2}\F_{(s)}(M)$ et
 $\theta' \, :\, M' \to \bigoplus_{s\ge 2}\F_{(s)}(M')$ qui pr\'eservent la
graduation totale venant de celle de $M$ et $M'$. Et soit, un morphisme
de dg-$\Sy$-bimodules $\Phi \, :\, \Po \to \Po'$  qui respecte la graduation
de $\F$ et la graduation totale. Ainsi, $\Phi$ induit un morphisme $\bar{\Phi}
\, :\, M=\F_{(1)}(M) \to M'=\
F_{(1)}(M')$.\\

Le morphisme $\Phi$ est un quasi-isomorphisme si et seulement si $\bar{\Phi}$
est un quasi-isomorphisme.
\end{thm}

\begin{deo}

Comme les morphismes $\theta$ et $\theta'$ pr\'eservent la graduation totale
venant de $M$ et $M'$, on peut d\'ej\`a remarquer que les d\'erivations
$d_\theta$ et $d_{\theta'}$ pr\'eservent aussi la graduation totale. Par
cons\'equent, les diff\'erentielles $\delta_\theta=\delta_M+d_\theta$ et
 $\delta_{\theta'}=\delta_{M'}+d_{\theta'}$ pr\'eservent cette graduation et
on peut donc d\'ecomposer les complexes $\F(M)=\bigoplus_{\rho \in \mathbb{N}}
\F(M)^ {(\rho)}$ et $\F(M')=\bigoplus_{\rho \in \mathbb{N}}
\F(M')^ {(\rho)}$ en fonction de cette graduation.

On introduit la filtration suivante
$$F_s(\F(M))=\bigoplus_{r\ge -s} \F_{(r)}(M). $$
Cette filtration est compatible avec la d\'ecomposition pr\'ec\'edente. Pour
des probl\`emes de convergence de suites spectrales, on s'interesse plut\^ot
au sous-complexe
$\F(M)^ {(\rho)}$ et \`a la filtration
$$F_s(\F(M)^{(\rho)})=\bigoplus_{r\ge -s} \F_{(r)}(M)^{(\rho)}.$$
Comme les dg-$\Sy$-bimodules $M$ et $M'$ sont de degr\'e
strictement positif pour la graduation par le poids, on a que
$\F_{(r)}(M)^{(\rho)}= \F_{(r)}(M')^{(\rho)}=0$ d\`es que
$r>\rho$, ainsi
$$F_s(\F(M)^{(\rho)})=\bigoplus_{-s \le r \le \rho} \F_{(r)}(M)^{(\rho)}.$$

On peut voir que la filtration $F_s$ est stable par la diff\'erentielle
$\delta_\theta=\delta_M+d_\theta$. En effet, on a $\delta_M(F_s)\subset F_s$ et
$d_\theta(F_s) \subset F_{s-1}$. Cette filtration induit donc une suite
spectrale $E^*_{s,\,t}$. Le premier terme de cette suite spectrale vaut
$$E^0_{s,\,t}=F_s(\F(M)^{(\rho)}_{s+t})/F_{s-1}(\F(M)^{(\rho)}_{s+t})
=\F_{-s}(M)^{(\rho)}_{s+t},  $$ et la diff\'erentielle $d^0$
correspond \`a la diff\'erentielle canonique issue de M :
$\delta_M$. De cette description et des propri\'et\'es de $\Phi$,
on tire que $E^0_{s,\,t}(\Phi)=\F_{(-s)}(\bar{\Phi})_{s+t}$.

On peut remarquer que la filtration $F_s$ sur $\F(M)^{(\rho)}$ est
born\'ee ($F_{-\rho-1}=0$ et $F_0=\F(M)^{(\rho)}$), ainsi par le th\'eor\`eme
classique de convergence des suites spectrales (\emph{cf.}
\cite{Weibel} 5.5.1)
on obtient que la suite spectrale $E^*_{s,\,t}$ converge vers l'homologie du
sous-complexe $\F(M)^{(\rho)}$.\\

D\'emontrons maintenant l'\'equivalence souhait\'ee:\\

($\Leftarrow $)   Si $\bar{\Phi}\, :\, M \to M'$ est un
quasi-isomorphisme,
comme le foncteur $\F$ est un foncteur exact (\emph{cf.}
proposition~\ref{Ffoncteurexact}), on obtient que
$$\F_{(-s)}(\bar{\Phi})_{s+t}=E^0_{s,\, t}(\Phi^{(\rho)}) \, : \,
E^0_{s,\,t}(\F(M)^{(\rho)}) \to E^0_{s,\,t}(\F(M')^{(\rho)})$$ est
un quasi-isomorphisme, c'est-\`a-dire $E^1(\Phi^{(\rho)})$ est un
isomorphisme. Enfin la convergence des suites spectrales en jeu,
montre que $\Phi^{(\rho)} \, :\, \F(M)^{(\rho)} \to
\F(M')^{(\rho)}$ est un quasi-isomorphisme. Ainsi, le morphisme
$\Phi$ est un quasi-isomorphisme.\\

($\Rightarrow$)  R\'eciproquement, si $\Phi$ est un quasi-isomorphisme, alors
chaque $\Phi^{(\rho)}$ est un quasi-isomorphisme. On va montrer par
r\'ecurrence sur l'entier $\rho$ que $\bar{\Phi}^{(\rho)}\, :\, M^{(\rho)}
\to M'^{(\rho)}$ est un quasi isomorphisme.

Pour $\rho=0$, on a $M^{(0)}=M'^{(0)}=0$, le morphisme $\bar{\Phi}^{(0)}$ est
donc un quasi-isomorphisme et la r\'ecurrence est fond\'ee.

Supposons que le r\'esultat soit vrai jusqu'\`a $\rho$ et montrons le pour
$\rho+1$. Dans ce cas, comme $E^0_{s,\,t}(\F(M)^{(\rho +1)})=
\F_{(-s)}(M)^{(\rho+1)}_{s+t}$,
 on voit que $E^1_{s,\, t}(\Phi^{(\rho +1)})$ est un
isomorphisme pour tout $t$, d\`es que $s<-1$.

On va \`a nouveau utiliser le c\^one de l'application $\Phi^{(\rho+1)}$. Pour
cela, on d\'efinit la m\^eme filtration sur $\cone(\Phi^{(\rho+1)})$ que
sur les $\F(M)^{(\rho+1)}$ :
$$F_s(\cone(\Phi^{(\rho+1)}))=F_s(\Sigma^{-1} \F(M)^{(\rho+1)})\oplus
F_s(\F(M')^{(\rho+1)}).$$

Cette filtration induit une suite spectrale qui v\'erifie
$$E^0_{s,\, *}(\cone(\Phi^{(\rho+1)}))=\cone(E^0_{s,\,*}(\Phi^{\rho+1})).$$

Par le m\^eme argument que celui de la d\'emonstration du lemme de
comparaison des $\Po$-modules quasi-libres, on a la suite exacte
longue
\begin{eqnarray*}
&&\xymatrix@C=20pt{\cdots \ar[r]& E^1_{s,\,t+1}(\F(M)^{(\rho+1)})
\ar[r] & E^1_{s,\,t+1}(\F(M')^{(\rho+1)}) \ar[r] &
E^1_{s,\,t}(\cone(\Phi^{(\rho+1)})) }\\
&&\xymatrix@C=20pt{\ar[r] &  E^1_{s,\,t}(\F(M)^{(\rho+1)})\ar[r]
& E^1_{s,\,t} (\F(M')^{(\rho+1)}) \ar[r] &
E^1_{s,\,t-1}(\cone(\Phi^{(\rho+1)}))\ar[r]& \cdots .}
\end{eqnarray*}

Comme $E^1_{s,\,t}(\Phi^{(\rho+1)})$ est un isomorphisme, lorsque $s<-1$, cette
suite exacte longue montre que $E^1_{s,\,t}(\cone(\Phi^{(\rho+1)})=0$ pour tout
$t$ d\`es que $s<-1$. Enfin, la forme de la filtration utilis\'ee donne que
$E^1_{s,\,t}(\cone(\Phi^{(\rho+1)}))=0$ pour tout $t$ si $s\ge0$. La suite
spectrale $E^*_{s,\,t}(\cone(\Phi^{(\rho+1)}))$ est donc d\'eg\'en\'er\'ee au
rang $E^1$ (seule la colonne $s=-1$ est non nulle). Ceci implique que
$$E^\infty_{-1,\,t}(\cone(\Phi^{(\rho+1)}))=
E^1_{-1,\,t}(\cone(\Phi^{(\rho+1)})).$$
En appliquant le
th\'eor\`eme classique de converge des suites spectrales, on a que
la suite spectrale $E^*_{s,\,t}(\cone(\Phi^{(\rho+1)}))$ converge
vers l'homologie du  c\^one de $\Phi^{(\rho+1)}$ qui est nulle
puisque $\Phi^{(\rho+1)}$ est un quasi-isomorphisme. On a donc que
$E^1_{-1,\,t}(\cone(\Phi^{(\rho+1)}))=0$, pour tout $t$, ce qui
donne, une fois r\'einject\'e dans la suite exacte longue, que
$$ E^0_{-1,\,t}(\Phi^{(\rho+1)})= \bar{\Phi}^{(\rho+1)} \ : \ M^{(\rho+1)} \to
 M'^{(\rho+1)}.$$
est un quasi-isomorphisme. D'o\`u le r\'esultat escompt\'e.
$\cqfd$
\end{deo}

\textsc{Remarque :} Ce th\'eor\`eme g\'en\'eralise aux
prop\'erades un th\'eor\`eme de B. Fresse \cite{Fresse} pour les
op\'erades ``connexes'' (c'est-\`a-dire des op\'erades $\Po$
telles que $\Po(0)=\Po(1)$). Cette hypoth\`ese technique de
connexit\'e est l\`a pour assurer la convergence de la suite
spectrale dans la d\'emonstration. Le probl\`eme est que le
r\'esultat de B. Fresse ne s'applique pas aux alg\`ebres, alors
que le th\'eor\`eme que nous proposons ici s'applique aux
alg\`ebres, aux op\'erades et aux prop\'erades. La seule
restriction vient du fait qu'il faille prendre des objets
gradu\'es par un poids, ce qu'il est le cas de tous les mono\"\i des quadratiques
rencontr\'es ici.  \\

On peut \'etendre ce th\'eor\`eme au cadre des PROPs.

\begin{thm}[Lemme de comparaison des PROPs quasi-libres]
\label{lemmecomparaisonPROPs}
Soient $M$ et $M'$ deux dg-$\Sy$-bimodules gradu\'es par un poids et
de poids au moins $1$. Soient
$\Po$ et $\Po'$ deux PROPs quasi-libres de la forme $\Po=S_\otimes(\F(M))$
et
$\Po'=S_\otimes(\F(M'))$, munis de d\'erivations $d_\theta$ et $d_{\theta'}$
provenant de
morphismes $\theta \, :\, M \to \bigoplus_{s\ge 2}\F_{(s)}(M)$ et
 $\theta' \, :\, M' \to \bigoplus_{s\ge 2}\F_{(s)}(M')$ qui pr\'eservent la
graduation totale venant de celle de $M$ et $M'$. Et soit, un morphisme
de dg-$\Sy$-bimodules $\Phi \, :\, \Po \to \Po'$  qui respecte la graduation
de $S_\otimes(\F)$ et la graduation totale. Ainsi, $\Phi$ induit un morphisme
$\bar{\Phi}
\, :\, M=\F_{(1)}(M) \to M'=\
F_{(1)}(M')$.\\

Le morphisme $\Phi$ est un quasi-isomorphisme si et seulement si $\bar{\Phi}$
est un quasi-isomorphisme.
\end{thm}

\begin{deo}
Dans un sens, on se sert du fait que le foncteur $S_\otimes$ est
un foncteur exact. La d\'emonstration de l'implication
r\'eciproque se montre de la m\^eme mani\`ere que dans le cas des
prop\'erades. $\cqfd$
\end{deo}

Nous utilisons ces th\'eor\`emes au chapitre $7$ lors des d\'emonstrations
des principaux r\'esultats de ce papier.

\chapter{Bar construction simpliciale}

\thispagestyle{empty}

On donne dans ce chapitre une autre g\'en\'eralisation de la bar
construction des alg\`ebres associatives dans le cadre mono\"\i
dal des prop\'erades diff\'erentielles. On montre le m\^eme type
de r\'esultat que pour la bar construction (diff\'erentielle) du
chapitre $4$ (acyclicit\'e de la bar construction augment\'ee
notamment). Puis, on construit explicitement un morphisme qui
devrait r\'ealiser un quasi-isomorphisme entre la bar construction
et la bar construction simpliciale.

\section{D\'efinitions et premi\`eres propri\'et\'es}

Dans toute cat\'egorie mono\"\i dale $(\mathcal{A},\, \Box,\,
I)$, on peut construire un complexe
simplicial \`a partir d'un mono\"\i de $M$ en consid\'erant les objets
$M^{\Box n}$,
pour tout entier $n$. Dans la cat\'egorie des $k$-modules, on obtient sur
la bar construction des alg\`ebres associatives. Dans le cas des monades, on
retrouve la bar construction des triples de J. Beck \cite{Beck}. On applique
cette construction \`a la cat\'egorie mono\"\i dale des $\Sy$-bimodules
diff\'erentiels.

\subsection{Bar construction simpliciale \`a coefficients}

\begin{dei}[Bar construction simpliciale \`a coefficients]
\index{bar construction simpliciale \`a coefficients}

Soit $\Po$ une prop\'erade diff\'erentielle.
Soient $(L,\, l)$ et $(R,\, r)$ deux $\Po$-modules diff\'erentiels
\`a droite et \`a gauche.
On pose
$$\mathcal{C}_n(L,\, \Po,\, R)=L\boxtimes_c \underbrace{\Po\boxtimes_c \cdots
\boxtimes_c \Po}_n \boxtimes_c R.$$ On d\'efinit les faces $d_i\,
:\, \mathcal{C}_n(L,\,\Po,\, R) \to \mathcal{C}_{n-1}(L,\,\Po,\,
R)$ par \vspace{6pt}

\begin{tabular}{l}
$\bullet$ l'action \`a droite $l$ si $i=0$, \\

$\bullet$ la composition $\mu$ de la $i^{\textrm{\`eme}}$ ligne, compos\'ee
d'\'el\'ements de $\Po$,
avec la  ${(i+1)}^{\textrm{\`eme}}$,\\

$\bullet$ l'action \`a gauche $l$ si $i=n$.

\end{tabular}\\

Les d\'eg\'en\'er\'escences $s_j \, : \, \mathcal{C}_n(L,\, \Po,\, R)
\to \mathcal{C}_{n+1}(L,\, \Po,\, R)$ sont donn\'ees par l'insertion de
l'unit\'e $\eta$ de la prop\'erade : $L\boxtimes_c \Po^{\boxtimes j}
\boxtimes_c \eta \boxtimes_c
\Po^{\boxtimes_c n-j} \boxtimes_c R$.\\

On munit $\mathcal{C}(L,\, \Po,\, R)$ de la diff\'erentielle
$d_\mathcal{C}$, somme des diff\'erentielles canoniques avec la
diff\'e\-ren\-ti\-elle simpliciale :
$$ d_\mathcal{C}=\delta_L + \delta_\Po +\delta_R +\sum_{i=0}^n
(-1)^{i+1} d_i.$$
Ce complexe est appel\'e \emph{bar construction simpliciale de
$\Po$ \`a coefficients dans $L$ et $R$}.
\end{dei}

\textsc{Remarque :} Lorsque l'on travaille dans la cat\'egorie non
diff\'erentielle des $\Sy$-bimodules, cette construction est
simpliciale au sens strict du terme.\\

On d\'efinit un \emph{morphisme d'augmentation} $\varepsilon \,  :\,
\mathcal{C}(L,\,
\Po,\, R) \to L{\boxtimes_c}_\Po R$ gr\^ace \`a la projection canonique
$$ \mathcal{C}_0(L,\, \Po,\, R)=L\boxtimes_c R \to L{\boxtimes_c}_\Po R. $$

Comme dans la cadre de la bar construction, on a ici une notion de
bar construction simpliciale r\'eduite.

\begin{dei}[Bar construction simpliciale r\'eduite]
\index{bar construction simpliciale r\'eduite}

La bar construction simpliciale \`a coefficients triviaux
$\mathcal{C}(I,\, \Po ,\, I)$ est
appel\'ee \emph{bar construction simpliciale r\'eduite}.
On la note $\bar{\mathcal{C}}(\Po)$.
\end{dei}

\subsection{Bar construction simpliciale augment\'ee}

\begin{dei}[Bar construction simpliciale augment\'ee]
\index{bar construction simpliciale augment\'ee}

Les complexes de cha\^\i nes $\mathcal{C}(I,\, \Po,\, \Po)=
\bar{\mathcal{C}}(\Po)\boxtimes_c \Po$ et $ \mathcal{C}(\Po,\,
\Po,\, I)=\Po \boxtimes_c \bar{\mathcal{C}}(\Po)$ sont appel\'es
\emph{bar constructions simpliciales augment\'ees \`a droite et
\`a gauche}.
\end{dei}

\begin{pro}
Pour toute prop\'erade $\Po$ et tout $\Po$-module \`a droite $L$, la bar
construction simpliciale $\mathcal{C}(L,\, \Po,\, \Po)$ est un $\Po$-module
quasi-libre
analytique \`a droite. Et, le morphisme d'augmentation
$$\varepsilon \ : \ \mathcal{C}(L,\, \Po,\, \Po) \to L $$
est un quasi-isomorphisme de $\Po$-modules simpliciaux \`a droite.
\end{pro}

On a \'evidement le m\^eme r\'esultat \`a gauche pour tout $\Po$-module $R$.

\begin{deo}
La forme de la diff\'erentielle sur $\mathcal{C}(L,\, \Po,\, \Po)=L
\boxtimes_c \bar{\mathcal{C}} \boxtimes_c
\Po$ montre qu'il s'agit bien d'un $\Po$-module
quasi-libre. Le c\^ot\'e analytique de la construction vient du foncteur $\Po
\to L\boxtimes_c \Po$ qui est analytique.

Le reste de la d\'emonstration est le m\^eme que dans le cas des
alg\`ebres associatives unitaires (\emph{cf.} \cite{Cartan}). On
introduit une d\'eg\'en\'er\'escence suppl\'ementaire
$$ s_{n+1}
\, : \, L \boxtimes_c \Po^{\boxtimes_c n} \boxtimes_c \Po \simeq L
\boxtimes_c \Po^{\boxtimes_c n+1}\boxtimes_c I
\xrightarrow{L\boxtimes_c \Po^{\boxtimes_c n+1} \boxtimes_c \eta}
L \boxtimes_c \Po^{\boxtimes_c n+1} \boxtimes_c \Po, $$ qui induit
une homotopie contractante. $\cqfd$
\end{deo}

\begin{cor}
Les bar constructions simpliciales augment\'ees \`a gauche et \`a droite sont
acycliques.
\end{cor}

\subsection{Bar construction normalis\'ee}

Comme pour tout module simplicial, on r\'eduit note \'etude au complexe
normalis\'e,
complexe de cha\^\i nes quotient du complexe de d\'epart mais ayant la m\^eme
homologie.

\begin{dei}[Bar construction normalis\'ee]
\index{bar construction normalis\'ee}

La \emph{bar construction normalis\'ee} correspond au quotient de la
bar construction simpliciale par les images des d\'eg\'en\'er\'escences.
On pose
$$ \N_n(L,\, \Po,\, R)=\coker\left( L\boxtimes_c \Po^{\boxtimes_c (n-1)}
\boxtimes_c R \xrightarrow{\sum_{i=0}^{n-1} L\boxtimes_c
\Po^{\boxtimes_c i} \boxtimes_c \eta \boxtimes_c \Po^{\boxtimes_c
(n-i-1)} \boxtimes_c R} L \boxtimes_c \Po^{\boxtimes_c n}
\boxtimes_c R  \right). $$ Le complexe de cha\^\i nes $\N(I,\,
\Po,\, I)$, not\'e $\bar{\N}(\Po)$, est appel\'e \emph{bar
construction normalis\'ee r\'eduite}.
\end{dei}

\textsc{Remarque :} Nous avons vu pr\'ec\'edemment que la bar
construction $\B(L,\, \Po,\, R)$ se repr\'esentait avec des
graphes et que la cod\'erivation $d_\theta$ correspondait \`a la
notion g\'en\'eralis\'ee d'\emph{edge contraction}, qui revient
\`a composer les paires d'op\'erations adjacentes. Au contraire,
la bar construction simpliciale se repr\'esente par des graphes
\`a niveaux et les faces $d_i$ correspondent \`a des compositions
entre deux niveaux d'op\'erations.

\section{Morphisme d'\'echelonnement}

Le \emph{morphisme d'\'echelonnement} est un morphisme injectif
entre la bar construction et la bar construction simpliciale qui
induit un isomorphisme en homologie.

\subsection{D\'efinition}

\begin{figure}[h]
$$\xymatrix{& \ar[d] &  \ar[d]&  \\
 \ar@{--}[r] & *+[F-,]{\nu_1} \ar@{--}[r] \ar@/_1pc/[d] \ar[dr]   &
*+[F-,]{\nu_2} \ar@/^1pc/[d]  \ar[dl] | \hole  \ar@{--}[r] & \\
 \ar@{--}[r] & *+[F-,]{\mu_1}  \ar[d]  \ar@{--}[r] &
*+[F-,]{\mu_2}  \ar[d]  \ar@{--}[r] &  \\
 & & &  } $$
\caption{Un exemple de $g_\xi$.}  \label{xi}
\end{figure}
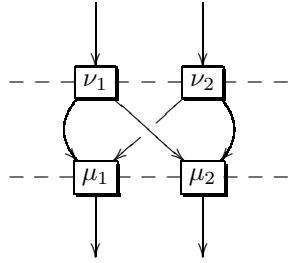

Soit $\xi$ un \'el\'ement de $\bar{\B}_{(n)}(\Po)=\F^c_{(n)}(\Sigma \oPo)$
repr\'esent\'e par un graphe $g_\xi$ \`a $n$ sommets (\emph{cf.} figure~\ref{xi}).

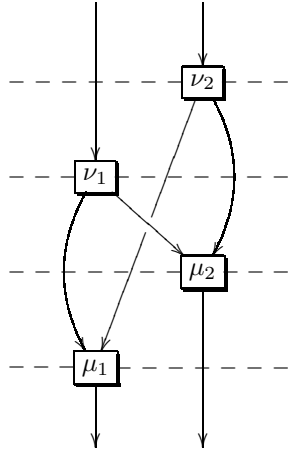
\begin{figure}[h]
$$\xymatrix{& \ar[dd] &  \ar[d]&  \\
 \ar@{--}[rr] &    &
*+[F-,]{\nu_2} \ar@/^1pc/[dd]  \ar[dddl] | \hole  \ar@{--}[r] & \\
 \ar@{--}[r] & *+[F-,]{\nu_1} \ar@/_1pc/[dd]    \ar[dr] \ar@{--}[rr] &
& \\
 \ar@{--}[rr] & &
*+[F-,]{\mu_2}  \ar[dd]  \ar@{--}[r] &  \\
 \ar@{--}[r] & *+[F-,]{\mu_1}  \ar[d]  \ar@{--}[rr] &
 &  \\
 & & &  } $$
\caption{Un exemple de $g\in \mathcal{G}_4(\xi)$.} \label{echel}
\end{figure}

L'\'el\'ement $\xi$ vient d'un graphe $g_r$, \`a $r$ niveaux
pass\'e au quotient par la relation $\equiv$ (\emph{cf.} chapitre
$3$ section $3$). Rappelons que la relation $\equiv$ revient \`a
changer une op\'eration de niveau, au signe pr\`es. Gr\^ace \`a
cette relation d'\'equivalence, on peut repr\'esenter $\xi$ avec
au moins un graphe \`a $n$ niveaux et $n$ sommets (c'est-\`a-dire
un sommet par niveau, \emph{cf.} figure~\ref{echel}). Posons,
$\mathcal{G}_n(\xi)$ l'ensemble des graphes \`a $n$ niveaux avec
un sommet par niveau, qui redonne $g_\xi$ apr\`es passage au
quotient par la relation $\equiv$.

Soit $g$ un graphe de $\mathcal{G}_n(\xi)$. L'\'el\'ement $\xi$ est
enti\`erement d\'etermin\'e par $g$ et par la suite des op\'erations
$\Sigma p_1 \otimes \cdots \otimes \Sigma p_n$ qui indicent les sommets de $g$.
Pour pouvoir associer \`a $\xi$ un \'el\'ement de $\bar{\N}(\Po)$, il faut
d\'esuspendre les op\'erations de $\xi$. Ceci fait appara\^\i tre des signes et
 on pose
$$g(\xi)=(-1)^{\sum_{i=1}^n(n-i)|p_i|} g(p_1\otimes \cdots \otimes p_n),$$
o\`u $g(p_1\otimes \cdots \otimes p_n)$ correspond au graphe $g$
dont les sommets sont indic\'es par les op\'erations $p_1,\, \ldots,\, p_n$. Ce
signe est obtenu en faisant passer toutes les suspensions \`a gauche.

\begin{dei}[Morphisme d'\'echelonnement]
\index{morphisme d'\'echelonnement}

On d\'efinit le \emph{morphisme d'\'echelonnement} $$e \ : \
\bar{\B}(\Po) \to \bar{\N}(\Po)$$ par $e(\xi)=\sum_{g\in
\mathcal{G}_n(\xi)} g(\xi)$, pour $\xi$ dans
$\bar{\B}_{(n)}(\Po)$.
\end{dei}

\textsc{Remarque :} Dans la d\'efinition de $\xi$ on a pris garde de respecter
les r\`egles de signes dues aux commutations d'\'el\'ements de $\Sigma \oPo$.
Gr\^ace \`a ceci, le morphisme d'\'echelonnement $e$ est bien d\'efini.\\

Comme la d\'efinition de $e$ fait intervenir des permutations de suspensions
vers la gauche, on peut \'etendre naturellement le morphisme d'\'echelonnement
aux bar constructions \`a coefficients \`a droite dans un $\Po$-module
diff\'erentiel $R$ :
$$e \ : \  \B(I,\, \Po,\, R) \to \N(I,\, \Po,\, R).$$

\subsection{Propri\'et\'es homologiques}

Nous esp\'erons montrer que le morphisme d'\'echelonnement est un
quasi-isomorphisme de dg-$\Sy$-bimodules.

\begin{lem}
Soit $\Po$ une prop\'erade diff\'erentielle augment\'ee.
Soient $R$ un $\Po$-module diff\'erentiel \`a
gauche.

Le morphisme d'echelonnement induit un morphisme injectif de
dg-$\Sy$-bimodules
$$ e \ : \ \B_{(n)}(I,\, \Po,\,R) \to \Sigma^n \N_n(I,\, \Po,\,R).$$
\end{lem}

\begin{deo}
Notons $d_\B$ la diff\'erentielle de $\B(I,\, \Po ,\, R)$. Elle est
compos\'ee de $4$ termes :
$$ d_\B = \delta_\Po +\delta_R +d_\theta +d_{\theta_R}.$$
Les deux premiers correspondent aux diff\'erentielles canoniques
issues de $\Po$ et de $R$. Le morphisme $d_\theta$ est la
cod\'erivation qui provient du produit partiel sur $\Po$. Quant \`a
$d_{\theta_R}$, il vient de l'action d'un seul \'el\'ement de $\oPo$
sur $R$.

Appelons $d_\N$ la diff\'erentielle sur $\Sigma ^n \N_n(I,\, \Po,\,R)$. Elle est
aussi compos\'ee de $4$ termes. Sur un \'el\'ement $\tau'=\Sigma^n \tau$ de
$\Sigma^n \N_n(I,\,
\Po,\, R)$, elle s'\'ecrit :
$$ d_\N(\tau') = (-1)^n\Sigma^n \delta_\Po(\tau) +(-1)^n\Sigma^n \delta_R(\tau) +
\sum_{i=1}^{n-1}(-1)^{i+1} \Sigma^{n-1} d_i(\tau)  +(-1)^{n+1} \Sigma^{n-1}
d_n(\tau).$$

Montrons que $d_\N\circ e(\xi)= e\circ d_\B(\xi)$.
\begin{itemize}
\item La commutativit\'e des diff\'erentielles canoniques
$$( \delta_\Po +\delta_R)\circ e(\xi)= e\circ
( \delta_\Po +\delta_R)(\xi)$$
vient des bons choix dans les r\`egles de signes et du respect des suspensions.
(Les calculs sont
du m\^eme type que ceux des chapitres pr\'ec\'edents.)

\item La cod\'erivation $d_\theta$ revient \`a composer les op\'erations
de $\Sigma \oPo$ indi\c cant les couples de sommets adjacents du graphe
repr\'esentant $\xi$. Et, $d_i\circ e$ correspond  \`a composer deux
\'etages avec au total deux op\'erations de $\oPo$. Il faut
distinguer deux cas. Posons
$\xi=X\otimes \Sigma p \otimes \Sigma q \otimes Y$.

\begin{enumerate}

\item Si les op\'erations $\Sigma p$ et $\Sigma q$ sont reli\'es
par au moins une branche, alors la composante de $\sum_{i=1}^{n-1}
(-1)^{i+1}d_i \circ e(\xi)$ faisant intervenir la composition de
$p$ avec $q$ est de la forme $$\sum \varepsilon(-1)^{j+2}
\Sigma^{n-1} X'\otimes \mu(p\otimes q) \otimes Y',$$ o\`u
$X'=p_1\otimes \cdots \otimes p_j$, $Y'=q_1\otimes \cdots \otimes
q_{n-j-2} \otimes \rho_1\otimes \cdots \otimes \rho_r$ et
$$\varepsilon= (-1)^{\sum_{i=1}^j
(n-i)|p_i|+(n-j-1)|p|+(n-j-2)|q|+\sum_{k=1}^{n-j-2}
(n-j-k-2)|q_k|}.$$ Et la composant de $d_\theta(\xi)$ faisant
intervenir la composition de $\Sigma p$ avec $\Sigma q$ est de la
forme
$$(-1)^{|X|+|p|} X\otimes \Sigma \mu(p\otimes q) \otimes Y.$$
L'image par $e$ d'un tel \'el\'ement donne
\begin{eqnarray*}
&& \sum (-1)^{|X|+|p|}\varepsilon(-1)^{|X'|+|p|}\Sigma^{n-1}X'\otimes
\mu(p\otimes q) \otimes Y'=\\
&& \sum \varepsilon (-1)^j \Sigma^{n-1}X'\otimes
\mu(p\otimes q) \otimes Y'.
\end{eqnarray*}

\item Si les op\'erations $\Sigma p$ et $\Sigma q$ ne sont pas
reli\'ees, alors il n'y a aucune composante de $d_\theta$ qui fait
intervenir une composition entre $\Sigma p$ et $\Sigma q$. Et,
la composante de $\sum_{i=0}^n  (-1)^{i+1} d_i\circ e(\xi)$ qui vient de
la composition des deux
\'etages o\`u se trouvent $p$ et $q$ est la somme de deux termes :
\begin{eqnarray*}
&& \sum (-1)^j \Sigma^n d_{j+1} \left(\varepsilon X'\otimes p \otimes q \otimes Y' +
\varepsilon (-1)^{(|p|+1)(|q|+1)+|p|+|q|}
X'\otimes q \otimes p \otimes Y' \right)=\\
&&  \sum (-1)^j \Sigma^n\left( \varepsilon X'\otimes p \otimes q \otimes Y'
+ \varepsilon (-1)^{(|p|+1)(|q|+1) +|p|+|q|+|p||q|} X'\otimes
p\otimes q \otimes Y' \right) =0.
\end{eqnarray*}
\end{enumerate}

\noindent Comme l'image de $\xi$ par $\sum_{i=1}^{n-1} (-1)^{i+1}d_i \circ e - e
\circ d_\theta$ est une somme de termes de la forme $(1)$ ou $(2)$, on
conclut de l'\'etude pr\'ec\'edente la commutativit\'e voulue.

\item Le morphisme $d_{\theta_R}$ revient \`a extraire une op\'eration
$\oPo$ par
le haut et \`a la faire agir par la gauche sur $R$. Or, $e$
consiste, entre autre, \`a ne placer qu'une seule op\'eration $\oPo$ sur
la deuxieme ligne (celle en dessous de $R$) et $d_n$ la fait agir sur
les \'el\'ements de $R$. Il ne reste plus qu'\`a v\'erifier que les signes
correspondent. Posons $\xi=X\otimes \Sigma  p_n \otimes \rho_1 \otimes \cdots
\otimes \rho_m$ o\`u $X=\Sigma p_1 \otimes \Sigma p_{n-1}$. Alors, la composante
de $d_{\theta_R}(\xi)$ faisant intervenir l'action de $p_n$ sur les
\'el\'ements $\rho_1, \ldots ,\, \rho_m$ de $R$ est de la forme
$(-1)^{|X|} X\otimes r(p_n\otimes(\rho_1 \otimes \cdots \otimes \rho_m))$.
Et son image par $e$ donne
\begin{eqnarray*}
&& \sum (-1)^{|X|} \varepsilon (-1)^{|X|+n-1} \Sigma^{n-1} X' \otimes
r(p_n\otimes(\rho_1 \otimes \cdots \otimes \rho_m))   \\
&& \sum (-1)^{n+1} \varepsilon \Sigma^{n-1} X' \otimes
r(p_n\otimes(\rho_1 \otimes \cdots \otimes \rho_m)),
\end{eqnarray*}
o\`u $\varepsilon=\sum_{i=1}^n (n-i)p_i$.
De la m\^eme mani\`ere, la composante de $(-1)^{n+1}d_n \circ e (\xi)$ faisant
intervenir l'action de
$p_n$ sur les \'el\'ements $\rho_1,\, \ldots ,\,  \rho_m$ de $R$ vaut
$$(-1)^{n+1} \sum \varepsilon \Sigma^{n-1} X'\otimes
r(p_n\otimes(\rho_1 \otimes \cdots \otimes \rho_m)).$$
\end{itemize}

Montrons maintenant l'injectivit\'e de $e$. Soit $\xi$ un
\'el\'ement de $\bar{\B}_{(n)}(\Po)=\F_{(n)}^c(\Sigma \oPo)$. Par
d\'efinition de la coprop\'erade colibre, on sait que
$\xi=\xi_1+\cdots +\xi_r$ o\`u chaque $\xi_i$ une somme finie
d'\'el\'ements de $\F_{(n)}^c(\Sigma \oPo)$ qui viennent de
l'indi\c cage des sommets d'un m\^eme graphe $g_{\xi_i}$. On
introduit un morphisme $\pi \, : \, \bar{\N}_n(\Po) \to
\bar{\B}_{(n)}(\Po)$. A un \'el\'ement $\tau$ de $
\bar{\N}_n(\Po)$ repr\'esent\'e par un graphe \`a $n$ niveaux, on
associe l'\'el\'ement correspondant $\tau'$ qui vient de la
suspension des op\'erations de chaque ligne. On introduit ici le
m\^eme signe $\varepsilon$ que celui qui d\'efinit $e$. L'objet
$\pi(\tau)$ est donn\'e par la classe d'\'equivalence de
$\varepsilon \tau'$ pour la relation $\equiv$. On a alors $\pi
\circ e(\xi)=n_1\xi_1+\cdots n_r\xi_r$ o\`u les $n_i$ sont des
entiers non nuls. Ainsi, l'\'equation $e(\xi)=0$ impose
$n_1\xi_1+\cdots n_r\xi_r=0$. Par identification des graphes
sous-jacents, on obtient $\xi_i=0$, pour tout $i$, d'o\`u $\xi=0$.
$\cqfd$
\end{deo}

\textsc{Remarque :} La suspension $\Sigma^n$ est l\`a pour faire commuter les
diff\'erentielles
canoniques avec le morphisme $e$.

\begin{cor}
Dans le cas o\`u $\Po$ et $R$ sont des $\Sy$-bimodules, c'est-\`a-dire de
diff\'erentielle nulle et concentr\'es en degr\'e $0$, le morphisme
d'\'echelonnement
$$ e \ : \ \B(I,\, \Po ,\, R) \to \N(I,\, \Po ,\, R) $$
est un morphisme injectif de dg-$\Sy$-bimodules.
\end{cor}

\begin{conj}
\label{echelonnement} Soit $\Po$ une prop\'erade gradu\'ee par un
poids.

Le morphisme d'\'echelonnement
$$ e \ : \ \bar{\B}(\Po) \to \bar{\N}(\Po) $$
est un quasi-isomorphisme.
\end{conj}

En r\'esum\'e, on esp\`ere montrer que la bar construction
r\'eduite est un sous-complexe de cha\^\i nes de la bar
construction normalis\'ee r\'eduite qui donne la m\^eme homologie.

\textsc{Remarque : } Cette conjecture est une g\'en\'eralisation
aux prop\'erades d'un th\'eor\`eme de B. Fresse \cite{Fresse} pour
les op\'erades.

\chapter{Dualit\'e de Koszul}

\thispagestyle{empty}

A l'aide des r\'esultats des chapitres pr\'ec\'edents, on peut
conclure l'\'etude de la dualit\'e de Koszul des prop\'erades et
des PROPs.

 Pour cela, on commence par d\'efinir les notions de
duale de Koszul d'une prop\'erade (respectivement d'une
coprop\'erade) et d'un PROP (respectivement d'un coPROP). On
montre ensuite que cette duale est une coprop\'erade (coPROP)
quadratique (respectivement une prop\'erade (PROP) quadratique),
puis on fait le lien avec les constructions classiques d'alg\`ebre
et d'op\'erade duales donn\'ees par S. Priddy dans \cite{Priddy}
et par V. Ginzburg et M. M. Kapranov dans \cite{GK}. 

On d\'efinit
la notion de prop\'erade (respectivement PROP) de Koszul et on
montre que le \emph{mod\`ele minimal} de cette prop\'erade est
donn\'e par la cobar construction sur la coprop\'erade duale.
On introduit un petit complexe, appel\'e complexe de Koszul,
dont l'acyclicit\'e est un crit\`ere qui permet de d\'eterminer si
la prop\'erade (respectivement le PROP) est de Koszul ou non.
Puis, on montre qu'un PROP $\Po$ est de Koszul si et seulement si
la prop\'erade associ\'ee $U_c(\Po)$ est de Koszul. Cette
derni\`ere proposition montre que pour \'etudier un PROP, il
suffit d'\'etudier la prop\'erade qui lui est associ\'ee. 

Enfin, on montre que les prop\'erades construites \`a partir de 
deux prop\'erade de Koszul et d'une \emph{loi de remplacement} sont 
de Koszul. En appliquant ce r\'esultat, on montre que la prop\'erade de 
big\`ebres de Lie et celle des big\`ebres de Hopf infinit\'esimales sont de Koszul, 
et on donne leurs duales.

\section{Dual de Koszul}

\index{dual de Koszul}
 \label{constructionduale} Pour une
prop\'erade gradu\'ee par un poids $\Po$, on d\'efinit sa duale de
Koszul comme une sous-coprop\'erade de la bar construction
r\'eduite sur $\Po$. De m\^eme pour une coprop\'erade gradu\'ee
par un poids $\mathcal{C}$, on d\'efinit sa duale de Koszul comme
une prop\'erade quotient de la cobar construction r\'eduite sur
$\mathcal{C}$.\\

Soit $\Po$ une prop\'erade diff\'erentielle gradu\'ee par un poids
augment\'ee. Cette graduation induit une graduation totale
$(\rho)$ sur la bar construction r\'eduite
$\bar{\B}(\Po)=\F^c(\Sigma \oPo)$. Cette derni\`ere est compatible
avec la d\'ecomposition en fonction du nombre de sommets $(s)$ de
la coprop\'erade colibre
$$\bar{\B}_{(s)}(\Po)=\bigoplus_{\rho \in \mathbb{N}}
\bar{\B}_{(s)}(\Po)^{(\rho)}.$$ On rappelle que la cod\'erivation
$d_\theta$, issue du produit partiel sur $\Po$, consiste \`a
composer les paires de sommets adjacents, soit
$$d_\theta\big( \bar{\B}_{(s)}(\Po)^{(\rho)}\big) \subset
\big( \bar{\B}_{(s-1)}(\Po)^{(\rho)}\big).$$

Dans le cas o\`u $\Po$ est connexe ($\Po^{(0)}=I$), la bar construction
a la forme suivante :

\begin{lem}
Soit $\Po$ une prop\'erade diff\'erentielle gradu\'ee par un poids
connexe. On a alors les \'egalit\'es
$$\left\{
\begin{array}{l}
\bar{\B}_{(\rho)}(\Po)^{(\rho)}=\F^c_{(\rho)}(\Sigma \oPo^{(1)}), \\
\bar{\B}_{(s)}(\Po)^{(\rho)}=0 \quad \textrm{si} \quad s>\rho.
\end{array} \right.$$
\end{lem}

\textsc{Remarque :} On a le m\^eme r\'esultat pour la cobar
construction r\'eduite sur une coprop\'erade diff\'erentielle
gradu\'ee par un poids connexe.

\begin{dei}[Duale de Koszul d'une prop\'erade]
Soit $\Po$ une prop\'erade diff\'erentielle  gradu\'ee par un
poids connexe. On d\'efinit la \emph{duale de Koszul de $\Po$} par
le $\Sy$-bimodule diff\'erentiel gradu\'e par un poids
$${\Po^{\ac}}_{(\rho)}=H_{(\rho)}\big( \bar{\B}_*(\Po)^{(\rho)},\,
d_\theta \big). $$
\end{dei}

Le lemme pr\'ec\'edent donne la forme du complexe
$\big( \bar{\B}_*(\Po)^{(\rho)},\,
d_\theta \big)$ :
$$\xymatrix{\cdots \ar[r]^{d_\theta} & 0 \ar[r]^(0.3){d_\theta} &
\bar{\B}_{(\rho)}(\Po)^{(\rho)} \ar[r]^(0.45){d_\theta} &
\bar{\B}_{(\rho-1)}(\Po)^{(\rho)} \ar[r]^(0.6){d_\theta}& \cdots \
.}$$

Ce qui donne l'\'egalit\'e
$${\Po^{\ac}}_{(\rho)}=\textrm{ker}\big(d_\theta
\, : \, \bar{\B}_{(\rho)}(\Po)^{(\rho)} \to
\bar{\B}_{(\rho-1)}(\Po)^{(\rho)} \big).$$

Cette \'ecriture montre bien que ${\Po^{\ac}}_{(\rho)}$ est un
$\Sy$-bimodule diff\'erentiel gradu\'e par un poids, o\`u la
diff\'erentielle est induite par celle de $\Po$ \`a savoir
$\delta_\Po$.

Si la prop\'erade $\Po$ est concentr\'ee en degr\'e $0$, on a

$$\big(\bar{\B}_{(s)}(\Po)^{(\rho)}   \big)_d =\left\{
\begin{array}{ll}
\bar{\B}_{(s)}(\Po)^{(\rho)}  & \textrm{si} \quad d=s, \\
0 & \textrm{sinon}.
\end{array} \right.$$

La coprop\'erade duale n'est pas concentr\'ee en degr\'e $0$ (et
sa diff\'erentielle n'est \`a priori pas nulle) et v\'erifie
$$\big( {\Po^{\ac}}_{(\rho)}\big)_d =\left\{\begin{array}{ll}
 {\Po^{\ac}}_{(\rho)}& \textrm{si} \quad d=\rho, \\
0 & \textrm{sinon}.
\end{array} \right. $$

De la m\^eme mani\`ere, on d\'efinit la duale d'une coprop\'erade.

\begin{dei}[Duale de Koszul d'une coprop\'erade]
Soit $\mathcal{C}$ une coprop\'erade di\-ff\'e\-ren\-tielle
gradu\'ee par un poids connexe. On d\'efinit la \emph{duale de
Koszul de $\mathcal{C}$} par le $\Sy$-bimodule diff\'erentiel
naturellement gradu\'e
$${\mathcal{C}^{\ac}}_{(\rho)}=H_{(\rho)}\big( \bar{\B}^c_*
(\mathcal{C})^{(\rho)},\,d_{\theta'} \big). $$
\end{dei}

La cobar construction munie de la d\'erivation $d_{\theta'}$ est
un complexe de la forme suivante :

$$\xymatrix{\cdots  \ar[r]^(0.4){d_{\theta'}} &
\bar{\B}_{(\rho-1)}(\mathcal{C})^{(\rho)}
\ar[r]^(0.55){d_{\theta'}} &
\bar{\B}_{(\rho)}(\mathcal{C})^{(\rho)}
\ar[r]^(0.65){d_{\theta'}}& 0
 \ .}$$

La duale d'une coprop\'erade v\'erifie donc l'\'egalit\'e

$${\mathcal{C}^{\ac}}_{(\rho)}=\textrm{coker}\big(d_{\theta'}
\, : \, \bar{\B}^c_{(\rho-1)}(\mathcal{C})^{(\rho)} \to
\bar{\B}^c_{(\rho)}(\mathcal{C})^{(\rho)} \big). $$

Le module ${\mathcal{C}^{\ac}}_{(\rho)}$ est donc un
$\Sy$-bimodule diff\'erentiel gradu\'e par un poids, o\`u la
diff\'erentielle est induite par celle de $\mathcal{C}$ \`a savoir
$\delta_\mathcal{C}$.

\begin{pro}
Soit $\Po$ une prop\'erade diff\'erentielle gradu\'ee par un poids
connexe. La duale de Koszul de $\Po$, not\'ee $\Po^{\ac}$, est une
sous-coprop\'erade diff\'erentielle gradu\'ee par un
poids de $\F^c(\Sigma \Po^{(1)})$.\\

Soit  $\mathcal{C}$ une coprop\'erade diff\'erentielle gradu\'ee
par un poids connexe. La duale de Koszul de $\mathcal{C}$,
$\mathcal{C}^{\ac}$, est une prop\'erade diff\'erentielle
gradu\'ee par un poids quotient de $\F(\Sigma^{-1}
\mathcal{C}^{(1)})$.
\end{pro}

\begin{deo}
Pour toute prop\'erade $\Po$ diff\'erentielle gradu\'ee par un
poids connexe, comme
${\Po^{\ac}}_{(\rho)}=\textrm{ker}\big(d_\theta \, : \,
\bar{\B}_{(\rho)}(\Po)^{(\rho)} \to
\bar{\B}_{(\rho-1)}(\Po)^{(\rho)} \big)$, on a imm\'ediatemment
que ${\Po^{\ac}}_{(\rho)}$ est un sous-dg-$\Sy$-bimodule
naturellement gradu\'e de $\F^c_{(\rho)}(\Sigma \Po^{(1)})$. Il
reste \`a montrer que $\Po^{\ac}$ est stable par le coproduit
$\Delta$ de $\F^c(\Sigma \Po^{(1)})$. Soit $\xi$ un \'el\'ement de
${\Po^{\ac}}_{(\rho)}$, c'est-\`a-dire $\xi\in
\F_{(\rho)}^c(\Sigma\Po^{(1)})$ et $d_\theta(\xi)=0$. Posons
$$\Delta(\xi)=\sum_\Xi (\xi^1_1,\, \ldots ,\, \xi^1_{a_1})\, \sigma
\, (\xi^2_1,\, \ldots ,\, \xi^2_{a_2}), $$
o\`u la somme porte sur une famille $\Xi$ de graphes \`a deux
\'etages, avec  $\xi_i^j \in \F_{(s_{ij})}^c(\Sigma \Po^{(1)})$. Le
fait que $d_\theta$ soit
une cod\'erivation signifie que $d_\theta \circ \Delta (\xi)=
\Delta \circ d_\theta (\xi)$. On a donc que
\begin{eqnarray*}
d_\theta \circ \Delta (\xi)&=&
\sum_\Xi  \Big( \sum_{k=1}^{a_1} \pm (\xi^1_1,\, \ldots ,\,
d_\theta(\xi^1_k),\, \ldots ,\,
\xi^1_{a_1})\, \sigma
\, (\xi^2_1,\, \ldots ,\, \xi^2_{a_2}) \\
&+&  \sum_{k=1}^{a_2}
   \pm (\xi^1_1,\, \ldots ,\,
\xi^1_{a_1})\, \sigma
\, (\xi^2_1,\, \ldots ,\, d_\theta(\xi^2_k) \, \ldots ,\,\xi^2_{a_2})
\Big)=0,
\end{eqnarray*}
o\`u les $d_\theta(\xi^j_i)$ appartiennent \`a $\F^c_{(s_{ij}-1)}
(\Sigma (\Po^{(1)}\oplus\Po^{(2)}))$ avec un seul sommet indic\'e par
$\Sigma \Po^{(2)}$.
Par identification, on a donc, pour tout $i,\,j$, que
$d_\theta(\xi^j_i)=0$.\\

De la m\^eme mani\`ere, on voit que, pour toute coprop\'erade
diff\'erentielle gradu\'ee par un poids connexe $\mathcal{C}$, sa
duale de Koszul ${\mathcal{C}^{\ac}}$ est un quotient de
$\F(\Sigma^{-1}\mathcal{C}^{(1)})$. Il reste donc \`a montrer que
le produit $\mu$ sur $\F(\Sigma^{-1} \mathcal{C}^{(1)})$ passe \`a
ce quotient. Pour cela, on consid\`ere le produit
$$\mu \big( (c_1^1,\, \ldots ,\, d_{\theta'}(c),\, \ldots ,\,
c^1_{a_1}) \, \sigma \,
(c^2_1,\, \ldots ,\, c^2_{a_2}) \big), $$
o\`u les $c^j_i$ appartiennent \`a $\F_{(s_{ij})}(\Sigma^{-1}
\mathcal{C}^{(1)})$ et $c$ \`a $\F_{(s)}(\Sigma^{-1} (\mathcal{C}^{(1)}
\oplus \mathcal{C}^{(2)}))$ avec un seul sommet indic\'e par
$\mathcal{C}^{(2)}$. Comme les $c^j_i$ sont des \'el\'ements
de $\F_{(s_{ij})}(\Sigma^{-1}
\mathcal{C}^{(1)})$, on a $d_{\theta'}(c^j_i)=0$. Et, le fait que
$d_{\theta'}$ soit une d\'erivation donne ici que
$$\mu \big( (c_1^1,\, \ldots ,\, d_{\theta'}(c),\, \ldots ,\,
c^1_{a_1}) \, \sigma \, (c^2_1,\, \ldots ,\, c^2_{a_2}) \big)=
d_{\theta'}\big( \mu \big(  (c_1^1,\, \ldots ,\, c,\, \ldots ,\,
c^1_{a_1}) \, \sigma \, (c^2_1,\, \ldots ,\, c^2_{a_2})  \big)
\big).$$ Ainsi, $\mu \big( (c_1^1,\, \ldots ,\, d_{\theta'}(c),\,
\ldots ,\, c^1_{a_1}) \, \sigma \, (c^2_1,\, \ldots ,\, c^2_{a_2})
\big) $ est nulle dans le conoyau de $d_{\theta'}$, d'o\`u le
r\'esultat.

 $\cqfd$
\end{deo}

On d\'efinit les m\^emes objets dans le cadre des PROPs.

\begin{dei}[Dual de Koszul d'un PROP]
Soit $\Po$ un PROP diff\'erentiel augment\'e gradu\'e par un poids et connexe.
Son \emph{dual de Koszul} est donn\'e par le $\Sy$-bimodule diff\'erentiel
gradu\'e par un poids
$$\Po^{\ac}_{(\rho)} = H_{(\rho)}\left( \bar{\B}_*(\Po)^{(\rho)},\, d_\theta
\right).$$
\end{dei}

On a aussi une notion de dualit\'e pour un coPROP.

\begin{dei}[Dual de Koszul d'un coPROP]
Soit $\mathcal{C}$ un coPROP diff\'erentiel gradu\'e par un poids et connexe.
Son \emph{dual de Koszul} est donn\'e par le $\Sy$-bimodule diff\'erentiel gradu\'e
par un poids
$${\mathcal{C}^{\ac}}_{(\rho)}=H_{(\rho)}\big( \bar{\B}^c_*
(\mathcal{C})^{(\rho)},\,d_{\theta'} \big). $$
\end{dei}

On peut relier les notions de dualit\'e au niveau des PROPs et coPROPs avec
celles des prop\'erades et coprop\'erades.

\begin{pro}
\label{dualPROPpropérade} Soit $\Po$ un PROP diff\'erentiel
augment\'e gradu\'e par un poids et connexe. On rappelle que
$U_c(\Po)$ repr\'esente la prop\'erade associ\'ee \`a $\Po$. Le
dual $\Po^{\ac}$ de $\Po$ est isomorphe en tant que coPROP
diff\'erentiel gradu\'e par un poids \`a
$S_\otimes(U_c(\Po)^{\ac})$ et il s'agit d'un sous-coPROP de
$S_\otimes(\F^c(\Sigma \Po^{(1)}))$.
\end{pro}

\begin{deo}
Tout repose sur l'isomorphisme de coPROPs diff\'erentiels gradu\'es par un poids
$$ \bar{\B}(\Po)=S_\otimes(\bar{\B}(U_c(\Po))).$$
$\cqfd$
\end{deo}

\begin{pro}
\label{dualcoPROPpropérade} Soit $\mathcal{C}$ un coPROP
diff\'erentiel augment\'e gradu\'e par un poids et connexe. On
rappelle que $U_c(\mathcal{C})$ repr\'esente la coprop\'erade
associ\'ee \`a $\mathcal{C}$. Le dual $\mathcal{C}^{\ac}$ de
$\mathcal{C}$ est isomorphe en tant que PROP diff\'erentiel
gradu\`e par un poids \`a $S_\otimes(U_c(\mathcal{C})^{\ac})$ et
il s'agit d'un PROP quotient de $S_\otimes(\F(\Sigma^{-1}
\mathcal{C}^{(1)}))$.
\end{pro}

\begin{deo}
La d\'emonstration repose encore sur l'isomorphisme de PROPs diff\'erentiels
gradu\'es par un poids
$$\bar{\B}^c(\mathcal{C})=S_\otimes(\bar{\B}^c(U_c(\mathcal{C}))).$$
$\cqfd$
\end{deo}

\section{Prop\'erade quadratique }

Nous montrons ici que le construction duale est une construction
quadratique.\\

Nous avons vu au chapitre $2$ qu'une prop\'erade quadratique
\'etait une prop\'erade de la forme $\F(V)/(R)$ o\`u $R$
appartenait \`a $\F_{(2)}(V)$. Ainsi, l'id\'eal $(R)$ est
homog\`ene pour la graduation de $\F(V)$ en fonction du nombre de
sommets. Une prop\'erade quadratique est donc gradu\'ee par un
poids (\emph{cf.} chapitre $2$ section $3$).

R\'eciproquement, toute prop\'erade quadratique est enti\`erement
d\'etermin\'ee par sa graduation $\Po=\bigoplus_{\rho \in
\mathbb{N}} \Po_{(\rho)}$, le dg-$\Sy$-bimodule $\Po_{(1)}$ et le
dg-$\Sy$-bimodule $\Po_{(2)}$ car on retrouve $R$ via la formule
$R=\textrm{ker}\big(\F_{(2)}(\Po_{(1)}) \xrightarrow{\mu}
\Po_{(2)} \big)$.

\begin{lem}
\label{dualequadratique} Soit $\mathcal{C}$ une coprop\'erade
diff\'erentielle gradu\'ee par un poids connexe. Alors
$\mathcal{C}^{\ac}$ est une prop\'erade quadratique d\'etermin\'ee
par les relations
$$ {\mathcal{C}^{\ac}}_{(1)}=\Sigma^{-1} \mathcal{C}^{(1)}
\quad \textrm{et} \quad  {\mathcal{C}^{\ac}}_{(2)}=
\coker \big(\theta' \, :\, \Sigma^{-1}\mathcal{C}^{(2)}
\to \F_{(2)}(\Sigma^{-1} \mathcal{C}^{(1)})\big).$$
\end{lem}

\begin{deo}
Posons $R=\textrm{ker}\big(\F_{(2)}(\mathcal{C}^{\ac}_{(1)}) \to
\mathcal{C}_{(2)}^{\ac}\big)=\textrm{im}(\theta'_{\Sigma^{-1}\mathcal{C}^{(2)}})$.
Par d\'efinition, $\mathcal{C}^{\ac}_{(\rho)}$ est le quotient de
$\F_{(\rho)}(\Sigma^{-1} \mathcal{C}^{(1)})$ par l'image de
$d_{\theta'}$ sur $\bar{\B}^c_{(\rho-1)}(\mathcal{C}^{(\rho)})$.
Or, cette image correspond aux graphes \`a $\rho$ sommets dont au
moins un couple de sommets adjacents est l'image d'un \'el\'ement
de $\Sigma^{-1} \mathcal{C}^{(2)}$ via $\theta'$. Ceci correspond
bien \`a la partie de degr\'e $\rho$ de l'id\'eal libre engendr\'e
par $R$. Il en r\'esulte $\mathcal{C}^{\ac}=\F(\Sigma^{-1}
\mathcal{C}^{(1)})/(R)$. $\cqfd$
\end{deo}

\begin{cor}
Soit $\mathcal{C}$ un coPROP diff\'erentiel gradu\'e par un poids et connexe.
Alors son PROP dual $\mathcal{C}^{\ac}$ est un PROP quadratique d\'etermin\'e
par les m\^emes relations que pr\'ec\'edemment.
\end{cor}

\begin{deo}
De la proposition~\ref{dualcoPROPpropérade}, on a
$\mathcal{C}^{\ac}= S_\otimes (U_c(\mathcal{C})^{\ac})$. Par le
lemme pr\'ec\'edent, on sait que $U_c(\mathcal{C})^{\ac}$ est une
prop\'erade quadratique. Et, la
proposition~\ref{Lambdapropéradequotient} permet de conclure que
$S_\otimes(U_c(\mathcal{C}^{\ac}))$ est un PROP quadratique.
$\cqfd$
\end{deo}

\section{Prop\'erades et r\'esolution de Koszul}

On donne ici les d\'efinitions de \emph{prop\'erade} et de
\emph{coprop\'erade de Koszul} (respectivement de \emph{PROP} et
de \emph{coPROP de Koszul}). Lorsqu'une prop\'erade $\Po$, de
diff\'erentielle nulle, est de Koszul, alors $\Po$ est
quadratique, sa duale est encore de Koszul et sa biduale est
isomorphe \`a $\Po$. En outre, on montre qu'une prop\'erade $\Po$
est de Koszul si et seulement si la cobar construction r\'eduite
sur la duale $\Po^{\ac}$ est une r\'esolution de $\Po$.

\subsection{D\'efinitions}

\begin{dei}[Prop\'erade de Koszul]
\index{prop\'erade de Koszul}

Soit $\Po$ une prop\'erade diff\'erentielle gradu\'ee par un poids
connexe. On dit que $\Po$ est une \emph{prop\'erade de Koszul} si
l'inclusion $\Po^{\ac}\hookrightarrow \bar{\B}(\Po)$ est un
quasi-isomorphisme.
\end{dei}

De la m\^eme mani\`ere, on a la d\'efinition de coprop\'erade de
Koszul.

\begin{dei}[Coprop\'erade de Koszul]
\index{coprop\'erade de Koszul}

Soit $\mathcal{C}$ une coprop\'erade diff\'erentielle gradu\'ee
par un poids connexe. On dit que $\mathcal{C}$ est une
\emph{coprop\'erade de Koszul} si la projection
$\bar{\B}^c(\mathcal{C}) \twoheadrightarrow \mathcal{C}^{\ac}$ est
un quasi-isomorphisme.
\end{dei}

\begin{pro}
Si $\Po$ est une prop\'erade gradu\'ee par un poids connexe de
Koszul, alors sa duale $\Po^{\ac}$ est une coprop\'erade de Koszul
et ${\Po^{\ac}}^{\ac}=\Po$.
\end{pro}

\begin{deo}
La prop\'erade $\Po$ est  concentr\'ee en degr\'e $0$
($\delta_\Po=0$ et $\Po_0=\Po$). Dans ce cas, $\Po^{\ac}_{(\rho)}$
est un $\Sy$-bimodule homog\`ene de degr\'e homologique $\rho$ et
les \'el\'ements de degr\'e homologique nul de
$\bar{\B}^c(\Po^{\ac})^{(\rho)}$ correspondent aux \'el\'ements de
$\bar{\B}^c_{(\rho)} ({\Po^{\ac}})^{(\rho)}$. Ceci montre que
$$H_0\big(\bar{\B}^c_{(\rho)}(\Po^{\ac})^{(\rho)} \big)=H_\rho\big(
\bar{\B}^c_*(\Po^{\ac})^{(\rho)},\,
d_{\theta'} \big)={\Po^{\ac}}_{(\rho)}^{\ac}.$$

La prop\'erade $\Po$ est de Koszul, c'est-\`a-dire que l'inclusion
$\Po^{\ac}\hookrightarrow \bar{\B}(\Po)$ est un
quasi-isomorphisme. En utilisant le lemme de comparaison des
prop\'erades quasi-libres
(th\'eor\`eme~\ref{lemmecomparaisonproperades}), on montre que le
morphisme induit $\bar{\B}^c(\Po^{\ac})\to
\bar{\B}^c(\bar{\B}(\Po))$ est lui aussi un quasi-isomorphisme.
(On travaille dans le cadre des $\Sy$-bimodules gradu\'es par un
poids connexes et toutes les hypoth\`eses du th\'eor\`eme sont
v\'erifi\'ees). Comme la construction bar-cobar est une
r\'esolution de $\Po$ (th\'eor\`eme~\ref{barcobarresolution}), la
cobar construction $\bar{\B}^c(\Po^{\ac})$ est quasi-isomorphe \`a
$\Po$. Et comme le $\Sy$-bimodule $\Po$ est homog\`ene de degr\'e
$0$, on a
$$ H_*\big(\bar{\B}^c(\Po^{\ac})^{(\rho)}\big) =
\left\{ \begin{array}{ll}
\Po^{(\rho)} & \textrm{si} \quad *=0,\\
0 & \textrm{sinon}.
\end{array} \right.$$

Ces deux r\'esultats sur l'homologie de $\bar{\B}^c(\Po^{\ac})$
mis bout \`a bout montrent que $\Po^{\ac}$ est de Koszul et que
${\Po^{\ac}}^{\ac}=\Po$.$\cqfd$
\end{deo}

\begin{cor}
Soit $\Po$ une prop\'erade gradu\'ee par un poids connexe. Si
$\Po$ est de Koszul alors $\Po$ est n\'ecessairement quadratique.
\end{cor}

\begin{deo}
Si $\Po$ est de Koszul, par la proposition pr\'ec\'edente, on sait
que $\Po={\Po^{\ac}}^{\ac}$. Et le lemme~\ref{dualequadratique}
montre que ${\Po^{\ac}}^{\ac}$ est quadratique.$\cqfd$
\end{deo}

On peut donner les m\^emes d\'efintions dans le cas des PROPs.

\begin{dei}[PROP de Koszul]
\index{PROP de Koszul}

Soit $\Po$ un PROP diff\'erentiel gradu\'e par un poids et
connexe.
On dit que $\Po$ est un \emph{PROP de Koszul} si l'inclusion
$\Po^{\ac}\hookrightarrow \bar{\B}(\Po)$ est un quasi-isomorphisme.
\end{dei}

\begin{dei}[CoPROP de Koszul]
\index{coPROP de Koszul}

Soit $\mathcal{C}$ un coPROP diff\'erentiel
gradu\'e par un poids et connexe.
On dit que $\mathcal{C}$ est un \emph{coPROP de Koszul} si
la projection $\bar{\B}^c(\mathcal{C}) \twoheadrightarrow
\mathcal{C}^{\ac}$ est un quasi-isomorphisme.
\end{dei}

\begin{pro}
\label{PROPpropéradeKoszul}
Soit $\Po$ un PROP diff\'erentiel augment\'e gradu\'e par un poids et connexe.
Le PROP $\Po$ est de Kosuzl si et seulement si la prop\'erade $U_c(\Po)$ est
de Koszul.
\end{pro}

\begin{deo}
$\ $

$(\Rightarrow)$ Si $\Po$ est un PROP de Koszul, cela signifie que
le morphisme $\Po^{\ac}\to \bar{\B}(\Po)$ est un
quasi-isomorphisme. La proposition~\ref{dualPROPpropérade} montre
que $\Po^{\ac}=S_\otimes (U_c(\Po)^{\ac})$ et la
proposition~\ref{Lambdabar} donne l'isomorphisme
$\bar{\B}(\Po)=S_\otimes(\bar{\B}(U_c(\Po))$. De ces propositions,
on tire que le morphisme $S_\otimes(U_c(\Po)^{\ac})\to S_\otimes
(\bar{\B}(U_c(\Po)))$ est un quasi-isomorphisme. En outre, nous
avons vu que la diff\'erentielle sur la bar construction
respectait les graphes connexes et que ce morphisme correspondait
\`a l'image du morphisme $U_c(\Po)^{\ac} \to \bar{\B}(U_c(\Po)$
via le foncteur $S_\otimes$. Cela donne que ce quasi-morphisme se
d\'ecompose en une somme directe de quasi-isomorphismes $
(S_\otimes)_{(n)}(U_c(\Po)^{\ac})\to (S_\otimes)_{(n)}
(\bar{\B}(U_c(\Po)))$. En particulier, on a pour $n=1$ que
$U_c(\Po)^{\ac} \to \bar{\B}(U_c(\Po))$ est un quasi-isomorphisme,
c'est-\`a-dire que $U_c(\Po)$ est une prop\'erade de Koszul.\\

$(\Leftarrow)$ R\'eciproquement, si $U_c(\Po)$ est une prop\'erade
de Koszul, comme le foncteur $S_\otimes$ est un foncteur exacte,
on a que $\Po$ est un PROP de Koszul. $\cqfd$
\end{deo}

On a le m\^eme r\'esultat au niveau des coPROPs et des coprop\'erades.

\subsection{R\'esolution de Koszul et mod\`ele minimal}

Comme nous l'avons vu pr\'ec\'edemment, l'inclusion $\Po^{\ac}
\hookrightarrow \bar{\B}(\Po)$ induit un morphisme
$\bar{\B}^c(\Po^{\ac}) \to \bar{\B}^c(\bar{\B}(\Po))$. Et en
composant ce morphisme avec le morphisme
$\zeta\, :\,  \bar{\B}^c(\bar{\B}(\Po))
\to \Po$ de la r\'esolution bar-cobar
(\emph{cf.} th\'eor\`eme~\ref{barcobarresolution}), on obtient
un morphisme de prop\'erades diff\'erentielles gradu\'ees par un poids
de la forme
$$ \bar{\B}^c(\Po^{\ac}) \to \Po.$$

\begin{thm}[R\'esolution de Koszul]
\index{r\'esolution de Koszul}

\label{resolutionKoszul}
Soit $\Po$ une prop\'erade diff\'erentielle gradu\'ee par un poids et
connexe.
La prop\'erade $\Po$ est de Koszul si et seulement si le morphisme
de prop\'erades diff\'erentielles gradu\'ees par un poids
$ \bar{\B}^c(\Po^{\ac}) \to \Po$ est un quasi-isomorphisme.
\end{thm}

\begin{deo}
Nous sommes dans le cadre gradu\'e par un poids, et toutes les
prop\'erades en jeu sont connexes. On peut donc appliquer le lemme
de comparaison des prop\'erades quasi-libres.\\

$(\Rightarrow)$ Si $\Po$ est
de Koszul, cela signifie que $\Po^{\ac}\hookrightarrow \bar{\B}(\Po)$
est un quasi-isomorphisme. Donc $ \bar{\B}^c(\Po^{\ac}) \to
 \bar{\B}^c(\bar{\B}(\Po))$ est aussi un quasi-isomorphisme. Il en
va de m\^eme pour $ \bar{\B}^c(\Po^{\ac}) \to \Po$, gr\^ace \`a la
r\'esolution bar-cobar (th\'eor\`eme~\ref{barcobarresolution}).\\

$(\Leftarrow)$ R\'eciproquement, supposons que le morphisme $\bar{\B}^c
(\Po^{\ac}) \to \Po $ soit un quasi-isomorphisme. Comme la
construction bar-cobar est un quasi-isomorphisme, le morphisme
$ \bar{\B}^c(\Po^{\ac}) \to
 \bar{\B}^c(\bar{\B}(\Po))$  est aussi un quasi-isomorphime. Et, on
conclut en utilisant le lemme de comparaison des prop\'erades
quasi-libres. $\cqfd$
\end{deo}

En reconduisant les m\^emes arguments dans le cas des PROPs, on montre
le th\'eor\`eme analogue.

\begin{thm}[R\'esolution de Koszul d'un PROP]
Soit $\Po$ un PROP diff\'erentiel gradu\'e par un poids et
connexe.
Le PROP $\Po$ est de Koszul si et seulement si le morphisme
de PROPs diff\'erentiels gradu\'es par un poids
$ \bar{\B}^c(\Po^{\ac}) \to \Po$ est un quasi-isomorphisme.
\end{thm}

\begin{deo}
La d\'emonstration est exactement la m\^eme. On utilise la
r\'esolution bar-cobar donn\'ee dans le cadre des PROPs par la
proposition~\ref{barcobarresolutionPROP} et le lemme de
comparaison des PROPs quasi-libres. $\cqfd$
\end{deo}

Ces th\'eor\`emes l\'egitiment la d\'efinition suivante :

\begin{dei}[R\'esolution de Koszul]
Lorsque $\Po$ est de Koszul (prop\'erade ou PROP), la r\'esolution
$\bar{\B}^c(\Po^{\ac})
\to \Po$
est appell\'ee \emph{r\'esolution de Koszul}.
\end{dei}

\begin{cor}
Lorsque $\Po$ est une prop\'erade ou un PROP de Koszul de
diff\'erentielle nulle, alors $\bar{\B}^c(\Po^{\ac})$ est le
\emph{mod\`ele minimal} de $\Po$.\index{mod\`ele minimal}
\end{cor}

\begin{deo}
Nous sommes en pr\'esence d'une r\'esolution quasi-libre
$\bar{\B}^c(\Po^{\ac})=\F^c(\Sigma^{-1} \oPo^{\ac}) \to \Po$ dont
la diff\'erentielle $d_{\theta'}$ est quadratique puisque, par
d\'efinition de la d\'erivation $d_{\theta'}$, on a
$d_{\theta'}(\Sigma^{-1} \oPo^{\ac})\subset \F_{(2)}(\Sigma^{-1}
\oPo^{\ac})$. $\cqfd$
\end{deo}

\begin{dei}[$\Po$-g\`ebre \`a homotopie pr\`es]
\index{$\Po$-g\`ebre \`a homotopie pr\`es}

Lorsque $\Po$ est une prop\'erade ou un PROP de Koszul, on appelle
\emph{$\Po$-g\`ebre \`a homotopie pr\`es} toute g\`ebre sur
$\bar{\B}^c(\Po^{\ac})$. On parle aussi de
\emph{$\Po_\infty$-g\`ebre}.
\end{dei}

Cette notion g\'en\'eralise celle d'alg\`ebre \`a homotopie pr\`es
(sur une op\'erade).

\section{Complexe de Koszul}

Il est \'equivalent et aussi difficile de montrer qu'une
prop\'erade (ou qu'un PROP) $\Po$ est de Koszul ($\Po^{\ac}\to
\bar{\B}(\Po)$ quasi-isomorphisme), \`a partir de la d\'efinition,
que de montrer que la cobar construction sur $\Po^{\ac}$ est une
r\'esolution de $\Po$ ($\bar{\B}^c(\Po^{\ac} \to \Po$
quasi-isomorphisme). On introduit donc un petit complexe dont
l'acyclicit\'e est un crit\`ere pour d\'eterminer si $\Po$ est de
Koszul et donc pour avoir le mod\`ele minimal sur $\Po$.

\subsection{Complexe de Koszul \`a coefficients}

De la m\^eme mani\`ere que l'on avait d\'efini la bar construction
\`a coefficients (\emph{cf.} chapitre $4$ section $2.2$) et la bar
construction simpliciale \`a coefficients (\emph{cf.} chapitre $6$
section $1$), on peut d\'efinir le \emph{complexe de Koszul \`a
coefficients}.

\begin{dei}[Complexe de Koszul \`a coefficients]
\index{complexe de Koszul \`a coefficients}

Soit $\Po$ une prop\'erade diff\'erentielle gradu\'ee par un poids
connexe et soient $L$ et $R$ deux modules (\`a droite et \`a
gauche) sur $\Po$. On appelle \emph{complexe de Koszul \`a
coefficients dans les modules $L$ et $R$}, le complexe d\'efini
sur le $\Sy$-bimodule $L\boxtimes_c\Po^{\ac}\boxtimes_c R$ par la
diff\'erentielle $d$, somme des trois termes suivants :
\vspace{6pt}

\begin{enumerate}
\item la diff\'erentielle canonique $\delta_\Po$ induite par celle de
$\Po$,

\item le morphisme homog\`ene $d_{\theta_L}$ de degr\'e $-1$ qui
vient de la structure de $\Po$-comodule \`a gauche de $\Po^{\ac}$ :
$$ \theta_l\ : \ \Po^{\ac} \xrightarrow{\Delta}  \Po^{\ac}
\boxtimes_c \Po^{\ac} \twoheadrightarrow (I \oplus
\underbrace{\oPo^{\ac}}_1)  \boxtimes_c  \Po^{\ac} \to (I\oplus
\underbrace{\Po^{(1)}}_1)\boxtimes_c \Po^{\ac},$$

\item le morphisme homog\`ene $d_{\theta_R}$ de degr\'e $-1$ qui
vient de la structure de $\Po$-comodule \`a droite de $\Po^{\ac}$ :
$$ \theta_r\ : \ \Po^{\ac} \xrightarrow{\Delta} \Po^{\ac}
\boxtimes_c \Po^{\ac} \twoheadrightarrow \Po^{\ac} \boxtimes_c (I
\oplus \underbrace{\oPo^{\ac}}_1) \to \Po^{\ac}\boxtimes_c
(I\oplus \underbrace{\Po^{(1)}}_1).$$
\end{enumerate}

On le note $\mathcal{K}(L,\, \Po,\, R)$.
\end{dei}

Tout comme $\Po^{\ac}$ s'injecte dans la bar construction
$\bar{\B}(\Po)$, le complexe de Koszul \`a co\'efficients est un
sous-complexe de la bar construction \`a coefficients.

\begin{pro}
Soit $\Po$ une prop\'erade diff\'erentielle gradu\'ee
par un poids connexe et soient $L$ et $R$ deux modules (\`a droite et \`a gauche)
sur $\Po$.
Le $\Sy$-bimodule diff\'erentiel $\mathcal{K}(L,\, \Po, \, R)=L\boxtimes_c
\Po^{\ac} \boxtimes_c R$ est un
sous-complexe de la bar construction \`a coefficients
$\B(L,\, \Po,\, R)=L\boxtimes_c \bar{\B}(\Po)_c \boxtimes R$.
\end{pro}

\begin{deo}
Tout repose sur le fait que la diff\'erentielle de la bar
construction est d\'efinie \`a partir du coproduit $\Delta$ sur
$\bar{\B}(\Po)=\F^c(\Sigma \oPo)$ et que la diff\'erentielle du
complexe de Koszul est aussi d\'efinie \` a partir du coproduit
sur $\Po^{\ac}$. Comme  $\Po^{\ac}$ est une sous-coprop\'erade
diff\'erentielle de $\bar{\B}(\Po)$ les diff\'erentielles
correspondent. $\cqfd$
\end{deo}

L'inclusion $\Po^{\ac}\hookrightarrow \bar{\B}(\Po)$
induit un morphisme de dg-$\Sy$-bimodules
$$L\boxtimes\Po^{\ac}\boxtimes R \hookrightarrow \B(L,\, \Po,\, R).$$

De la m\^eme mani\`ere, on d\'efinit le complexe de Koszul \`a coefficients
pour un PROP.

\begin{dei}[Complexe de Koszul \`a coefficients d'un PROP]
Soit $\Po$ un PROP diff\'erentiel gradu\'e par un poids
connexe et soient $L$ et $R$ deux modules (\`a droite et \`a gauche)
sur $\Po$.
On appelle \emph{complexe de Koszul
\`a coefficients dans les modules $L$ et $R$}, le complexe d\'efini
sur le $\Sy$-bimodule $L\boxtimes\Po^{\ac}\boxtimes R$
par la diff\'erentielle $d$ somme des trois m\^emes termes que dans le cas
des prop\'erades.
\end{dei}

\begin{pro}
Le complexe de Koszul \`a coefficents $\mathcal{K}(L,\, \Po,\, R)=L \boxtimes
\Po^{\ac} \boxtimes R$ sur un PROP
$\Po$ est un sous-complexe de
la bar construction \`a coefficients $\B(L,\,\Po,\, R)=L \boxtimes
\bar{\B}(\Po) \boxtimes R$.
\end{pro}

\subsection{Complexe de Koszul et mod\`ele minimal}

Tout comme nous avions \'etudi\'e une bar construction
particuli\`ere, la bar construction augment\'ee
$\B(I,\, \Po,\, \Po)=\bar{\B}(\Po)\boxtimes \Po$ et la bar
construction normalis\'ee augment\'ee $\N(I,\, \Po,\,\Po)=\bar{\N}
(\Po) \boxtimes \Po$,
on consid\`ere
ici le complexe \'equivalent au niveau des complexes de Koszul.

\begin{dei}[Complexe de Koszul]
\index{complexe de Koszul}

On appelle \emph{complexe de Koszul} le complexe
$\mathcal{K}(I,\, \Po,\, \Po)=\Po^{\ac}\boxtimes_c \Po$
(et $\Po^{\ac}\boxtimes
\Po$ dans le cas des PROPs).
\end{dei}

La diff\'erentielle de ce complexe est d\'efinie par le morphisme
$d_{\theta_r}$ pr\'ec\'edent plus \'eventuellement la
diff\'erentielle canonique $\delta_\Po$ induite par celle de
$\Po$. Remarquons que le morphisme $d_{\theta_r}$ revient \`a
extraire une op\'eration $\nu\in\Po^{\ac}_{(1)}$ de $\Po^{\ac}$
par le haut, \`a identifier cette op\'eration $\nu$ comme une
op\'eration de $\Po_{(1)}$ (puisque $\Po^{\ac}_{(1)}=\Po_{(1)}$)
et finalement \`a la composer dans $\Po$. Ce qui se r\'esume par
les diagrammes donn\'es par J.-L. Loday dans \cite{Loday1}.\\

Le principal th\'eor\`eme de cette th\`ese est le crit\`ere
suivant.

\begin{thm}[Crit\`ere de Koszul]
\index{crit\`ere de Koszul}

\label{criteredeKoszul} Soit $\Po$ une prop\'erade
diff\'erentielle gradu\'ee par un poids et connexe. Les
propositions suivantes sont \'equivalentes \vspace{6pt}

\begin{tabular}{ll}
$(1)$ & $\Po$ est de Koszul (l'inclusion $\Po^{\ac}\hookrightarrow \bar{\B}(\Po)$
est un quasi-isomorphisme)\\

$(2)$  &Le complexe de Koszul $\Po^{\ac}\boxtimes_c \Po$ est acyclique.\\

$(2')$ &Le complexe de Koszul $\Po \boxtimes_c \Po^{\ac}$ est acyclique.\\

$(3)$ &Le morphisme de prop\'erades diff\'erentielles gradu\'ees par un poids \\
& $\bar{\B}^c(\Po^{\ac})\to \Po $ est un quasi-isomorphisme.
\end{tabular}
\end{thm}

\begin{deo}
Nous avons d\'ej\`a vu l'\'equivalence $(1)\iff(3)$ (\emph{cf.}
th\'eor\`eme~\ref{resolutionKoszul}).

Par le lemme de comparaison des $\Po$-modules quasi-libre (\`a droite),
l'assertion $(1)$ est \'equivalente au fait que le morphisme
$\Po^{\ac}\boxtimes_c \Po \to \bar{\B}(\Po)\boxtimes_c \Po$
soit un quasi-isomorphisme. L'acyclicit\'e de la bar construction
augment\'ee (\emph{cf.} th\'eor\`eme~\ref{acyclicitébaraugmentée})
permet de voir que la
prop\'erade $\Po$ est de Koszul $(1)$ si et seulement si
le complexe de Koszul $\Po^{\ac}\boxtimes_c \Po$ est acyclique $(2)$.

On proc\`ede de la m\^eme mani\`ere, avec le lemme de comparaison
des $\Po$-modules quasi-libres \`a gauche, pour montrer
l'\'equivalence $(1)\iff(2')$. $\cqfd$
\end{deo}

\begin{thm}[Crit\`ere de Koszul pour les  PROPs]
Soit $\Po$ un PROP diff\'erentiel gradu\'e par un poids et
connexe. Les propositions suivantes sont \'equivalentes
\vspace{6pt}

\begin{tabular}{ll}
$(1)$ & $\Po$ est de Koszul (l'inclusion $\Po^{\ac}\hookrightarrow
\bar{\B}(\Po)$
est un quasi-isomorphisme)\\

$(2)$  &Le complexe de Koszul $\Po^{\ac}\boxtimes \Po$ est acyclique.\\

$(2')$ &Le complexe de Koszul $\Po \boxtimes \Po^{\ac}$ est acyclique.\\

$(3)$ &Le morphisme de PROPs diff\'erentiels gradu\'es par un poids\\
& $\bar{\B}^c(\Po^{\ac})\to \Po $ est un quasi-isomorphisme.
\end{tabular}
\end{thm}

\begin{deo}
La d\'emonstration est exactement la m\^eme. On utilise ici
l'acyclicit\'e de la bar construction augment\'ee sur un PROP et
les versions PROPiques des lemmes de comparaison. $\cqfd$
\end{deo}

\begin{pro}
On a un isomorphisme de $\Sy$-bimodules diff\'erentiels gradu\'es par un poids
$$\Po^{\ac}\boxtimes \Po= S_\otimes(U_c(\Po)^{\ac}\boxtimes_c \Po).$$
\end{pro}

Cette proposition ainsi que le
proposition~\ref{PROPpropéradeKoszul} justifient que lorsque l'on
veut montrer que la cobar construction sur le coPROP dual fournit
le mod\`ele minimal d'un PROP $\Po$, il suffit de prouver
l'acyclicit\'e du complexe de Koszul associ\'e \`a la prop\'erade
d\'efinie par $\Po$, \`a savoir $\Po^{\ac}\boxtimes_c \Po$. La
notion de prop\'erade fournit le bon cadre d'\'etude pour la
dualit\'e de Kosuzl des PROPs. Toute la th\'eorie d\'evelopp\'ee
ici montre que l'information essentielle d'un PROP quadratique \`a
relations connexes est pr\'esente dans la prop\'erade associ\'ee
et que c'est plut\^ot sur elle qu'il faut travailler en pratique.

\begin{conj}
Soit $\Po$ une prop\'erade de Koszul, la composition suivante est un
quasi-isomorphisme
$$\Po^{\ac}\to \bar{\B}(\Po) \to \bar{\N}(\Po).$$
\end{conj}

Lorsqu'une prop\'erade $\Po$ est de Koszul, ce corollaire devrait
permettre de calculer l'homologie de sa bar construction
normalis\'ee r\'eduite. Elle doit correspondre \`a la duale de
Koszul $\Po^{\ac}$. Cette propri\'et\'e a permis \`a B. Fresse
\cite{Fresse} de montrer que l'homologie du poset des partitions
\'etait donn\'ee par l'op\'erade $\Li$ en interpr\'etant les
modules simpliciaux engendr\'es par le poset des partitions comme
la bar construction normalis\'ee r\'eduite de l'op\'erade $\C$. En
utilisant cette m\'ethode, nous avons calcul\'e , dans \cite{BV},
l'homologie d'autres types de posets en les reliant aux op\'erades
$\mathcal{A}s$, $\mathcal{P}erm$ et $\mathcal{D}ias$. Ainsi ces
homologies sont donn\'ees par les op\'erades duales, \`a savoir
$\mathcal{A}s$, $\mathcal{P}relie$ et $\mathcal{D}end$.

\subsection{Lien avec les th\'eories classiques (alg\`ebres et op\'erades)}

La  notion de duale de Kozsul, telle que nous l'avons d\'efinie
ici, est en fait une coprop\'erade ou un coPROP (\emph{cf.}
\ref{constructionduale}), alors que la duale de Koszul d'une
alg\`ebre est une alg\`ebre et que la duale de Koszul d'une
op\'erade est une op\'erade (\emph{cf.} \cite{Priddy} et
\cite{GK}). Pour retrouver ces constructions classiques, il suffit
de consid\'erer la duale lin\'eaire, duale de Czech, de
$\Po^{\ac}$.

\begin{dei}[Dual de Czech d'un $\Sy$-bimodule]
\index{dual de Czech d'un $\Sy$-bimodule}

Soit $\Po$ un $\Sy$-bimodule, on d\'efinit le dual de Czech
$\Po^\vee$ par le $\Sy$-bimodule $\Po^\vee = \bigoplus_{\rho,\,
m,\, n} \Po_{(\rho)}^\vee(m,\, n)$, o\`u
$$\Po_{(\rho)}^\vee (m,\,n)=sgn_{\Sy_m}\otimes_k \Po_{(\rho)}(m,\, n)^* \otimes_k
sgn_{\Sy_n}.$$
\end{dei}

Le dual de Czech revient \`a consid\'erer le dual lin\'eaire tordu
par les repr\'esentations signatures.

\begin{lem}
Soit $(\Co,\, \Delta,\, \varepsilon)$ une coprop\'erade
(repectivement un coPROP) gradu\'ee par un poids telle que les
modules $\Co_{(\rho)}(m,\,n)$ soient de dimension finie sur $k$,
pour tout $m$, $n$ et $\rho$. Alors, le $\Sy$-bimodule $\Co^\vee$
est naturellement muni d'un structure de prop\'erade
(respectivement PROP) gradu\'ee par un poids.
\end{lem}

\begin{deo}
La comultiplication $\Delta$ se d\'ecompose avec le poids en
$\Delta=\bigoplus_{\rho \in \mathbb{N}} \Delta_{(\rho)}$. Pour
d\'efinir la multiplication $\mu$ sur $\Co^\vee$, on dualise
lin\'eairement chaque
$$\Delta_{(\rho)}(m,\,n) \, :\,
\Co_{(\rho)}(m,\,n)\to (\Co\boxtimes_c\Co)_{(\rho)}(m,\,n).$$

Plus pr\'ecisement, on a :
\begin{eqnarray*}
(\Po\boxtimes_c \Po)_{(\rho)}(m,\, n)&=&\bigoplus_{g\in
\mathcal{G_2^c}(m,\, n)} \left( \bigotimes_{\nu \, \textrm{sommet
de} \, g \atop \sum \rho_\nu =\rho}
\Co_{(\rho_\nu)}^\vee(|Out(\nu)|,\, |In(\nu)|)
\right)\Bigg/_\approx \\
&=&\bigoplus_{g\in \mathcal{G}_2^c(m,\, n)} \left( \left(
\bigotimes_{\nu \, \textrm{sommet de} \, g \atop \sum \rho_\nu
=\rho} \Co_{(\rho_\nu)}(|Out(\nu)|,\, |In(\nu)|)
\right)\Bigg/_\approx\right)^\vee,
\end{eqnarray*}
en utilisant l'hypoth\`ese sur la dimension des $\Co_{(\rho)}(m,\,
n)$ et en identifiant invariants et coinvariants (nous
travaillons sur un corps de caract\'eristique nulle). On d\'efinit
alors $\mu_{(\rho)}$ par la composition :
\begin{eqnarray*}
(\Po\boxtimes_c \Po)_{(\rho)}(m,\, n) &=&\bigoplus_{g\in
\mathcal{G}^2(m,\, n)} \left( \left( \bigotimes_{\nu \,
\textrm{sommet de} \, g \atop \sum \rho_\nu =\rho}
\Co_{(\rho_\nu)}(|Out(\nu)|,\, |In(\nu)|)
\right)\Bigg/_\approx\right)^\vee \\ &\longrightarrow & \left(
\bigoplus_{g\in \mathcal{G}^2(m,\, n)}  \left( \bigotimes_{\nu \,
\textrm{sommet de} \, g \atop \sum \rho_\nu =\rho}
\Co_{(\rho_\nu)}(|Out(\nu)|,\, |In(\nu)|)
\right)\Bigg/_\approx\right)^\vee\\
&=&(\Co\boxtimes_c\Co)_{(\rho)}^\vee(m,\,n) \xrightarrow{ ^t\Delta_{(\rho)}}
\Co^\vee_{(\rho)}(m,\, n)=\Po_{(\rho)}(m,\, n).
\end{eqnarray*}

La coassociativit\'e de la comultiplication $\Delta$ induit l'associativit\'e
de la mutliplication $\mu$. Et la counit\'e de $\Co$
$\varepsilon \, :\, \Co \to I$
donne, par passage au dual, l'unit\'e de $\Po$ :
$\varepsilon \, : \, I \to \Po$.

Dans le cas PROPique, on dualise, de la m\^eme mani\`ere, la
d\'econcat\'eantion horizontale. $\cqfd$
\end{deo}

\begin{pro}
Soit $\Po$ une prop\'erade (respectivement un PROP) gradu\'ee par un poids
(par exemple quadratique). Posons $V=\Po^{(1)}$, le $\Sy$-bimodule engendr\'e
par les \'el\'ements de poids $1$ de $\Po$.

S'il existe deux entiers $M$ et $N$ tels que $V(m,\, n)=0$ lorsque $m>M$ ou
$n>N$ et si les modules $V(m,\, n)$ sont tous de dimension finie sur
$k$, alors la cobar construction $\bar{\B}^c(\Po)$ sur $\Po$ et la duale de
Koszul $\Po^{\ac}$ v\'erifient les hypoth\`eses du lemme pr\'ec\'edent.
\end{pro}

\begin{deo}
Comme $\Po^{\ac}$ est une sous-coprop\'erade (sous-coPROP) de la cobar
construction $\bar{\B}^c(\Po)$ sur $\Po$, il suffit de d\'emontrer
que les modules $\bar{\B}^c(\Po)_{(\rho)}(m,\,n)$  sont de dimension finie
sur $k$.

Nous avons vu que
$\bar{\B}^c(\Po)_{(\rho)}(m,\,n)=\F^c_{(\rho)}(V)(m,\, n)$. Dans
le cas o\` u $V$ v\'erifient les hypoth\`eses de la proposition,
ce dernier module est donn\'e par une somme sur l'ensemble des
graphes \`a $\rho$ sommets et tels que chaque sommet poss\`ede au
plus $N$ entr\'ees et $M$ sorties. Cet ensemble \'etant fini et
les modules $V(m,\, n)$ \'etant de dimension finie, on a le
r\'esultat escompt\'e. $\cqfd$
\end{deo}

\begin{cor}
Pour toute prop\'erade (respectivement tout PROP) quadratique engendr\'ee
par un $\Sy$-bimodule
$V$ dont la somme des dimensions $\sum_{m,\,n} \dim_k V(m,\,n)$ est finie,
le dual de Czech ${\Po^{\ac}}^\vee$ de la duale de Koszul de $\Po$ est
muni d'une strucutre naturelle de prop\'erade (respectivement de PROP).

De plus, si $\Po$ est de la forme $\F(V)/(R)$, alors la prop\'erade
${\Po^{\ac}}^\vee$ est quadratique et de la forme  ${\Po^{\ac}}^\vee
=\F(\Sigma V)/(\Sigma^2 R^\perp)$. On note cette prop\'erade (ou le PROP
associ\'e) $\Po^!$.
\end{cor}

On peut remarquer que la prop\'erade (ou le PROP) $\Po^!$ est enti\`erement
d\'etermin\'ee par $\Po^{\ac}_{(1)}=\Sigma V$ et $\Po^{\ac}_{(2)}=\Sigma^2 R$.

Rappelons que dans le cas d'une alg\`ebre quadratique
$A=T(V)/(R)$, S. Priddy d\'efinit l'alg\`ebre duale $A^!$ par
$T(V^*)/(R^\perp)$, lorsque $V$ est de dimension finie sur $k$. De
m\^eme, V. Ginzburg et M. M. Kapranov construisent la duale d'une
op\'erade quadratique $\Po=\F(V)/(R)$ en posant
$\Po^!=\F(V^\vee)/(R^\perp)$, encore une fois lorsque $V$
est de dimension finie.\\

Dans le cas de la dimension finie, les d\'efinitions conceptuelles
des objets duaux donn\'ees ici coincident
avec les d\'efinitions classiques, \`a suspension pr\`es
(\emph{cf.} \cite{BGS} et
\cite{Fresse}). En particulier, la duale d'une alg\`ebre $A^{\ac}_{(n)}$
 est isomorphe \`a $\Sigma^n {A^!}_n^*$ et la duale d'une op\'erade
$\Po^{\ac}_{(n)}$ est isomorphe \`a $\Sigma^n{\Po^!}_n^\vee$. \\

Les complexes de
Koszul et les r\'esolutions venant de la cobar construction introduits
ici correspondent, dans le cas des alg\`ebres, \`a ceux donn\'es
dans \cite{Priddy} et, dans le cas des op\'erades, \`a ceux
de \cite{GK}.

Comme nous n'avons aucune hypoth\`ese sur les dimensions
de modules en jeu, les notions introduites ici sont plus g\'en\'erales.

\section{Exemples}

La derni\`ere difficult\'e est de pouvoir montrer que le complexe
de Koszul d'une prop\'erade (un PROP) est acyclique dans des cas
concrets, comme celui des big\`ebres de Lie et celui des
big\`ebres de Hopf infinit\'esimales par exemple. Pour cela, nous
demontrons une proposition affirmant que les prop\'erades
construites suivant un sch\'ema particulier sont de Koszul. Cette
section est une g\'en\'eralisation aux prop\'erades des m\'ethodes
de T. Fox, M. Markl \cite{FM} et W. L. Gan \cite{Gan}.

\subsection{Loi de remplacement}
\index{loi de remplacement}

Soit $\Po$ une prop\'erade quadratique de le forme
$$\Po=\F(V,\, W)/(R\oplus D \oplus S),$$ o\`u $R\subset \F_{(2)}(V)$, $S\subset
\F_{(2)}(W)$ et
$$D\subset (I\oplus \underbrace{W}_1)\boxtimes_c
(I\oplus \underbrace{V}_1)  \bigoplus
(I\oplus \underbrace{V}_1)\boxtimes_c
(I\oplus \underbrace{W}_1).$$
Les deux couples de $\Sy$-bimodules $(V,\, R )$ et $(W,\, S)$ induisent
des prop\'erades que l'on note $A=\F(V)/(R)$ et $B=\F(W)/(S)$.

\begin{dei}[Loi de remplacement]
Soit $\lambda$ un morphisme de $\Sy$-bimodules
$$\lambda \ : \   (I\oplus \underbrace{W}_1)\boxtimes_c
(I\oplus \underbrace{V}_1)  \to
(I\oplus \underbrace{V}_1)\boxtimes_c
(I\oplus \underbrace{W}_1).$$
Lorsque le $\Sy$-bimodule $D$ est d\'efini comme l'image
de
$$(id,\, -\lambda) :\  (I\oplus \underbrace{W}_1)\boxtimes_c
(I\oplus \underbrace{V}_1) \to  (I\oplus \underbrace{W}_1)\boxtimes_c
(I\oplus \underbrace{V}_1)  \bigoplus
(I\oplus \underbrace{V}_1)\boxtimes_c
(I\oplus \underbrace{W}_1),$$
on dit que $\lambda$ est une \emph{loi de remplacement} et on note
$D$ par $D_\lambda$.
\end{dei}

Les exemples que nous traitons ici sont tous de cette forme. Ils proviennent m\^eme
d'un ``m\'elange'' (\emph{cf.} \cite{FM}) de deux op\'erades.

\begin{dei}[$\Sy$-bimodule oppos\'e]
\index{$\Sy$-bimodule oppos\'e}

A partir d'un $\Sy$-bimodule $\Po$, on d\'efinit un $\Sy$-bimodule \emph{oppos\'e}
$\Po^{op}$ par
$$\Po^{op}(m,\, n)=\Po(n,\, m) .$$
\end{dei}

Prendre l'oppos\'e d'un $\Sy$-bimodule revient \`a inverser le sens de parcours du
$\Sy$-bimodule. Soient $\Po$ et $\Qo$ deux $\Sy$-bimodules, on retrouve cette remarque
au niveau du produit $\Po\boxtimes_c \Qo$, car on a
$$(\Po\boxtimes_c \Qo )^{op}=\Qo^{op} \boxtimes_c \Po^{op}.$$
Lorsque le $\Sy$-bimodule $\Po$ est muni d'une structure de prop\'erade (ou de PROP),
le $\Sy$-bimodule oppos\'e $\Po^{op}$ est lui aussi muni d'une structure de
prop\'erade (ou de PROP).  \\

Les prop\'erades donn\'ees en exemple ici (big\`ebres de Lie,
big\`ebres de Hopf infinit\'esimales) sont engendr\'ees uniquement
par des op\'erations ($n$ entr\'ees et une sortie) et des
coop\'erations (une entr\'ee et $m$ sorties). On peut les \'ecrire
sous la forme $\F(V\oplus W)/(R\oplus D_\lambda \oplus S)$, o\`u
$V$ repr\'esente les op\'erations g\'en\'eratrices ($V(m,\, n)=0$
si $m>1$) et o\`u $W$ repr\'esente les coop\'erations
g\'en\'eratrices ($W(m,\, n)=0$ si $n>1$). La prop\'erade
$A=\F(V)/(R)$ est alors une op\'erade et la prop\'erade $B^{op}=
\F(W^{op})/(S^{op})$ est aussi une op\'erade. Dans ce cas
particulier, les lois de remplacement sont de la forme
$$\lambda \ : \   W\otimes_k V \to (I\oplus \underbrace{V}_1)\boxtimes_c
(I\oplus \underbrace{W}_1).$$

Pour la prop\'erade associ\'ee aux big\`ebres de Lie (\emph{cf.}
chapitre 2 section~\ref{BiLie}), on a $A=B^{op}=\Li$ et la loi de
remplacement $\lambda$ est donn\'ee par
$$\lambda \ : \
\vcenter{\xymatrix@M=0pt@R=6pt@C=4pt{1 & & 2 \\\ar@{-}[dr] &  &\ar@{-}[dl]  \\  &\ar@{-}[d] & \\
& \ar@{-}[dl]\ar@{-}[dr]& \\ & & \\ 1 & & 2}} \mapsto
 \vcenter{\xymatrix@M=0pt@R=6pt@C=2pt{ & 1 & & 2 \\ & \ar@{-}[d] &  &\ar@{-}[d]  \\  & \ar@{-}[dl]
 \ar@{-}[dr] &  & \ar@{-}[dl]\\
\ar@{-}[d]&  & \ar@{-}[d] & \\ & & & \\ 1 & & 2 &}} -
 \vcenter{\xymatrix@M=0pt@R=6pt@C=2pt{ & 2 & & 1 \\ & \ar@{-}[d] &  &\ar@{-}[d]  \\  & \ar@{-}[dl]
 \ar@{-}[dr] &  & \ar@{-}[dl]\\
\ar@{-}[d]&  & \ar@{-}[d] & \\ & & & \\ 1 & & 2 &}} +
\vcenter{\xymatrix@M=0pt@R=6pt@C=2pt{1 & & 2 & \\ \ar@{-}[d] &
&\ar@{-}[d] &  \\ \ar@{-}[dr] &
& \ar@{-}[dl] \ar@{-}[dr]& \\
& \ar@{-}[d]&  & \ar@{-}[d] \\ & & & \\ & 1 & & 2}}-
\vcenter{\xymatrix@M=0pt@R=6pt@C=2pt{2 & & 1 & \\ \ar@{-}[d] &
&\ar@{-}[d] &  \\ \ar@{-}[dr] &
& \ar@{-}[dl] \ar@{-}[dr]& \\
& \ar@{-}[d]&  & \ar@{-}[d] \\ & & & \\ & 1 & & 2}}.$$
Dans le cas
des big\`ebres de Hopf infinit\'esimiales
 (\emph{cf.} chapitre 2 section~\ref{InfBi}), les
op\'erades $A$ et $B^{op}$ correspondent \`a l'op\'erade $\A$ des
alg\`ebres associatives et la loi de remplacement vient de
$$\lambda \ : \
\vcenter{\xymatrix@M=0pt@R=6pt@C=6pt{\ar@{-}[dr] &  &\ar@{-}[dl]  \\  &\ar@{-}[d] & \\
& \ar@{-}[dl]\ar@{-}[dr]& \\ & & }} \mapsto
 \vcenter{\xymatrix@M=0pt@R=6pt@C=6pt{ & \ar@{-}[d] &  &\ar@{-}[d]  \\  & \ar@{-}[dl]
 \ar@{-}[dr] &  & \ar@{-}[dl]\\
\ar@{-}[d]&  & \ar@{-}[d] & \\ & & & }}  +
\vcenter{\xymatrix@M=0pt@R=6pt@C=6pt{ \ar@{-}[d] & &\ar@{-}[d] &
\\ \ar@{-}[dr] &
& \ar@{-}[dl] \ar@{-}[dr]& \\
& \ar@{-}[d]&  & \ar@{-}[d] \\ & & & }} \ .$$

Ces deux lois de remplacement permettent de permuter verticalement
op\'erations et coop\'erations.

\begin{dei}[Loi de remplacement compatible]
On dit qu'une loi de remplacement $\lambda$ est compatible avec les 
realtions $R$ et $S$ si les deux morphismes suivants sont injectifs 
$$\left\{ 
\begin{array}{l}
\underbrace{A}_{(1)}\boxtimes_c \underbrace{B}_{(2)} \to \Po\\
\underbrace{A}_{(2)}\boxtimes_c \underbrace{B}_{(1)} \to \Po.
\end{array}
\right. $$
\end{dei}

\begin{lem}
\label{FormeProperade} Toute prop\'erade de la forme $\Po=\F(V,\,
W)/(R\oplus D_\lambda \oplus S)$ d\'efinie par une loi de
remplacement $\lambda$ compatible v\'erifie l'isomorphisme de $\Sy$-bimodules
$$\Po \cong A\boxtimes_c B .$$
\end{lem}

\begin{deo}
L'hypoth\`ese de compatibilit\'e de la loi de remplacement $\lambda$ permet 
de montrer que le morphisme $A\boxtimes_c B \to \Po$ est injectif 
(\emph{cf.} \cite{FM}).

Ensuite, on montre que ce morphisme est surjectif. Pour cela, \`a tout 
\'el\'ement de $\Po$, on choisit un
repr\'esentant dans $\F(V\oplus W)$. Ce dernier s'\'ecrit comme un
somme finie de graphes indic\'es par des op\'erations de $V$ et
des coop\'erations de $W$. Pour chacun des graphes, on fixe
arbitrairement un sommet par niveau (\emph{cf.} chapitre 6) et on
permute op\'erations et coop\'erations pour placer l'ensemble des
coop\'erations au dessus de celui des op\'erations. L'\'el\'ement
ainsi obtenu appartient \`a $\F(V)\boxtimes_c\F(W)$ et, une fois
project\'e dans $A\boxtimes_c B$, il fournit l'ant\'ec\'edent
voulu. $\cqfd$
\end{deo}

Dans les deux exemples pr\'ec\'edents ce lemme montre que l'on
peut \'ecrire tout \'el\'ement de $\Po$ comme somme d'\'el\'ements
de $A\boxtimes_c B$, c'est-\`a-dire en mettant toutes les
coop\'erations en haut et toutes les op\'erations en bas. Au
niveau des big\`ebres de Lie, ce r\'esultat \'etait d\'ej\`a
pr\'esent dans \cite{EE} (section 6.4).

\subsection{Duale de Koszul d'une prop\'erade donn\'ee par une loi 
de remplacement}

\begin{pro}
\label{FormeProperadeDuale} Soit $\Po$ une prop\'erade de la forme
$\Po=\F(V,\, W)/(R\oplus D_\lambda \oplus S)$ d\'efinie par une
loi de remplacement $\lambda$ compatible avec les relations $R$ et 
$S$ et telle que la somme totale des
dimensions de $V$ et $W$ sur $k$, $\sum_{m,\, n} \dim_k (V\oplus
W)(m,\, n)$.

La prop\'erade duale est alors donn\'ee par
$$\Po^!=\F( \Sigma W^\vee \oplus \Sigma V^\vee)/(\Sigma^2 S^\perp
 \oplus \Sigma^2
D_{ ^t\lambda} \oplus  \Sigma^2 R^\perp),$$ c'\'est-\`a-dire par
la loi de remplacement $^t\lambda$. Et la coprop\'erade duale
v\'erifie l'isomorphisme de $\Sy$-bimodules
$$\Po^{\ac}\cong B^{\ac}\boxtimes_c A^{\ac}.$$
\end{pro}

\begin{deo}
L'hypoth\`ese sur la dimension des $\Sy$-bimodules g\'en\'erateurs
$V$ et $W$ de $\Po$, implique que la dimension de
$\F_{(\rho)}(\Sigma V^\vee \oplus \Sigma W^\vee)(m,\ n)$ est
finie, pour tous les entiers $\rho$, $m$ et $n$. On en conclut que
les dimensions de $\Po^{\ac}_{(\rho)}(m,\,n)$ et de
$\Po^!_{(\rho)}(m,\,n)$ sont aussi finies, d'o\`u ${\Po^!}^\vee
\cong \Po^{\ac}$.

On \'ecrit $R^\perp$ l'orthogonal de $R$ dans $\F_{(2)}(V^\vee)$,
$S^\perp$ l'orthogonal de $S$ dans $\F_{(2)}(W^\vee)$ et
$D_\lambda^\perp$ l'orthogonal de $D_\lambda$ dans $(I\oplus
\underbrace{W^\vee}_1)\boxtimes_c (I\oplus \underbrace{V^\vee}_1)
\bigoplus (I\oplus \underbrace{V^\vee}_1)\boxtimes_c (I\oplus
\underbrace{W^\vee}_1)$. Alors, la prop\'erade duale $\Po^!$ est
donn\'ee par
$$\Po^!=\F(\Sigma V^\vee \oplus \Sigma W^\vee)/(\Sigma^2 R^\perp \oplus \Sigma^2
D_\lambda^\perp \oplus \Sigma^2 S^\perp).$$ Il suffit maintenant
de remarquer que le $\Sy$-bimodule $D_\lambda^\perp$ correspond
\`a l'image du morphisme
$$(id,- {^t\lambda}) \ : \  (I\oplus \underbrace{V^\vee}_1)\boxtimes_c
(I\oplus \underbrace{W^\vee}_1) \to (I\oplus
\underbrace{V^\vee}_1)\boxtimes_c (I\oplus \underbrace{W^\vee}_1)
\bigoplus (I\oplus \underbrace{W^\vee}_1)\boxtimes_c (I\oplus
\underbrace{V^\vee}_1),$$ c'est-\`a-dire que
$D_\lambda^\perp=D_{^t\lambda}$.

En appliquant le lemme pr\'ec\'edent \`a la prop\'erade $\Po^!$,
on obtient que $\Po^!\cong B^!\boxtimes_c A^!$ puis
$\Po^{\ac}\cong B^{\ac} \boxtimes_c A^{\ac}$ en dualisant
lin\'eairement. $\cqfd$
\end{deo}

Les deux exemples donn\'es par les prop\'erades des big\`ebres de
Lie et  des big\`ebres de Hopf infinit\'esimales v\'erifient les
hypoth\`eses de cette proposition. Elles sont toutes les deux
engendr\'ees par un nombre fini de g\'en\'erateurs. Comme dans les
deux cas, aucun \'el\'ement de la forme
$\vcenter{\xymatrix@M=0pt@R=6pt@C=6pt{& \ar@{-}[d] &   \\
&\ar@{-}[dl]
\ar@{-}[dr] & \\
\ar@{-}[dr]& & \ar@{-}[dl]\\ &\ar@{-}[d] & \\ & & }}$ n'apparait
dans les relations, on les retrouve parmi les relations de la
prop\'erade duale (au sens de Koszul).

\begin{cor}
$ $
\begin{enumerate}
\item La prop\'erade duale de Koszul $\BLi^!$ de celle des
big\`ebres de Lie est donn\'ee par
$$\BLi^!=\F(V)/(R),$$
o\`u $V=\vcenter{\xymatrix@M=0pt@R=6pt@C=4pt{{\scriptstyle 1} & &
{\scriptstyle 2}
\\\ar@{-}[dr] &
&\ar@{-}[dl]  \\  &\ar@{-}[d] & \\
& & \\ & & }} \oplus \vcenter{\xymatrix@M=0pt@R=6pt@C=4pt{  & \ar@{-}[d]&  \\
&\ar@{-}[dl] \ar@{-}[dr]& \\ & & \\{\scriptstyle 1}&
&{\scriptstyle 2} }}$, c'est-\`a-dire une op\'eration commutative
et une coop\'eration cocommutative, et
\begin{eqnarray*}
R&=&\left\{ \begin{array}{l} \left(
\vcenter{\xymatrix@M=0pt@R=6pt@C=4pt{{\scriptstyle 1} & &
{\scriptstyle 2}& & {\scriptstyle 3}\\
\ar@{-}[dr] &
&\ar@{-}[dl] & & \ar@{-}[dl]  \\
& \ar@{-}[dr] & &\ar@{-}[dl]  & \\
& &\ar@{-}[d] & & \\
& & \\ & & }} - \vcenter{\xymatrix@M=0pt@R=6pt@C=4pt{{\scriptstyle
1} & &
{\scriptstyle 2}& & {\scriptstyle 3}\\
\ar@{-}[dr] &
&\ar@{-}[dr] & & \ar@{-}[dl]  \\
& \ar@{-}[dr] & &\ar@{-}[dl]  & \\
& &\ar@{-}[d] & & \\
& & \\ & & }} \right).k[\Sy_3]
 \\
\oplus \ k[\Sy_3].\left( \vcenter{\xymatrix@M=0pt@R=6pt@C=4pt{ & & & & \\
 & &\ar@{-}[d] & & \\
& & \ar@{-}[dl]  \ar@{-}[dr] & & \\
& \ar@{-}[dl]  \ar@{-}[dr] & &\ar@{-}[dr] & \\
& & & &\\ {\scriptstyle 1} & & {\scriptstyle 2}& &{\scriptstyle 3}
}} -
\vcenter{\xymatrix@M=0pt@R=6pt@C=4pt{ & & & & \\
 & &\ar@{-}[d] & & \\
& & \ar@{-}[dl]  \ar@{-}[dr] & & \\
& \ar@{-}[dl]   & &\ar@{-}[dl] \ar@{-}[dr] & \\
& & & &\\ {\scriptstyle 1} & & {\scriptstyle 2}& &{\scriptstyle 3}
}} \right)
\\ \ \\
\oplus \ \vcenter{\xymatrix@M=0pt@R=6pt@C=4pt{{\scriptstyle 1} & &
{\scriptstyle 2}
\\\ar@{-}[dr] &  &\ar@{-}[dl]  \\  &\ar@{-}[d] & \\
& \ar@{-}[dl]\ar@{-}[dr]& \\ & & \\ {\scriptstyle 1} & &
{\scriptstyle 2}}} -
 \vcenter{\xymatrix@M=0pt@R=6pt@C=2pt{ & {\scriptstyle 1} & &{\scriptstyle 2}
 \\ & \ar@{-}[d] &  &\ar@{-}[d]  \\  & \ar@{-}[dl]
 \ar@{-}[dr] &  & \ar@{-}[dl]\\
\ar@{-}[d]&  & \ar@{-}[d] & \\ & & & \\ {\scriptstyle 1} & &
{\scriptstyle 2} &}} \oplus
\vcenter{\xymatrix@M=0pt@R=6pt@C=4pt{{\scriptstyle 1} & &
{\scriptstyle 2}
\\\ar@{-}[dr] &  &\ar@{-}[dl]  \\  &\ar@{-}[d] & \\
& \ar@{-}[dl]\ar@{-}[dr]& \\ & & \\ {\scriptstyle 1} & &
{\scriptstyle 2}}} -
 \vcenter{\xymatrix@M=0pt@R=6pt@C=2pt{ & {\scriptstyle 2} & & {\scriptstyle 1}
  \\ & \ar@{-}[d] &  &\ar@{-}[d]  \\  & \ar@{-}[dl]
 \ar@{-}[dr] &  & \ar@{-}[dl]\\
\ar@{-}[d]&  & \ar@{-}[d] & \\ & & & \\ {\scriptstyle 1} & &
{\scriptstyle 2} &}} \oplus
\vcenter{\xymatrix@M=0pt@R=6pt@C=4pt{{\scriptstyle 1} & &
{\scriptstyle 2}
\\\ar@{-}[dr] &  &\ar@{-}[dl]  \\  &\ar@{-}[d] & \\
& \ar@{-}[dl]\ar@{-}[dr]& \\ & & \\ {\scriptstyle 1} & &
{\scriptstyle 2}}} -
\vcenter{\xymatrix@M=0pt@R=6pt@C=2pt{{\scriptstyle 1} & & {\scriptstyle 2} & \\
\ar@{-}[d] & &\ar@{-}[d] &  \\ \ar@{-}[dr] &
& \ar@{-}[dl] \ar@{-}[dr]& \\
& \ar@{-}[d]&  & \ar@{-}[d] \\ & & & \\ &{\scriptstyle 1} &
&{\scriptstyle 2}}} \oplus
\vcenter{\xymatrix@M=0pt@R=6pt@C=4pt{{\scriptstyle 1} & &
{\scriptstyle 2}
\\\ar@{-}[dr] &  &\ar@{-}[dl]  \\  &\ar@{-}[d] & \\
& \ar@{-}[dl]\ar@{-}[dr]& \\ & & \\ {\scriptstyle 1} & &
{\scriptstyle 2}}} -
\vcenter{\xymatrix@M=0pt@R=6pt@C=2pt{{\scriptstyle 2} & &
{\scriptstyle 1} &
\\ \ar@{-}[d] & &\ar@{-}[d] &  \\ \ar@{-}[dr] &
& \ar@{-}[dl] \ar@{-}[dr]& \\
& \ar@{-}[d]&  & \ar@{-}[d] \\ & & & \\ & {\scriptstyle 1} & &
{\scriptstyle 2}}} \\ \ \\
\oplus \ \vcenter{\xymatrix@M=0pt@R=7pt@C=7pt{& \ar@{-}[d] &   \\
&\ar@{-}[dl]
\ar@{-}[dr] & \\
\ar@{-}[dr]& & \ar@{-}[dl]\\ &\ar@{-}[d] & \\ & & }} \ .
\end{array}
\right. \end{eqnarray*}

Elle correspond \`a la diop\'erade des alg\`ebres de Frobenius
commutatives unitaires. Au niveau des $\Sy$-bimodules, on a
$$\BLi^!(m,\, n)=k.$$

\item La prop\'erade duale de Koszul $\IBi^!$ de celle des
big\`ebres de Hopf infinit\'esimales est donn\'ee par
$$\IBi^!=\F(V)/(R),$$
o\`u $V=\vcenter{\xymatrix@M=0pt@R=6pt@C=6pt{\ar@{-}[dr] &
&\ar@{-}[dl]  \\  &\ar@{-}[d] & \\
& & \\ & & }} \oplus \vcenter{\xymatrix@M=0pt@R=6pt@C=6pt{  & \ar@{-}[d]&  \\
&\ar@{-}[dl] \ar@{-}[dr]& \\ & &  \\ & &}}$  et
\begin{eqnarray*}
R&=&\left\{ \begin{array}{l} \vcenter{\xymatrix@M=0pt@R=6pt@C=6pt{
\ar@{-}[dr] &
&\ar@{-}[dl] & & \ar@{-}[dl]  \\
& \ar@{-}[dr] & &\ar@{-}[dl]  & \\
& &\ar@{-}[d] & & \\
& & \\ & & }} - \vcenter{\xymatrix@M=0pt@R=6pt@C=6pt{
 \ar@{-}[dr] &
&\ar@{-}[dr] & & \ar@{-}[dl]  \\
& \ar@{-}[dr] & &\ar@{-}[dl]  & \\
& &\ar@{-}[d] & & \\
& & \\ & & }}
 \\
\oplus \  \vcenter{\xymatrix@M=0pt@R=6pt@C=6pt{ & & & & \\
 & &\ar@{-}[d] & & \\
& & \ar@{-}[dl]  \ar@{-}[dr] & & \\
& \ar@{-}[dl]  \ar@{-}[dr] & &\ar@{-}[dr] & \\
& & & & }} -
\vcenter{\xymatrix@M=0pt@R=6pt@C=6pt{ & & & & \\
 & &\ar@{-}[d] & & \\
& & \ar@{-}[dl]  \ar@{-}[dr] & & \\
& \ar@{-}[dl]   & &\ar@{-}[dl] \ar@{-}[dr] & \\
& & & & }}
\\ \ \\
\oplus \  \vcenter{\xymatrix@M=0pt@R=6pt@C=6pt{\ar@{-}[dr] &  &\ar@{-}[dl]  \\  &\ar@{-}[d] & \\
& \ar@{-}[dl]\ar@{-}[dr]& \\ & & }} -
 \vcenter{\xymatrix@M=0pt@R=6pt@C=6pt{ & \ar@{-}[d] &  &\ar@{-}[d]  \\  & \ar@{-}[dl]
 \ar@{-}[dr] &  & \ar@{-}[dl]\\
\ar@{-}[d]&  & \ar@{-}[d] & \\ & & & }}  \oplus
\vcenter{\xymatrix@M=0pt@R=6pt@C=6pt{\ar@{-}[dr] &  &\ar@{-}[dl]  \\  &\ar@{-}[d] & \\
& \ar@{-}[dl]\ar@{-}[dr]& \\ & & }} -
\vcenter{\xymatrix@M=0pt@R=6pt@C=6pt{\ar@{-}[d] & &\ar@{-}[d] &
\\ \ar@{-}[dr] &
& \ar@{-}[dl] \ar@{-}[dr]& \\
& \ar@{-}[d]&  & \ar@{-}[d] \\ & & & }} \\ \ \\
\oplus \ 
\vcenter{\xymatrix@M=0pt@R=7pt@C=7pt{& \ar@{-}[d] &   \\
&\ar@{-}[dl]
\ar@{-}[dr] & \\
\ar@{-}[dr]& & \ar@{-}[dl]\\ &\ar@{-}[d] & \\ & & }} 
\oplus 
\vcenter{\xymatrix@R=1pt@C=1pt{ &\ar@{-}[d] & & &\ar@{-}[d] \\
 & *=0{}\ar@{-}[dl]\ar@{-}[dr] & & & *=0{}\ar@{-}[d] \\
 *=0{}\ar@{-}[ddrr] & & *=0{}\ar@{-}[ddll] | \hole  & & *=0{}\ar@{-}[d] \\
 & & & & *=0{}\ar@{-}[d] \\
 *=0{}\ar@{-}[d] & & *=0{}\ar@{-}[dr] & & *=0{}\ar@{-}[dl] \\
 *=0{}\ar@{-}[d] & & & *=0{}\ar@{-}[d] & \\
 & & & & }}
\oplus
\vcenter{\xymatrix@R=1pt@C=1pt{ *=0{}\ar@{-}[dddd]  & & & *=0{}\ar@{-}[d]  & \\
 & & & *=0{}\ar@{-}[dr] \ar@{-}[dl]  & \\
 & & *=0{}\ar@{-}[ddrr]  & & *=0{}\ar@{-}[ddll] | \hole  \\
 & & & & \\
 *=0{}\ar@{-}[dr] & & *=0{}\ar@{-}[dl] & & *=0{}\ar@{-}[d]  \\
 &  *=0{}\ar@{-}[d]& & & *=0{}\ar@{-}[d]  \\
 & & & &  }}
\oplus
\vcenter{\xymatrix@R=1pt@C=1pt{ *=0{}\ar@{-}[dd] &  &  &  *=0{}\ar@{-}[d] &   \\
 &  &  &  *=0{}\ar@{-}[dr]\ar@{-}[dl] &   \\
 *=0{}\ar@{-}[ddrr] &  & *=0{}\ar@{-}[ddll] | \hole  &  &   *=0{}\ar@{-}[dddd] \\
 &  &  &  &   \\
 *=0{}\ar@{-}[dr] &  &  *=0{}\ar@{-}[dl] &  &   \\
 &  *=0{}\ar@{-}[d] &  &  &   \\
 &  &  &  &   }}
\oplus
\vcenter{\xymatrix@R=1pt@C=1pt{ &  *=0{}\ar@{-}[d] &  & &   *=0{}\ar@{-}[dd] \\
 &  *=0{}\ar@{-}[dr]\ar@{-}[dl] &  &  &   \\
 *=0{}\ar@{-}[dddd] &  & *=0{}\ar@{-}[ddrr]  &  &  *=0{}\ar@{-}[ddll] | \hole  \\
 &  &  &  &   \\
 &  &  *=0{}\ar@{-}[dr] &  & *=0{}\ar@{-}[dl]   \\
 &  &  &  *=0{}\ar@{-}[d] &   \\
 &  &  &  &   }}

\ .

\end{array}
\right. \end{eqnarray*}

Au niveau des $\Sy$-bimodules, on a
$$\IBi^!(m,\, n)=k[\Sy_m]\otimes_k k[\Sy_n].$$
\end{enumerate}
\end{cor}

\subsection{Complexe de Koszul d'une prop\'erade donn\'ee par une loi de remplacement}

\begin{pro}
Soit $\Po$ une prop\'erade de la forme $\Po=\F(V,\, W)/(R\oplus
D_\lambda \oplus S)$ d\'efinie par une loi de remplacement
$\lambda$, telle que la somme totale des dimensions de $V$ et $W$
sur $k$, $\sum_{m,\, n} \dim_k (V\oplus W)(m,\, n)$, soit finie.
On d\'efinit les deux prop\'erades $A$ et $B$ par $A=\F(V)/(R)$ et
$B=\F(W)/(S)$. On suppose que $W$ est un $\Sy$-bimodule de
degr\'e homologique nul.\\
Si $A$ et $B$ sont des prop\'erades de Koszul, alors $\Po$ est
aussi une prop\'erade de Koszul.
\end{pro}

\begin{deo}
En appliquant le lemme~\ref{FormeProperade} et la
proposition~\ref{FormeProperadeDuale}, on montre que le complexe
de Koszul de $\Po$ est de la forme
$$\Po^{\ac}\boxtimes_c \Po =(B^{\ac}\boxtimes_c A^{\ac})\boxtimes_c(A\boxtimes_c B)=
B^{\ac}\boxtimes_c (A^{\ac}\boxtimes_c A )\boxtimes_c B.$$ On
introduit la filtration suivante du complexe de Koszul :
$F_n(\Po^{\ac}\boxtimes_c\Po)$ correspond au sous-$\Sy$-bimodule
de $B^{\ac}\boxtimes_c (A^{\ac}\boxtimes_c A )\boxtimes_c B$
engendr\'e par les graphes \`a $4$ niveaux pr\'esentant au plus
$n$ sommets sur le quatri\`eme niveau, c'est-\`a-dire celui
indic\'e par des \'el\'ements de $B^{\ac}$. Cette filtration est
stable par la diff\'erentielle du complexe de Koszul, elle induit
donc une suite spectrale not\'ee $E^*_{p,\, q}$. Le premier terme
de cette suite spectrale $E^0_{p,\, q}$ est compos\'e des
\'el\'ements de $B^{\ac}\boxtimes_c (A^{\ac}\boxtimes_c A
)\boxtimes_c B$ de degr\'e homologique $p+q$ et qui s\'ecrivent
avec exactement $p$ \'el\'ements de $B^{\ac}$. Et la
diff\'erentielle $d_0$ correspond \`a la diff\'erentielle de
Koszul de la prop\'erade $A$. Comme $A$ est une prop\'erade de
Koszul, on a $E^1_{p,\,q}=B^{\ac}\boxtimes_c B$ et $d_1$ est la
diff\'erentielle du complexe de Koszul de la prop\'erade $B$.
Comme celle-ci est acyclique, la suite spectrale est
d\'eg\'en\'er\'ee en $E^2$. Plus pr\'ecisement, on a
$$E^2_{p,\,q}=\left\{
\begin{array}{ll}
I & \textrm{si} \ p=0 \ \textrm{et} \ q=0, \\
0 & \textrm{sinon}.
\end{array}\right. $$
La filtration est exhaustive et born\'ee inf\'erieurement, on peut
appliquer \`a  $E^r_{p,\, q}$ le th\'eor\`eme classique de
convergence des suites spectrales (\emph{cf.} \cite{Weibel}
5.5.1). On obtient que la suite spectrale converge vers
l'homologie du complexe de Koszul de $\Po$. Ce complexe est donc
acyclique et $\Po$ est une prop\'erade de Koszul. $\cqfd$
\end{deo}

\begin{cor}
Les prop\'erades des big\`ebres de Lie $\BLi$ et des big\`ebres de
Hopf infinit\'esimales $\IBi$ sont de Koszul.
\end{cor}

\begin{deo}
Dans le cas $\BLi$, la prop\'erade $A$ est l'op\'erade de Koszul
des alg\`ebres de Lie $\Li$ et la prop\'erade $B$ est l'oppos\'ee
de $\Li$, $B=\Li^{op}$, qui est aussi de Koszul.

Dans le cas $\IBi$, la prop\'erade $A$ est l'op\'erade de Koszul
des alg\`ebres associatives $\A$ et la prop\'erade $B$ est
l'oppos\'ee de $\A$, $B=\A^{op}$, qui est aussi de Koszul. $\cqfd$
\end{deo}

\textsc{application :} La cobar construction sur la coprop\'erade
$\BLi^{\ac}$ est une r\'esolution de la prop\'erade $\BLi$. Si on
interpr\`ete cela en termes de cohomologie des graphes, on
retrouve les r\'esultats de M. Markl et A. A. Voronov \cite{MV}.
La cohomologie des graphes "commutatifs" connexes est \'egal \`a
la prop\'erade $\BLi$. Dans le cas de la prop\'erade $IBi$, on
trouve que la cohomologie des graphes ``ribbon'' connexes est
\'egale \`a la prop\'erade $\IBi$.

\chapter{S\'eries de Poincar\'e}

\thispagestyle{empty} 
\index{s\'erie de Poincar\'e}
 On
g\'en\'eralise ici une d\'emarche d\'ej\`a utilis\'ee dans le
cadre des alg\`ebres associatives (\emph{cf.} Y. Manin
\cite{Manin}) et des op\'erades binaires quadratiques (\emph{cf.}
V. Ginzburg et M.M. Kapranov \cite{GK}).

Soit $\Po$ une prop\'erade quadratique (issue \'eventuellement
d'un PROP). A $\Po$ (respectivement \`a $\Po^{\ac}$), on associe
une s\'erie de Poincar\'e $f_\Po$ (respectivement $f_{\Po^{\ac}}$)
dont les coefficients viennent de dimensions des modules
$\Po_{(\rho)}(m,\, n)$ (respectivement $\Po^{\ac}_{(\rho)}(m,\,
n)$) lorsque ces quantit\'es sont finies.

Au chapitre pr\'ec\'edent, nous avons vu que $\Po$ est de Koszul si et
seulement si le complexe de Koszul $\Po^{\ac}\boxtimes_c \Po$ est acyclique.
Dans ce cas, la caract\'eristique d'Euler-Poincar\'e du complexe de Koszul
est nulle et induit une relation fonctionnelle entre $f_\Po$ et
$f_{\Po^{\ac}}$.

Dans une premi\`ere section nous d\'efinissons les s\'eries de
Poincar\'e dans le cadre g\'en\'eral des $\Sy$-bimodules et nous
d\'emontrons le th\'eor\`eme v\'erifi\'e par les s\'eries
associ\'ees aux prop\'erades de Koszul. Dans la suite, nous
appliquons ce r\'esultat aux cas particuliers des alg\`ebres
associatives, des op\'erades binaires et des op\'erades
quadratiques non n\'ecessairement binaires. Dans les deux premiers
cas, nous retrouvons les \'equations fonctionnelles donn\'ees par
\cite{Loday1} et \cite{GK} respectivement. Le dernier cas est
nouveau. Et nous donnons un exemple d'application.

\section{S\'eries de Poincar\'e des prop\'erades}

On se place ici dans la cat\'egorie mono\"\i dale des
$\Sy$-bimodules gradu\'es par un poids munie du produit
$\boxtimes_c$ (\emph{cf.} chapitre $2$ section $1.3$). On rappelle
que toute prop\'erade quadratique est gradu\'ee par un poids (qui
vient du nombre de sommet des graphes d\'ecrivant la prop\'erade
libre) et que le complexe de Koszul se d\'ecompose de la mani\`ere
suivante :

\begin{pro}
\label{decompositioncomplexeKoszul} Le complexe de Koszul
associ\'e \`a une prop\'erade gradu\'ee par un poids
(\'even\-tu\-el\-le\-ment quadratique) $\Po$ se d\'ecompose en
somme directe de sous-complexes indic\'es par le nombre
d'``entr\'ees et de sorties'' (bigraduation naturelle du
$\Sy$-bimodule $\Po$) et par la graduation totale venant du poids.
Soit
$$ \K = \bigoplus_{m,\, n,\, d \, \ge 0} \K_{(d)}(m,\, n), $$
o\`u $\K_{(d)}(m,\, n)$  correspond \`a
$$ 0 \to \Po_{(d)}^{\ac}(m,\, n) \to \underbrace{\Po^{\ac}}_{(d-1)}\boxtimes_c
\underbrace{\Po}_{(1)}(m,\, n) \to \cdots \to
\underbrace{\Po^{\ac}}_{(1)}\boxtimes_c
\underbrace{\Po}_{(d-1)}(m,\, n) \to \Po_{(d)}(m,\, n) \to 0.$$
\end{pro}

\textsc{Remarque :} On a la m\^eme proposition pour le complexe de
Koszul $\mathcal{K}'$, version sym\'etrique de $\mathcal{K}$, o\`u
$\mathcal{K}'_{(d)}(m,\, n)=\left( \Po \boxtimes_c \Po^{\ac}
\right)_{(d)}$.

\begin{lem}
\label{lemmedim} Soit $V$ un $\Sy$-bimodule tel que
$\bigoplus_{m,\, n \in \mathbb{N}^*} V(m, \,n)$ soit de dimension
finie sur $k$.
Alors la dimension de $\F_{(d)}(V)(m,\,n)$ est finie. \\
Soit $\Po$ une prop\'erade  engendr\'ee par un tel module $V$
($\Po=\F(V)/(R)$). Alors la dimension de $\Po_{(d)}(m,\, n)$ est
finie.
\end{lem}

\begin{deo}
L'hyopth\`ese que le module $\bigoplus_{m,\, n \in \mathbb{N}^*}
V(m, \,n)$ soit de dimension finie sur $k$ implique qu'il existe
$M$ et $N$ dans $\mathbb{N}$ tels que $V(m,\, n)=0$ si $m\ge M$ ou
$n\ge N$. Comme il n'existe qu'un nombre fini de graphes connexes
avec moins de $d$ sommets et tels que chaque sommet admette un
nombre d'entr\'ees inf\'erieur \`a $M$ et un nombre de sorties
inf\'erieur \`a $N$, on obtient le r\'esultat voulu. $\cqfd$
\end{deo}

Le principal th\'eor\`eme \'enonc\'e dans cette th\`ese, le
crit\`ere de Koszul (\emph{cf.}
th\'e\-o\-r\`eme~\ref{criteredeKoszul}), affirme qu'une
prop\'erade $\Po$ est de Koszul si et seulement si les complexes
de Koszul $\mathcal{K}_{(d)}(m,\, n)$ sont acycliques pour $d\ge
1$ (ce qui est \'equivalent \`a l'acyclicit\'e des complexes de
Koszul
$\mathcal{K}_{(d)}'(m,\, n)$ pour $d\ge 1$).\\

Dans ce cas, la caract\'eristique d'Euler-Poincar\'e des complexes
de Koszul s'annule pour donner la formule
$$ \sum_{k=0}^{d}(-1)^k dim (\underbrace{\Po^{\ac}}_{(k)} \boxtimes_c
\underbrace{\Po}_{(d-k)})(m,\, n)=0.$$

Dans le cas o\`u la prop\'erade quadratique $\Po$ est engendr\'ee
par un espace de g\'en\'erateurs de dimension totale finie, on
peut appliquer le lemme~\ref{lemmedim} qui montre que tous les
modules consid\'er\'es dans la formule pr\'ec\'edente sont de
dimension finie. On peut donc la d\'evelopper de la mani\`ere
suivante :
\begin{eqnarray}
& & \sum_{k=0}^{d}(-1)^k dim (\underbrace{\Po^{\ac}}_{(k)}
\boxtimes_c
\underbrace{\Po}_{(d-k)})(m,\, n)= \nonumber \\
& &\sum_\Xi \sharp \mathcal{S}_c^{\bar{k},\, \bar{j}}
\frac{n!}{\oi !\, \oj !}\, \frac{m!}{\ok ! \, \ol !} \,
dim\Po_{o_1}^{\ac} (l_1,\, k_1)\ldots dim\Po_{o_b}^{\ac}(l_b,\,
k_b). \, dim\Po_{q_1}(j_1,\, i_1) \ldots dim\Po_{q_a}(j_a,\,
i_a),\nonumber
\end{eqnarray}
o\`u la somme $\Xi$ court sur les $n$-uplets $\oi$, $\oj$, $\ok$,
$\ol$, $\bar{o}$ et $\bar{q}$ tels que $|\oi|=n$, $|\oj|=|\ok|$,
$|\ol|=m$, $|\bar{o}|=k$ et $|\bar{q}|=d-k$.

\begin{dei}[S\'erie de Poincar\'e associ\'ee \`a un $\Sy$-bimodule]
A un $\Sy$-bimodule gradu\'e par un poids $\Po$ r\'eduit
(c'est-\`a-dire $\Po(m,\,n)=0$ d\`es que $m=0$ ou $n=0$), on
associe la \emph{s\'erie de Poincar\'e}
$$f_\Po(y,\, x,\, z)=\sum_{m,n \ge 1 \atop d\ge 0} \frac{\textrm{dim}\Po_{(d)}(m,\, n)}{m!\,
n!}y^mx^nz^d,$$ lorsque les modules $\Po_{(d)}(m,\, n)$ sont tous
de dimension finie sur $k$.
\end{dei}

Si on pose $$\Psi\big(g(y,\, X,\, z),\, f(Y,\, x,\,
z')\big)=\sum_{\Xi'} \sharp S_c^{\bar{k},\, \bar{j}}
\prod_{\beta=1}^b \frac{1}{k_\beta!}  \frac{\partial^{k_\beta
}g}{{\partial X}^{k_\beta}}(y,\, 0,\, z)
\prod_{\alpha=1}^a\frac{1}{j_\alpha!}\frac{\partial^{j_\alpha}f}{{\partial
Y}^{j_\alpha}}(0,\, x,\, z'),
$$
o\`u la somme $\Xi'$ court sur les $n$-uplets $\ok$ et $\oj$ tels
que $|\ok|=|\j|$. On a alors la proposition suivante :

\begin{thm}
\label{seriedepoincare}
 Toute prop\'erade $\Po$  de Koszul
engendr\'ee par un espace de g\'en\'erateurs de dimension finie
v\'erifie
$$\Psi\big(f_{\Po^{\ac}}(y,\, X,\, -z),\, f_\Po(Y,\, x,\, z)\big)=xy.$$
\end{thm}

\begin{deo}
On a
\begin{eqnarray*}
&& \Psi\big(f_{\Po^{\ac}}(y,\, X,\, -z),\, f_\Po(Y,\, x,\, z)\big) = \\
&&  \sum_{\Xi'}  \sharp S_c^{\bar{k},\, \bar{j}} \prod_{\beta=1}^b
\frac{1}{k_\beta!}  \frac{\partial^{k_\beta
}f_{\Po^{\ac}}}{{\partial X}^{k_\beta}}(y,\, 0,\, -z)
\prod_{\alpha=1}^a\frac{1}{j_\alpha!}\frac{\partial^{j_\alpha}f_\Po}{{\partial
Y}^{j_\alpha}}(0,\, x,\, z)  \\
& &\sum_{m,\, n \ge 1 \atop d\ge 0} \Big( \underbrace{\sum_{k=0}^d
(-1)^k \Big(\sum_\Xi \sharp S_c^{\bar{k},\, \bar{j}}
\prod_{\beta=1}^b \frac{dim\Po_{q_\beta}^{\ac}(l_\beta,\,
k_\beta)}{l_\beta!\, k_\beta!} \prod_{\alpha=1}^a
\frac{dim\Po_{o_\alpha}(j_\alpha,\, i_\alpha)}{j_\alpha!\,
i_\alpha!} \Big)}_{=0 \quad \textrm{pour} \quad d\ge 1} \Big) y^mx^nz^d    \\
& &=xy.
\end{eqnarray*}
$\cqfd$
\end{deo}

\textsc{Remarque :} Comme le complexe de Koszul $\mathcal{K}$ est
acyclique si et seulement si le complexe $\mathcal{K}'$ est
acyclique, dans le cas o\`u la prop\'erade $\Po$ est de Koszul, on
a la formule sym\'etrique
$$ \Psi\big(f_\Po (y,\, X,\, -z),\, f_{\Po^{\ac}}(Y,\, x,\, z)\big)=xy.$$

On va maintenant appliquer ce th\'eor\`eme aux sous-cat\'egories
pleine de $\Sy$-biMod que sont $k$-Mod et $\Sy$-Mod.

\section{S\'eries de Poincar\'e des alg\`ebres quadratiques}

On se place ici dans la sous-cat\'egorie mono\"\i dale pleine
$(\textrm{gr-Mod},\, \otimes_k, \, k)$ de $($gr-$\Sy$-biMod,
$\boxtimes_c,\, I)$. Dans cette sous-cat\'egorie mono\"\i dale,
une prop\'erade (un mono\"\i de) correspond \`a la notion
classique l'alg\`ebre associative gradu\'ee. Dans ce cas, la
prop\'erade (alg\`ebre) libre sur un espace $V$ correspond \`a
l'alg\`ebre tensorielle $T(V)$. Et une prop\'erade quadratique est
une alg\`ebre quadratique $A$ de la forme $A=T(V)/(R)$, o\`u $R\in
T_2(V)=V^{\otimes 2}$. Sa duale de Koszul est donn\'ee par
$A^{\ac}=\big(T(V^*)/(R^\perp)\big)^*= (R^\perp)^\perp$. Si
l'alg\`ebre $A$ est engendr\'ee par un module $V$ concentr\'e en
degr\'e $0$, alors $A^{\ac}_n=A^{\ac}_{(n)}$. On se place ici dans
le cas o\`u $V$ est un espace de dimension finie afin de
v\'erifier les hypoth\`eses du lemme~\ref{lemmedim}.

\begin{pro}
Le \textit{complexe de Koszul} a ici la forme suivante :
$$ \K_{(n)}\quad : \quad 0 \to A^{\ac}_{(n)} \to A^{\ac}_{(n-1)}\otimes A_{(1)}
\to \cdots \to
A^{\ac}_{(1)}\otimes A_{(n-1)}\to A_{(n)}\to 0 .$$
\end{pro}

Dans \cite{Loday1}, une alg\`ebre est dite de Koszul si les complexes
$\mathcal{K}_n$ sont acycliques pour $n\ge 1$. Cette d\'efinition est
\'equivalente \`a celle donn\'ee ici (\emph{cf.}
th\'eor\`eme~\ref{criteredeKoszul}).

\begin{dei}[S\'erie de Poincar\'e associ\'ee \`a un $k$-module]
La \textit{s\'erie de Poincar\'e} associ\'ee \`a un $k$-module
gradu\'e par un poids $A$ est donn\'ee par
$$f_A(x)=\sum_{n\ge 0} dim(A_{(n)})x^n,$$
si tous les modules $A_{(n)}$ sont de dimension finie.
\end{dei}

\textsc{Remarque :} Cette d\'efinition correspond bien \`a celle
donn\'ee dans la section pr\'ec\'edente. En consid\'erant $A$
comme un $\Sy$-bimodule, on a la s\'erie de Poincar\'e
$$f_A(y,\, x,\, z)=\sum_{n\ge 0} dim (A_{(n)}) z^n \, xy.$$
Les deux d\'efinitions sont reli\'ees par la formule
$$f_A(z)=f_A(1,\, 1,\, z). $$

Comme corollaire du th\'eor\`eme~\ref{seriedepoincare}, on a la
formule suivante dans le cadre des alg\`ebres associatives.

\begin{pro}
Si $A$ est une alg\`ebre de Koszul, alors les s\'eries de
Poincar\'e associ\'ees v\'erifient l'equation fonctionnelle :
$$f_A(x).f_{A^{\ac}}(-x)=1.$$
\end{pro}

\begin{deo}
En appliquant le th\'eor\`eme~\ref{seriedepoincare}, on a ici
\begin{eqnarray*}
f_A(z).f_{A^{\ac}}(-z) &=& \Psi\left( f_{A^{\ac}}(1,\, X,\, z), \,
f_A (1,\, Y,\, -z) \right) \\
&=& 1.
\end{eqnarray*}
$\cqfd$
\end{deo}

Dans la litt\'erature, on trouve plut\^ot la notion d'alg\`ebre
duale d\'efinie par $A^!=T(V^*)/(R^\perp)$. Nous avons vu \`a la
 fin du chapitre pr\'ec\'edent que ${A^!}^*=A^{\ac}$. On a alors l'\'egalit\'e
$f_{A^!}(x)=f_{A^{\ac}}(x)$, d'o\`u
$$f_A(x).f_{A^!}(-x)=1.$$

\textsc{Exemples :} L'exemple le plus connu vient du complexe
original de Koszul construit \`a partir de l'alg\`ebre
sym\'etrique $S(V)$ et de sa duale $\Lambda(V)$ (\emph{cf.} J.-L.
Koszul \cite{Koszul}).

\section{S\'eries de Poincar\'e des op\'erades binaires}

Nous reprenons ici le raisonnement effectu\'e par V. Ginzburg et
M. M. Kapranov
dans \cite{GK} pour les op\'erades binaires quadratiques. \\

On se place maintenant dans la cat\'egorie mono\"\i dale
$(\Sy\textrm{-Mod},\, \circ,\, I)$ des $\Sy$-modules munie du
produit mono\"\i dal
$$\Po\circ \Qo \
(n)=\bigoplus_{1\leq k\leq n \atop i_1+\cdots +i_k=n}
\Po(k)\otimes_{\Sy_k}\Qo(i_1)\otimes \cdots \otimes \Qo(i_k).$$
 Par d\'efinition, une op\'erade $(\Po,\,
\gamma,\, \eta)$ est un mono\"\i de dans cette cat\'egorie
mono\"\i dale. L'op\'erade libre $\F(V)$ sur $V$ peut \^etre
repr\'esent\'e par les arbres o\`u chaque sommet \`a $k$ entr\'ees
est indic\'e par une op\'eration de $V(k)$. Une op\'erade
quadratique est une op\'erade de la forme : $\Po=\F(V)/(R)$, avec
$R\in \F_{(2)}(V)$. Sa duale, au sens de Ginzburg et Kapranov, est
d\'efinie par $\Po^!=\F(V^\vee)/(R^\bot)$. Ici, nous avons
d\'efini sa coduale par $\Po^{\ac}$ qui correspond \`a
${\Po^!}^\vee$.

\subsection{Cas binaire g\'en\'eral}

Dans le cas o\`u $V$ est un $\Sy_2$-module (de dimension finie),
c'est-\`a-dire compos\'e uniquement d'op\'erations binaires, on
parle d'op\'erade \textit{binaire} (toutes les op\'erations sont
engendr\'ees par des op\'erations \`a deux variables). La
graduation en poids est directement li\'ee \`a la graduation
donn\'ee par le $\Sy$-module $\Po$ par la formule
$\Po_{(n-1)}=\Po(n)$.\index{op\'erade binaire}

\begin{pro}
Dans le cadre des op\'erades binaires quadratiques, le
\textit{complexe de Koszul} correspond \`a
$$\K_{(n-1)}(n) \quad : \quad 0 \to \Po^{\ac}(n) \to (\Po^{\ac}(n-1)\circ
\Po)(n) \to \cdots \to (\Po^{\ac}(2)\circ \Po)(n) \to \Po(n) \to
0.$$
\end{pro}

Dans le cas o\`u l'op\'erade $\Po$ est de Koszul, la caract\'erisque
d'Euler-Poincar\'e s'annule pour donner
: $$\sum_{k=1}^n (-1)^k dim((\Po^{\ac}(k)\circ \Po)(n))=0.$$ Or,
on a
$$dim((\Po^{\ac}(k)\circ
\Po)(n))=\sum_{i_1+\cdots +i_k=n}\frac{n!\, dim\Po^{\ac}(k)\,
dim\Po(i_1) \ldots  dim\Po(i_k)}{k!\, i_1! \ldots  i_k!}.$$ Ceci
pousse \`a la d\'efintion suivante :

\begin{dei}[S\'erie de Poincar\'e associ\'ee \`a une op\'erade
binaire quadratique]
A une op\'erade binaire quadratique $\Po$, on associe la
\textit{s\'erie de Poincar\'e} $$f_\Po(x)=\sum_{n\ge 1} (-1)^n
\frac{dim\Po(n)}{n!} x^n.$$
\end{dei}

\textsc{Remarque :} Cette d\'efinition apparait comme un cas particulier
des s\'eries de Poincar\'e associ\'ees aux prop\'erades. En effet, si on
consid\`ere $\Po$ comme une prop\'erade, on obitent
$$f_\Po(y,\, x,\, z)=\sum_{n\ge 1} \frac{\textrm{dim} \Po(n)}{n !}y\, x^n z^{n-1}.$$
Et les deux s\'eries sont reli\'ees par la formule
$$f_\Po(x)=f_\Po(1,\, -x,\, 1).$$

\begin{pro}
\label{EqOpBin} Si $\Po$ est une op\'erade binaire quadratique
de Koszul, les s\'eries de Poincar\'e v\'erifient l'\'equation fonctionnelle
$$f_{\Po^{\ac}}(f_\Po(x))=x.$$
\end{pro}

\begin{deo}
On a
\begin{eqnarray*}
f_{\Po^{\ac}}(f_\Po(x)) &=& f_{\Po^{\ac}}\left(1,\, -f_\Po(1,\, -x,\, 1),\, 1
 \right) \\
&=& \Psi\left(f_{\Po^{\ac}}(1,\, X,\, -1),\, f_\Po(Y,\, x,\, 1)
\right)=x
\end{eqnarray*}
$\cqfd$
\end{deo}

\textsc{Exemples :}
\begin{itemize}
\item Les op\'erades $\C$ et $\Li$ codant les alg\`ebres
commutatives et
les alg\`ebres de Lie sont des op\'erades binaires,
quadratiques et de Koszul.
Elles sont duales l'une de l'autre. Comme $\C(n)=k$,
on a $f_{\C}(x)=e^{-x}-1$
et comme $\textrm{dim}\ \Li(n)=(n-1)!$, on a $f_{\Li}(x)=-\ln(1+x)$.
Ces deux
s\'eries v\'erifient bien la relation pr\'ec\'edente.
(\emph{cf.} \cite{GK})

\item Un autre exemple est donn\'e par les op\'erades
$\mathcal{L}eib$ et
$\mathcal{Z}inb$ repr\'esentant les alg\`ebres de Leibniz
et les alg\`ebres Zinbiel
(\emph{cf.} \cite{Loday2}). Ces deux op\'erades sont binaires,
quadratiques et de
Koszul. Elles sont aussi duales l'une de l'autre.
On a $\mathcal{L}eib(n)=k[\Sy_n]$, d'o\`u
$f_{\mathcal{L}eib}(x)=\frac{-x}{1+x}$. De m\^eme, on sait
que $\mathcal{Z}inb(n)=k[\Sy_n]$, d'o\`u on tire
$f_{\mathcal{Z}inb}(x)=\frac{-x}{1+x}$. L'\'equation fonctionnelle
est encore v\'erifi\'ee ici.
\end{itemize}

\subsection{Cas non-sym\'etrique}

\index{op\'erade non-sym\'etrique}
Si les op\'erades en question
peuvent \^etre d\'ecrites sans l'action du groupe sym\'etrique, on
parle d'op\'erades \emph{non-sym\'etriques} ou
\emph{non-$\Sigma$}. Soit $\Po$ une telle op\'erade, alors le
$\Sy_n$-module $\Po(n)$ est un $\Sy_n$-module libre que l'on note
$\Po'(n) \otimes_k k\lbrack \Sy_n \rbrack$. Dans ce cas, l'objet
sous-jacent \`a l'op\'erade est juste le $k$-module gradu\'e
$\Po'$. La s\'erie de Poincar\'e devient $f_\Po(x)=\sum_{n\ge
1}(-1)^n \textrm{dim}\, \Po'(n) x^n$
et le r\'esultat pr\'ec\'edent reste vrai.\\

\textsc{Exemples :}
\begin{itemize}

\item L'exemple le plus c\'el\'ebre est celui de l'op\'erade $\A$
des alg\`ebres associatives. Cette op\'erade est une op\'erade non
sym\'etrique binaire quadratique et de Koszul. Elle est autoduale.
Et, sa s\'erie de Poincar\'e vaut $f_{\A}(x)=\frac{-x}{1+x}$.

\item Un autre exemple de telles op\'erades est donn\'e par les
op\'erades $\mathcal{D}ias$ et $\mathcal{D}end$ repr\'esentant les
dig\`ebres et les alg\`ebres dendriformes (\emph{cf.}
\cite{Loday2}). Le $\Sy_n$-module $\mathcal{D}ias(n)$ correspond
\`a $n$ copies de $k[\Sy_n]$, d'o\`u
$f_{\mathcal{D}ias}(x)=\frac{-x}{(1+x)^2}$. Quant au
$\Sy_n$-module $\mathcal{D}end(n)$, il est isomorphe \`a une somme
directe indic\'ee par les arbres binaires planaires \`a $n$
sommets. Le cardinal de cet exemple est \'egal au nombre de
Catalan $c_n$. On a donc
$f_{\mathcal{D}end}(x)=\frac{-1-2x+\sqrt{1+4x}}{2x}$.

\end{itemize}

\section{S\'eries de Poincar\'e des op\'erades quadratiques}

Le cas pr\'ec\'edent \'etait d\'ej\`a connu et donn\'e dans
 \cite{GK}. Cette article d\'ecrit la th\'eorie de Koszul des
op\'erades binaires quadratiques (comme $\A$, $\C$, $\Li$ etc...).
Dans ce cas particulier, on s'est fortement servi de la relation
$\Po_{(n-1)}=\Po(n)$ pour pouvoir \'ecrire l'\'equation
fonctionnelle v\'erifi\'ee par les s\'eries de Poincar\'e. Pour
d\'epasser cette difficult\'e dans le cas g\'en\'eral, il faut
introduire le poids dans la d\'efinition de la s\'erie de
Poincar\'e. Ceci apparait naturellement lorsque l'on \'ecrit la
s\'erie de Poincar\'e d'une op\'erade gradu\'e par un poids comme
une prop\'erade.

\subsection{Cas g\'en\'eral}

On reformule la propositon~\ref{decompositioncomplexeKoszul} dans
le cadre des op\'erades.

\begin{pro} Le complexe de Koszul se d\'ecompose en somme directe en
fonction du "nombre de feuilles" (graduation du $\Sy$-module) et
de la graduation totale venant du le poids
$$\mathcal{K}=\bigoplus_{d,\, n \ge 0} \mathcal{K}_{(d)}(n)$$ o\`u
$\mathcal{K}_{(d)}(n)$ correspond \`a :
$$  \K_{(d)}(n)  \quad : \quad 0 \to \Po_{(d)}^{\ac}(n) \to \Po_{(d-1)}^{\ac}\circ
\underbrace{\Po}_{(1)}(n) \to \cdots \to \Po_{(1)}^{\ac}\circ
\underbrace{\Po}_{(d-1)}(n) \to \Po_{(d)}(n) \to 0.$$
\end{pro}
On a le m\^eme r\'esultat pour le complexe de Koszul
$\mathcal{K}'$ o\`u
$$ \K'_{(d)}(n)  \quad : \quad 0 \to \Po_{(d)}^{\ac}(n) \to \Po_{(1)}(n)\circ
\underbrace{\Po^{\ac}}_{(d-1)}(n) \to \cdots \to \Po_{(d-1)}\circ
\underbrace{\Po^{\ac}}_{(1)}(n) \to \Po_{(d)}(n) \to 0.$$

De la m\^eme mani\`ere que nous avions d\'efini une s\'erie de
Poincar\'e pour les $\Sy$-bimodules gradu\'es pas un poids, on
peut le faire pour un $\Sy$-module gradu\'e par un poids.

\begin{dei}[S\'erie de Poincar\'e associ\'ee \`a un
$\Sy$-module gradu\'e par un poids]

A tout $\Sy$-module gradu\'e par un poids $\Po$, on associe la
\emph{s\'erie de Poincar\'e}
$$f_\Po(x,\, y)=\sum_{d\ge 0,\, n\ge
1}\frac{dim\Po_{(d)}(n)}{n !}\, y^d\, x^n.$$
\end{dei}

\textsc{Remarque :} Cette d\'efintion est bien un cas particulier
de s\'erie de Poincar\'e associ\'ee \`a un $\Sy$-bimodule. En
effet, ces deux d\'efinitions v\'erifient la relation
$$f_\Po(y,\, x,\ z)=y f_\Po(x,\, z).$$

Dans ce cas particulier des op\'erades l'\'equation fonctionnelle
donn\'ee au th\'eor\`eme~\ref{seriedepoincare} se simplifie.

\begin{pro}
\label{EqOp}
 Soit une op\'erade $\Po$ quadratique et de Koszul, alors les
 s\'eries de Poincar\'e v\'erifient les \'equations
$$f_{\Po^{\ac}}(
f_{\Po}(x,\, y),\, -y)=x \quad \mathrm{et} \quad
 f_{\Po}( f_{\Po^{\ac}}(x,\, y),\, -y)=x.$$
\end{pro}

\begin{deo}
D'apr\`es le th\'eor\`eme~\ref{seriedepoincare}, on a
\begin{eqnarray*}
f_{\Po^{\ac}}( f_{\Po}(x,\, y),\, -y) &=& \Psi\left(
f_{\Po^{\ac}}(1,\, X,\, -y),\, f_\Po(Y,\, x,\, y) \right) \\
&=& x
\end{eqnarray*}
$\cqfd$
\end{deo}

\subsection{Cas non-sym\'etrique}

Comme dans le cas binaire, si $\Po$ est une op\'erade
non-sym\'etrique, on peut simplifier la s\'erie de Poincar\'e en
posant
$$f_\Po(x,\, y)=\sum_{d,\, n} dim\Po'_{(d)}(n)\, y^d x^n$$
 et la proposition pr\'ec\'edente reste vraie.

\subsection{Exemple : Cas d'une op\'erade libre engendr\'ee
par une op\'eration $n$-aire pour tout $n$}
 On consid\`ere
l'op\'erade non sym\'etrique libre $\Po$ engendr\'ee par l'espace
$V=k\{\,
\vcenter{\xymatrix@M=0pt@R=6pt@C=6pt{\ar@{-}[dr] &  &\ar@{-}[dl]  \\  &\ar@{-}[d] & \\
& &}} ,\,
\vcenter{\xymatrix@M=0pt@R=6pt@C=6pt{\ar@{-}[dr] & \ar@{-}[d] &\ar@{-}[dl]  \\
&\ar@{-}[d] & \\  & &}}
 ,\, \ldots \}$. Donc $\Po=\F(V)$ et
$\Po^{\ac}=k\oplus V$.\\

Remarquons que $\F(V)(n)$ est de dimension finie alors que
l'espace des g\'en\'erateurs $V$ est de dimension infinie.

\begin{pro}
L'op\'erade $\F(V)$ est de Koszul.
\end{pro}

\begin{deo}
\rm Il s'agit de montrer que $\K_{(d)}(n)$ est acyclique pour
$n\ge 2$. Comme $\Po^{\ac}=k\oplus V$, ce complexe se r\'eduit \`a
:

$$  \K_{(d)}(n)  \quad : \quad 0 \to 0 \to  \cdots \to 0 \to \Po_{(1)}^{\ac}\circ
\underbrace{\Po}_{(d-1)}(n) \to \Po_{(d)}(n) \to 0.$$ Or,
l'op\'erade $\Po$ \'etant sans relation, le morphisme $d : V\circ
\underbrace{\Po}_{(d-1)}(n) \to \Po_{(d)}(n)$ est un isomorphisme.
D'o\`u le r\'esultat. $\cqfd$
\end{deo}

De la forme de $\Po^{\ac}$ on tire la s\'erie de Poincar\'e
associ\'ee,
$$ f_{\Po^{\ac}}(x,\, y)= \sum_{d\ge 0,\, n\ge 1} dim\Po_{(d)}^{\ac}(n)x^n y^d=x+\sum_{n\ge
2} x^n y=x+y \frac{x^2}{1-x}.$$ D'apr\`es la
proposition~\ref{EqOp}, on obtient
$$f_\Po(x,\, y)-y\frac{f^2_\Po(x,\, y)}{1-f_\Po(x,\, y)}=x.$$
Ce qui implique
$$(y+1)f_\Po^2(x,\, y)-(1+x)f_\Po(x,\, y)+x=0.$$

Soit $P_n(x)$ le polyn\^ome de Poincar\'e du polytope du Stasheff
de dimension $n$, aussi appel\'e associa\`edre et not\'e $K^n$ ou
$K_{n+2}$. Ce polyn\^ome s'\'ecrit $P_n(y)=\sum_{k=0}^n \sharp
Cel_{k}^{n}.y^k$ o\`u $Cel_{k}^{n}$ repr\'esente l'ensembles des
cellules de dimension $k$ du polytope de Stasheff de dimension
$n$. On pose, $f_K(x,\, y)=\sum_{n\ge 0}P_n(y)x^n$ la s\'erie
g\'en\'eratrice associ\'ee \`a ces polyn\^omes.

Les cellules de dimension $k$ de $K^n$ peuvent \^etre incid\'ees
par les arbres planaires \`a $n+2$ feuilles et $n+1-k$ sommets.
Cette bijection implique que $\sharp Cel_{k}^{n}=dim
\Po_{n+1-k}(n+2)$. Ce qui donne au niveau des s\'eries
g\'en\'eratrices :

\begin{eqnarray}
f_\Po(x,\, y) &=& x+ \sum_{n\ge 2} \sum_{k=1}^{n-1} dim \Po_k^n \, y^kx^n \nonumber \\
&=& x+\sum_{n\ge2} \Big( \sum_{k=1}^{n-1} \sharp
Cel_{n-2-(k-1)}^{l-2}\, y^k \Big) x^n
\nonumber \\
&=& x+yx^2\sum_{n\ge 0} \Big( \sum_{k=0}^{n}\sharp Cel_{n-k}^{n}\,
y^k \Big) x^n
\nonumber \\
&=& x +yx^2\sum_{n\ge 0} P_n\Big(\frac{1}{y}\Big)\, (xy)^n =
x+yx^2f_K\Big(xy,\, \frac{1}{y}\Big). \nonumber
\end{eqnarray}

En utilisant l'equation v\'erifi\'ee par $f_\Po$, on en trouve une
pour $f_K$

$$ ((1+y)x^2)f_K^2(x,\, y)+(-1+(2+y)x)f_K(x,\, y)+1=0.$$

Ce qui donne la formule suivante :

\begin{pro}
La s\'erie g\'en\'eratrice associ\'ee aux polytopes de Stasheff
v\'erifie
$$ f_K(x,\, y)=\frac{1+(2+y)x-\sqrt{1-2(2+y)x+y^2x^2}}{2(1+y)x^2}.$$
\end{pro}

\textsc{Remarque :} On retrouve ici le m\^eme formule d\'ej\`a que
celle de J.-L. Loday et M. Ronco dans \cite{JLLTri}. La m\'ethode
utilis\'ee dans cette article est la m\^eme qu'ici mais repose sur
l'op\'erade de Koszul des trig\`ebres dendriformes. Notons que
cette \'egalit\'e se d\'emontre aussi de mani\`ere purement
combinatoire \`a l'aide de la formule de r\'ecurrence qui donne le
nombre d'arbres planaire \`a $n$ sommets en fonction du nombre
d'arbres planaires \`a $k$ sommets, avec $k \le n$.

\newpage

\thispagestyle{empty}

\newpage
\fancyhead{}
\fancyhead[RO]{\MakeUppercase{Index}}
\fancyhead[LE]{\MakeUppercase{Index}}
\printindex

\newpage

\fancyhead{}
\fancyhead[RO]{\MakeUppercase{Table des mati\`eres}}
\fancyhead[LE]{\MakeUppercase{Table des mati\`eres}}
\thispagestyle{empty}
\tableofcontents

\end{document}